\documentclass{amsart}
\usepackage[textsize=scriptsize]{todonotes}

\usepackage{etoolbox}
\usepackage[margin=1in,marginparwidth=0.8in]{geometry}
\usepackage{enumerate}

\usepackage{amsmath}
\usepackage{amssymb}
\usepackage{amsthm}
\usepackage{amscd}
\usepackage[pdfusetitle,unicode,hidelinks]{hyperref}
\usepackage{bbm}
\usepackage{etoolbox}
\usepackage{stmaryrd}

\usepackage[utf8]{inputenc}
\usepackage[T1]{fontenc}

\usepackage[textsize=scriptsize]{todonotes}
\newcommand{\NB}[1]{\todo[color=gray!40]{#1}}
\newcommand{\TODO}[1]{\todo[color=red]{#1}}

\newcounter{proposition}

\newtheorem{proposition}[equation]{Proposition}
\newtheorem{corollary}[equation]{Corollary}
\newtheorem{lemma}[equation]{Lemma}
\newtheorem{theorem}[equation]{Theorem}

\newtheorem*{conjecture*}{Conjecture}
\newtheorem*{theorem*}{Theorem}
\newtheorem*{corollary*}{Corollary}
\newtheorem*{proposition*}{Proposition}
\newtheorem*{lemma*}{Lemma}
\theoremstyle{definition}
\newtheorem{definition}[equation]{Definition}
\newtheorem{construction}[equation]{Construction}

\newtheorem*{definition*}{Definition}
\newtheorem*{construction*}{Construction}
\theoremstyle{remark}
\newtheorem{remark}[equation]{Remark}
\newtheorem*{remark*}{Remark}

\newtheorem{example}[equation]{Example}
\newtheorem*{example*}{Example}

\newcommand{\fpsr}[1]{\llbracket #1 \rrbracket}
\newcommand{\id}{\operatorname{id}}
\newcommand{\Z}{\mathbb{Z}}
\def\C{\mathbb C}
\newcommand{\N}{\mathbb{N}}
\newcommand{\Q}{\mathbb{Q}}
\newcommand{\F}{\mathbb{F}}

\let\scr=\mathcal
\let\bb=\mathbb
\newcommand{\Gm}{{\mathbb{G}_m}}

\def\A{\bb A}
\def\P{\bb P}
\def\R{\bb R}
\def\1{\mathbbm{1}}

\newcommand{\eff}{{\mathrm{eff}}}
\newcommand{\veff}{{\text{veff}}}
\newcommand{\DM}{\mathcal{DM}}
\newcommand{\SH}{\mathcal{SH}}

\DeclareMathOperator*{\colim}{colim}

\let\lim=\relax
\DeclareMathOperator*{\lim}{lim}
\def\Map{\mathrm{Map}}
\def\iMap{\ul{\mathrm{Map}}}
\def\map{\mathrm{map}}
\def\CAlg{\mathrm{CAlg}}
\def\CMon{\mathrm{CMon}}
\def\NAlg{\mathrm{NAlg}}

\def\PSh{\mathcal{P}}
\def\Ind{\mathrm{Ind}}
\def\Span{\mathrm{Span}}
\def\Cat{\mathcal{C}\mathrm{at}{}}
\def\Spc{\mathcal{S}\mathrm{pc}{}}

\def\Fun{\mathrm{Fun}}

\newcommand{\Spec}{\mathrm{Spec}}
\newcommand{\gp}{\mathrm{gp}}

\newcommand{\wequi}{\simeq}
\newcommand{\Mod}{\mathcal{M}\mathrm{od}}
\def\adj{\leftrightarrows}

\DeclareRobustCommand{\ul}{\underline}
\newcommand{\heart}{\heartsuit}

\newcommand{\Hom}{\operatorname{Hom}}

\def\op{\mathrm{op}}

\let\cat=\mathrm
\def\Sm{{\cat{S}\mathrm{m}}}

\def\Sch{\cat{S}\mathrm{ch}{}}

\def\Nis{\mathrm{Nis}}
\def\Zar{\mathrm{Zar}}

\def\mot{\mathrm{mot}}
\newcommand{\et}{{\acute{e}t}}

\newcommand{\fr}{\mathrm{fr}}

\def\ph{\mathord-}

\newcommand{\Th}{\mathrm{Th}}

\newcommand{\Shv}{\cat{S}\mathrm{hv}}
\newcommand{\cof}{\mathrm{cof}}
\newcommand{\fib}{\mathrm{fib}}

\newcommand{\GL}{\mathrm{GL}}
\newcommand{\SL}{\mathrm{SL}}

\newcommand{\MGL}{\mathrm{MGL}}

\def\Alg{\mathrm{Alg}}
\def\cAlg{\mathrm{cAlg}}

\usepackage{xstring}
\usepackage{xr}
\newcommand{\myexternaldocument}[2][] {{
\let\nl\newlabel
\def\nlxx##1##2##3##4##5##6{\nl{##1}{{#1##2}{}{}{}{}}}
\renewcommand\newlabel[2]{\IfBeginWith{##1}{tocindent}{}{\nlxx{##1}##2}}
\externaldocument{#2}
}}

\newtoggle{final}
\toggletrue{final}

\iftoggle{final} {
\renewcommand{\NB}[1]{}
\renewcommand{\TODO}[1]{}
\renewcommand{\todo}[1]{}
\newcommand{\gpt}[1]{}
}{ 
}

\usepackage{comment}
\usepackage{array,longtable}
\usepackage{tikz}
\usetikzlibrary{matrix,arrows}
\usepackage{tikz-cd}

\numberwithin{proposition}{section}
\numberwithin{equation}{section}

\def\fet{\mathrm{f\acute et}}
\def\all{\mathrm{all}}
\def\cell{\mathrm{cell}}
\def\tate{\mathrm{tate}}
\def\NSym{\mathrm{NSym}}
\def\LSym{\mathrm{LSym}}
\def\Sym{\mathrm{Sym}}
\def\NMon{\mathrm{NMon}}
\def\cNMon{\mathrm{cNMon}}
\def\NHopf{\mathrm{NHopf}}
\def\Hopf{\mathrm{Hopf}}
\def\CHopf{\mathrm{CHopf}}
\def\dHopf{\mathrm{dHopf}}
\def\Hopf{\mathrm{Hopf}}
\def\cNAlg{\mathrm{cNAlg}}
\def\bNAlg{\mathrm{bNAlg}}
\def\cCAlg{\mathrm{cCAlg}}
\def\Ext{\mathrm{Ext}}
\newcommand{\cMod}{\mathrm{cMod}}
\def\cn{\mathrm{cn}}
\def\DAlg{\mathrm{DAlg}}
\def\imap{\ul{\mathrm{map}}}
\newcommand{\comp}{{{\kern -.5pt}\wedge}}
\newcommand{\Sp}{\mathrm{Sp}}
\newcommand{\SW}{\mathrm{SW}}
\newcommand{\Ga}{{\mathbb{G}_a}}
\newcommand{\Gr}[1]{{\mathrm{Gr}_{#1}}}
\def\Free{\mathrm{F}}
\def\Cofree{\mathrm{C}}
\def\FFree{\mathrm{FFree}}
\def\Perf{\mathrm{Perf}}

\def\dHFp{dH\F_p}
\def\hotimes{\odot}
\def\hLoops{L}
\def\h{\mathrm h}
\newcommand{\Vect}{\mathrm{Vect}}
\newcommand{\CB}{\mathrm{CB}}
\newcommand{\Tot}{\mathrm{Tot}}
\newcommand{\naive}{\mathrm{naive}}

\title{Motivic Eilenberg--Mac Lane spaces at the characteristic}
\date{\today}

\author{Tom Bachmann}
\author{Michael J. Hopkins}

\begin{document}
\maketitle

\begin{abstract}
We provide a novel, homotopy theoretic model for the motivic Eilenberg--Mac Lane spaces $K(\Z_{(p)}(n), 2n) \in \Spc(k)$ over fields of characteristic $p$.
Using this, we determine their motives and, consequently, determine the motivic Steenrod algebra at the characteristic.
Applications include the Hopkins--Morel isomorphism over arbitrary base schemes as well as the computation of the $p$-completed motivic stable stems and algebraic cobordism groups of fields of characteristic $p$.

Along the way we expand the theory of normed algebraic structures in motivic homotopy theory, compare them to derived algebra at the characteristic, and produce a functor computing the motivic homology of the $\P^1$-loop spaces of (certain) motivic spaces in terms of the motives.
\end{abstract}

\tableofcontents

\section{Introduction}
\subsection{The motivic Steenrod algebra}
Cohomology operations have been investigated by topologists for a long time.
The question is straightforward to formulate, and we do it here motivically: working over a field $k$ and with $\F_p$-coefficients for a fixed prime $p$, what are (all) the natural operations among the motivic cohomology groups $H^{**}(X, \F_p)$ for smooth $k$-schemes $X$?
The classical analog has been investigated by Steenrod \cite{steenrod1952reduced}, Serre \cite{Serre:em}, Cartan \cite{Cartan:em} and many others.
It turns out that there are certain operations called \emph{Steenrod power operations}, every other natural operation can be expressed in terms of these, and all the universal relations among them are known.
Attempting to obtain similar results for motivic cohomology has been a long-running effort.
Voevodsky constructed motivic versions of the Steenrod power operations as long as $k$ is of characteristic $\ne p$ \cite{voevodsky2003reduced} and proved that if $k$ has characteristic zero these generate all the bistable operations, and determined all the relations \cite{voevodsky2010motivic}.
Using these tools, he established the Milnor and Bloch--Kato conjectures \cite{voevodsky2003motivic,voevodsky-BK} and received a Fields Medal for his work.\NB{Chronology: the Fields Medal was awarded in 2002, before the published Bloch--Kato proof.}
Somewhat later, Hoyois--Kelly--{\O}stv{\ae}r generalized Voevodsky's results to fields of characteristic $\ne p$.
However, fields of characteristic $p$ have seemed to mostly defy attack.
Frankland--Spitzweck \cite{frankland-spitzweck} have managed to at least produce the expected operations; more recently Annala--Elmanto \cite{annala-elmanto} proved that they satisfy the expected relations. Whether these account for \emph{all} operations remained open.
The importance of this problem was highlighted in Lurie's 2026 ICM address \cite[Remark~9.5]{lurie-prismatic}, where he identified a related motivic question---also resolved by our computation---as an obstacle to his conjectural \emph{prismatic stable homotopy theory}.

In this article we solve this problem in full: we determine all the natural operations in motivic cohomology.

\begin{remark}
A follow-up question is what makes this problem harder in characteristic $p$ than in characteristic zero.
Slightly curiously, given the fact that Voevodsky's results famously used resolution of singularities, the lack thereof in positive characteristic does \emph{not} seem to be the primary stumbling block.
Indeed no results like ours can be found in the literature, even conditional on resolution of singularities.
Instead the principal difficulty in characteristic $p$ seems to be that cycles (which are usually used to model motivic cohomology) are just more complicated.
While in characteristic zero cycles on a smooth scheme can be captured by schemes (namely, symmetric powers), in characteristic $p$ no such geometric model is known.
Instead we provide a new model for motivic cohomology at the characteristic not utilizing cycles, which additionally lends itself toward easy determination of the operations.
We ultimately rely on a universal characterization of motivic cohomology inside the motivic stable homotopy category, conjectured by Voevodsky \cite[Conjecture~10]{voevodsky-slice-filtration} and established by Levine \cite[Theorem~10.5.1]{levine2008homotopy}.
\end{remark}

\subsection{Applications} \label{subsec:intro-app}
There are many open problems in motivic homotopy theory which are well-known to be equivalent to determining the motivic Steenrod algebra at the characteristic.

We list some of them here, without much claim to originality.
See \S\ref{sec:applications} for more details and examples.

The dual motivic Steenrod algebra can be determined over any base (and hence the same is true for the motivic Steenrod algebra).
\begin{theorem}[see Theorem \ref{thm:dual-steenrod}]
Let $S$ be any scheme.
The motivic dual Steenrod algebra over $S$ has the expected form.
In other words, the canonical map \[ \bigoplus_{I \in B}\Sigma^{p(I),q(I)} H\F_p \to H\F_p \otimes H\F_p \in \SH(S) \] constructed in \cite[Theorem 10.26]{spitzweck2012commutative} is an equivalence.
\end{theorem}

The so-called Hopkins--Morel isomorphism holds over any base.
\begin{theorem}[see Theorem \ref{thm:hopkins-morel}]
Let $S$ be any scheme.
Then $\MGL/(t_1, t_2, \dots) \wequi H\Z$.
\end{theorem}

The $p$-completed motivic stable stems of fields of characteristic $p$ can be determined fully.
\begin{theorem}[see Theorem \ref{thm:stable-stems}]
Let $k$ be a field of characteristic $p$.
Then $\pi_{**} \1_p^\comp \wequi \Ext(BP_{2*}BP) \otimes_{\Z[\eta]} K_*^{MW}(k)_p^\comp$.
\end{theorem}

The same is true for algebraic cobordism.
\begin{theorem}[see Theorem \ref{thm:MGL-htpy}]
Let $k$ be a field of characteristic $p$.
Then $\pi_{**} \MGL_p^\comp \wequi K_*^M(k)_p^\comp \otimes L_*$.
\end{theorem}

\subsection{Ingredients}
Our main effort is to construct a novel model for the motivic Eilenberg--Mac Lane spaces $K(\Z_{(p)}(n), 2n) \in \Spc(k)$ (representing the motivic cohomology groups $H^{2n,n}(\ph, \Z_{(p)})$ in the homotopy category of motivic spaces), where $k$ is a field of characteristic $p$.
Instead of utilizing the geometry of cycles, it is homotopy theoretic in nature and captures the fact that motivic cohomology has both additive and multiplicative transfers along finite étale morphisms.
In other words, given $f: X \to Y \in \Sm_k$ finite étale of degree $d$, there are operations \[ f_*: H^{2n,n}(X, \Z) \to H^{2n,n}(Y, \Z) \text{ and } f_\otimes: H^{2n,n}(X, \Z) \to H^{2nd,nd}(Y, \Z) \] (the latter of which is not even additive!).
They satisfy a number of relations, coherently for all $f$ at once.
We essentially prove that the $\{K(\Z_{(p)}(n), 2n)\}_n$ are the universal family of spaces representing a cohomology theory with such operations, which is in weight $1$ given by the usual motivic cohomology (i.e., line bundles).
In order to make this precise, we employ the theory of \emph{norms} in motivic homotopy theory, developed by the first author and Hoyois in \cite{bachmann-norms}.

There are some further important ingredients.
In order to establish a new model for the spaces $K(\Z_{(p)}(n), 2n)$, we clearly need to appeal to \emph{some} other model.
Here, instead of using the cycle-based models of Bloch \cite{bloch1986algebraic} and Suslin--Voevodsky \cite[Definition~3.1]{suslin2000bloch}, we use the universal characterization of motivic cohomology conjectured by Voevodsky \cite[Conjecture~10]{voevodsky-slice-filtration} and established (in the positive characteristic case we need) by Levine \cite[Theorem~10.5.1]{levine2008homotopy}.
In order to show that our proposed model satisfies Levine's universal property, we need to show that its spaces assemble into a $\P^1$-$\Omega$-spectrum.
In other words, we need some control over (certain) $\P^1$-loop spaces.
By the first author's conservativity result \cite[Corollary~4.15]{bachmann-moving2}, it will suffice to control the homology of these $\P^1$-loop spaces.
Our starting point is the first author's affine Grassmannian model for $\Gm$-loops on $\SL_n$ and the associated homology computation in characteristic zero \cite[Theorem~20 and Corollary~25]{bachmann-grassmann}; we extend the computation to positive characteristic below.

\subsection{Main result} \label{subsec:intro-main}
We now state our main result precisely.

Let $k$ be a perfect field of characteristic $p$.
The category $\NMon(\Spc(k))$ of normed monoids (see \S\ref{sec:normed-monoids}) carries a transported norm structure (\S\ref{subsec:Spc-fet}), inducing normed functors \[ \Spc(k)_* \xrightarrow{T} \NMon(\Spc(k)) \xrightarrow{\bar\Sigma^\infty} \SH(k), \] with $T$ left adjoint to the forgetful functor and $\bar\Sigma^\infty \circ T \wequi \Sigma^\infty$.
The motivic Eilenberg--Mac Lane spectrum $H\Z_{(p)} \in \NAlg(\SH(k))$ is periodizable, i.e., $\bigoplus_n \Sigma^{2n,n}H\Z_{(p)} \wequi H\Z_{(p)} \otimes_{\MGL} \mathrm{PMGL}$ lifts to $\NAlg(\SH(k)^\Z)$ \cite[Example~16.34]{bachmann-norms}.
The functor $\bar\Sigma^\infty$ being normed, its right adjoint $\bar\Omega^\infty$ is lax normed and consequently we obtain \[ \bar\Omega^\infty \Sigma^{2\bullet,\bullet} H\Z_{(p)} \in \NAlg(\NMon(\Spc(k))^\Z). \]
Let us write $\bar\Omega^\infty(H\Z_{(p)})^\bullet$ for the constituent normed monoids, so that $\bar\Omega^\infty(H\Z_{(p)})^n = K(\Z_{(p)}(n),2n)$.
For $X \in \NMon(\Spc(k))$ and $d \in \Z$ write $\NSym^{\Z\otimes}(X,d)$ for the free object in $\NAlg(\NMon(\Spc(k))^\Z)$ on $X$ placed in degree $d$.
Since the norm structure on $\NMon(\Spc(k))^\Z$ satisfies the distributivity law, the usual formula for free normed objects holds \cite[Example~13.19]{bachmann-norms}.
In particular, if $d=1$ we get \[ \NSym^{\Z\otimes}(X,1)^n = D_n^\otimes(X), \] where for $n \ge 0$, $D_n^\otimes$ denotes the motivic extended power \cite[Definition~5.6]{bachmann-colimits} for the normed category $\NMon(\Spc(k))$ and for $n<0$ we have $D_n^\otimes = *$.
Our main result is as follows.
\begin{theorem}[see Theorem \ref{thm:mainmain}] \label{thm:intro-main}
Let $k$ be a field of characteristic $p$.
The canonical map \[ \NSym^{\Z\otimes}(K(\Z_{(p)}(1),2),1)^\bullet \to \bar\Omega^\infty(H\Z_{(p)})^\bullet \] is an equivalence for $\bullet \ge 1$.\NB{what about $\bullet=0$?}
\end{theorem}

\begin{remark}
The above description of $K(\Z_{(p)}(n), 2n)$ in terms of norms is firmly a characteristic $p$ phenomenon.
For example, in classical topology, the analogous claim would be that $D_n(\Sigma^{2} H\Z_{(p)}) \wequi \Sigma^{2n} H\Z_{(p)}$ (where $D_n$ refers to extended powers in the category of connected commutative monoids, i.e., connected spectra), which is easily seen to be false.
Some further comments:
\begin{enumerate}
\item In the counterexample above, instead of taking extended powers in spectra, one might take them in $H\Z_{(p)}$-modules.
  The analogous statement remains wrong.
\item The topological situation does work out rationally: we indeed have $D_n(\Sigma^2 H\Q) \wequi \Sigma^{2n} H\Q$.
  The analogy seems to extend: in ways we don't fully understand yet, $p$-local homotopy theory in characteristic $p$ seems to behave a lot like rational homotopy theory in characteristic zero.
\item It seems likely that many other spaces and spectra of interest have descriptions purely in terms of norms in characteristic $p$.
  For example, $\MGL_{(p)} \in \SH(\F_p)$ is the free normed algebra on the pointed spectrum $\MGL_{(p)}(1) = \Th_{B\GL_1}(\gamma - \scr O)$ ($\wequi \Sigma^{-2,-1} \Sigma^\infty \P^\infty_{(p)}$).
  (There is a tautological map, which is an equivalence in homology by the formulas for free normed algebras at the characteristic, i.e., Proposition \ref{prop:NAlg-DAlg}.)
\end{enumerate}
\end{remark}

\subsection{Summary of the proof of Theorem \ref{thm:intro-main}} \label{subsec:into-summary}
Write \[ b_0 \in \pi_{2,1} \NSym^{\Z\otimes}(K(\Z_{(p)}(1),2),1)^1 \wequi H^{2,1}(\P^1, \Z_{(p)}) \] for the tautological class.
Multiplication with $b_0$ thus induces a map \begin{gather*} S^{2,1} \wedge \NSym^{\Z\otimes}(K(\Z_{(p)}(1),2),1)^n \to \\ \NSym^{\Z\otimes}(K(\Z_{(p)}(1),2),1)^1 \wedge \NSym^{\Z\otimes}(K(\Z_{(p)}(1),2),1)^n \to \NSym^{\Z\otimes}(K(\Z_{(p)}(1),2),1)^{n+1}. \end{gather*}
We shall argue that these maps make $\NSym^{\Z\otimes}(K(\Z_{(p)}(1),2),1)^\bullet$ into an $\Omega^{2,1}$-spectrum (in positive degrees), i.e., that the adjoint maps \begin{equation} \label{eq:Omega-spectrum} \NSym^{\Z\otimes}(K(\Z_{(p)}(1),2),1)^n \to \Omega^{2,1} \NSym^{\Z\otimes}(K(\Z_{(p)}(1),2),1)^{n+1} \in \Spc(k)_* \end{equation} are equivalences.
Let us denote by $K \in \SH(k)$ the corresponding spectrum, so that we obtain a map $K \to H\Z_{(p)}$.
It will suffice to show this is an equivalence.
By construction, \[ \NSym^{\Z\otimes}(K(\Z_{(p)}(1),2),1)^1 \wequi \Omega^\infty \Sigma^{2,1} K \to \Omega^\infty \Sigma^{2,1}H\Z_{(p)} \wequi K(\Z_{(p)}(1),2) \] is an equivalence.
Thus in order to prove Theorem \ref{thm:intro-main}, it is enough to show that $K, H\Z_{(p)} \in \Sigma^{-2,-1} \SH(k)^\veff$.
But we even have $K \in \SH(k)^\veff$ (essentially because $K^1$ is very $1$-effective and $K^n$ is more or less an $n$-th tensor power of $K^1$; see Proposition \ref{prop:Dn-tens-eff}), and the same holds for $H\Z_{(p)}$ by \cite[Lemma~13.7]{bachmann-norms}.

It thus remains to establish that the map in \eqref{eq:Omega-spectrum} is an equivalence for $n \ge 1$.
We will argue that it is enough to show that this map of motivic spaces induces an equivalence on motives with $\Z_{(p)}$-coefficients.
It is fairly straightforward to compute the homology of the left hand side.
Our main effort will thus be to determine the homology of the $\P^1$-loop space on the right hand side.
We will access the homology of $\Gm$-loops on the spheres $S^{2n-1,n}$ by studying a model via affine Grassmannians, and then propagate this information.

\begin{remark}
The way we ``propagate the information'' may be of independent interest.
If $M \in \NMon(\Spc(k))$ is a normed monoid, then $H\F_p[M] \in \Mod_{H\F_p}(k)$ is a \emph{normed Hopf algebra}.
Such an object has an underlying derived algebra (in the sense of \cite[Definition~4.2.22]{raksit2020hochschild}) encoding the multiplication; see \S\ref{subsec:normed-and-derived}.
We can thus extract the \emph{derived indecomposables} $QH\F_p[M]$ of the derived algebra underlying $H\F_p[M]$.
By dualizing, we similarly extract the \emph{derived primitives} $PH\F_p[M]$.
Now, in strong analogy with the classical situation (see, e.g., \cite[Lemma~6.3.6]{neisendorfer-algebraic}), given an appropriate normed monoid $M$, we construct a map \[ \Sigma^{2,1}QH\F_p[\Omega^{2,1}M] \to PH\F_p[M]. \]
We show that it is an equivalence in favorable cases (see Theorem \ref{thm:homology-of-loops}).
In other words, in these favorable cases, we can read off the derived indecomposables of the homology of a $\P^1$-loop space from the derived primitives of the homology of the original space.
Since taking derived indecomposables is a conservative functor in many cases, this information is often sufficient.
\end{remark}

\begin{remark}
Because they are buried in technicalities, let us spell out the following two key points which make the argument work.
\begin{enumerate}
\item We do not (try to) utilize any model of $H\Z_{(p)}$ involving cycles.
  In characteristic $p$, any such attempt seems to run into insurmountable (to us) problems specific to characteristic $p$ geometry, i.e., inseparable field extensions and Frobenius.
  Instead we use Levine's characterization of $H\Z_{(p)}$ via effectivity.
\item $p$-locally in characteristic $p$, (free) normed objects behave essentially like ``derived'' structures, e.g., any normed $H\Z_{(p)}$-algebra has an underlying derived $\Z_{(p)}$-algebra \cite{bachmann-linearity}.
  Thus, while étale norms in characteristic $p$ are deficient for describing power operations and loop space theory, they do provide us with additional expressivity to construct objects with no analog away from the characteristic.
\end{enumerate}
Expanding on (2), much of this work is inspired by the slogan that $p$-local motivic homotopy theory in characteristic $p$ behaves like motivic homotopy theory ``over cofiber $\tau$'' in characteristic zero.
For example, the formula $D_n^\otimes(K(\Z_{(p)}(1),2)) \wequi K(\Z_{(p)}(n),2n)$ is directly inspired by an analogous formula in terms of derived Hopf algebras, which we obtained for the homology modulo $\tau$ of motivic Eilenberg--Mac Lane spaces in characteristic $0$ \cite[Theorem~3.16]{bachmann-wilson}.
\end{remark}

\subsection{Invitation: the motivic dual Steenrod algebra at the characteristic}
It is fairly clear that Theorem \ref{thm:intro-main} implies all the applications mentioned in \S\ref{subsec:intro-app}.
Since it is illustrative of our methods, here we show how to compute the motivic dual Steenrod algebra at the characteristic.
\begin{corollary} \label{cor:intro-steenrod}
Let $k$ be a field of characteristic $p$.
We have \[ H\Z \otimes H\F_p \wequi H\F_p[\tau_1,\tau_2, \dots, \xi_1, \xi_2, \dots]/(\tau_i^2 = 0) \in \SH(k). \]
Here by the right hand side we mean a free $H\F_p$-module with basis $\tau_1^{\epsilon_1} \tau_2^{\epsilon_2} \dots \xi_1^{n_1} \xi_2^{n_2} \dots$, where $\epsilon_i \in \{0,1\}$ and $n_i \ge 0$.
We have $|\tau_i| = (2p^i-1,p^i-1)$ and $|\xi_i| = (2p^i-2,p^i-1)$.
\end{corollary}
\begin{proof}[Proof sketch.]
By base change, we may assume that $k$ is perfect.
Since $\bar\Sigma^\infty$ induces $\NMon(\Spc(k))[(\P^1)^{-1}] \wequi \SH(k)$ \cite[Theorem~1.1]{bachmann-MGM}, we see that \[ H\Z_{(p)} \wequi \NSym(\bar\Sigma^\infty K(\Z_{(p)}(1),2))/(b_0=1) \] and so \[ H\Z \otimes H\F_p \wequi H\Z_{(p)} \otimes H\F_p \wequi \NSym_{H\F_p}(H\F_p \otimes \bar\Sigma^\infty K(\Z_{(p)}(1),2))/(b_0=1). \]
(Here by setting $b_0=1$ we mean inverting $b_0$ and then passing to the degree zero part of the $\Z$-grading.)
Recall that for cellular $H\F_p$-modules over $k$, the monad $\NSym$ identifies with the derived algebra monad for the Chow $t$-structure (see Proposition \ref{prop:NAlg-DAlg}).
We shall show that $H\F_p \otimes \bar\Sigma^\infty K(\Z_{(p)}(1),2)$ is a free module on classes $a_i$, $b_j$ for $i \ge 1$, $j \ge 0$ with $|a_i| = (2p^i+1,p^i)$ and $|b_i| = (2p^i,p^i)$.
Since the $b_i$ have Chow degree $0$ whereas the $a_i$ have Chow degree $1$, this implies that the free normed algebra (i.e., the free derived algebra) on these classes is polynomial on the $b_i$ and exterior on the $a_i$, and the result follows by setting $\tau_i = a_i/b_0$, $\xi_i = b_i/b_0$.

We claim that for $X \in \NMon(\Spc(k))$ with cellular motive, \[ H\F_p \otimes \bar\Sigma^\infty X \wequi Q(H\F_p \otimes \Sigma^\infty_+ X), \] where $Q$ denotes the derived indecomposables (of an animated ring).
To see this, resolve $X$ by iterated application of the free normed monoid functor, observe that the claim holds when $X$ is a free normed monoid, and that derived indecomposables commute with colimits.
We thus need to compute the derived indecomposables of $H\F_p[K(\Z_{(p)}(1),2)]$.
Since $K(\Z(1),2) \wequi \P^\infty$, we see that \[ (H\F_p)_{**} K(\Z_{(p)}(1),2) \wequi \Gamma_{\F_p}[b_0] \wequi \F_p[b_0, b_1, \dots]/(b_i^p = 0). \]
This shows that as a derived algebra, $H\F_p[K(\Z_{(p)}(1),2)]$ is equivalent to the tensor product of elementary algebras $E(b_i)$, each given by a pushout
\begin{equation*}
\begin{CD}
\F_p[c_{i+1}] @>{c_{i+1} \mapsto b_i^p}>> \F_p[b_i] \\
@VVV  @VVV \\
\F_p @>>> E(b_i).
\end{CD}
\end{equation*}
From this we read off that the derived indecomposables of $E(b_i)$ are free on a class in Chow degree $0$ corresponding to $b_i$, and a class in Chow degree one coming from $c_{i+1}$, i.e., corresponding to $a_{i+1}$.
This concludes the proof.
\end{proof}

\subsection{Organization}
Let us briefly record the organization of this article.
Each section has its own introduction with more details.
We begin in \S\ref{sec:prelim} with some preliminary recollections and observations.
The next sections deal with the theory of norms.
First in \S\ref{sec:norm-hopf-1}, we introduce normed Hopf algebras and develop some basic properties, mainly with an eye towards our application.
Second in \S\ref{sec:normed-monoids}, we consider the special case of normed monoids in $\Spc(k)$.
Third in \S\ref{sec:normed-hopf-charp}, we establish some properties of normed Hopf algebras specific to characteristic $p$---these are essentially extensions of the ideas from \cite{bachmann-linearity}.
At this point we have established all necessary prerequisites, and in \S\ref{sec:main} we can put everything together and prove our main results.
We first compute the homology of $\Omega_\Gm JS^{2n,n}$ by analyzing motivic analogs of Toda's fiber sequences using techniques specific to characteristic $p$.
Then we show how to control the homology of the $\P^1$-loops of cellular $2$-effective normed monoids by combining the previous computation with the normed Hopf algebra theory developed in the previous sections.
Using this, we establish our new model for $K(\Z_{(p)}(n),2n)$ as outlined above.
Finally in \S\ref{sec:applications} we derive some applications.

We also include an appendix, in which we record the definition and important properties of the \emph{derived Ravenel--Wilson Hopf ring}.

\subsection{Notation}\label{subsec:notation}
Throughout, $p$ denotes a prime.  Unless otherwise indicated, $k$ denotes a perfect field, $S$ a base scheme, and $R$ a coefficient ring.
In the parts of the paper concerned with the characteristic, $\operatorname{char}(k)=p$, and $R$ will usually be $\F_p$, $\Z/p^\nu$, or $\Z_{(p)}$.
For a homogeneous class $x$, we write $|x|=(a,b)$ for its motivic bidegree and we call $a - 2b$ its \emph{Chow degree}.
Classical Eilenberg--Mac Lane spectra are written without an $H$; thus, for example, $\Mod_{\F_p}$ denotes ordinary $\F_p$-module spectra, whereas $H\F_p$ denotes the motivic Eilenberg--Mac Lane spectrum.
We freely use the language of $\infty$-categories, as set out in \cite{lurie-htt,lurie-ha}.

\medskip
\noindent\textbf{Category theory notation.}
\renewcommand{\arraystretch}{1.15}
\begin{longtable}{@{}>{\raggedright\arraybackslash}p{0.22\textwidth}>{\raggedright\arraybackslash}p{0.50\textwidth}>{\raggedright\arraybackslash}p{0.20\textwidth}@{}}
\textbf{Notation} & \textbf{Meaning} & \textbf{Reference} \\
\hline
\endfirsthead
\textbf{Notation} & \textbf{Meaning} & \textbf{Reference} \\
\hline
\endhead
$\Map(X,Y)$ & mapping space between objects in an $\infty$-category \\
$\map(X,Y)$ & mapping spectrum in a stable $\infty$-category \\
$[X,Y]$ & $\pi_0\Map(X,Y)$ \\
$\imap(X,Y)$, $\iMap(X,Y)$ & internal mapping objects in symmetric monoidal $\infty$-categories \\
$\Map_*(X,Y), \iMap_*(X,Y)$ & notation to emphasize working in pointed objects \\
$D(\ph)$ & (weak) dual in a symmetric monoidal category, $D(\ph) = \imap(\ph, \1)$ \\
\end{longtable}

\medskip
\noindent\textbf{Ambient motivic notation.}
\renewcommand{\arraystretch}{1.15}
\begin{longtable}{@{}>{\raggedright\arraybackslash}p{0.22\textwidth}>{\raggedright\arraybackslash}p{0.50\textwidth}>{\raggedright\arraybackslash}p{0.20\textwidth}@{}}
\textbf{Notation} & \textbf{Meaning} & \textbf{Reference} \\
\hline
\endfirsthead
\textbf{Notation} & \textbf{Meaning} & \textbf{Reference} \\
\hline
\endhead
$k,S,p,R$ & A perfect base field, a general base scheme, a prime, and a coefficient ring, respectively. &  \\
$\Spc(S),\Spc(S)_*$; $\SH(S)$ & Motivic spaces, pointed motivic spaces, and motivic spectra over $S$. & \S\ref{subsec:mot-spc}; \cite[\S2.2, \S4.1]{bachmann-norms} \\
$\SH^{S^1}(k)$ & The $S^1$-stable motivic category $\Sp(\Spc(k))$. & \\
$S^{a,b}$; $\Sigma^{a,b}$; $\Omega^{a,b}$ & The motivic sphere of bidegree $(a,b)$, suspension by it, and the corresponding loop functor.  Thus $\Omega^{1,1}=\Omega_{\Gm}$ and $\Omega^{2,1}=\Omega_{\P^1}$. &  \\
$\ul\pi_i(X)_j$, $\pi_{**}E$ & Motivic homotopy sheaves and bigraded homotopy groups.  For a homogeneous class, $|x|=(a,b)$ records its bidegree. & \\
$E_{\ge n},E_{\le n}$ & Homotopy $t$-structure truncation notation. & \cite[\S\S4.3, 5.2]{morel-trieste} \\
$X_{(p)}$ & $p$-localization of a motivic space (or related object). & \S\ref{subsec:p-loc} \\
$HR$ & The motivic Eilenberg--Mac Lane spectrum associated with $R$; in particular $H\Z$, $H\F_p$, and $H\Z_{(p)}$. & \cite{spitzweck2012commutative} \\
$\Mod_A(S)$ & The $\infty$-category of $A$-modules in $\SH(S)$. &  \\
$A[X]$ & $A\otimes\Sigma^\infty_+X$ for a motivic space $X$. & \\
$\overline{A[X]}$ & $A\otimes\Sigma^\infty X$ for pointed $X$. & \\
$JX,J_mX$ & The James construction and its $m$-th filtration stage. & \S\ref{subsub:toda}; \cite[Construction~5.1]{devalapurkar2019james} \\
$\CB^\bullet(L,A,R)$; $\CB(L,A,R)$ & The cobar cosimplicial object associated with a left $A$-comodule $L$ and a right $A$-comodule $R$, and its totalization. & \S\ref{subsub:eilenberg-moore}; \cite[Theorem~4.4.2.8]{lurie-ha} \\
\end{longtable}

\medskip
\noindent\textbf{Cellularity and effectivity.}
\begin{longtable}{@{}>{\raggedright\arraybackslash}p{0.22\textwidth}>{\raggedright\arraybackslash}p{0.50\textwidth}>{\raggedright\arraybackslash}p{0.20\textwidth}@{}}
\textbf{Notation} & \textbf{Meaning} & \textbf{Reference} \\
\hline
\endfirsthead
\textbf{Notation} & \textbf{Meaning} & \textbf{Reference} \\
\hline
\endhead
$\scr C^{\cell i\eff}_{\ge n}$, $\scr C^{\cell i\eff}$, $\scr C^{\cell}_{\ge n}$; $\scr C^{\cell}$ & Various cellular subcategories. & Definition~\ref{def:cellular-categories} \\
$\Mod_A(S)^{\eff}$; $\Mod_A(S)^{\eff\perp}$ & The subcategory of effective $A$-modules and its right orthogonal. & \S\ref{subsec:duality} \\
$\Mod_A(S)^\fet$, $\Mod_A(S)^{\fet\pm}$ & Various subcategories generated by finite étale $S$-schemes. & \S\ref{subsec:duality} \\
``quickly converging'' & An $\N$-indexed diagram whose successive cofibers eventually lie arbitrarily far in the effective, respectively coeffective, direction; the same terminology is used for sums, products, and limits. & Definition~\ref{def:quick-conv} \\
$\Mod_R^\Z$ & Evenly $\Z$-graded ordinary $R$-module spectra. & \S\ref{subsub:charp} \\
\end{longtable}

\medskip
\noindent\textbf{Normed algebraic structures.}
\begin{longtable}{@{}>{\raggedright\arraybackslash}p{0.22\textwidth}>{\raggedright\arraybackslash}p{0.50\textwidth}>{\raggedright\arraybackslash}p{0.20\textwidth}@{}}
\textbf{Notation} & \textbf{Meaning} & \textbf{Reference} \\
\hline
\endfirsthead
\textbf{Notation} & \textbf{Meaning} & \textbf{Reference} \\
\hline
\endhead
$\NAlg(\scr C)$, $\NMon(\scr C)$ & Normed algebras in a normed category $\scr C$, called normed monoids in the cartesian case. & \S\ref{subsec:normed-cat}; \cite{bachmann-norms} \\
$\cNAlg(\scr C)$, $\cNMon(\scr C)$ & Co-normed coalgebras, by definition $\NAlg(\scr C^\op)^\op$; co-normed comonoids in the cocartesian case. & \S\ref{subsec:normed-cat} \\
$\bNAlg(\scr C)$ & Normed bialgebras, $\NMon(\cNAlg(\scr C))$. & \S\ref{sec:norm-hopf-1} \\
$\NHopf(\scr C)$ & Normed Hopf algebras, $\bNAlg(\scr C)^{\gp}$. & \S\ref{sec:norm-hopf-1} \\
$T$ & The free normed monoid functor $\Spc(k)_*\to\NMon(\Spc(k))$. & \S\ref{subsub:free-normed} \\
$\NSym$ & The free normed algebra monad. & \S\ref{subsec:normed-cat} \\
$D_n^{\scr C}$; $D_n$ & The $n$-th normed extended power in $\scr C$; the superscript is omitted when the normed category is clear. & Equation~\eqref{eq:formula-for-dn}; \cite[Definition~5.6]{bachmann-colimits} \\
$D_n^\wedge$ & The extended power for the normed category of pointed motivic spaces. & \S\ref{subsec:Spc-fet} \\
$D_n^\otimes$ & The extended power for the transported norm structure on normed monoids. & \S\ref{subsec:Spc-fet} \\
$\hotimes$ & The symmetric monoidal product on the relevant normed/derived Hopf category. & Proposition~\ref{prop:NHopf-stability} \\
$\Free$; $\Cofree$ & The free and cofree normed Hopf-algebra constructions used in the characteristic-$p$ theory. & Proposition~\ref{prop:NHopf-adj} \\
\end{longtable}

\medskip
\noindent\textbf{Derived algebra and Hopf algebra notation.}
\begin{longtable}{@{}>{\raggedright\arraybackslash}p{0.22\textwidth}>{\raggedright\arraybackslash}p{0.50\textwidth}>{\raggedright\arraybackslash}p{0.20\textwidth}@{}}
\textbf{Notation} & \textbf{Meaning} & \textbf{Reference} \\
\hline
\endfirsthead
\textbf{Notation} & \textbf{Meaning} & \textbf{Reference} \\
\hline
\endhead
$\DAlg(\scr C)$ & Derived algebras in $\scr C$. & \S\ref{subsec:normed-and-derived}; \cite[\S4]{raksit2020hochschild} \\
$\LSym$; $\LSym_n$ & The free derived algebra monad and its $n$-th homogeneous piece. & \S\ref{subsec:normed-and-derived}; \cite[\S4]{raksit2020hochschild} \\
$D_n^\hotimes$ & The extended power for the transported norm structure on normed Hopf algebras. & \S\ref{sec:normed-hopf-charp} \\
$\LSym_n^\hotimes$ & The corresponding derived symmetric-power operation in the derived Hopf category. & \\
$C, B$ & co/bar construction for (co)augmented (co)algebras & \\
$\A^1(n)$ & The object $\Cofree\,\Sigma^{2n,n}HR$. Not to be confused with the scheme $\A^1$. & Definition~\ref{def:dhopfbp} \\
$\Hopf_R^{bp}$ & Bipolynomial Hopf algebras. & Definition~\ref{def:dhopfbp} \\
$\dHopf_R^{bp}$ & Derived bipolynomial Hopf algebras. & Definition~\ref{def:dhopfbp} \\
$\dHopf_{\F_p}^{cn}$ & The derived category of connected, evenly graded $\F_p$-Hopf algebras. & \S\ref{subsec:derived-hopf-alg} \\
$W(n)$ & The $p$-typical Witt Hopf algebra whose initial generator has degree $n$. & \S\ref{subsub:witt} \\
$F,V$ & Frobenius and Verschiebung on the Witt Hopf algebras. & \S\ref{subsub:witt} \\
$w_i$ & The Witt polynomial primitives. & \S\ref{subsub:witt} \\
$W(n)_{\le d}$ & $\fib(V^{d+1}:W(n)\to W(p^{d+1}n))$. & \S\ref{subsub:witt} \\
$W(n)^{\le d}$ & $\cof(F^{d+1}:W(p^{d+1}n)\to W(n))$. & \S\ref{subsub:witt} \\
$\Ga(n)$ & $W(n)^{\le0}$. & \S\ref{subsub:witt} \\
$Q$ & Derived indecomposables. & \S\ref{subsec:more-dHopf}; Remark~\ref{rmk:P-Q-explanation} \\
$P$ & Derived primitives. & \S\ref{subsec:more-dHopf}; Remark~\ref{rmk:P-Q-explanation} \\
$\hLoops^{2,1}$ & The ``algebraic loop'' functor: the relevant right adjoint to tensoring with the degree-$(2,1)$ generator. & Proposition~\ref{prop:alg-loops-defined}; \S\ref{subsec:P1-loops} \\
$\dHFp[\ph]$ & Variant of $H\F_p[\ph]$ landing in derived Hopf algebras. & \S\ref{subsec:P1-loops} \\
\end{longtable}
\renewcommand{\arraystretch}{1}

\subsection{Acknowledgements} The problems resolved in this article
have occupied us (and other mathematicians working in motivic homotopy
theory) for many years.  It would be infeasible to list everyone who
has influenced our work on them.  Nonetheless we want to express our
special gratitude to and acknowledge the deep influence of: Aravind
Asok, Robert Burklund, Elden Elmanto, Jeremy Hahn, Marc Hoyois, Jacob
Lurie, Marc Levine, and Markus Spitzweck.  These results were announced by
the first author in a lecture at
SLMath\footnote{\url{https://www.slmath.org/workshops/1183/schedules/40286}}.  We are
posting this manuscript so that it will be available for a followup
workshop.  We wish to thank SLMath for support during the
preparation of this manuscript and for hosting the workshops.

\section{Preliminaries} \label{sec:prelim}
We record some preliminary observations, also using the opportunity to fix some terminology that will be used throughout the article.
In \S\ref{subsec:mot-spc} we recall basic facts about motivic spaces.
Next in \S\ref{subsec:cellularity} we discuss notions of cellularity in stable categories.
Then in \S\ref{subsub:eilenberg-moore} we construct a motivic analog of the Eilenberg--Moore spectral sequence.
We immediately utilize this in \S\ref{subsub:toda} to construct motivic versions of Toda's classical fibrations, linking certain $p$-local loop spaces of spheres.
Finally in \S\ref{sec:grass} we revisit the affine Grassmannian model for $\Omega_\Gm$ and extend the computation of its homology to positive characteristic.

Throughout, let $k$ be a perfect field.

\subsection{Motivic spaces} \label{subsec:mot-spc}
We denote the $\infty$-category of motivic spaces over $k$ \cite[\S2.2]{bachmann-norms} by $\Spc(k)$.

\subsubsection{$p$-localization} \label{subsec:p-loc}
We write $X \mapsto X_{(p)}$ for the $p$-localization functor on $\Spc(k)$ \cite[\S4]{AFHLocalization}.
We only apply this to nilpotent spaces, in which case we have \cite[Theorem~4.3.9(1)]{AFHLocalization} \[ \ul\pi_*(X_{(p)}) \wequi \ul\pi_*(X)_{(p)}. \]
Since $\ul\pi_i \Omega^{1,1} X \wequi \ul \pi_i(X)_{-1}$ \cite[Theorem~3.1.19(2)]{AFHLocalization} and contraction is an exact operation \cite[Proposition~3.1.18(3)]{AFHLocalization}, this in particular implies that (for nilpotent spaces) \[ \Omega^{1,1}(X_{(p)}) \wequi (\Omega^{1,1} X)_{(p)}. \]

\subsubsection{Effectivity}
\begin{definition} \label{def:n-effective-space}
We call $X \in \Spc(k)$ \emph{$n$-effective} if it is weakly $S^{n,n}$-cellular in the sense of \cite[Definition 3.1.3]{freudenthal} \emph{and connected}.
We call $X$ \emph{very $n$-effective} if it is weakly $S^{2n,n}$-cellular and connected.
\end{definition}
We find that the additional connectivity assumption in our definitions over those of \cite{freudenthal} simplifies some of our arguments.

\begin{remark}\label{rmk:effective-operations}
Note that among connected spaces, $n$-effective and very $n$-effective ones are stable under cofiber-extension and colimits (being defined as those which become contractible in some Bousfield localization).
For many more properties, see \cite[\S3 and \S4]{freudenthal}.
\end{remark}

\begin{example} \label{ex:Poo-effective}
$\P^n$ and $\P^\infty$ are very $1$-effective.
This follows via Remark \ref{rmk:effective-operations} from the cofiber sequences $\P^{n-1} \to \P^n \to S^{2n,n}$.
\end{example}

\subsubsection{Limit-determinedness}
\begin{definition} \label{def:limit-determined}
We call $X \in \Spc(k)$ \emph{limit-determined} if it lies in the closure under limits of $\Omega^\infty(\SH(k)) \subset \Spc(k)$.
\end{definition}

\begin{example} \label{ex:resolvable-spaces}
If $X \in \Spc(k)$ is $2$-effective and nilpotent, it is limit-determined by \cite[Amplification~4.1.12]{freudenthal}.
\end{example}

\begin{example} \label{ex:loops-limit-determined}
If $X \in \Spc(k)_*$ is limit-determined then so is $\Omega^{p,q} X$ (and more generally so is $\ul\Map_*(Y,X)$ for any $Y \in \Spc(k)_*$).
\end{example}

\begin{lemma}\label{lemm:limit-determined-M-test}
Let $f: X \to Y \in \Spc(k)$ be a map between limit-determined spaces.
Then $f$ is an equivalence in the following cases.
\begin{enumerate}
\item $\Sigma^\infty_+ f$ is an equivalence.
\item $cd_2(k) < \infty$ and $H\Z[f]$ is an equivalence.
\item -1 is a sum of squares in $k$, $p$ is a prime, if $p=2$, additionally $cd_2(k) < \infty$, 
$X, Y$ nilpotent and $p$-local, and $H\Z_{(p)}[f]$ is an equivalence.
\end{enumerate}
\end{lemma}
\begin{proof}
(1) Let $\scr C \subset \Spc(k)$ be the class of spaces $T$ such that $\Map(f,T)$ is an equivalence.
By Yoneda, it suffices to show that $X, Y \in \scr C$.
Since $\scr C$ is closed under limits, it contains all limit-determined spaces, whence the result.

(2) By (1) it will suffice to show that, under the stated assumption on $k$, any map $f: E \to F \in \SH(k)^\veff$ with $f \otimes H\Z$ an equivalence is an equivalence.
By \cite[Theorem~7.3.4 and \S5.2]{mantovani2018localizations}, the totalization of the cosimplicial object $f \otimes H\Z^{\otimes \bullet +1}$ is the $\eta$-completion of $f$.
The claim thus follows from $\eta$-completeness of $E,F$ \cite[Lemma 6.43]{bachmann-resolution}.

(3) Since $X, Y$ are $p$-local and nilpotent, $\Sigma^\infty X, \Sigma^\infty Y$ are $p$-local \cite[Proposition~2.2.16]{freudenthal}.
It follows that $H\Z_{(p)}[f]$ is an equivalence if and only if $H\Z[f]$ is.
Arguing as in (2), it thus suffices to show that any object in $\SH(k)^\veff_{(p)}$ is $\eta$-complete.
If $p=2$, we have $cd_2(k) < \infty$ by assumption, and the claim reduces to (2).
If $p \ne 2$ then actually $\eta = 0$ on any object of $\SH(k)^\veff_{(p)}$, $k$ being unorderable.
(Recall that $\SH(k)[1/2] \wequi \SH(k)[1/2]^+ \times \SH(k)[1/2]^-$ with $\eta=0$ on the $+$ part, and the $-$ part vanishing if $k$ is unorderable \cite[\S2.7.3 and Corollary~2.19]{bachmann-eta}.)
\end{proof}

\subsection{Cellularity} \label{subsec:cellularity}
We will eventually be able to control the $\Gm$-loops of certain spaces made out of spheres.
For this and other reasons, we will have to contend with various ``cellular subcategories'', i.e., subcategories generated under colimits (and perhaps other operations) by (certain) spheres.
We observe here some very general properties of this kind of construction.

\begin{definition} \label{def:cellular-categories}
Let $\scr C$ be a stable, presentable $\infty$-category under $\Spc(k)_*$.
For $i \ge 0, n \in \Z$ we define $\scr C^{\cell i \eff}_{\ge n} \subset \scr C$ as the subcategory generated under colimits and extensions by $\Sigma^n S^{j,j}$ for $j \ge i$.
We write $\scr C^{\cell i \eff}$ for the subcategory generated under colimits by $\scr C^{\cell i \eff}_{\ge n}$ for all $n \in \Z$.
If $\scr C$ is under $\SH(k)$, so that negative spheres make sense, we also allow $i < 0$ in the above definitions.
We then also write $\scr C^{\cell}_{\ge n}$ for the subcategory generated under colimits and extensions by $S^{n+i,i}$ for $i \in \Z$, and $\scr C^{\cell}$ for the subcategory generated under colimits by all the $S^{i,j}$ for $i, j \in \Z$.

We call the objects in these categories \emph{($n$-connective, $i$-effective) cellular}.
\end{definition}

\begin{remark} \label{rmk:t-str}
These are nonnegative parts of $t$-structures, by \cite[Proposition~1.4.4.11]{lurie-ha}.
\end{remark}

\begin{lemma} \label{lemm:cell-crit}
Let $\scr C$ be a presentable, stable $\infty$-category.
Let $\scr S \subset \scr C$ be a (small) family of compact objects.
Write $\scr C_{\ge \scr S+n}$ for the subcategory generated under colimits and extensions by $\Sigma^n \scr S$.

\begin{enumerate}
\item Let $E \in \scr C$.
  Then $E \in \scr C_{\ge \scr S+n}$ if and only if there exists a (possibly transfinite) family of cofiber sequences $S_i \to E_i \to E_{i+1}$ with $S_i \in \scr C_{\ge \scr S+n}$, such that $E_0 = E$ and $\colim_i E_i = 0$.

\item Suppose furthermore that for each $S, T \in \scr S$ and $i > 0$, $[S, \Sigma^i T] = 0$.
  Then given $E \in \scr C_{\ge \scr S + n}$ there exists a sequence of $S_i, E_i$ as above with $S_i$ a sum of terms of the form $\Sigma^{i+n} S$ for $S \in \scr S$, and we have $E_i \in \scr C_{\ge \scr S+n+i}$.

\item In the situation of (2), if $E$ lies in the localizing subcategory of $\scr C$ generated by $\scr S$, then $E \in \scr C_{\ge \scr S + n}$ if and only if $\map(S, E)$ is $n$-connective for every $S \in \scr S$.
\end{enumerate}
\end{lemma}
\begin{proof}
As in Remark \ref{rmk:t-str}, $\scr C_{\ge \scr S+n}$ is the nonnegative part of a $t$-structure.

(1)
Necessity is clear (take $S_0 = E$, $S_i = 0$ for $i>0$), so we show sufficiency.
Write $\tau$ for the left adjoint truncation functor in the $t$-structure corresponding to $\scr C_{\ge \scr S+n}$.
We must show that $\tau E = 0$.
The cofiber sequences show that $\tau E_i \wequi \tau E_{i+1}$ and so the filtered system $\{\tau E_i\}_i$ consists of equivalences.
It follows that $0 = \tau \colim_i E_i \wequi \colim_i \tau E_i \wequi \tau E$, as needed.

(2)
Let $E \in \scr C_{\ge \scr S + n}$.
Set \[ S_0 = \bigoplus_{S \in \scr S, \Sigma^n S \to E} \Sigma^n S, \] so that $S_0 \to E$ is surjective on $[\Sigma^n S, \ph]$ for every $S \in \scr S$.
Then $E_1 = \cof(S_0 \to E)$ has the property that $[\Sigma^n S, E_1] = 0$ for all $S \in \scr S$.
Then define $S_1$ by the analogous sum of maps $\Sigma^{n+1}S \to E_1$ to obtain $E_2$, and so on.
Repeating this process we obtain a colimiting object $E'$ with $[\Sigma^i S, E'] = 0$ for all $i \ge n$.
But also $[\Sigma^i S, F] = 0$ for any $i < n$ and any object $F \in \scr C_{\ge \scr S +n}$ (i.e., $\map(\Sigma^n S, F)$ is connective, a property which is stable under colimits and extensions and holds for the generators $\Sigma^n T$, $T \in \scr S$).
Thus $[\Sigma^i S, E'] = 0$ for all $S \in \scr S, i \in \Z$ and thus $E'=0$ (since $E$ lies in localizing subcategory generated by $\scr S$; now use \cite[Lemma~2.2.1]{schwede-shipley2003modules}).
This shows that $E_i$ can be killed by attaching cells (of connectivity $\ge n+i$), and so, using (1), we deduce that in fact $E_i \in \scr C_{\ge \scr S + n+i}$.

(3) Under the stated assumption, we can also run the construction of (2) and thus appeal to (1).
\end{proof}

We now specialize to our most common case.
\begin{corollary} \label{cor:module-cellular-criterion}
Let $A \in \SH(k)_{\ge 0}$.
\begin{enumerate}
\item Let $E \in \Mod_A(k)$.
  Then $E \in \Mod_A(k)^{\cell j \eff}_{\ge n}$ if and only if there exists a (possibly transfinite) family of cofiber sequences $S_i \to E_i \to E_{i+1}$ with $S_i$ a sum of objects of the form $\Sigma^{c+w,w} A$ where $w \ge j$, $c \ge n$, such that $E_0 = E$ and $\colim_i E_i = 0$.
\item Moreover given $E \in \Mod_A(k)^{\cell j \eff}_{\ge n}$ one can ensure that $E_i \in \Mod_A(k)^{\cell j \eff}_{\ge n+i}$.
\item Also, $\Mod_A(k)^{\cell j \eff}_{\ge n} = \Mod_A(k)^{\cell j \eff} \cap \Mod_A(k)_{\ge n}$.
\item If furthermore $A \wequi s_0(A)$, then $\Mod_A(k)^{\cell j \eff} \wequi \Mod_A(k)^{\cell} \cap \Mod_A(k)^{\eff}(j)$.
\end{enumerate}
\end{corollary}
\begin{proof}
For (1) and (2), apply Lemma \ref{lemm:cell-crit} to the compact free $A$-modules $\Sigma^{w,w}A$, $w\ge j$.
For $i>0$, the free/forgetful adjunction identifies
\[
[\Sigma^{w,w}A,\Sigma^{i+v,v}A]_A
\wequi \ul\pi_{-i}(A)_{v-w}(k)=0,
\]
where the vanishing follows from $A\in\SH(k)_{\ge0}$.
For (3), one containment is clear, and the other is immediate from Lemma \ref{lemm:cell-crit}.

(4) In this case, the $A$-effective cover functor $f_n^A: \Mod_A(k) \to \Sigma^{2n,n}\Mod_A(k)^\eff$ preserves cellular objects (check on generators), and so if $E \in \Mod_A(k)^{\cell} \cap \Sigma^{2j,j}\Mod_A(k)^{\eff}$ then $E \wequi f_j^A E \in \Mod_A(k)^{\cell j \eff}$.
The other containment is clear.
\end{proof}

\subsubsection{Characteristic $p$ case} \label{subsub:charp}
We now discuss what is essentially an extended example.
Let $k$ be a perfect field of characteristic $p$.

Let $R = \Z/p^\nu$.
By our assumption on $k$, we have $\pi_{**}H\Z/p^\nu \wequi \Z/p^\nu$ \cite[Theorem~1.1]{geisser2000k}, and thus by \cite[Theorem~5.1.1]{schwede-shipley2003modules} we find \[ \Mod_{HR}(k)^\cell \wequi \Mod_{R}^\Z. \]
Here $\Mod_{R}$ denotes the category of (ordinary) $R$-module spectra, and $(\ph)^\Z$ denotes the category of (evenly) $\Z$-graded objects.
Our equivalence sends $\Sigma^{2n,n}HR$ to $R$ placed in topological degree $0$ and external degree $n$, i.e., the standard $t$-structure on $\Mod_{R}^\Z$ corresponds to the Chow $t$-structure on $\Mod_{HR}(k)^\cell$.
We will use this identification (between cellular modules over motivic cohomology and ordinary graded modules) frequently in the rest of the article.

\begin{remark}
Write $R\{n\} \in \Mod_R^\Z$ for $R$ placed in homotopical degree $0$ and external degree $n$.
Our choice above, of identifying $R\{n\}$ with $\Sigma^{2n,n}HR$ instead of (say) $\Sigma^{n,n}HR$ is, at this point, somewhat arbitrary (of course, if we want a symmetric monoidal equivalence, after making a choice for $n=1$ all others have been fixed).
It will eventually turn out to be the most convenient one (see, e.g., Proposition \ref{prop:NAlg-DAlg}), largely because on the motivic side, the Thom twist $\Sigma^{2n,n}$ plays a distinguished role.
\end{remark}

\begin{remark} \label{rmk:finite-field-motivic-cohomology}
Let $k$ be a finite field of characteristic $p$.
Its rational motivic cohomology is concentrated in bidegree $(0,0)$ \cite[Corollary~IV.1.13 and Example~VI.4.9]{weibel-k-book}.
Together with the preceding mod-$p$ calculation, this implies $H^{a,b}(k,\Z_{(p)})=0$ for $(a,b)\ne(0,0)$.
Hence we also obtain
\[
\Mod_{H\Z_{(p)}}(k)^\cell
\wequi \Mod_{\Z_{(p)}}^\Z.
\]
\end{remark}

\subsection{Eilenberg--Moore} \label{subsub:eilenberg-moore}
We want to be able to access the homology of a fiber of a map between appropriate motivic spaces.
For this, we will utilize a motivic generalization of the classical Eilenberg--Moore spectral sequence \cite{eilenberg-moore1962limits}.

We begin with some formalities.
Given a symmetric monoidal $\infty$-category $\scr C$, by a \emph{coalgebra} we mean an object of $\cAlg(\scr C) = \Alg(\scr C^\op)^\op$.
By a \emph{left comodule} for $A$ we mean an object of $(\mathrm{RMod}_{A^\op})^\op$.
We define right comodules similarly.
Given a coalgebra $A$, a left comodule $L$ and a right comodule $R$, we can form the cosimplicial diagram (by applying \cite[Construction~4.4.2.7]{lurie-ha} in $\scr C^\op$)
\[
\CB^\bullet(L,A,R)
=
\left[
L\otimes R
\;\substack{
  \xrightarrow{\hspace*{1.5em}}\\[-0.7ex]
  \xleftarrow{\hspace*{1.5em}}\\[-0.7ex]
  \xrightarrow{\hspace*{1.5em}}
}\;
L\otimes A\otimes R
\;\substack{
  \xrightarrow{\hspace*{1.5em}}\\[-0.7ex]
  \xleftarrow{\hspace*{1.5em}}\\[-0.7ex]
  \xrightarrow{\hspace*{1.5em}}\\[-0.7ex]
  \xleftarrow{\hspace*{1.5em}}\\[-0.7ex]
  \xrightarrow{\hspace*{1.5em}}
}\;
L\otimes A^{\otimes 2}\otimes R
\;\cdots
\right];
\]
We write $\CB(L,A,R)$ for its totalization.

\begin{remark}
If the tensor product in $\scr C$ preserved totalizations, taking the limit in $\CB^\bullet(L, A, R)$ would define a well-behaved cotensor product of bicomodules (e.g. satisfying associativity) \cite[Theorem~4.4.2.8]{lurie-ha}.
Since this preservation of totalizations is rather atypical, extra caution must be exercised.
\end{remark}

\begin{example}
If $\scr C$ is cartesian symmetric monoidal, then $\cAlg(\scr C) \to \scr C$ is an equivalence \cite[Corollary~2.4.3.10]{lurie-ha}, i.e., every object is automatically a coalgebra.
Any map $X \to Y \in \scr C$ thus makes $X$ into a comodule over $Y$, and given a span $X \to Y \leftarrow Z$, we have $\CB(X, Y, Z) \wequi X \times_Y Z$ (this is the standard way of expressing a pullback as a cosifted limit of finite products \cite[Theorem~4.4.2.8]{lurie-ha}).
\end{example}

\begin{proposition} \label{prop:EM}
Let $\scr X$ be an $\infty$-topos, $\scr D$ a Postnikov complete, presentable, prestable, symmetric monoidal $\infty$-category and $F: \scr X \to \scr D$ a cocontinuous, symmetric monoidal functor.
Given a span $X \to Y \leftarrow Z \in \scr X$ with $Y$ simply connected\NB{or connected and acts nilpotently on the fibers?} and $X, Z$ connected, we have \[ F(X \times_Y Z) \wequi \CB(FX, FY, FZ). \]
\end{proposition}
\begin{proof}
Set $G = \Omega Y$ acting on $X' = \fib(X \to Y)$.
Write $\scr X$ as an accessible localization \[ L: \scr P \adj \scr X: i \] of a presheaf topos $\scr P$ \cite[Definition~6.1.0.4]{lurie-htt}.
Set $\tilde G = (iG)_{\ge 1}$ which acts on $\tilde X' := (iX')_{\ge 1}$.
Put $\tilde X = \tilde X'_{h\tilde G} \in \scr P$.
Construct $\tilde Z$ similarly, and set $\tilde Y = B\tilde G$.
We obtain a span $\tilde X \to \tilde Y \leftarrow \tilde Z$.
Applying $L$ we recover the original span and so, replacing $\scr X$ by $\scr P$, we may assume that $\scr X$ is a presheaf $\infty$-topos.

We now continue in this setting.
Since $\scr D$ is prestable, we may factor $F$ as \[ \scr X \xrightarrow{\Sigma^\infty_+} \Sp(\scr X)_{\ge 0} \xrightarrow{F'} \scr D. \]
Let us put $\scr D_0 = \Sp(\scr D)$ and write $F_0: \Sp(\scr X) \to \scr D_0$ for the functor induced by $F'$.
It will suffice to show that both $\Sigma^\infty_+$ and $F_0$ commute with the totalization of $\CB^\bullet(X, Y, Z)$.
For $\Sigma^\infty_+$, since $\scr X$ is a presheaf topos, we need only check the claim sectionwise, i.e., reduce to the case $\scr X = \Spc$.
Since $Y$ is simply connected, the integral homology spectral sequence of $\CB^\bullet(X,Y,Z)$ converges strongly \cite[Theorem~4.1]{bousfield1987homology}.
Hence $\Sigma^\infty_+$ commutes with its totalization \cite[Corollary~A.3 and Proposition~A.6]{hess-shipley-cohochschild}.
Now set \[
C_n = \Tot_n \CB^\bullet(
  \Sigma^\infty_+ X,\Sigma^\infty_+ Y,\Sigma^\infty_+ Z)
\quad\text{and}\quad C = \lim_n C_n.
\]
Since $\Tot_n$ is essentially a finite limit, it is preserved by $F_0$ and so we need only prove that $F_0 C \wequi \lim_n F_0 C_n$.
Put $\overline A = \fib(\Sigma^\infty_+ Y \to \1)$.
By the dual of \cite[Remark~1.2.4.3]{lurie-ha}, the totalization layers are
\[
\fib(C_r \to C_{r-1}) \wequi
\Omega^r\bigl(
  \Sigma^\infty_+ X \otimes \overline A^{\otimes r}
  \otimes \Sigma^\infty_+ Z
\bigr),
\]
since the codegeneracies apply the counit to the factors $\Sigma^\infty_+Y$.
As $Y$ is simply connected, $\overline A$ is $2$-connective, so the $r$-th layer is $r$-connective.
Taking the inverse limit of the fibers of $C_m \to C_n$ therefore shows that $E_n := \fib(C \to C_n)$ is $n$-connective (recall that we are working in a presheaf topos), and hence $F_0 E_n \in \Sigma^n\scr D$.
By Postnikov completeness, and since limits commute with limits, we have \[ \lim_n F_0E_n \wequi \lim_n\lim_m\tau_{\le m}F_0E_n \wequi \lim_m\lim_n\tau_{\le m}F_0E_n = 0.  \]
This concludes the proof.
\end{proof}

\begin{corollary} \label{cor:EM}
Let $k$ be a field and $A \in \CAlg(\SH(k)_{\ge 0})$.
Given a fiber sequence \[ F \to E \to B \in \Spc(k)_* \] with $B$ simply connected and $E$ connected, we have \[ A[F] \wequi \CB(A[E], A[B], A). \]
\end{corollary}
\begin{proof}
Immediate from Proposition \ref{prop:EM} applied to $\scr X = \Shv_\Nis(\Sm_k)$, $\scr D = \Mod_A(\SH(k))_{\ge 0}$ and $F= A \otimes \Sigma^\infty_+(\ph) \circ L_\mot$.
\end{proof}

We briefly return to the case of a general symmetric monoidal $\infty$-category $\scr C$ and coalgebra $A$.
The forgetful functor from comodules over $A$ to $\scr C$ admits a right adjoint of the form $X \mapsto A \otimes X$ (the relevant argument \cite[Corollary~4.2.4.8]{lurie-ha} for the opposite categories does not require tensor products to be compatible with geometric realizations).
We will employ the following somewhat technical fact.
\begin{lemma} \label{lemm:cofree-cotor}
Let $X \in \scr C$ be strongly dualizable and $A \to B \in \cAlg(\scr C)$.
Then \[ \CB(X \otimes B, B, A) \wequi X \otimes A \] as $A$-comodules.
\end{lemma}
\begin{proof}
We have the equivalence of diagrams $\CB^\bullet(X \otimes B, B, A) \wequi X \otimes \CB^\bullet(B, B, A)$.
Since $X$ is strongly dualizable, tensoring with it is right adjoint to tensoring with its dual \cite[Lemma~4.6.1.6]{lurie-ha}, and hence commutes with all limits. Consequently, $\CB(X \otimes B, B, A) \wequi X \otimes \CB(B,B,A)$.
But the latter cobar construction admits an extra codegeneracy, and so totalizes to $A$ as needed \cite[Remark~4.7.2.3]{lurie-ha}.
\end{proof}

\subsection{Toda fibrations} \label{subsub:toda}
Given a pointed (ordinary) space $X$, we have the \emph{James construction} $JX$, which is a filtered space \[ * = J_0X \to X=J_1 X \to J_2 X \to \dots \] with $J_nX/J_{n-1}X \wequi X^{\wedge n}$.
By working sectionwise, this immediately extends to presheaf topoi, and hence by using the embedding into a presheaf topos (and relocalizing) we may extend the functors $J, J_n$ to $\Spc(k)_*$.

\begin{remark} \label{rmk:J-omegasigma}
A key point is that, if $X \in \Spc(k)_*$ is connected, then $JX \wequi \Omega\Sigma X$ \cite[Theorem~2.4.2(1)]{asok2017simplicial}.
\end{remark}

\begin{remark} \label{rmk:Jn-functor}
For a pointed ordinary space $X$, we have pushout squares \cite[Construction~5.1]{devalapurkar2019james}
\begin{equation*}
\begin{CD}
J_{n-2} X @>>> J_{n-2} X \times X \\
@VVV               @VVV \\
J_{n-1} X @>>> A_n(X)
\end{CD}
\end{equation*}
and
\begin{equation*}
\begin{CD}
A_n(X) @>>> J_{n-1} X \times X \\
@VVV               @VVV \\
J_{n-1} X @>>> J_n X.
\end{CD}
\end{equation*}
This implies that the functor \[ J_n: \Shv_\Nis(\Sm_k)_* \to \Shv_\Nis(\Sm_k)_* \] preserves motivic equivalences and hence induces the functor $J_n: \Spc(k)_* \to \Spc(k)_*$.
\end{remark}

\begin{example} \label{ex:J-effective}
If $X \in \Spc(k)_*$ is (very) $n$-effective then $JX$ is also (very) $n$-effective; the same holds for $J_n(X)$ (and $A_n(X)$).
This follows from the above pushout squares and stability of effective objects under finite products \cite[Lemma 3.1.8]{freudenthal} and colimits (see Remark \ref{rmk:effective-operations}).
\end{example}

Here is the main result of this subsection.
\begin{proposition} \label{prop:Toda-fibrations}
Let $k$ be an unorderable\NB{is this necessary?} field and $n \ge 2$.
We have $p$-local fiber sequences \[ S^{2n-1,n} \to \Omega J_{p-1} S^{2n,n} \xrightarrow{H} \Omega S^{2np-1,np} \] and \[ J_{p-1}S^{2n,n} \to J S^{2n,n} \xrightarrow{H} JS^{2np,np}. \]
\end{proposition}

We begin by constructing the Hopf maps appearing in the fibration sequences.
We can more generally do this for all $n > 0$, without assumptions on $k$.
We start with the second fiber sequence, where the construction is easier.
\begin{construction} \label{cons:H-2}
Let $X \in \Spc(k)_*$ be connected.
The map $H: \Omega \Sigma X \to \Omega \Sigma X^{\wedge p}$ is the adjoint of the projection map \[ \Sigma \Omega \Sigma X \wequi \Sigma(X \vee X^{\wedge 2} \vee \dots) \to \Sigma X^{\wedge p}, \] where the equivalence is \cite[\S4.2]{devalapurkar2019james}.
Since $X$ is connected, this takes the form $H: JX \to J(X^{\wedge p})$ (see Remark \ref{rmk:J-omegasigma}).
\end{construction}

Now we treat the first fiber sequence.
We employ the following elementary fact.
\begin{lemma} \label{lem:cone-elem}
Let $\scr X$ be an $\infty$-category in which products (with a fixed object) preserve pushouts (e.g., one with universal colimits).
Suppose given $f: K \to A \in \scr X_*$ and $X \in \scr X_*$ and set $B = A \cup_f CK$.
Then \[ (B \amalg_A (A \times \Sigma X)) \cup C(K * X) \wequi B \times \Sigma X. \]
\end{lemma}
\begin{proof}
Since products preserve pushouts, we have (a) $K \times \Sigma X \wequi K \amalg_{K \times X} K$ and (b) $B \times \Sigma X \wequi \Sigma X \amalg_{K \times \Sigma X} (A \times \Sigma X)$.
Consider the pushout squares
\begin{equation*}
\begin{CD}
X \amalg_{K \times X} K @>>> * \amalg_{K \times X} K @. \quad @. (K \times \Sigma X) \amalg_K * @>>> (A \times \Sigma X) \amalg_A B \\
@VVV                         @V{\alpha}VV                         @.        @V{\beta}VV                                       @VVV \\
* @>>>                        \Sigma X \amalg_* *    @.       @.  \Sigma X \amalg_* * @>>> P.
\end{CD}
\end{equation*}
Observation (a) identifies $\alpha$ and $\beta$, whereas (b) identifies $P$ with $B \times \Sigma X$.
This concludes the proof.
\end{proof}

Now let $X \in \Spc(k)_*$.
We claim that \[ J_n \Sigma X = J_{n-1} \Sigma X \cup C \Sigma^{n-1} X^{\wedge n}. \]
This is clear for $n=1$ (since $J_1 \Sigma X = \Sigma X = J_0 \Sigma X \cup CX$).
Inductively apply Lemma \ref{lem:cone-elem} with $K=\Sigma^{n-2} X^{\wedge n-1}$ to see that in the diagram below, the left hand square is a pushout
\begin{equation*}
\begin{CD}
K * X @>>> A_n(\Sigma X) @>>> J_{n-1} \Sigma X \\
@VVV          @VVV               @VVV \\
* @>>> J_{n-1} \Sigma X  \times \Sigma X @>>> J_n \Sigma X.
\end{CD}
\end{equation*}
Since the right hand square is a pushout by Remark \ref{rmk:Jn-functor}, so is the outer one by pasting, as needed.
From this we obtain the coaction or ``pinch'' map \[ J_n(\Sigma X) \to (\Sigma X)^{\wedge n} \vee J_n(\Sigma X). \]
\begin{construction} \label{cons:H-1}
If $p=2$, the map we are looking for is another instance of Construction \ref{cons:H-2}.
So we now assume $p$ odd.
The map $\Omega J_{p-1} S^{2n,n} \xrightarrow{H} \Omega S^{2np-1,np}$ is obtained as the composite
\begin{gather*} \Omega J_{p-1}S^{2n,n} \xrightarrow{(1)} \Omega(J_{p-1}S^{2n,n} \vee S^{2n(p-1),n(p-1)}) \\
  \stackrel{(2)}{\wequi} \Omega(J_{p-1}S^{2n,n}) \times \Omega(S^{2n(p-1),n(p-1)}) \times \Omega \Sigma (\Omega J_{p-1}S^{2n,n} \wedge \Omega S^{2n(p-1),n(p-1)}) \\
  \xrightarrow{(3)} \Omega \Sigma (\Omega J_{p-1}S^{2n,n} \wedge \Omega S^{2n(p-1),n(p-1)}) \xrightarrow{c} \Omega(\Omega J_{p-1}S^{2n,n} \wedge S^{2n(p-1),n(p-1)}) \\
  \wequi \Omega(\Omega J_{p-1}S^{2n,n} \wedge S^2 \wedge S^{2n(p-1)-2,n(p-1)}) \xrightarrow{c} \Omega(J_{p-1}S^{2n,n} \wedge S^1 \wedge S^{2n(p-1)-2,n(p-1)}) \\
  \stackrel{(4)}{\wequi} \Omega(\bigvee_{i=1}^{p-1} S^{2ni,ni} \wedge S^1 \wedge S^{2n(p-1)-2,n(p-1)}) \\
  \xrightarrow{(5)} \Omega (S^{2n,n} \wedge S^1 \wedge S^{2n(p-1)-2,n(p-1)}) \wequi \Omega S^{2np-1,np}. \end{gather*}
Here (1) is $\Omega$ applied to the pinch map above, (2) is the Hilton--Milnor splitting \cite[Theorem~3.1]{devalapurkar2019james}, (3) is the projection, $c$ are counit maps $\Sigma \Omega \to \id$, (4) is the splitting of $\Sigma J_{p-1}$ \cite[Proposition~5.8]{devalapurkar2019james}, and (5) is projection to the summand $\Sigma S^{2n,n}$ in the splitting.
(We have used that $p > 2$ and $n > 0$ to split off $S^2$ from $S^{2n(p-1),n(p-1)}$.)
\end{construction}

\begin{proof}[Proof of Proposition \ref{prop:Toda-fibrations}.]
By base change we may assume that $k$ is perfect and $cd_2(k) < \infty$.

\textbf{Second fiber sequence.}
The composite $J_{p-1}S^{2n,n} \to J S^{2n,n} \to JS^{2np,np}$ is null, being adjoint to \[ \Sigma\left(\bigvee_{1 \le i < p} S^{2ni,ni} \to \bigvee_{i \ge 1} S^{2ni,ni} \to S^{2np,np}\right). \]
We consequently obtain a map $J_{p-1}S^{2n,n} \to F := \fib(J S^{2n,n} \to JS^{2np,np})$ which we must show is a $p$-local equivalence.
The spaces $J_{p-1} S^{2n,n}$, $J S^{2n,n}$ and $JS^{2np,np}$ are nilpotent, and are limit-determined by Examples \ref{ex:resolvable-spaces} and \ref{ex:J-effective}.
Hence so is $F$.
Consequently by Lemma \ref{lemm:limit-determined-M-test}, it suffices to show that the map is inverted by $H\Z_{(p)}[\ph]$.
We claim that the stable map $JS^{2n,n} \to J_{p-1}S^{2n,n}$ exhibits $H\Z_{(p)}[JS^{2n,n}]$ as the cofree $H\Z_{(p)}[JS^{2np,np}]$-comodule $H\Z_{(p)}[JS^{2np,np}] \otimes H\Z_{(p)}[J_{p-1}S^{2n,n}]$.
Assuming this, by Corollary \ref{cor:EM} and Lemma \ref{lemm:cofree-cotor} we see that $H\Z_{(p)}[F] \wequi H\Z_{(p)}[J_{p-1}S^{2n,n}]$ via the stable map $F \to J S^{2n,n} \to J_{p-1} S^{2n,n}$, and we conclude.

To prove the claim, we note that all the involved motives are split Tate and degreewise finite, so we may check the statement in cohomology in Chow degree zero.
Note that $H\Z_{(p)}[JS^{2n,n}]$ is the free pointed $\scr E_1$-algebra on $H\Z_{(p)}[S^{2n,n}]$, or equivalently the free $\scr E_1$-algebra on $\overline{H\Z_{(p)}[S^{2n,n}]}$, and receives a coalgebra map from $H\Z_{(p)}[S^{2n,n}]$.
Thus its homotopy is of the form $H\Z_{(p)2*,*}[t_{2n,n}]$ (where $t_{2n,n}$ corresponds to $S^{2n,n}$) with Hopf algebra structure $\Delta(t_{2n,n}) = t_{2n,n} \otimes 1 + 1 \otimes t_{2n,n}$ (because this is true in $S^{2n,n}$).
Thus $H^{2*,*}(JS^{2n,n}, \Z_{(p)}) \wequi \Gamma_{\Z_{(p)}}(x_{2n,n})$ is a divided power algebra.
Similarly $H^{2*,*}(JS^{2np,np}, \Z_{(p)}) \wequi \Gamma_{\Z_{(p)}}(y_{2np,np})$.
By construction, the map $H^{2*,*}(JS^{2np,np}, \Z_{(p)}) \to H^{2*,*}(JS^{2n,n}, \Z_{(p)})$ sends $y_{2np,np}$ to a generator in the appropriate degree, i.e., $x_{2n,n}^{(p)}$.
Since $\Z_{(p)}$ is torsion-free, the map must be compatible with the divided power operations, and so is determined completely.
The map $H^{2*,*}(J_{p-1}S^{2n,n}, \Z_{(p)}) \to H^{2*,*}(JS^{2n,n},\Z_{(p)})$ is (by construction) the inclusion of $x_{2n,n}^{(i)}$ for $i < p$.
In the end, then, the map \[ H^{2*,*}(JS^{2np,np}, \Z_{(p)}) \otimes H^{2*,*}(J_{p-1}S^{2n,n}, \Z_{(p)}) \to H^{2*,*}(JS^{2n,n}, \Z_{(p)}) \] identifies with the $\Gamma_{\Z_{(p)}}(y_{2np,np})$-module map \[ \Gamma_{\Z_{(p)}}(y_{2np,np}) \otimes \Z_{(p)}\{1, x_{2n,n}, x_{2n,n}^{(2)}, \dots, x_{2n,n}^{(p-1)}\} \to \Gamma_{\Z_{(p)}}(x_{2n,n}), \] and the claim is equivalent to the assertion that the elements $(x^{(p)}_{2n,n})^{(i)}x_{2n,n}^{(j)}$ with $i \ge 0$, $0 \le j < p$ form a basis of $\Gamma_{\Z_{(p)}}(x_{2n,n})$.
This is well-known \cite[Proposition~5.1.7]{neisendorfer-algebraic}.

\textbf{First fiber sequence.}
The arguments are analogous, but somewhat more elaborate.
We treat here the case $k=\F_p$.
By base change, this immediately implies the case of any field of characteristic $p$.
Essentially the same argument also works in characteristic $\ne p$, but a bit more care has to be taken.\todo{should we reformulate the proposition to only claim this for fields of characteristic $p$?}

The map $S^{2n-1,n} \to \Omega J_{p-1} S^{2n,n}$ is the inclusion of the bottom cell.
The composite with $H$ vanishes for degree reasons.\NB{If $n=1$, $p=2$, this would not work. But then the source is $\Gm$ which we shouldn't $p$-localize anyway.}
We hence obtain $S^{2n-1,n} \to F := \fib(\Omega J_{p-1} S^{2n,n} \to \Omega S^{2np-1,np})$, and we must show as before that this is inverted by $H\Z_{(p)}[\ph]$.
We claim that the stable map $\Omega J_{p-1} S^{2n,n} \to S^{2n-1,n}$ (loop evaluation followed by projection to the bottom cell) exhibits $H\Z_{(p)}[\Omega J_{p-1} S^{2n,n}]$ as the cofree $H\Z_{(p)}[\Omega S^{2np-1,np}]$-comodule $H\Z_{(p)}[\Omega S^{2np-1,np}] \otimes H\Z_{(p)}[S^{2n-1,n}]$.
Assuming this, we conclude as before.
We know by Eilenberg--Moore that $H\Z_{(p)}[\Omega J_{p-1} S^{2n,n}]$ is the cobar construction of $H\Z_{(p)}[J_{p-1}S^{2n,n}]$.
The latter is a finite cellular motive, and as seen in the proof of Proposition \ref{prop:EM}, the corresponding limit converges rapidly.
In particular, the formation of this limit commutes with taking duals, and the limit is reflexive (in the sense of equalling its double dual).
It will thus suffice to consider cohomology.
We will simplify our considerations by using a fact established later: since $DH\Z_{(p)}[J_{p-1}S^{2n,n}]$ is a cellular normed algebra, it is by Proposition \ref{prop:NAlg-DAlg} in fact a derived algebra.
Recalling that $H^{**}(\F_p, \Z_{(p)}) \wequi \Z_{(p)}$ (see Remark \ref{rmk:finite-field-motivic-cohomology}) we get $H^{**}(J_{p-1} S^{2n,n}, \Z_{(p)}) \wequi \Z_{(p)}[z_{2n,n}]/z_{2n,n}^p$, and so as a derived algebra, $DH\Z_{(p)}[J_{p-1}S^{2n,n}] \wequi \LSym(z_{2n,n})/z_{2n,n}^p$ (where by the right hand side we mean the evident pushout in derived algebras).
By Corollary \ref{cor:EM} and the compatibility with duality explained above, this implies that \[ DH\Z_{(p)}[\Omega J_{p-1}S^{2n,n}] \wequi B(\LSym(z_{2n,n})/z_{2n,n}^p) \wequi \LSym(z_{2n-1,n},y_{2np-2,np}) \] and so\footnote{Over other fields, this conclusion would need to be reached by a more careful spectral sequence argument.} \[ H^{**}(\Omega J_{p-1}S^{2n,n}, \Z_{(p)}) \wequi E(z_{2n-1,n}) \otimes \Gamma_{\Z_{(p)}}(y_{2np-2,np}). \]
It remains to show that the Hopf map $H: \Omega J_{p-1} S^{2n,n} \to \Omega S^{2np-1,np}$ sends the bottom cohomology generator $\bar y_{2np-2,np}$ of the target to a $p$-local unit multiple of $y_{2np-2,np}$.
We may as well check the dual statement in homology.
Put $R=\Z_{(p)}$, $A=J_{p-1}S^{2n,n}$, and abbreviate $z=z_{2n,n}$ and $y=y_{2np-2,np}$.
Let $x \in H_{2n,n}(JS^{2n,n};R)$ be the standard generator and put $b_i=x^i/i!$ for $0\le i<p$.
Then $b_i$ is dual to $z^i$ and
\[
  \Delta b_i=\sum_{j=0}^i b_j\otimes b_{i-j}.
\]
One verifies that the expression
\begin{equation} \label{eq:Toda-cobar-generator}
  \zeta=\sum_{i=1}^{p-1}[b_i\mid b_{p-i}]
  \in H_{2np-2,np}(\Omega A;R)
\end{equation}
is a cycle generating the indicated homology group (it pairs to $1$ with the relation class $[z^{p-1}\mid z]$ representing $y$ in the bar complex).

Suppose first that $p$ is odd and put $T=S^{2n(p-1),n(p-1)}$.
Let $t$ be the top homology generator of $T$, oriented so that the top-cell collapse sends $x^{p-1}$ to $t$.
For the pinch map $\rho:A\to A\vee T$ of Construction \ref{cons:H-1}, we then have
\[
  \rho_*(b_i)=b_i\quad (i<p-1),\qquad
  \rho_*(b_{p-1})=b_{p-1}+\frac{t}{(p-1)!}.
\]
Consequently, writing $\iota:\Omega A\to\Omega(A\vee T)$ for the pure-factor inclusion (i.e., not the looped pinch map), the induced cobar map gives
\[
  (\Omega\rho)_*(\zeta)=\iota_*(\zeta)
  +\frac{1}{(p-1)!}\bigl([b_1\mid t]+[t\mid b_1]\bigr).
\]
One verifies\NB{insert?} that the mixed Hilton--Milnor projection kills $\iota_*(\zeta)$ and sends $[b_1\mid t]+[t\mid b_1]$ to the bottom mixed generator, up to sign.
The counits and the projection to the bottom James summand in Construction \ref{cons:H-1} identify that generator with the bottom generator of $\Omega S^{2np-1,np}$, again up to sign.
For $p=2$, the James--Hopf map of Construction \ref{cons:H-2} sends the length-two class $\zeta=[b_1\mid b_1]$ to the bottom target generator, by the defining projection in the James splitting.
Thus in both cases, if $\bar\zeta$ denotes the bottom homology generator of the target, we have
\begin{equation} \label{eq:Toda-Hopf-bottom}
  H_*(\zeta)=\pm\frac{1}{(p-1)!}\bar\zeta.
\end{equation}
In particular, the coefficient is a unit in $R$, as needed.
\end{proof}

Let us also record the following boundary map.
Rotating the second fiber sequence and looping, we obtain a map $\Omega^2 JS^{2np,np} \to \Omega J_{p-1} S^{2n,n}$.
Composing with the Hopf map from the first fiber sequence, we obtain $\Omega^2 JS^{2np,np} \to JS^{2np-2,np}$.
The bottom cell on each side is $S^{2np-2,np}$.
\begin{lemma} \label{lem:boundary-p-map}
The bottom cell map above is of degree $p$ (up to a $p$-local unit).
\end{lemma}
\begin{proof}
As in the proof of Proposition \ref{prop:Toda-fibrations}, we work over $k=\F_p$ with coefficients $R=\Z_{(p)}$; base change then gives the characteristic $p$ case.
The characteristic $\ne p$ case can be argued similarly but more carefully (or one can use étale realization to reduce to topology).
Retain $x$, $b_i$, and $\zeta$ from \eqref{eq:Toda-cobar-generator}.
Write $q:JS^{2n,n}\to JS^{2np,np}$ for the second Hopf map, normalize the bottom target class $v$ by $q_*(x^p)=v$, and let $\beta$ be its double desuspension.
We compute the looped boundary $\partial_*(\beta)$, as follows.
First use $Y^\bullet=\CB^\bullet(JS^{2n,n},JS^{2np,np},*)$ to model the fiber, and then apply the loop cobar construction.
The boundary class is represented by $1\otimes v$ in inner and outer cobar degrees $(1,1)$.
The element $x^p$ in degrees $(0,1)$ has inner differential $1\otimes v$, since its coaction is $x^p\otimes1+1\otimes v$, and outer differential given by its reduced coproduct.
Thus the total differential identifies the boundary class, up to sign, with
\[
  \partial_*(\beta)
  =\pm\sum_{i=1}^{p-1}\binom pi[x^i\mid x^{p-i}]
  =\pm p!\sum_{i=1}^{p-1}[b_i\mid b_{p-i}]
  =\pm p!\,\zeta.
\]
All factors in this representative lie in $H_{2*,*}(J_{p-1}S^{2n,n};R)$; the cofree-comodule calculation in Proposition \ref{prop:Toda-fibrations} identifies the inner cobar cohomology with this fiber coalgebra, so the representative computes the actual boundary without further corrections.
Finally, \eqref{eq:Toda-Hopf-bottom} gives
\[
  H_*\partial_*(\beta)=\pm p\,\bar\zeta,
\]
which proves the assertion.
\end{proof}

\subsection{Affine Grassmannians revisited} \label{sec:grass}
In \cite{bachmann-grassmann}, the first author showed that the motivic loop space $\Omega_\Gm \SL_n$ is modeled by an affine Grassmannian.
Moreover, under the assumption of resolution of singularities, he exhibited a cell structure and used it to compute the homology of $\Omega_\Gm \SL_n$.
The aim of this section is to revisit the proofs of \emph{loc. cit.}, modifying them slightly to show that resolution of singularities can be avoided.
(Specifically, instead of using the stratification of the affine Grassmannian by finite dimensional \emph{singular} closed subvarieties, we employ the opposite stratification by finite codimensional open subschemes.)

\begin{definition}
The affine Grassmannian \cite[\S2]{faltings-loop-groups} of $\SL_n$ is an fppf sheaf \[ \Gr{\SL_n}: \Sch^\op \to Set \] which exhibits an fppf sheafification of the presheaf $X \mapsto \SL_n(X\fpsr{t}[t^{-1}])/\SL_n(X\fpsr{t})$.
\end{definition}

\begin{lemma}
\begin{enumerate}
\item The sheaf $\Gr{\SL_n}$ is isomorphic to a filtered colimit of sheaves represented by projective $\Z$-schemes.
\item The canonical map \[ a_\Zar \SL_n((\ph)[t,t^{-1}])/\SL_n((\ph)[t]) \to \Gr{\SL_n} \] is an isomorphism.
\end{enumerate}
\end{lemma}
\begin{proof}
For (1), see \cite[\S2]{faltings-loop-groups}.
Now we prove (2).
Write $L\SL_n$ for the sheaf of groups $A \mapsto \SL_n(A\fpsr{t}[t^{-1}])$ and consider the map $p: L\SL_n \to \Gr{\SL_n}$.
We shall show that there is an open cover $\{U_\alpha\}$ of $\Gr{\SL_n}$ and for each $\alpha$ a map $U_\alpha \to \SL_n((\ph)[t,t^{-1}])$ such that the composition with $\SL_n((\ph)[t,t^{-1}]) \to L\SL_n$ is a section of $p$.
This shows that $\Gr{\SL_n}$ is the Zariski quotient $a_\Zar L\SL_n/\SL_n((\ph)\fpsr{t})$, and also that $\SL_n((\ph)[t,t^{-1}]) \to \Gr{\SL_n}$ is (Zariski-locally) surjective.
Since $\SL_n(A\fpsr{t}) \cap \SL_n(A[t,t^{-1}]) = \SL_n(A[t])$, the result follows.

It remains to construct the open cover $\{U_\alpha\}$ and the sections.
For this, see Lemma \ref{lemm:aff-grass}(3) below: take $U_\alpha$ of the form $\mu\Sigma_0$, which identifies with $\mu L^{--}\SL_n \subset \SL_n((\ph)[t,t^{-1}])$ (and admits a section via this identification).
\end{proof}

We write $L^-\SL_n$ for the sheaf of groups $\SL_n((\ph)[t^{-1}])$ and $L^+\SL_n$ for the sheaf of groups $\SL_n((\ph)\fpsr{t})$, both of which act on $\Gr{\SL_n}$.
We write $N^-\SL_n \subset L^-\SL_n$ and $N^+ \SL_n \subset L^+ \SL_n$ for the subsheaves of matrices whose constant term is upper, respectively lower, unitriangular.
Let $L^{--}\SL_n \subset L^-\SL_n$ denote the kernel of the reduction mod $t^{-1}$.
Recall the lattice of homomorphisms $\Gm \to \SL_n$ landing in diagonal matrices; we denote such homomorphisms generically by $\mu$.
The constant homomorphism at the identity matrix is denoted $0$.
Note that $\mu \in \SL_n(\Z[t,t^{-1}])$ and thus we obtain a point $[\mu] \in \Gr{\SL_n}$.

\begin{definition}
We denote by $\Sigma_\mu$ the $N^-\SL_n$-orbit of $[\mu]$, and by $C_\mu$ the $N^+\SL_n$-orbit.
\end{definition}
These can be constructed as locally closed, sub-ind-schemes of $\Gr{\SL_n}$ \cite[\S2,\S3]{faltings-loop-groups}.
\begin{lemma} \label{lemm:aff-grass}
\begin{enumerate}
\item $\Sigma_0$ is open (and $C_0 = [0]$ is closed).
\item For each $\mu$, the open subset $\mu\Sigma_0$ contains both $\Sigma_\mu$ and $C_\mu$.
\item We have $\Sigma_0 \wequi L^{--}\SL_n$.
\item Both the $\Sigma_\mu$ and the $C_\mu$ define locally closed stratifications of $\Gr{\SL_n}$.
\item There is an isomorphism \[ \mu\Sigma_0 \wequi \Sigma_\mu \times C_\mu, \] compatible with the embedding of $\Sigma_\mu$ and $C_\mu$ in $\mu\Sigma_0$.
\item The closure of $C_\mu$ is the union of those $C_\lambda$ with $\lambda \le \mu$, whereas the closure of $\Sigma_\mu$ is the union of those $\Sigma_\lambda$ with $\lambda \ge \mu$ in the Bruhat order.\NB{or its opposite?}
\item $\Sigma_\mu$ is $\A^1$-contractible, and $C_\mu \wequi \A^{d(\mu)}$ where $d(\mu)$ is the length of $\mu$ \cite[p. 46]{faltings-loop-groups}.
\end{enumerate}
\end{lemma}
\begin{proof}
(1) to (6) as well as the second half of (7) (i.e., $C_\mu \wequi \A^{d(\mu)}$) can be extracted from \cite[\S3]{faltings-loop-groups}.\NB{unfortunately this is not immediate, since Faltings writes about the flag variety, not the Grassmannian}
To prove contractibility of $\Sigma_\mu$, we can e.g. note by (5) that it is a retract of $\mu\Sigma_0$ to reduce to proving that $\Sigma_0$ is contractible.
By (3), $\Sigma_0$ identifies with the variety of $n \times n$ matrices with values in $A[t^{-1}]$ of determinant $1$ and constant term the identity matrix.
This is $\A^1$-contractible by scaling $t^{-1}$ down to zero.
\end{proof}

Now choose a total ordering $\mu_0, \mu_1, \dots$ of the $\mu$ refining the Bruhat order.
(Thus $\mu_0 = 0$, $\mu_1$ is minimal among those $\mu \ne 0$, $\mu_2$ is minimal among those $\mu \ne 0,\mu_1$, and so on.)
Write $\Sigma_r=\Sigma_{\mu_r}$ and $C_r=C_{\mu_r}$ for the individual strata.
Let $U_r=\bigcup_{i=0}^r\Sigma_i$, which is open by Lemma \ref{lemm:aff-grass}(6).
(To see that $U_r$ is open, by induction we need only prove that $\Sigma_r$ is open in the complement of $U_{r-1}$.
It is given by the intersection of the complements of the closures of the finitely many $\Sigma_\mu$ with $\mu > \mu_r$ minimal.)

\begin{lemma} \label{lemm:aff-grass-cells}
As Zariski sheaves on all schemes, we have $\Gr{\SL_n} \wequi \colim_r U_r$, and for $r\ge1$ there are pushout squares (of Zariski sheaves)
\begin{equation*}
\begin{CD}
\Sigma_r \times (C_r \setminus 0) @>>> U_{r-1} \\
@VVV         @VVV \\
\mu_r \Sigma_0 @>>> U_r.
\end{CD}
\end{equation*}
\end{lemma}
Note that the topological analog of the above decomposition was used in \cite[Proof of (8.6.6)]{pressley1984loop}.
\begin{proof}
The first statement is clear.

Now we prove the second statement.
By Lemma \ref{lemm:aff-grass}(5), $\mu_r \Sigma_0 \subset N^-\SL_n C_{\mu_r}$, which is the union of those $N^-\SL_n$-orbits meeting $C_{\mu_r}$.
By \cite[p. 48]{faltings-loop-groups} this happens only for orbits around $\mu_i$, $i \le r$.
This shows that $\mu_r \Sigma_0 \subset U_r$.
Consequently by Lemma \ref{lemm:aff-grass}(2), $U_r = U_{r-1} \cup \mu_r \Sigma_0$ is an open cover, and to obtain the desired pushout square we need only show that $\mu_r \Sigma_0 \cap U_{r-1} \wequi \Sigma_r \times (C_r \setminus 0)$.
Since $U_{r-1} = U_r \setminus \Sigma_r$, this is immediate from Lemma \ref{lemm:aff-grass}(5).
\end{proof}

For any scheme $S$, restricting the presheaf $\Gr{\SL_n}$ to $\Sm_S$ and then motivically localizing, we obtain a motivic space $\Gr{\SL_n}_S \in \Spc(S)$.
Similarly for $U_{r,S}$.
Recall also the \emph{generating complex} $\P^{n-1} \to \Gr{\SL_n}$, e.g. from \cite[\S5, Example~(i), and \S6]{littig2011generating}.
We have arrived at the main result of this section.
\begin{proposition} \label{prop:aff-grass}
\begin{enumerate}
\item We have $\Gr{\SL_n}_S \wequi \colim_r U_{r,S} \in \Spc(S)$, $U_0$ is contractible, and $r\ge1$ there are pushout squares
\begin{equation*}
\begin{CD}
\A^{d(\mu_r)} \setminus 0 @>>> U_{r-1,S} \\
@VVV         @VVV \\
* @>>> U_{r,S}.
\end{CD}
\end{equation*}
\item For $f: T \to S \in \Sch$, the canonical map $f^* \Gr{\SL_n}_S \to \Gr{\SL_n}_T$ is an equivalence in $\Spc(T)$.
\item If $S$ is the spectrum of a field or $\Z$, then $\Gr{\SL_n}_S \wequi \Omega_\Gm \SL_{n,S} \in \Spc(S)$.
\item We have \[ H\Z[\Gr{\SL_n}_S] \wequi H\Z[b_1, b_2, \dots, b_{n-1}] \in \h\Mod_{H\Z}(S), \] where $b_i \in H\Z_{2i,i}\Gr{\SL_n}_S$ is the image of the corresponding class in the generating complex.
  The coalgebra structure is \[ \Delta(b_i) = \sum_{i=j+k} b_j \otimes b_k. \]
\end{enumerate}
\end{proposition}
\begin{proof}
(1) is immediate from Lemma \ref{lemm:aff-grass-cells} and Lemma \ref{lemm:aff-grass}(7).
(2) follows directly from (1).
(3) was proved for fields in \cite[Theorem 15]{bachmann-grassmann}, and the argument over $\Z$ is the same.
To prove (4), it suffices by base change to work over $\Z$.
Note first that the cell decomposition implies by induction that \[ H\Z[\Gr{\SL_n}] \wequi \bigoplus_r \Sigma^{2d(\mu_r),d(\mu_r)}H\Z. \]
The monomials in the $b_i$ provide homology classes in the same degrees as the above generators, by comparison with the analogous topological result \cite[Proposition~10.27.6]{neisendorfer-algebraic}.
It will thus suffice to show that the polynomial ring on the $b_i$ is a summand of $H\Z[\Gr{\SL_n}]$.
By construction, the map $\SL_n \to \SL_{n+1}$ induces a map $\Gr{\SL_n} \to \Gr{\SL_{n+1}}$ compatible with the cell structures, so that $H\Z[\Gr{\SL_n}] \to H\Z[\Gr{\SL_{n+1}}]$ is a split injection.
It is thus enough to show that the homology of \[ \Gr{\SL_\infty} \wequi \Omega_\Gm \SL_\infty \] has the required form.
By motivic Bott periodicity \cite[\S4, Theorem~3.13]{A1-homotopy-theory}, we have $\Omega_\Gm\SL_\infty \wequi \Omega_{\P^1}B\SL \wequi B\GL$, and the result follows from, e.g., \cite[Proposition~6.2(ii)(b)]{naumann2009motivic}.
\end{proof}


Let us record already here our most important consequence.
\begin{corollary} \label{cor:Omega11-S2n-1}
Let $k$ be a field (or $\Z$) and $n \ge 2$.
Then $\Omega^{1,1}S^{2n-1,n}$ is strongly cellular (in the sense of being built by iterated colimits out of motivic spheres $S^{p,q}$), so in particular $H\Z[\Omega^{1,1}S^{2n-1,n}]$ is cellular.
We have $H\Z_{**}\Omega^{1,1}S^{2n-1,n} \wequi H\Z_{**}[b]$ as coalgebras, where $|b|=(2n-2,n-1)$ and $b$ is primitive.
\end{corollary}
\begin{proof}
We have the fiber sequence $S^{2n-1,n} \wequi \A^n \setminus 0 \to B\SL_{n-1} \to B\SL_n$ \cite[Theorems~2.2.5 and 3.3.1, proof of Corollary~4.2.6]{asok2015affine}.
Applying $\Omega^{1,1}$ and rotating we find that $\Omega^{1,1}S^{2n-1,n} \wequi \Omega^{1,1}\SL_n \sslash \Omega^{1,1}\SL_{n-1}$, and so the same holds after applying $H\Z[\ph]$.
We now use Proposition \ref{prop:aff-grass}.
We immediately deduce that $\Omega^{1,1}S^{2n-1,n}$ is strongly cellular, and so in particular this holds for the motive.
We also see that \[ H\Z_{**} \Omega^{1,1}S^{2n-1,n} \wequi H\Z_{**}[b_1, \dots, b_{n-1}] \otimes^L_{H\Z_{**}[b_1, \dots, b_{n-2}]} H\Z_{**} \wequi H\Z_{**}[b_{n-1}], \] i.e., the tor spectral sequence degenerates.
Setting $b=b_{n-1}$ we obtain the desired identification of $H\Z_{**} \Omega^{1,1}S^{2n-1,n}$ (as an $H\Z_{**}$-module).
We immediately get $\Delta(b) = 1 \otimes b + b \otimes 1$ from the coalgebra structure on $H\Z_{**}\Omega^{1,1}\SL_n$.
\end{proof}

\section{Normed Hopf algebras: general theory} \label{sec:norm-hopf-1}
In classical topology, it is well-known that the homology of any space (with field coefficients) is a coalgebra, the homology of a connected $H$-space is a Hopf algebra\NB{nonassociative...}, and that these additional structures on the homology are often extremely useful (see, e.g., \cite[\S8.9]{milnor-moore}).
Throughout this article, we will utilize this philosophy, but enhanced in two important ways.
\begin{enumerate}
\item Instead of looking at the separate homology groups, we look at the entire chain complex.
  This carries the structure of a coalgebra up to \emph{coherent} homotopies, i.e., an $\scr E_\infty$-coalgebra.
  If the space is a connected commutative monoid up to coherent homotopy, then its chains even form an $\scr E_\infty$-Hopf algebra.
\item Instead of spaces we look at motivic spaces, and there the usual algebraic structures can be upgraded using \emph{norms}.
  So in the end, we will look at normed monoids (in motivic spaces) and want to argue that their chains form a \emph{normed Hopf algebra}.
\end{enumerate}
The aim of this section is to justify this claim (see Example \ref{ex:homology-hopf-algebra}), and also to develop some of the theory of normed Hopf algebras, so that we can actually manipulate these objects.

\subsection{Some normed category theory}
In the rest of the article, we will need various fragments of a theory of normed categories which has not been developed in detail.

\subsubsection{Normed categories and normed algebras} \label{subsec:normed-cat}
We begin by recalling some salient definitions and results.

Following \cite{bachmann-norms}, by a \emph{normed category over $S$} we mean a functor $\scr C^\otimes: \Span(\Sch_S, \all,\fet) \to \Cat$ which preserves finite products.
For each $X \in \Sch_S$ we obtain an $\infty$-category $\scr C^\otimes(X)$, which we typically just denote by $\scr C(X)$.
In fact, we will often call $\scr C(S)$ a normed category, leaving the extension to a functor on $\Span(\Sch_S, \all, \fet)$ implicit.

\begin{example}
The category $\Spc(k)$ is a normed category, where for $p: X \to Y \in \Sch_S$ finite étale, the norm functor $p_\otimes: \Spc(X) \to \Spc(Y)$ is given by the functor $p_*$ right adjoint to the pullback $p^*$ \cite[Proposition~2.11 and \S6.1]{bachmann-norms}.
We call this the \emph{cartesian norm structure}.
(See \S\ref{subsub:cocart-normed} for more about cartesian norm structures.)
\end{example}

By a \emph{normed functor} we mean a natural transformation in $\Fun(\Span(\Sch_S, \all,\fet), \Cat)$.
\begin{example}
The categories $\Spc(S)_*$ and $\SH(S)$ admit norm structures such that \[ \Spc(S) \xrightarrow{(\ph)_+} \Spc(S)_* \xrightarrow{\Sigma^\infty} \SH(S) \] upgrade to normed functors \cite[\S6.1]{bachmann-norms}.
\end{example}

Given a normed category $\scr C(S)$, we have the notion of a \emph{normed algebra} in $\scr C(S)$, defined as an appropriate partially cocartesian section of $\scr C(S)$ over $\Span(\Sm_S^{fp}, \all, \fet)$.
We denote the category of normed algebras by $\NAlg(\scr C(S))$.
There is a forgetful functor $\NAlg(\scr C(S)) \to \scr C(S)$ which, under fairly mild assumptions, admits a left adjoint $\NSym$.
In the special case that the norm structure is cartesian (e.g. on $\Spc(k)$), we speak of \emph{normed monoids} and denote the category by $\NMon(\scr C(k))$.
We write $\cNAlg(\scr C(S)) := \NAlg(\scr C(S)^\op)^\op$ for the category of \emph{co-normed coalgebras}.
In the cocartesian case we denote this by $\cNMon(\scr C(S))$.

\begin{example} \label{ex:HR-normed}
Many rings in $\SH(S)$ upgrade to normed spectra.
Most salient for us are the motivic Eilenberg--Mac Lane spectra $HR$ for $R$ a $\Z$-algebra.
When $R=\Z$ (or $\Z/m$) this is proved in \cite[Theorem 13.9]{bachmann-norms}; the same proof applies in general (note that $H\Z$ is a $\Z$-algebra, so $HR := H\Z \otimes_{\Z} R$ is defined).
\end{example}

Given any normed category $\scr C$ (with sufficiently many colimits), we have a notion of \emph{normed extended power}.
These are functors \[ D_n^{\scr C}: \scr C(k) \to \scr C(k) \] which can be defined by the explicit colimit \cite[Definition~5.6 and Example~5.8]{bachmann-colimits} \begin{equation} \label{eq:formula-for-dn} D_n^{\scr C}(E) = \colim_{X \xrightarrow{p} Y \xrightarrow{f} k} f_\sharp p_\otimes(E|_X). \end{equation}
Here the colimit runs over all smooth $k$-schemes $Y$ and finite étale morphisms $p: X \to Y$ of constant degree $n$.
Under certain assumptions (a key one being that the so-called ``distributivity law'' \cite[\S5.2]{bachmann-norms} holds in $\scr C$) there is an equivalence \cite[Theorem~3.10]{bachmann-MGM} \begin{equation} \label{eq:formula-for-nsym} \NSym(E) \wequi \amalg_{n \ge 0} D_n^{\scr C}(E). \end{equation}
When no confusion can arise, we denote the extended powers just by $D_n$.
A convenient consequence of \eqref{eq:formula-for-nsym} (and the fact that $\NSym(A \oplus B) \wequi \NSym(A) \otimes \NSym(B)$) is that \begin{equation} \label{eq:formula-for-Dn-sum} D_n(E \amalg F) \wequi \amalg_{n = a + b} D_a(E) \otimes D_b(F). \end{equation}

\begin{definition}
We call a normed category $\scr C^\otimes(\ph)$ \emph{presentably normed} if:
\begin{enumerate}
\item Each $\scr C(X)$ is presentable.
\item Each of the pullback functors $f^*$ preserves colimits.
\item Whenever $f$ is smooth, $f^*$ has a left adjoint $f_\sharp$.
\item Each of the norms $p_\otimes$ preserves sifted colimits.
\item The smooth base change law holds, i.e., for every cartesian square
  \[
  \begin{CD}
  X' @>{g'}>> X \\
  @V{f'}VV @VV{f}V \\
  Y' @>{g}>> Y
  \end{CD}
  \]
  with $f$ smooth, the exchange map $f'_\sharp g'^* \to g^*f_\sharp$ is an equivalence.
\item Norms along finite étale morphisms distribute over $f_\sharp$ for $f$ smooth, in the sense of \cite[\S5.2]{bachmann-norms}.
\end{enumerate}
\end{definition}
\begin{remark}
This definition differs from the one in \cite[Definition 6.5]{bachmann-norms} slightly: we require $f_\sharp$ for smooth maps (not just finite étale) and we require the extended distributivity law.
The definition is rigged so that the proof of Lemma \ref{lemm:norms-for-nhopf} below works.
\end{remark}

Finally let us recall the following standard technique for constructing new normed categories by passing to modules.
\begin{lemma}
Let $\scr C$ be presentably normed over $S$ and $A \in \NAlg(\scr C(S))$.
Then $X \mapsto \Mod_{A_X}(\scr C(X))$ is a presentably normed category.
\end{lemma}
\begin{proof}
The norm structure is constructed by the same argument as in \cite[Proposition 7.6(4)]{bachmann-norms}.
The functors $f_\sharp$ (for $f$ smooth) in $\Mod_A$ are computed in $\scr C$.
The norms in $\Mod_A$ are given by $p_\otimes^A(\ph) = p_\otimes(\ph) \otimes_{p_\otimes(A)}A$.
From these two assertions one easily deduces that the distributivity law holds.
\end{proof}

\subsubsection{(Co)cartesian normed categories} \label{subsub:cocart-normed}
We call a norm structure \emph{(co)cartesian} if the norms are given by right (left) adjoints to the restriction.
This is in fact a property, not structure: given a finite étale morphism $p$, there is a $2$-step procedure (starting with the clopen map $\Delta_p$ in place of $p$) for constructing a canonical map $p_\otimes \to p_*$ (respectively $p_\sharp \to p_\otimes$) which is then required to be an equivalence.
In fact, the forgetful functor from (co)cartesian normed categories to $\Fun(\Sch_S^{\op}, \Cat_\infty)$ is fully faithful by \cite[Theorem 3.3]{dauser2024uniqueness}.

Here is another way of characterizing (co)cartesian norm structures.
\begin{lemma} \label{lemm:transported-normed-structure-technical}
There is an adjunction \[ (\ph)^{\amalg}: \{ \scr C \in \Fun^\times(\Sch_S^{fp,\op}, \Cat_\infty) \mid \scr C \text{ admits indexed colimits} \} \adj \Fun^{\times}(\Span(\Sch_S^{fp}, \all, \fet), \Cat_\infty)^{lax}: \NAlg(\ph). \]
Here by ``$\scr C$ admits indexed colimits'' we mean that for every finite étale morphism $p: X \to Y \in \Sch_S^{fp}$ the functor $p^*: \scr C(Y) \to \scr C(X)$ admits a left adjoint $p_\sharp$, and that these functors are compatible with base change.
By $\Fun(\dots)^{lax}$ we mean functors and lax natural transformations, i.e., the right hand side denotes normed categories and lax normed functors.
Finally $\scr C^\amalg$ means the cocartesian norm structure on $\scr C$ (we assert its (functorial) existence, uniqueness is clear by what we said above).

A dual adjunction involving cartesian norm structures and normed coalgebras also holds.
\end{lemma}
\begin{remark} \label{rmk:NAlg-alternative}
Given a normed category $\scr C(S)$, the category $\NAlg(\scr C(S))$ is defined in terms of partially cocartesian sections of $\Span(\Sm_S, \all, \fet)$.
One could alternatively consider the category $\NAlg^{ex}(\scr C(S))$, defined in terms of partially cocartesian sections of $\Span(\Sch_S^{fp}, \all, \fet)$.
As explained in \cite[\S B]{bachmann-splitting}, as long as $\scr C$ satisfies Zariski descent, the forgetful functor $\NAlg^{ex}(\scr C(S)) \to \NAlg(\scr C(S))$ is an equivalence.
Both $\NAlg(\scr C(S))$ and $\NAlg^{ex}(\scr C(S))$ have convenient properties.
In order not to overburden notation, from now on we silently\todo{is this a good idea?} assume that they are equivalent (e.g., $\scr C$ satisfies Zariski descent).
\end{remark}
\begin{proof}[Proof of Lemma \ref{lemm:transported-normed-structure-technical}.]
As explained in Remark \ref{rmk:NAlg-alternative}, we can replace $\NAlg$ by $\NAlg^{ex}$.
With this, we can apply the machinery of \cite{LLP}.\footnote{The main point is that for $T \in \Sch_S^{fp}$, $\NAlg^{ex}(\scr C(T))$ is defined in terms of the slice category $(\Sch_{S}^{fp})_{/T} \wequi \Sch_T^{fp}$, whereas $\NAlg(\scr C(T))$ is defined in terms of $\Sm_T \ne (\Sm_S)_{/T}$.}
Specifically, Lemma 2.31 and Definition 2.28 (and Lemma 2.16(3)) of \emph{loc. cit.} supply our desired adjunction.

The case of cartesian norm structures is dealt with by passing to opposites.
\end{proof}

\begin{corollary} \label{cor:forgetful-normed}
Given a normed category $\scr C(S)$, the forgetful functor $\cNAlg(\scr C(S)) \to \scr C(S)$ upgrades to a normed functor, where the source is given the cartesian norm structure.
Similarly for $\NAlg(\scr C(S))$ and the cocartesian norm structure.
\end{corollary}
\begin{proof}
Lemma \ref{lemm:transported-normed-structure-technical} supplies a(n) (op)lax normed functor.
That it is strongly normed (i.e., the left adjoint to pullback on $\NAlg(\scr C(S))$ is given by the norm) follows from \cite[Theorem 8.5]{bachmann-norms}.
\end{proof}

Here is a situation where passage to (co)normed (co)algebras has no effect.
\begin{lemma} \label{lemm:cNAlg-cartesian-normed}
If $\scr C(S)$ is cartesian normed, then $\cNAlg(\scr C(S)) \to \scr C(S)$ is an equivalence.
Similarly if $\scr C(S)$ is cocartesian normed, then $\NAlg(\scr C(S)) \wequi \scr C(S)$.
\end{lemma}
\begin{proof}
We prove the $\NAlg$ case, $\cNAlg$ is dealt with similarly.
By Lemma \ref{lemm:transported-normed-structure-technical}, we can construct a functor $\scr C(S) \to \NAlg(\scr C(S))$ by instead exhibiting a lax normed functor $\scr C(S)^\amalg \to \scr C(S)$.
By assumption, we can choose for this the identity.
By construction, the composite $\scr C(S) \to \NAlg(\scr C(S)) \to \scr C(S)$ is the identity.
Conversely, to show that $\NAlg(\scr C(S)) \to \scr C(S) \to \NAlg(\scr C(S))$ is the identity, by adjunction we must show that the lax normed functor $\NAlg(\scr C(S))^\amalg \to \scr C(S)^\amalg \wequi \scr C(S)$ is the canonical forgetful functor (where the first functor is $(\ph)^\amalg$ applied to the forgetful functor, and the second one witnesses that $\scr C$ is cocartesian normed).
Since all three categories are cocartesian normed and the functors are normed, it suffices to check this ignoring the norms, i.e., restricting to functors on $\Sch_S^{fp}$ \cite[Theorem 3.3]{dauser2024uniqueness}.
But this holds by construction.
\end{proof}

\subsubsection{Coalgebras and bialgebras}
Now we establish some basic properties of colimits in the category of conormed coalgebras and co/limits in normed algebras.
\begin{lemma} \label{lemm:cNAlg-basics}
Let $\scr C$ be a normed category.
\begin{enumerate}
\item Assume that each $\scr C(X)$ has colimits of shape $K$, and that such are preserved by pullback along smooth morphisms.
  Then $\cNAlg(\scr C(S))$ has colimits of shape $K$, and their formation commutes with the forgetful functor to $\scr C(S)$.
\item Let $p: X \to Y$ be a morphism such that $p_\sharp$ exists in $\scr C$ and is compatible with base change.
  Then the forgetful functor $\cNAlg(\scr C(S)) \to \scr C(S)$ commutes with $p_\sharp$.
\item If $\scr C$ is presentably normed, then $\cNAlg(\scr C(S))$ is cartesian presentably normed.
\end{enumerate}
\end{lemma}
\begin{proof}
(1) The argument of \cite[Proposition 7.6(3)]{bachmann-norms} goes through unchanged.

(2) The argument of \cite[Theorem 8.2]{bachmann-norms} goes through unchanged.

(3) By (1), $\cNAlg(\scr C(S))$ has all small colimits, and pullbacks in $\cNAlg(\scr C(S))$ preserve them (because this is true for $\scr C$).
Similarly norms preserve sifted colimits in $\cNAlg(\scr C(S))$.
Moreover, by (1) and (2), the (natural) distributivity transformation commutes with the forgetful functor.
Consequently the distributivity law for $\cNAlg(\scr C(S))$ follows from that of $\scr C(S)$.
It thus remains to show that $\cNAlg(\scr C(S))$ is accessible.
This follows from \cite[Proposition 5.4.7.11]{lurie-htt} (applied with $\scr C$ the accessible $\infty$-categories).
\end{proof}

\begin{lemma} \label{lemm:NAlg-basics}
If $\scr C$ is presentably normed, then $\NAlg(\scr C(S))$ is presentable.
  The forgetful functor $\NAlg(\scr C(S)) \to \scr C(S)$ preserves limits and sifted colimits; the coproduct in $\NAlg(\scr C(S))$ is the tensor product.
\end{lemma}
\begin{proof}
Same argument as in \cite[Theorem 8.2]{bachmann-norms}.
\end{proof}

\begin{remark}
Beware that Lemmas \ref{lemm:cNAlg-basics} and \ref{lemm:NAlg-basics} cannot be straightforward duals of one another.
For example, $\NAlg(\scr C(S))$ typically is not cocartesian presentably normed (just like how the category of abelian groups \emph{with the direct sum operation} is not presentably symmetric monoidal).
\end{remark}

\begin{definition}
We put \[ \bNAlg(\scr C(S)) = \NMon(\cNAlg(\scr C(S))) \text{ and } \NHopf(\scr C(S)) = \bNAlg(\scr C(S))^\gp. \]
\end{definition}

Here is how colimits work in categories of normed bialgebras.
\begin{lemma} \label{lemm:colim-in-bNAlg}
If $\scr C$ is presentably normed then $\bNAlg(\scr C(S))$ admits colimits preserved by the forgetful functor $\bNAlg(\scr C(S)) \to \NAlg(\scr C(S))$.
\end{lemma}
\begin{proof}
Since $\cNAlg(\scr C(S))$ is presentably normed by Lemma \ref{lemm:cNAlg-basics}(3), we may apply Lemma \ref{lemm:NAlg-basics} to deduce that $\bNAlg(\scr C(S))$ is presentable.
Moreover, the coproduct in $\bNAlg(\scr C(S))$ is the tensor product in $\cNAlg(\scr C(S))$, i.e., in $\scr C(S)$, and sifted colimits in $\bNAlg(\scr C(S))$ are computed in $\cNAlg(\scr C(S))$, i.e., in $\scr C(S)$.
The same facts about coproducts and sifted colimits apply to $\NAlg(\scr C(S))$ (by the same reference).
We thus conclude since all colimits are built using finite coproducts and sifted colimits \cite[Lemma~2.8]{bachmann-norms}.
\end{proof}

Let us observe the following exchange property.
\begin{lemma} \label{lemm:bNAlg-exchange}
There is a canonical equivalence $\bNAlg(\scr C(S)) \wequi \bNAlg(\scr C(S)^\op)^\op$.
\end{lemma}
\begin{proof}
In other words, $\NMon(\cNAlg(\scr C(S))) \wequi \cNMon(\NAlg(\scr C(S)))$.
Let us write $\scr S \subset \Fun^\times((\Sch_S^{fp})^\op, \Cat_\infty)$ for the subcategory on those functors (``motivic categories'') admitting indexed limits and colimits (i.e. functors $p_\sharp$ and $p_*$ for $p$ finite étale) stable under base change, and which \emph{coincide} (i.e. $p_* \wequi p_\sharp$), and transformations (``motivic functors'') which commute with $p_\sharp$ and (equivalently) $p_*$.
Note that both $\NMon(\cNAlg(\scr C(S)))$ and $\cNMon(\NAlg(\scr C(S)))$ lie in $\scr S$, with $p_\sharp = p_\otimes = p_*$ (use Corollary \ref{cor:forgetful-normed} and Lemma \ref{lemm:cNAlg-basics}(2)).
We now show that the presheaves on $\scr S$ represented by the two categories agree.
Thus let $X \in \scr S$.
Using Lemma \ref{lemm:transported-normed-structure-technical} we see that $\Map_{\scr S}(X, \NMon(\cNAlg(\scr C(S))))$ is the same as normed functors $X^\amalg \to \cNAlg(\scr C(S))$.
But $X^\amalg \wequi X^\times$ (i.e. the cartesian and cocartesian norm structures agree), and so (reading the lemma in reverse) this is the same as functors $X \to \cNAlg(\cNAlg(\scr C(S)))$ preserving indexed limits.
By Lemma \ref{lemm:cNAlg-cartesian-normed}, $\cNAlg(\cNAlg(\scr C(S))) \wequi \cNAlg(\scr C(S))$ and so, applying Lemma \ref{lemm:transported-normed-structure-technical} one further time, we see that this is the same as normed functors $X^\oplus \to \scr C$.
Here we write $X^\oplus$ for the norm structure which is equivalently given by the cartesian or the cocartesian one, to emphasize the symmetry.
Essentially the same argument applies to $\cNMon(\NAlg(\scr C(S)))$.
\end{proof}

\subsubsection{Transported norm structures}
Finally we see that categories of normed monoids admit norms which arise as linearizations of the original ones.
\begin{lemma} \label{lemm:norms-for-nhopf}
If $\scr C(S)$ is cartesian presentably normed, then $\NMon(\scr C(S))$ admits a functorial (in $\scr C$) presentably normed structure making $\NSym: \cNAlg(\scr C(S)) \to \NMon(\scr C(S))$ normed.
\end{lemma}
\begin{proof}[Proof sketch.]\todo{...}
In order to not run into representability issues for Weil restrictions, in this proof we will only work with affine schemes.
As long as $\scr C$ satisfies Zariski descent, this is not an issue.
In general, the statement should be formulated slightly differently, but following the convention from Remark \ref{rmk:NAlg-alternative}, we ignore this issue.

Begin by choosing, for each $X\in\Sch_S^{fp}$, a small full subcategory $\scr C_0(X)\subset\scr C(X)$ containing all $\kappa$-compact objects for a sufficiently large regular cardinal $\kappa$ such that $\scr C(X)$ is $\kappa$-accessible.
Enlarge these subcategories so that they are closed under finite coproducts and the functors $f^*$, $p_\otimes=p_*$ ($p$ finite étale), and $f_\sharp$ ($f$ smooth).
Each $\scr C_0(X)$ is then dense in $\scr C(X)$ \cite[Proposition~5.3.5.11]{lurie-htt}.
Let $\scr E(X)$ be the Grothendieck construction of $\scr C_0|_{\Sm_X}$, i.e., the category of pairs $(Y \to X, E)$ with $Y \to X$ smooth and $E \in \scr C_0(Y)$.
This can be made functorial in $X$, where if $f: X' \to X$ is any morphism, then $f^*(Y, E) = (Y \times_X X', (f \times_X Y)^* E)$.
We obtain this way a functor $\scr E: \Sch_S^{fp,\op} \to \Cat_\infty$.
If $f$ is smooth, then $f^*$ (in $\scr E$) has a left adjoint given by $f_\#(Y, E) = (f_\sharp Y, E)$.
If $f$ is finite étale, there is also a right adjoint satisfying $f_*(Y, E) = (R_f Y, q_*e^* E)$.
Here $R_f Y \to X$ is the Weil restriction, and $Y \xleftarrow{e} R_f Y \times_X X' \xrightarrow{q} R_f Y$ is the exponential diagram \cite[\S5.2, diagram~(5.7)]{bachmann-norms}.
One verifies that $\scr E(S)$ is a cartesian normed category.
Passing to presheaves sectionwise, we define $\scr F(X) = \PSh_\Sigma(\scr E(X))$; this admits a cartesian norm structure and one verifies that $\scr E(S) \to \scr F(S)$ is normed.
The functor $\scr E(S) \to \scr C(S)$ sending $(p: Y \to S, E)$ to $p_\sharp E$ extends to a (cartesian) presentably normed functor $L: \scr F(S) \to \scr C(S)$.
($L$ being normed requires the distributivity law for $\scr C$, in the extended form for smooth morphisms.)
The functor $L$ admits a right adjoint given by \[ (RM)((Y, E)) = \Map_{\scr C(Y)}(E, M|_Y). \]
Note that the restriction of $R$ to $\scr C_0(S)$ is fully faithful: it is given by $RE = (S, E)$ and thus $LRE = E$.
Since $\scr C_0(S)$ is dense in $\scr C(S)$, the restricted Yoneda embedding $\scr C(S)\to\PSh_\Sigma(\scr C_0(S))$ is fully faithful.
Together with the previous sentence, this implies that $R$ is fully faithful, i.e., $L$ is a localization.
Let $W(X)$ be the class of morphisms inverted by $L(X)$.
Then by construction $W$ is stable under $f^*, f_\sharp$ (for $f$ smooth) and $f_*$ (for $f$ finite étale).
One verifies from this that we can obtain $\NMon(\scr C(S))$ by localizing $\NMon(\scr F(S))$ at the free monoids on the maps in $W$.

Now we come to the actual construction.
Call $(Y, E) \to (Z, F) \in \scr E(X)$ \emph{admissible} if $Y \to Z$ is finite étale and $E \to (Y \to Z)^* F$ is an equivalence.
One verifies that pullbacks of admissible maps exist and are admissible.
Consequently \[ \scr K(X) = \Span(\scr E(X), \mathrm{admissible}, \all) \] is defined.
Arguing as in \cite[Lemma 16.17, Remark 16.18]{bachmann-norms} we see that $\PSh_\Sigma(\scr K(S)) \wequi \NMon(\scr F(S))$, functorially in $S$.
The explicit description of the norms on $\scr E(S)$ shows that admissible maps are stable under norms.
Thus applying the span construction pointwise, we see that $\scr K(S)$ upgrades to a  normed category.
One verifies that the localization \[ \PSh_\Sigma(\scr K(S)) \wequi \NMon(\scr F(S)) \to \NMon(\scr C(S)) \] is compatible with the norms.
All claims follow.
\end{proof}

\begin{remark}
Consider the case $\scr C(S) = \Spc(S)$.
In the notation of the proof, let $\scr E'(X) \subset \scr E(X)$ be the subcategory where $E=*$ is always required to be the terminal object, and $\scr K'$ is defined in terms of $\scr E'$, not $\scr E$.
Then $\scr K'(S) = \Span(\Sm_S, \fet, \all)$ and we consequently obtain a normed functor $\PSh_\Sigma(\Span(\Sm_S, \fet, \all)) \to \PSh_\Sigma(\scr K(S)) \to \NMon(\Spc(S))$.
This is a normed upgrade of the usual identification, obtained from \cite[Remark 16.18]{bachmann-norms}.
This shows that the norm structure constructed by Lemma \ref{lemm:norms-for-nhopf} on $\NMon(\Spc(S))$ is the same as the one coming from the identification with $\Spc^\fet(S)$.
\end{remark}

\begin{remark} \label{rmk:smash-product-normed-str}
Let $\scr C(S)$ be cartesian presentably normed.
Then we can produce on $\scr C(S)_*$ a presentably normed structure (the smash product), and a factorization of normed functors $\NSym: \scr C(S) \xrightarrow{(\ph)_+} \scr C(S)_* \to \NMon(\scr C(S))$.
The argument is a minor modification of the proof of Lemma \ref{lemm:norms-for-nhopf}.
\end{remark}

The following Example is used all throughout the rest of this article.
\begin{example} \label{ex:homology-hopf-algebra}
Let $A \in \NAlg(\SH(k))$.
Applying the functor $\NHopf$ to the morphism of presentably normed categories $A[\ph]: \Spc(k) \to \Mod_A(k)$ we obtain normed functors \[ A[\ph]: \Spc(k) \wequi \cNAlg(\Spc(k)) \to \cNAlg(\Mod_A(k)) \] and \[ A[\ph]: \NMon(\Spc(k))^\gp \wequi \NHopf(\Spc(k)) \to \NHopf(\Mod_A(k)). \]
\end{example}

\subsection{Rational normed Hopf algebras}
In the presence of étale descent, normed versions of algebraic structures collapse to classical ones.
\begin{lemma} \label{lemm:Hopf-FEt}
Let $\scr C(S)$ be a normed category satisfying finite étale descent.
Then so do $\scr C^\op$ and $\NAlg(\scr C(S))$.

Moreover the canonical forgetful functors \begin{gather*} \NAlg(\scr C(S)) \to \CAlg(\scr C(S)), \quad \cNAlg(\scr C(S)) \to \cCAlg(\scr C(S)), \\ \bNAlg(\scr C(S)) \to \mathrm{bCAlg}(\scr C(S)) := \CMon(\cCAlg(\scr C(S))), \\ \text{and } \NHopf(\scr C) \to \CHopf(\scr C(S)) := \mathrm{bCAlg}(\scr C(S))^\gp \end{gather*} are equivalences.
\end{lemma}
\begin{proof}
The statement about forgetful functors follows from étale descent claims, by \cite[Corollary C.16]{bachmann-norms}.
The claim about $\scr C^\op$ is clear because passing to opposites is an automorphism of $\Cat_\infty$.
Using that $\NAlg(\scr C(S)) \wequi \CAlg(\scr C(S))$, the claim about $\NAlg(\scr C(S))$ reduces to showing that $\CAlg(\ph)$ preserves limits (whence sheaves), which holds because it is corepresentable \cite[Proposition 2.2.4.9]{lurie-ha}.
\end{proof}

\begin{corollary}
We have \[ \NHopf(\Mod_{H\Q}(k)) \wequi \CHopf(\Mod_{H\Q}(k)). \]
\end{corollary}
\begin{proof}
Immediate from Lemma \ref{lemm:Hopf-FEt}, using that $\Mod_{H\Q}(\ph)$ has étale descent \cite[Theorem~14.3.4(1)]{triangulated-mixed-motives}.
\end{proof}

\subsection{Duality} \label{subsec:duality}
For a scheme $S$ and $A \in \NAlg(\SH(S))$, write \[ \Mod_A(S)^\fet \subset \Mod_A(S) \] for the full subcategory generated under retracts, finite colimits and desuspensions by objects of the form $\Sigma^{2n,n}A[X]$, where $X/S$ is a finite étale scheme and $n \in \Z$.

We also want to define somewhat larger versions.
For this, let $\Mod_A(S)^\eff \subset \Mod_A(S)$ be the subcategory generated under colimits and desuspensions by $A[X]$ for $X/S$ smooth, and write $\Mod_A(S)^{\eff\perp}$ for its right orthogonal.
\begin{definition} \label{def:quick-conv}
We call an $\N$-indexed colimit \[ X_1 \to X_2 \to \dots \in \Mod_A(S) \] \emph{quickly converging} (in the $+$ or $-$ direction, to be precise) if for every integer $N$, all but finitely many of the cofibers $X_{i+1}/X_i$ lie in $\Sigma^{2N,N} \Mod_A(S)^\eff$ or $\Sigma^{2N,N} \Mod_A(S)^{\eff\perp}$.
We also apply this definition in other categories with a forgetful functor to $\Mod_A(S)$; we mean then an $\N$-indexed diagram such that the underlying diagram of modules is quickly converging.
We call an infinite sum $\bigoplus_i X_i \in \Mod_A(S)$ quickly converging if the $\N$-indexed diagram of partial sums is.
We define quickly converging limits and products similarly.
\end{definition}

\begin{definition}
We call an $\N$-indexed diagram $X_1 \to X_2 \to \dots$ split\footnote{perhaps splittable would be more appropriate} if there exist retractions $X_{i+1} \to X_i$ of $X_i \to X_{i+1}$ for each $i$, and a corresponding colimit \emph{split $\N$-indexed colimit}.
We proceed similarly for limits.
\end{definition}
\begin{remark} \label{rmk:split-immediate}
Let us point out the following immediate facts.
\begin{enumerate}
\item In a stable, idempotent complete category, any split $\N$-indexed colimit is equivalent to an $\N$-indexed sum.
  Similarly for split limits.
\item Split $\N$-indexed diagrams are preserved by any functor.
  The same holds for split $\N^{op}$-indexed diagrams.
\end{enumerate}
\end{remark}

With this notion in place, we define subcategories \[ \Mod_A(S)^{\fet+}, \Mod_A(S)^{\fet-} \subset \Mod_A(S) \] which are obtained from $\Mod_A(S)^{\fet}$ by additionally closing under quickly converging infinite sums (in the $+$ respectively $-$ direction).

\begin{definition} \label{def:wt-bounded}
We call $A$ \emph{weight-bounded}\NB{``and coeffective''?} if $A \in \Sigma^{2,1}\SH(S)^{\eff\perp}$ and, for every smooth scheme $X/S$ there exists an integer $N_X$ such that for any $Y/X$ finite étale and $M > N_X$ we have $\map(Y, \Sigma^{2M,M}A) = 0$.\NB{not minimal assumptions}
\end{definition}

\begin{example} \label{ex:wt-bounded-crit}
Let $S=\Spec(k)$, $k$ some perfect field.
By considering the descent spectral sequence and using unramifiedness of homotopy sheaves \cite[Chapter~2]{A1-alg-top}, we see\NB{details?} that $A$ is weight-bounded if, and only if, for every finitely generated field $K/k$ of transcendence degree $d$, there is an integer $b_d$ such that $\ul\pi_*(A)_j(K) = 0$ for $j \not\in[0,b_d]$.
\end{example}

\begin{example} \label{ex:charp-wt-bounded}
Let $S$ be the spectrum of a field $k$ of characteristic $p$ with $[k:k^p] < \infty$.
Then $H\F_p$ is weight-bounded by the Geisser--Levine computation of motivic cohomology in terms of logarithmic differential forms \cite[Theorem~1.1 and Proposition~3.1]{geisser2000k}:
in the notation of Example \ref{ex:wt-bounded-crit}, we may take $b_d = \log_p[k:k^p] + d$.
The same argument, using logarithmic de Rham--Witt forms, applies to $H\Z/p^\nu$.
\end{example}


\begin{example}
If $S$ is the spectrum of a field of characteristic $\ne p$ and with $cd_p(k) < \infty$ and containing a primitive $p$-th root, then $H\F_p/\tau$ is weight-bounded.
(Use Bloch--Kato.)
\end{example}

Here are some straightforward properties of quickly converging sums.
\begin{lemma} \label{lemm:quickly-converging-sum}
Let $A$ be weight-bounded.
\begin{enumerate}
\item The category $\Mod_A(S)^{\fet+} \cap \Sigma^{2N,N} \Mod_A(S)^{\eff}$ is generated (under finite colimits, desuspensions, and quickly converging sums) by objects of the form $\Sigma^{2i,i}A[Z]$ with $i \ge N$ and $Z/S$ finite étale.
\item Similarly $\Mod_A(S)^{\fet-} \cap \Sigma^{2N,N} \Mod_A(S)^{\eff\perp}$ is generated by $\Sigma^{2i,i}A[Z]$ with $i<N$.
\item Let $X \in \Sm_S$.
  For $Y/X$ finite étale and $E \in \Mod_A(S)^{\fet+} \cap \Sigma^{2N_X+2,N_X+1}\Mod_A(S)^{\eff}$ we have $\map(Y, E) = 0$.
\item For $\bigoplus_i E_i$ a quickly converging sum in $\Mod_A(S)^{\fet\pm}$, the canonical map $\bigoplus_i E_i \to \prod_i E_i$ is an equivalence.
\item For $E \in \Mod_A(S)^{\fet+} \cap \Sigma^{2N,N} \Mod_A(S)^{\eff}$ we have $DE \in \Mod_A(S)^{\fet-} \cap \Sigma^{2-2N,1-N} \Mod_A(S)^{\eff\perp}$, and vice versa, establishing an anti-equivalence between the categories.
\item For $E \in \Mod_A(S)^{\fet\pm}$ we have $DE \in \Mod_A(S)^{\fet\mp}$, establishing an anti-equivalence between the categories.
\end{enumerate}
\end{lemma}
\begin{proof}
(1)
We use the $A$-effective cover functor $f_n^A$ right adjoint to the inclusion $\Sigma^{2n,n}\Mod_A(S)^\eff \to \Mod_A(S)$.
Since $A$ is coeffective (by assumption) and motives of finite étale schemes are self-dual \cite[Proposition~2.4.31]{triangulated-mixed-motives}, we see that for $Z/S$ finite étale we have \[ f_n(A[Z]) = \begin{cases} 0 & \text{ for } n > 0 \\  A[Z] & \text{ for } n \le 0 \end{cases}. \]
Now if $E \in \Mod_A(S)^{\fet+} \cap \Sigma^{2N,N} \Mod_A(S)^{\eff}$ then $E \wequi f_N^A E$, so it suffices to show that $f_N^A$ takes $\Mod_A(S)^{\fet+}$ into $\Mod_A(S)^{\fet+} \cap \Sigma^{2N,N} \Mod_A(S)^{\eff}$.
Since $f_n^A$ preserves colimits and desuspensions, this reduces to the case $E = \Sigma^{2n,n} A[Z]$, which was treated above.

(2)
We similarly use the $A$-coeffective cover functor $f^n_A = \id/f_n^A$.
The proof is the same.

(3) As a functor of $E$, the expression preserves colimits and desuspensions, and so we reduce via (1) to the case $E=\Sigma^{2i,i} A[Z]$ for $Z/S$ finite étale and $i \ge N_X+1$.
Using self-duality of $A[Z]$ once more we may assume that $Z=*$.
The desired vanishing now holds by assumption of weight-boundedness of $A$.

(4) It suffices to show that, for $X/S$ smooth and $N \in \Z$ we have $\map(\Sigma^{2N,N}A[X], E_i) = 0$ for all but finitely many $i$.
By definition of quick convergence (Definition \ref{def:quick-conv}), for any $M$, all but finitely many $E_i$ lie in $\Sigma^{2M,M}\Mod_A(S)^\eff$ (in the $+$ case) or $\Sigma^{2M,M}\Mod_A(S)^{\eff\perp}$ (in the $-$ case).
In the $-$ case the vanishing thus holds by definition, and the $+$ case is treated in (3).

(5, 6) Write $\scr C \subset \Mod_A(S)^{\fet+}$ for the subcategory of objects $X$ such that $D(f_n^A X) \in \Mod_A(S)^{\fet-} \cap \Sigma^{2-2n,1-n} \Mod_A(S)^{\eff\perp}$.
This category is closed under desuspensions and finite colimits.
Note that if $X \in \scr C$ is $n$-effective, then $DX \in \Mod_A(S)^{\fet-} \cap \Sigma^{2-2n,1-n} \Mod_A(S)^{\eff\perp}$ (since $X \wequi f_n^A X$).
It follows that if $X = \bigoplus_i X_i$ is a quickly converging sum with $X_i \in \scr C$, then $DX$ is a quickly converging product, whence sum by (4).
The same argument applies to $f_n^A X$, showing that $\scr C$ is closed under quickly converging sums.
By the arguments from (1), $\Sigma^{2n,n} A[Z] \in \scr C$ for any $Z/S$ finite étale and $n \in \Z$.
Consequently $\Mod_A(S)^{\fet+} \subset \scr C$.
This proves the first halves of the containment assertions of (5) and (6).

The other halves of the containments are proved similarly, using $f^n_A$ in place of $f_n^A$.

To prove that we have equivalences, we need to show each time that $DDX \wequi X$.
Letting $\scr C$ be the subcategory of objects $X$ such that $DDf_n X \wequi f_nX$ for all $n$ (in the $+$ situation, with $f^n$ in the $-$ situation), we see again that $\scr C$ is closed under quickly converging sums.
It is also easily seen to be closed under desuspensions and finite colimits and to contain the generators, so these closure properties imply that it is the whole category.
This concludes the proof.
\end{proof}

This is the duality result we are really after.
\begin{proposition} \label{prop:duality}
Suppose that $A$ is orientable.

The category $\Mod_{A}(S)^{\fet}$ is closed under norms.
Passing to duals induces an autoequivalence (of normed categories on $\Sm_S$) of $\Mod_A(S)^{\fet}$.

If $A$ is in addition weight-bounded, then also $\Mod_{A}(S)^{\fet\pm}$ are closed under norms and passing to duals induces an equivalence of normed categories $\Mod_A(S)^{\fet+\op} \wequi \Mod_A(S)^{\fet-}$.
\end{proposition}
\begin{proof}
\textbf{Step 1.} We show that if $p: X \to Y \in \Sm_S$ is finite étale, then $p_\sharp$ preserves the categories $\Mod_A(\ph)^{\fet(\pm)}$.
Since $p_\sharp$ preserves colimits and desuspensions and, for $Z \to X$ finite étale, $p_\sharp(\Sigma^{2n,n} A_X[Z]) \wequi \Sigma^{2n,n}A_Y[p_\sharp Z]$, it will suffice to show that $p_\sharp$ preserves $n$-effective and $n$-coeffective modules.
For effectivity, this holds essentially by definition.
For coeffectivity, using that $p_\sharp \wequi p_*$ \cite[Theorem~2.4.50(2), (3)]{triangulated-mixed-motives}, it suffices to show that $p^*$ preserves $n$-effective objects, which is again clear.

\textbf{Step 2.} We show that $\Mod_A(\ph)^\fet$ is preserved by functors of the form $p_\otimes$.
We use induction on the degree of $p$.
If $p$ is a clopen immersion, the claim is clear.
Thus we may as well assume $p$ to be surjective and of constant degree.
By the distributivity law in the form \cite[Corollary 5.13]{bachmann-norms}, step 1 and the induction, if $M, N \in \Mod_A(X)^\fet$ with $p_\otimes(M), p_\otimes(N) \in \Mod_A(Y)^\fet$, then also $p_\otimes(M \oplus N) \in \Mod_A(Y)^\fet$.
From this we deduce that given a map $\alpha: M \to N$, $p_\otimes(\cof(\alpha)) \in \Mod_A(Y)^\fet$.
Indeed, we can write $\cof(\alpha)$ as the geometric realization of a simplicial objects $S_\bullet = N \oplus M^{\oplus \bullet}$.
This is a \emph{1-skeletal simplicial object}, and so, if $p$ has degree $d$, then $p_\otimes(S_\bullet)$ is $d$-skeletal by \cite[Corollary 5.26]{bachmann-norms} and \cite[Proposition 2.10]{k-theory-and-poly-functors}.
In particular, $p_\otimes(\cof \alpha) \wequi |p_\otimes(S_\bullet)| \in \Mod_A(Y)^\fet$, the geometric realization being a finite colimit.
Consequently the subcategory $\scr C$ of $\Mod_A(X)^\fet$ sent by $p_\otimes$ into $\Mod_A(Y)^\fet$ is closed under finite colimits.
It contains $\Sigma^{2n,n}A_X[Z]$ (since $p_\otimes(\Sigma^{2n,n}A_X[Z]) \wequi \Sigma^{2nd,nd}A_Y[R_p Z]$, by orientability and \cite[Remark 4.9]{bachmann-norms}), so in particular $\1$, and hence (by closure under finite colimits) $S^1$.
Since $p_\otimes$ preserves tensor products and strong duals (being symmetric monoidal) and $\Mod_A(Y)^\fet$ is closed under them, we see that $\scr C$ is closed under these operations as well.
In particular, $S^{-1} \in \scr C$ and $\scr C$ is closed under desuspensions.
Since closure under retracts is clear, we deduce that $\scr C = \Mod_A(X)^\fet$.

\textbf{Step 3.} We now show that, in any normed category with strong duals, taking duals upgrades to a normed functor.
For this, we observe that we have a functor (see Step 3' for details on the construction) \[ d: \CAlg(\Cat_\infty)^{dl} \to \Fun(\Delta^1, \Cat_\infty), \scr C \mapsto (D: \scr C^\op \to \scr C). \]
Here $\CAlg(\Cat_\infty)^{dl}$ denotes the category of symmetric monoidal $\infty$-categories in which every object has a strong dual (and symmetric monoidal functors).
The point is that any symmetric monoidal functor automatically preserves duals; see \cite[\S4.6.1]{lurie-ha} for details.
We can view a normed category with strong duals as a functor $\scr C: \Span(\Sm_S, \all, \fet) \to \CAlg(\Cat_\infty)^{dl}$ (use \cite[Proposition C.9]{bachmann-norms}).
Composing with $d$ constructs a normed functor $D: \scr C^\op \to \scr C$ as desired.

\textbf{Step 3'.} More generally, in any normed category with weak duals, taking weak duals upgrades to a lax normed functor.
To see this, let $\operatorname{Pair}(\scr C)$ denote the $\infty$-category of triples $(X,Y,e:X\otimes Y\to\1)$, with morphisms given by maps $u:X'\to X$ and $v:Y\to Y'$ together with a homotopy
\[ e\circ(u\otimes\id) \wequi e'\circ(\id\otimes v). \]
Formally, $\operatorname{Pair}(\scr C)$ is the cocartesian Grothendieck construction of the functor $\scr C^\op\to\Cat_\infty$ sending $X$ to the total $\infty$-category of the right fibration over $\scr C$ classified by the presheaf $Y\mapsto\Map_{\scr C}(X\otimes Y,\1)$; a map $u:X'\to X$ acts by precomposition with $u\otimes\id$.
There are projections
\[ \scr C^\op \xleftarrow{p_{\scr C}} \operatorname{Pair}(\scr C) \xrightarrow{q_{\scr C}} \scr C \]
selecting the first and second entries.
The universal property of weak duals gives a right adjoint $r_{\scr C}(X)=(X,DX,\mathrm{ev}_X)$ to $p_{\scr C}$.
Moreover, $q_{\scr C}r_{\scr C}\wequi D_{\scr C}$.
Since $\operatorname{Pair}(\scr C)$, $p_\scr C$ and $q_\scr C$ are natural in the symmetric monoidal category $\scr C$, the calculus of mates from \cite[Theorem~5.3.5]{hhln} applies to deduce that $r_{\scr C}$ and hence $D$ are lax natural in $\scr C$, as was to be shown.

\textbf{Step 4.} We have now shown that $\Mod_{A}(S)^{\fet}$ is a normed category and $D: \Mod_{A}(S)^{\fet\op} \to \Mod_{A}(S)^{\fet}$ is a normed functor.
It remains to verify that $D$ is an equivalence.
This no longer has anything to do with norms, and is in fact true in any category with strong duals (see, e.g., \cite[\S4.6.1]{lurie-ha}).

\textbf{Step 5.} We show that, if $\bigoplus_i M_i \in \Mod_A(X)^{\fet+}$ is a quickly converging sum, then $\colim_i p_\otimes(\bigoplus_{j \le i} M_j)$ is a quickly converging colimit.
We may assume that $\bigoplus_i M_i$ is $n_0$-effective, for some $n_0 \le 0$.
Again we may as well assume that $p$ is surjective of degree $|p|$.
We fix $i$ and set $N = \bigoplus_{j < i} M_j$, $N' = M_i$.
We assume that $N'$ is $n$-effective for some $n \ge 0$.
We shall show that $C := \cof(p_\otimes(N) \to p_\otimes(N \oplus N'))$ is $(n+n_0|p|)$-effective; this implies the quick convergence.
We again use the distributivity law in the form \cite[Corollary 5.13]{bachmann-norms}, which shows that \[ C \wequi c_\sharp(q_{l\otimes}e_l^*(N) \otimes q_{r\otimes}e_r^*(N')) \oplus p_\otimes(N'), \] for certain finite étale maps $c, e_l, q_l, e_r, q_r$ of degrees $\le |p|$; here $q_l$ and $q_r$ are surjective.
Functors of the form $c_\sharp$ and $e^*$ preserve effectivity, whereas functors of the form $q_\otimes$ send $r$-effective objects to $r|q|$-effective objects (since they preserve effective objects and send $\Sigma^{2r,r}\1$ to $\Sigma^{2r|q|,r|q|}\1$ by orientability and \cite[Remark 4.9]{bachmann-norms}).
The claim follows from this.

\textbf{Step 6.} We show that $\Mod_A(\ph)^{\fet+}$ is preserved by $p_\otimes$.
Let $\scr C \subset \Mod_A(X)^{\fet+}$ be the subcategory of objects mapped by $p_\otimes$ into $\Mod_A(Y)^{\fet+}$.
By arguing as in step 2, $\scr C$ is closed under finite colimits.
It is also closed under tensor products and contains $\Mod_A(X)^{\fet}$, so is closed under desuspensions.
Noting that $p_\otimes$ preserves split diagrams, and (quickly converging) split colimits are equivalent to (quickly converging) sums (see Remark \ref{rmk:split-immediate}), step 5 shows that $\scr C$ is closed under quickly converging sums.
It follows that $\scr C = \Mod_A(X)^{\fet+}$, as desired.

\textbf{Step 7.} We show that, for every $M \in  \Mod_A(X)^{\fet+}$, $p_\otimes(DM) \wequi Dp_\otimes(M)$.
To do so, let $\scr C$ be the subcategory of all such $M$.
It is closed under tensor products and contains $\Mod_A(X)^{\fet}$ by step 4, so is closed under desuspensions.
We next show that it is closed under finite sums.
Indeed, using once more the distributivity law, this reduces to showing that functors of the form $c_\sharp$ preserve duals.
This follows from the equivalences $Dc_\sharp \wequi c_* D \wequi c_\sharp D$ \cite[Theorem~2.4.50(2), (3), (5)]{triangulated-mixed-motives}.
From this we deduce that $\scr C$ is closed under cofibers.
Indeed, arguing as in step 2, we reduce to showing that $p_\otimes$ commutes with coskeletal totalizations.
By \cite[Proposition 2.15]{k-theory-and-poly-functors}, this will hold if $p_\otimes$ is excisive, which it is by \cite[Corollary 5.26]{bachmann-norms}.
It follows that $\scr C$ is closed under finite colimits.
We next show that $\scr C$ is closed under quickly converging sums.
As a preliminary observation, for any quickly converging, split $\N$-indexed system, the colimit along the forward maps coincides with the limit along the chosen retractions, by Lemma \ref{lemm:quickly-converging-sum}(4).
Let $M = \bigoplus_i M_i$ be such a sum, and $N_\bullet$ the split $\N$-indexed system of partial sums.
Then $p_\otimes(N_\bullet)$ is a quickly converging, split $\N$-indexed system, by step 5.
Consequently $Dp_\otimes(N_\bullet)$ is a quickly converging, split inverse system, with forward sections obtained by dualizing the retractions, by Lemma \ref{lemm:quickly-converging-sum}(5).
Since $\scr C$ is closed under finite sums, each term $N_\bullet$ lies in $\scr C$, and so $Dp_\otimes(N_\bullet) \wequi p_\otimes(DN_\bullet)$.
On the other hand $DM$ is the limit of $DN_\bullet$ along the duals of the inclusions, which coincides with the colimit along the duals of the retractions by the preliminary observation.
Thus $Dp_\otimes(M)$ is the limit of $p_\otimes(DN_\bullet)$ along the former maps, whereas $p_\otimes(DM)$ is its colimit along the latter maps.
They coincide by the preliminary observation.
It follows that $\scr C = \Mod_A(X)^{\fet+}$, as desired.

\textbf{Conclusion.}
Since every object of $\Mod_A(X)^{\fet-}$ is of the form $DM$ for some $M \in \Mod_A(X)^{\fet+}$ (Lemma \ref{lemm:quickly-converging-sum}(6)), steps 6 and 7 show that $\Mod_A(X)^{\fet-}$ is also preserved by $p_\otimes$.
It follows that both $\Mod_A(X)^{\fet\pm}$ are normed categories.
Using step 3', we upgrade $D: \Mod_A(S)^{\fet+\op} \to  \Mod_A(S)^{\fet-}$ to a lax normed functor.
It is an equivalence by Lemma \ref{lemm:quickly-converging-sum}(6), and strongly normed by step 7.
This was to be shown.
\end{proof}

\begin{lemma} \label{lemm:duals-cNAlg}
Let $A$ be orientable and weight-bounded.
Passing to duals induces an equivalence \[ \cNAlg(\Mod_A(S)^{\fet\pm}) \wequi \NAlg(\Mod_A(S)^{\fet\mp})^\op \] as well as \[ \bNAlg(\Mod_{A}(S)^{\fet\pm}) \wequi \bNAlg(\Mod_{A}(S)^{\fet\mp})^\op. \]
\end{lemma}
\begin{proof}
Immediate from Proposition \ref{prop:duality} and, for the statement about bialgebras, Lemma \ref{lemm:bNAlg-exchange}.
\end{proof}

\section{Normed monoids} \label{sec:normed-monoids}
Arguably the simplest category to apply the constructions of \S\ref{sec:norm-hopf-1} to is the category of motivic spaces $\Spc(k)$ itself.
This category being cartesian normed, normed Hopf algebras in $\Spc(k)$ are just grouplike \emph{normed monoids}.
Our goal in this section is to study the category $\NMon(\Spc(k))$ of normed monoids, and establish a variety of properties we will be relying on in the sequel.

\subsection{Spaces with finite étale transfers} \label{subsec:Spc-fet}
The category of motivic spaces with finite étale transfers $\Spc^\fet(k)$ is obtained by starting with $\Span(\Sm_k, \fet, \all)$ and applying the usual ``motivic machinery'' \cite[Definition~3.22]{bachmann-MGM}.
It turns out that $\NMon(\Spc(k)) \wequi \Spc^\fet(k)$ (see \cite[Remark 16.18]{bachmann-norms}), and so we can translate results about the latter category into the former.
For example, the motivic categories $\Spc(\ph)_*$ and $\Spc^\fet(\ph) \wequi \NMon(\Spc(\ph))$ admit norm structures making \[ \Spc(\ph) \xrightarrow{(\ph)_+} \Spc(\ph)_* \xrightarrow{T} \NMon(\Spc(\ph)) \] normed functors \cite[\S3.3.3, (4)]{bachmann-MGM}.
We denote the corresponding extended powers by $D_n^\wedge$ and $D_n^\otimes$, in contrast with the extended powers on $\Spc(\ph)$ denoted $D_n$.

\subsection{$p$-localization} \label{subsub:mon-ploc}
We write $\NMon(\Spc(k))_{(p)}$ for the subcategory of $p$-local normed monoids, i.e., those such that multiplication by numbers coprime to $p$ is an isomorphism.
The inclusion has a left adjoint, computed as usual as the colimit over multiplication by numbers coprime to $p$.
If $X \in \NMon(\Spc(k))$ is connected, then the $p$-localization map $X \to X_{(p)} \in \NMon(\Spc(k))$ is also a $p$-localization of the underlying motivic space, in the sense of \S\ref{subsec:p-loc}.
Indeed, by Eckmann--Hilton, multiplication by $n$ on $X$ induces multiplication by $n$ on $\ul\pi_*(X)$, and so $X \to X_{(p)}$ localizes the homotopy sheaves, as needed.

\begin{lemma} \label{lemm:Dnotimes-plocal}
The functors $D_n^\otimes$ (for $n \ge 1$) commute with the localization at any set of primes; in particular, they preserve the subcategory $\NMon(\Spc(k))_{(p)}^{\gp}$.
\end{lemma}
\begin{proof}
Rewriting the extended power using norms and applying Grothendieck's Galois theory \cite[\S10.2]{bachmann-norms}, we reduce to a statement in equivariant homotopy theory.
Using the theory of spectral Mackey functors \cite[Corollary~A.4.1]{nardin2016stability}, we translate from grouplike normed monoids to connective genuine $G$-spectra, and must now prove: given a finite group $G$ and a subgroup $H$ the functor $N^H_G$ preserves $p$-equivalences and $p$-local objects.
Since these properties can be checked on geometric fixed points (the latter forming a conservative family of cocontinuous functors \cite[Proposition~2.52]{HHR}), this follows from \cite[Proposition~2.57]{HHR}.
\end{proof}

\subsection{Stabilization}
Let $k$ be a perfect field.
Write $\NMon(\Spc(k))^{\gp\A^1} \subset \NMon(\Spc(k))$ for the category of normed monoids which are grouplike (i.e., $\ul\pi_0$ is a sheaf of groups, not just monoids) and such that $\ul\pi_0$ is strongly $\A^1$-invariant (in the sense of \cite{A1-alg-top}).

\begin{example} \label{ex:connected-gpa1}
Any connected normed monoid lies in $\NMon(\Spc(k))^{\gp\A^1}$.
\end{example}

The category $\SH^{S^1}(k) = \Sp(\Spc(k))$ admits a cartesian norm structure (see \S\ref{subsub:cocart-normed}).
The forgetful functor $\NMon(\SH^{S^1}(k)) \to \SH^{S^1}(k)$ preserves limits and sifted colimits \cite[Proposition~7.6(2) and Remark~7.7]{bachmann-norms}.
It also sends finite coproducts to cartesian products \cite[Remark~7.7]{bachmann-norms} \cite[Proposition~3.2.4.7]{lurie-ha}, i.e., coproducts, and so preserves all colimits.
In particular, $\NMon(\SH^{S^1}(k))$ is stable.

The functor $\Omega^\infty: \SH^{S^1}(k) \to \Spc(k)$ commutes with cartesian norms (being right adjoint to $\Sigma^\infty_+$, which commutes with pullbacks) and hence upgrades to a normed functor (see Lemma \ref{lemm:transported-normed-structure-technical} for details).
There is thus an induced functor $\Omega^\infty: \NMon(\SH^{S^1}(k)) \to \NMon(\Spc(k))$.

\begin{lemma} \label{lemm:NMon-prestable}
The functor $\Omega^\infty$ restricts to an equivalence between $\NMon(\SH^{S^1}(k)) \times_{\SH^{S^1}(k)} \SH^{S^1}(k)_{\ge 0}$ and $\NMon(\Spc(k))^{\gp\A^1}$.
The former subcategory is the nonnegative part of a $t$-structure on $\NMon(\SH^{S^1}(k))$.
\end{lemma}
\begin{proof}
Define $\scr C(X) \subset \SH^{S^1}(X)$ as the subcategory on those objects $E$ such that $f^* E$ is connective for every $f: Y \to X \in \Sm_k$.
Then $\scr C(k)$ is a normed subcategory of $\SH^{S^1}(k)$.
(Since functors $p_*$ for $p$ finite étale are compatible with pullbacks, it suffices to show that given $p: Y \to Z \in \Sm_k$ finite étale and $E \in \SH^{S^1}(Y)$ connective, also $p_*(E)$ is connective.
This is true since a finite étale extension of a henselian local ring is a finite product of henselian local rings \cite[Tag~04GH]{stacks-project}.)
Applying $\CMon(\ph)$ sectionwise to the cartesian normed functor $\Omega^\infty$, we obtain \[ \tilde\Omega^\infty: \SH^{S^1}(k) \wequi \CMon(\SH^{S^1}(k)) \to \CMon(\Spc(k)). \]
By \cite[Proposition~3.1.9]{EHKSY}, $\tilde\Omega^\infty$ restricted to $\scr C$ is fully faithful.
Using that $\NMon(\CMon(\Spc(k))) \wequi \NMon(\Spc(k))$ \cite[Proposition C.10]{bachmann-norms}, we have thus constructed a fully faithful functor \[ \NMon(\SH^{S^1}(k)) \times_{\SH^{S^1}(k)} \SH^{S^1}(k)_{\ge 0} \to \NMon(\Spc(k)). \]
Its essential image consists of those normed monoids whose underlying motivic space is of the form $\Omega^\infty E$ for some $E \in \SH^{S^1}(k)_{\ge 0}$, i.e., those in $\NMon(\Spc(k))^{\gp\A^1}$ \cite[Corollary~3.1.15]{EHKSY}.

The claim about being the nonnegative part of a $t$-structure follows from \cite[Proposition~1.4.4.11]{lurie-ha}, since $\NMon(\SH^{S^1}(k)) \to \SH^{S^1}(k)$ preserves limits and colimits, and so the stated subcategory is closed under colimits and extensions.
\end{proof}

\begin{remark}
The composite $\NMon(\Spc(k))^{\gp\A^1} \hookrightarrow \NMon(\SH^{S^1}(k)) \to \SH^{S^1}(k)_{\ge 0}$ of the inverse of the equivalence from Lemma \ref{lemm:NMon-prestable} and the forgetful functor is obtained as follows: given $X \in \NMon(\Spc(k))^{\gp\A^1}$, the commutative monoid structure allows us to construct deloopings $B^n X$, which are still $\A^1$-invariant by the assumptions on $X$.
Together these deloopings define the connective spectrum we are after.
\end{remark}

\subsection{Free normed monoids} \label{subsub:free-normed}
\begin{definition}
We denote the left adjoint of the forgetful functor $\NMon(\Spc(k)) \to \Spc(k)_*$ by $T$.
\end{definition}

Let $T_{\mathrm{cl}}:\Spc_* \to \CMon(\Spc)$ denote the free commutative monoid functor on pointed spaces, where the basepoint is sent to the unit \cite[Remark~2.1.3.10 and Corollary~3.1.3.7]{lurie-ha}.
We denote its arity filtration \NB{ref?} by $F_* T_{\mathrm{cl}}(A)$.
For a pointed étale sheaf $A$, define
\[
T^\et(A) := L_\et\bigl(U \mapsto T_{\mathrm{cl}}(A(U))\bigr),
\qquad
F_nT^\et(A) := L_\et\bigl(U \mapsto F_nT_{\mathrm{cl}}(A(U))\bigr).
\]
These constructions are functorial in all pointed maps.
In particular, their restrictions to $\Sm_{k+}$ are functorial in arbitrary pointed morphisms $X_+ \to Y_+$, not only those induced by morphisms $X \to Y$.

We now describe their values on $X_+$.
For an étale sheaf $X$, put
\[
D_n^\et(X) = L_\et(X^{\times n}_{h\Sigma_n}).
\]
Then
\[
T^\et(X_+) \wequi \coprod_{n\ge 0}^{\Shv_\et} D_n^\et(X),
\qquad
F_nT^\et(X_+) \wequi \coprod_{0\le j\le n}^{\Shv_\et} D_j^\et(X),
\]
pointed by the degree-zero summand.
These decompositions by arity are natural in unpointed maps $X \to Y$; an arbitrary pointed map $X_+ \to Y_+$ preserves the filtration but may lower arity.
Here the coproducts are taken in étale sheaves.

For $X \in \Sm_k$, the étale-orbit description of motivic extended powers \cite[Proposition~5.23]{bachmann-colimits}, together with the description of free normed monoids \cite[Lemma~3.24 and Corollary~3.28]{bachmann-MGM}, gives
\[
T(X_+) \wequi L_\mot T^\et(X_+).
\]
(Indeed, arity is locally constant, so the canonical maps from the corresponding presheaf coproducts to the displayed étale-sheaf coproducts are Zariski-local equivalences.
Thus either coproduct convention gives the same result after motivic localization.)

The functor
\[
\Sm_{k+} \longrightarrow \Fun(\N,\Spc(k)_*),
\qquad
A \longmapsto L_\mot F_*T^\et(A),
\]
is defined on all pointed morphisms.
Extending it by sifted colimits over $\PSh_\Sigma(\Sm_{k+}) \wequi \PSh_\Sigma(\Sm_k)_*$, we obtain a functor
\[
F_*T(\ph): \PSh_\Sigma(\Sm_k)_* \longrightarrow \Fun(\N,\Spc(k)_*).
\]
It is easily seen to preserve $\A^1$-homotopies, and arguing as in \cite[Proposition~2.1.5(iii)]{EHKSY} one sees that it also sends Nisnevich covering sieves (generated by one map) to Nisnevich equivalences.\NB{details?}
It consequently preserves motivic equivalences and induces
\[
F_*T(\ph): \Spc(k)_* \longrightarrow \Fun(\N,\Spc(k)_*),
\]
i.e., an arity filtration on any free normed monoid.

\begin{lemma} \label{lemm:T-filtration}
For $n \ge 1$, we have
\[
F_nTX/F_{n-1}TX \wequi D_n^\wedge(X).
\]
\end{lemma}
\begin{proof}
Naturally in $X \in \Sm_{k+}$, the above discussion gives \[ F_nTX/F_{n-1}TX \wequi L_\mot F_nT^\et(X)/F_{n-1}T^\et(X) \wequi D_n^\wedge(X). \]
The claim follows by extending by sifted colimits.
\end{proof}

Let us also recall the following.
\begin{lemma} \label{lemm:M-of-T}
For $X \in \Spc(k)_*$ and $A \in \NAlg(\SH(k))$ we have $A[TX] \wequi \NSym_A(\overline{A[X]})$.
\end{lemma}
\begin{proof}
See \cite[Corollary 3.39]{bachmann-MGM}.
\end{proof}

\subsection{Effectivity}
Let us call $X \in \NMon(\Spc(k))$ \emph{very $n$-effective} or $n$-effective if this holds for the underlying motivic space, in the sense of Definition \ref{def:n-effective-space}.

\begin{example} \label{ex:K1-very-1-eff}
$K(\Z_{(p)}(1), 2) \in \NMon(\Spc(k))$ is very $1$-effective.
Indeed very $n$-effective spaces are closed under colimits (Remark \ref{rmk:effective-operations}), so we need only verify this for $K(\Z(1),2) \wequi \P^\infty$, for which see Example \ref{ex:Poo-effective}.
\end{example}

\begin{lemma} \label{lemm:Dn-effectivity-Spc*}
Let $X \in \Spc(k)_*$ be very $n$-effective.
\begin{enumerate}
\item For $m \ge 1$, $D_m^\wedge(X)$ is very $nm$-effective.
\item $T(X)$ is very $n$-effective.
\end{enumerate}
\end{lemma}
\begin{proof}
(1) Write $L_n: \Spc(k)_* \to \Spc(k)_*$ for the localization defining very $n$-effective spaces.
Since $L_n X = *$ and $X$ is connected, by \cite[Construction 3.1.17]{freudenthal} there exists a (transfinite) sequence of cell attachments nullifying $X$.
From this we deduce that it suffices to show the following: if $X$ is very $n$-effective, $Y \in \Spc(k)_*$ and $\alpha: \Sigma^{2n,n}Y \to X$ is any map, then $D_m^\wedge(X) \to D_m^\wedge(\cof(\alpha))$ is an $L_{nm}$-equivalence.
Writing the cofiber as the geometric realization of the simplicial object $[q] \mapsto X \vee (\Sigma^{2n,n}Y)^{\vee q}$ and using the formula for $D_m^\wedge(A \vee B)$ from \eqref{eq:formula-for-Dn-sum} for this it is enough to show that $D_a^\wedge(\Sigma^{2n,n}Y) \wedge D_b^\wedge(X)$ is very $mn$-effective, for any $a,b \ge 0$ with $a+b=m$ and $b<m$.
By induction, $D_b^\wedge(X)$ is very $bn$-effective, and so using \cite[Theorem 3.1.20]{freudenthal} it suffices to show that $D_a^\wedge(\Sigma^{2n,n}Y)$ is very $an$-effective.
Writing $D_a^\wedge$ as a colimit of norms (see \eqref{eq:formula-for-dn}), using that norms are symmetric monoidal, and $p_\otimes(S^{2n,n}) \wequi Th(V)$ where $V = p_* \scr O^n$ \cite[Lemma 4.4]{bachmann-norms} we reduce to the following: given $Y \in \Sm_k$, $A \in \Spc(Y)_*$, $V$ a vector bundle of rank $na$ on $Y$, then $f_\sharp(Th(V) \wedge A)$ is very $na$-effective.
The problem is local on $Y$, so we may assume $V$ trivial.
But then $f_\sharp(Th(V) \wedge A) \wequi S^{2na,na} \wedge f_\sharp(A)$ is indeed very $na$-effective.

(2) This follows immediately from (1) and Lemma \ref{lemm:T-filtration} (together with stability of very $n$-effective objects under cofiber extensions and filtered colimits, i.e., Remark \ref{rmk:effective-operations}).
\end{proof}

\begin{proposition} \label{prop:Dn-tens-eff}
Let $X \in \NMon(\Spc(k))$ be very $n$-effective.
Then $D_m^\otimes(X)$ is very $nm$-effective for $m \ge 1$.
\end{proposition}
\begin{proof}
Since $\NMon(\Spc(k)) \to \Spc(k)_*$ is conservative and preserves sifted colimits, any normed monoid can be resolved as a sifted colimit of iterated free normed monoids in the usual way \cite[Proposition~4.7.3.14]{lurie-ha}.
By Lemma \ref{lemm:Dn-effectivity-Spc*}, applying $T$ iteratively preserves very $n$-effectivity, and so we may assume that $X=TY$ for some $Y \in \Spc(k)_*$ which is very $n$-effective.
Then $D_m^\otimes(X) \wequi D_m^\otimes(TY) \wequi T(D_m^\wedge(Y))$ (the functor $T: \Spc(k)_* \to \NMon(\Spc(k))$ being normed, essentially by construction) which has the desired effectivity by Lemma \ref{lemm:Dn-effectivity-Spc*}.
\end{proof}

We will make use of the following technical result later.
\begin{lemma} \label{lemm:2-eff-trick}
Let $X \in \NMon(\Spc(k))$ and $\alpha: Y \to X \in \Spc(k)_*$, with $X$ $1$-effective and $Y$ very $1$-effective.
\begin{enumerate}
\item Suppose that $\cof(\alpha)$ (computed in $\Spc(k)_*$) is $2$-effective.
  Then also the cofiber of $T\alpha: TY \to X \in \NMon(\Spc(k))$ is $2$-effective.
\item Suppose that $X$ is $p$-local and $\cof(\alpha)_{(p)}$ is simply-connected and $2$-effective.
  Then $\cof(T(Y)_{(p)} \to X) \in \NMon(\Spc(k))$ is $2$-effective.
\end{enumerate}
\end{lemma}
\begin{proof}
(1) Consider the commutative square
\begin{equation*}
\begin{CD}
Y @>{\id}>> Y \\
@VVV     @V{\alpha}VV \\
TY @>{T\alpha}>> X.
\end{CD}
\end{equation*}
We compute $\cof(T\alpha)/\cof(\id_Y)$ (with the first cofiber computed in normed monoids, and the second in pointed spaces) by writing each as a bar construction.
This yields \[ \cof(T\alpha) \wequi \cof(T\alpha)/\cof(\id_Y) \wequi |\cof(Y^{\vee \bullet} \vee Y \to (TY)^{\times \bullet} \times X)|, \] and so it suffices to show that each of the inner cofibers is $2$-effective.
Factoring the map as $Y^{\vee n} \vee Y \to (TY)^{\vee n} \vee X \to (TY)^{\times n} \times X$, the desired cofiber sits in a cofiber extension \[ (TY/Y)^{\vee n} \vee X/Y \to \cof \to \cof((TY)^{\vee n} \vee X \to (TY)^{\times n} \times X). \]
The first term is $2$-effective since $X/Y$ is by assumption, and $TY/Y$ is since it admits a filtration (the arity filtration of \S\ref{subsub:free-normed}) with $2$-effective subquotients (by Lemmas \ref{lemm:T-filtration} and \ref{lemm:Dn-effectivity-Spc*}(1)).
Filtering $A_1 \times \dots \times A_n$ by the number of non-basepoint coordinates gives a filtration of the cofiber of $A_1 \vee \dots \vee A_n \to A_1 \times \dots \times A_n$ with subquotients $\bigvee_{|I|=r}\bigwedge_{i\in I}A_i$ for $2\le r\le n$, so the remaining term is $2$-effective by Remark~\ref{rmk:effective-operations} and \cite[Theorem~3.1.20]{freudenthal}.

(2) The argument is essentially the same.
Since $X$ is connected so is $\cof(T\alpha)$, and hence the $p$-localization in normed monoids and spaces is the same (see \S\ref{subsub:mon-ploc}).
Consequently \[ \cof(T\alpha)_{(p)} \wequi |\cof(Y^{\vee \bullet} \vee Y \to (TY)^{\times \bullet} \times X)_{(p)}|_{(p)} \] and so, appealing to \cite[Proposition 4.3.4]{freudenthal}, it will suffice to show that the $p$-localization of each inner cofiber is $2$-effective and simply-connected.
The arity layers of $TY/Y$ are very $m$-effective for $m \ge 2$, hence simply-connected.
The smash-product layers in the wedge-to-product cofiber are simply-connected because their factors are connected.
Together with the hypothesis on $(X/Y)_{(p)}$, these observations justify the same cofiber-extension argument in the category of $p$-local spaces, using preservation of simple-connectedness under $p$-localization and \cite[Proposition 4.3.4]{freudenthal}.
\end{proof}

\subsection{Cellularity}
We use the categories \[ \NMon(\SH^{S^1}(k))^{\cell i \eff}_{\ge n} \subset \NMon(\SH^{S^1}(k)) \] arising as a special case of Definition \ref{def:cellular-categories}.
If $n \ge 1$, then we denote by $\NMon(\Spc(k))^{\cell i \eff}_{\ge n} \subset \NMon(\Spc(k))$ the corresponding subcategory, via Lemma \ref{lemm:NMon-prestable} and Example \ref{ex:connected-gpa1}.

\begin{corollary} \label{cor:NMon-cell-crit}
Let $E \in \NMon(\SH^{S^1}(k))$.
Then $E \in \NMon(\SH^{S^1}(k))^{\cell j \eff}_{\ge n}$ if and only if there exists a (possibly transfinite) family of cofiber sequences $S_i \to E_i \to E_{i+1}$ with $S_i \in \NMon(\SH^{S^1}(k))^{\cell j \eff}_{\ge n}$, such that $E_0=E$ and $\colim_i E_i=0$.
Moreover, given $E \in \NMon(\SH^{S^1}(k))^{\cell j \eff}_{\ge n}$, such a sequence of $E_i$ indexed on $\N$ can be found, with $S_i$ a sum of objects of the form $\Sigma^{n+i}T(S^{w,w})$ (where $w\ge j$) and $E_i \in \NMon(\SH^{S^1}(k))^{\cell j \eff}_{\ge n+i}$.
\end{corollary}
\begin{proof}
Apply Lemma \ref{lemm:cell-crit} to the compact free objects $T(S^{w,w})$, $w\ge j$.
By the free/forgetful adjunction, mapping spectra out of these objects are given by evaluating contracted homotopy sheaves of the underlying $S^1$-spectrum at $k$.
The underlying free objects are connective, so $[T(S^{w,w}),\Sigma^iT(S^{v,v})]=0$ for $i>0$, as needed.
\end{proof}

\begin{lemma} \label{lemm:motives-of-Dn}
Let $R$ be a $\Z$-algebra.
The endofunctor $D_n$ of $\Mod_{HR}(k)$ (arising from the normed structure of $HR$, see Example \ref{ex:HR-normed}) preserves cellular objects.
More specifically, \begin{gather*} D_n(\Mod_{HR}(k)^{\cell i \eff}) \subset \Mod_{HR}(k)^{\cell ni \eff}, \\ D_n(\Mod_{HR}(k)^{\cell i \eff}_{\ge i}) \subset \Mod_{HR}(k)^{\cell ni \eff}_{\ge ni}, \\ \text{and for $i,r \ge 0$, $n \ge 1$ } D_n(\Mod_{HR}(k)^{\cell i \eff}_{\ge r}) \subset \Mod_{HR}(k)^{\cell ni \eff}_{\ge r} \end{gather*}
\end{lemma}
\begin{proof}
We first show that $D_n(\Mod_{HR}(k)^{\cell 0 \eff}_{\ge 0}) \subset \Mod_{HR}(k)^{\cell 0 \eff}$.
Thus let $E \in \Mod_{HR}(k)^{\cell 0 \eff}_{\ge 0}$.
The inclusion $\Mod_{HR}(k)^{\cell 0 \eff} \to \Mod_{HR}(k)$ admits a right adjoint $\tau^0$; set $\tau = \id/\tau^0$.
By Corollary \ref{cor:module-cellular-criterion}, we can nullify $E$ by a sequence of cell attachments.
Let $S^{p,q} \to E \to E'$ (with $p \ge q$) be such a cell-attachment.
It will suffice to show that $\tau D_n(E) \wequi \tau D_n(E')$.
Rewriting the cofiber as a geometric realization, it is enough to show that $D_n(E \oplus S^{p,q} \oplus \dots \oplus S^{p,q}) \to D_n(E)$ is inverted by $\tau$.
Expanding $D_n$ of a sum into tensor products of $D_i$-s (using \eqref{eq:formula-for-nsym}) and using induction on $n$, we are reduced to showing that $D_n(S^{p,q}) \in \Mod_{HR}(k)^{\cell 0 \eff}$.
For $p=q=0$, this amounts to showing that $HR[B_\et\Sigma_n] \in \Mod_{HR}(k)^{\cell 0 \eff}$, which we will do at the end.
For $p = 2q > 0$, we have $D_n(\Sigma^{2q,q} HR) \wequi \Sigma^{2nq,nq}HR[B_\et\Sigma_n] \in \Mod_{HR}(k)^{\cell 0 \eff}$.
It will thus remain to show that if $D_n(S^{p,q}) \in \Mod_{HR}(k)^{\cell 0 \eff}$ the same is true for $S^{p+1,q}$ and $S^{p-1,q}$.
As above, we see that $D_n(S^{p,q} \oplus \dots \oplus S^{p,q}) \in \Mod_{HR}(k)^{\cell 0 \eff}$.
Writing $\Sigma S^{p,q}$ as a sifted colimit thus deals with $S^{p+1,q}$.
Since $D_n$ is polynomial, it commutes with finite totalizations \cite[Proposition 2.15]{k-theory-and-poly-functors}, and so $\Omega S^{p,q} \wequi S^{p-1,q}$ is also dealt with (note that $\Mod_{HR}(k)^{\cell 0 \eff}$ is stable under finite totalizations, those being finite limits).

Next we show that $D_n(\Mod_{HR}(k)^{\cell 0 \eff}_{\ge 0}) \subset \Mod_{HR}(k)^{\cell 0 \eff}_{\ge 0}$, i.e., if $E \in \Mod_{HR}(k)^{\cell 0 \eff}_{\ge 0}$ then $D_n(E)$ is connective (see Corollary \ref{cor:module-cellular-criterion}(3)).
Using the same strategy as above, we reduce to showing that $D_n(S^{p,q})$ is connective for $p \ge q$.
But this is clear since it is the suspension spectrum of a space (namely $D_n^\wedge(S^{p,q})$).

Now let $E \in \Mod_{HR}(k)^{\cell i \eff}_{\ge i}$.
Then $E = \Sigma^{2i,i}E'$ with $E' \in \Mod_{HR}(k)^{\cell 0 \eff}_{\ge 0}$ and so $D_n(E) \wequi \Sigma^{2ni,ni} D_n(E') \in \Mod_{HR}(k)^{\cell ni \eff}_{\ge ni}$ (see \cite[Proposition~5.37]{bachmann-colimits} for the equivalence).

Next we want to show that $D_n(\Mod_{HR}(k)^{\cell i \eff}) \subset \Mod_{HR}(k)^{\cell ni \eff}$.
By taking colimits and arguing by induction, we can assume that $D_n(\Mod_{HR}(k)^{\cell i \eff}_{\ge r}) \subset \Mod_{HR}(k)^{\cell ni \eff}$ and want to prove it with $r-1$ in place of $r$.
This we do by rewriting a desuspension as a finite totalization, as before.

Finally let $n \ge 1$, $i, r \ge 0$ and $E \in \Mod_{HR}(k)^{\cell i \eff}_{\ge r}$.
We want to show that $D_n(E) \in \Mod_{HR}(k)^{\cell ni \eff}_{\ge r}$, i.e., $D_n(E)$ is $ni$-effective and $r$-connective (see Corollary \ref{cor:module-cellular-criterion}(3,4)).
Effectivity was already established.
For connectivity, write $E \wequi \Sigma^r E'$ with $E' \in \SH(k)^\veff$.
Rewriting $D_n(E)$ as a colimit of norms, it will be enough to show that if $p: X \to Y \in \Sm_k$ is finite étale of positive degree, then $p_\otimes(S^1) \in \Sigma \SH(Y)^\veff$.
Relating this via Grothendieck's Galois theory to equivariant homotopy theory \cite[\S10.2]{bachmann-norms}, we reduce this to the following statement: if $G$ is a finite group acting on a finite nonempty set $Z$, then $S^Z$ (smash product of $S^1$ indexed by $Z$, viewed as a genuine $G$-space) is a suspension.
For this it suffices to show that the permutation representation $\R{Z}$ admits a trivial summand, which it does since it contains the invariant vector $\sum_{z \in Z} z$.

\textit{Proof that $HR[B_\et\Sigma_n]$ is cellular.}
We may assume that $R=\Z$, cellular objects being stable under base change.
Since $\tau$ preserves colimits it commutes with $p$-localization, and so it suffices to show that $H\Z_{(p)}[B_\et\Sigma_n]$ is cellular.
Suppose first $p \ne char(k)$.
The motive $H\F_p[B_\et\Sigma_p]$ is cellular by \cite[Theorem~6.16]{voevodsky2003reduced}.
Since also $H\Q[B_\et\Sigma_p] \wequi H\Q$, it follows that $H\Z_{(p)}[B_\et\Sigma_p]$ is cellular: indeed, its image under $\tau$ vanishes modulo $p$ and rationally, and a $p$-local object with these two properties is zero.
For $m<p$, the order of $\Sigma_m$ is invertible in $\Z_{(p)}$, so $H\Z_{(p)}[B_\et\Sigma_m] \wequi H\Z_{(p)}$ is cellular as well.
The above argument now already shows that $D_p$ preserves cellular objects.
Let $W_0=\{1\}$ and, for $a\ge1$, let $W_a\subset\Sigma_{p^a}$ be the iterated wreath product of $a$ copies of $\Sigma_p$ in its standard permutation action.
By \cite[Proposition~5.28 and Example~5.10]{bachmann-colimits}, the motive $H\Z_{(p)}[B_\et W_a]$ is obtained by iterating $D_p$ on $H\Z_{(p)}$, and hence is cellular.
Write $n=\sum_a n_ap^a$ with $0\le n_a<p$ and embed $H=\prod_a W_a^{n_a}$ into $\Sigma_n$ by its action on disjoint blocks.
Since $v_p(|W_a|)=1+p+\dots+p^{a-1}=v_p((p^a)!)$, the index of $H$ in $\Sigma_n$ is prime to $p$.
The motive $H\Z_{(p)}[B_\et H]$ is a tensor product of cellular motives, and $H\Z_{(p)}[B_\et\Sigma_n]$ is a retract of it by \cite[Proposition~5.39]{bachmann-colimits}.
This concludes the case $char(k) \ne p$.
Finally if $p = char(k)$ then $H\Z_{(p)}[B_\et\Sigma_n] \wequi H\Z_{(p)}$ \cite[Proposition~2.1(3)]{bachmann-linearity} is cellular.
\end{proof}

\begin{lemma} \label{lemm:T-sphere-cellular}\NB{more generally $\overline{HR[\ph]}$ preserves cellular categories}
For $n, c \ge 0$ we have $\overline{HR[T(S^{c+n,n})]} \in \Mod_{HR}(k)^{\cell n \eff}_{\ge c}$.
\end{lemma}
\begin{proof}
By Lemma \ref{lemm:M-of-T} we have \[ \overline{HR[T(S^{c+n,n})]} \wequi \bigoplus_{i \ge 1} D_i(\Sigma^{c+n,n}HR), \] which satisfies the required cellularity by Lemma \ref{lemm:motives-of-Dn}.
\end{proof}

\begin{lemma} \label{lem:cellularity}
Suppose that $k$ is unorderable and $cd_2(k) < \infty$.
Let $X \in \NMon(\Spc(k))^{\gp\A^1}$ be $2$-effective.
If $\overline{H\Z[X]} \in \Mod_{H\Z}(k)^{\cell j \eff}_{\ge n}$, then also $X \in \NMon(\Spc(k))^{\cell j \eff}_{\ge n}$.
\end{lemma}
\begin{proof}
We build a sequence of normed monoids $X=X_0 \to X_1 \to \dots$ such that $X_{i+1}$ is obtained from $X_i$ by coning off (free normed monoids on) cells, and $X_{i+1}$ satisfies the same assumptions as $X$ and is $(i+1)$-connective.
It follows that $\colim_i X_i = *$ and we conclude by Corollary \ref{cor:NMon-cell-crit}.
It will suffice to explain how to build $X_{i+1}$ given $X_i$.
By Corollary \ref{cor:module-cellular-criterion} we may assume that $\overline{H\Z[X_i]} \in \Mod_{H\Z}(k)^{\cell j \eff}_{\ge i}$, where $j \ge 2$ (and $i \ge n$), and furthermore there exists $S_i \in \Mod_{H\Z}(k)^{\cell j \eff}_{\ge i}$ a sum of cells, and a map $\alpha: S_i \to \overline{H\Z[X_i]}$ with cofiber in $\Mod_{H\Z}(k)^{\cell j \eff}_{\ge i+1}$.
Let $\tilde S_i \in \Spc(k)_*$ be the corresponding sum of cells.
We claim that there exists a map $\tilde \alpha: \tilde S_i \to X_i$ with $\overline{H\Z[\tilde\alpha]} \wequi \alpha$.
The map $\alpha$ corresponds to certain elements in $\ul\pi_i(\overline{H\Z[X_i]})_{-e}(k)$, and so it will suffice to show that $\ul\pi_i X_i \to \ul\pi_i \Sigma^\infty X_i \to \ul\pi_i(\overline{H\Z[X_i]})$ is surjective (contraction being an exact operation \cite[Proposition~3.1.18(3)]{AFHLocalization}).
Since $X_i$ is $i$-connective and $2$-effective, the first map is an isomorphism by \cite[Theorems~6.37 and 6.61]{A1-alg-top} and \cite[Theorem~2.2.30(1)]{freudenthal}, and the second map is given by modding out by $\eta$ \cite[Theorem~12]{bachmann-hurewicz}, whence surjective.
We have thus obtained a map $\tilde\alpha$ as claimed, inducing $T(\tilde \alpha): T(\tilde S_i) \to X_i$.
By construction, the composite \[ \overline{H\Z[\tilde S_i]} \to \overline{H\Z[T(\tilde S_i)]} \xrightarrow{\overline{H\Z[T\tilde \alpha]}} \overline{H\Z[X_i]} \] is $\alpha$ and so is surjective on $\ul\pi_i$.
Consequently the same is true for $\overline{H\Z[T\tilde \alpha]}$, and hence for $\Sigma^\infty T\tilde \alpha$ by \cite[Theorem 16]{bachmann-hurewicz} (the result there is stated only under the assumption that the exponential characteristic is invertible, but the relevant place in the proof also works without this assumption, by appeal to \cite[Proposition D.1]{bachmann-resolution}).
But as already remarked, $\ul\pi_i X_i \wequi \ul\pi_i \Sigma^\infty X_i$, and so, $T\tilde \alpha$ is surjective on $\ul\pi_i$.
This implies that $\cof(T\tilde \alpha)$ is $(i+1)$-connective, as needed.
\end{proof}

\begin{example} \label{ex:K1-cellular}
$\P^\infty \wequi K(\Z(1),2)$ and $K(\Z_{(p)}(1),2)$ are cellular.
Treating first $K(\Z(1),2)$, it suffices to show that $T(S^{2,1}) \to K(\Z(1),2)$ has cellular cofiber.
But this cofiber is $2$-effective by Lemma \ref{lemm:2-eff-trick}(1) (since $\P^\infty/S^{2,1}$ is), so it suffices to show that it has cellular motive (Lemma \ref{lem:cellularity}).
For this it suffices to consider the motive of $K(\Z(1),2)$, which is clear.
For $K(\Z_{(p)}(1),2)$, we can either repeat the argument (using Lemma \ref{lemm:2-eff-trick}(2)), or observe that cellular normed monoids are closed under colimits, and hence $p$-localizations.
\end{example}

\begin{remark} \label{rmk:cell-p-loc}
If $X$ is $p$-local, in Lemma \ref{lem:cellularity} it suffices to check that $\overline{H\Z_{(p)}[X]}$ is cellular.
Indeed $\overline{H\Z[X]}$ is also $p$-local \cite[Proposition~2.2.16]{freudenthal}, and so $\overline{H\Z[X]} \wequi \overline{H\Z_{(p)}[X]}$.
\end{remark}

\begin{proposition} \label{prop:Dn-cellular}
Suppose -1 is a sum of squares in $k$ and $cd_2(k) < \infty$.\NB{surely the first assumption is actually unnecessary. the second also is, but more effort here}
Let $X \in \NMon(\Spc(k))^{\cell 1 \eff}_{\ge 1}$.
Then $D_n^\otimes(X)$ is cellular.
\end{proposition}
\begin{proof}
For $n \le 1$ there is nothing to prove.
For $n \ge 2$, $D_n^\otimes(X)$ is $2$-effective by Proposition \ref{prop:Dn-tens-eff}.
We also have $H\Z[D_n^\otimes(X)] \wequi D_n^\hotimes(H\Z[X])$ and so, by appeal to Lemma \ref{lem:cellularity}, it suffices to prove the analogous statement for $\bNAlg(\Mod_{H\Z}(k))$.
Thus let $X \in \bNAlg(\Mod_{H\Z}(k))$ with cellular underlying module.
We can write $X \wequi |\NSym^{\bullet +1} X|$ and so $D_n^\hotimes(X) \wequi |D_n^\hotimes \NSym^{\bullet +1} X|$.
Degreewise we have $D_n^\hotimes \NSym^{d +1} X \wequi \NSym D_n \NSym^{d} X$.
Since also $\NSym \wequi \bigoplus_{i \ge 0} D_i$, it suffices to show that $\Mod_{H\Z}(k)^\cell$ is closed under $D_n$.
This was done in Lemma \ref{lemm:motives-of-Dn}.
\end{proof}

\subsection{Rational normed monoids}
\begin{proposition} \label{prop:rational-vs-h}
Let $k$ be a field.
The functor \[ H\Q \otimes \sigma^\infty: \SH^{S^1}(k)_\Q \to \Mod_{H\Q}(k)^\eff \] is a localization (namely in the $h$-topology).
\end{proposition}
\begin{proof}
By \cite[Theorems~14.2.9 and~16.1.4, Corollary~16.1.7(1)]{triangulated-mixed-motives}, $\Mod_{H\Q}(k) \wequi \DM(k,\Q)$.
Voevodsky's cancellation theorem \cite[Corollary~4.10]{voevodsky2002cancellation} (extended to imperfect fields by passage to the perfection, as in \cite[proof of Corollary~5.5.5(2)]{cisinski-deglise-etale-motives}) and the comparison with $qfh$-sheaves \cite[Theorem~10.5.6 and Proposition~10.5.12]{triangulated-mixed-motives} give
\[
  \omega^\infty H\Q[X]
  \wequi L_{\A^1}\Q_{\mathrm{tr}}X
  \wequi L_{\A^1}L_{qfh}\Q[X],
\]
where $\Q[X]$ denotes the free (pre)sheaf of $\Q$-vector spaces on $X \in \Sm_k$.
Since $H\Q[\ph]$ inverts $h$-equivalences \cite[Theorem~14.3.4(2)]{triangulated-mixed-motives}, and hence $qfh$-equivalences, we deduce that the counit
\[
  H\Q \otimes \sigma^\infty \omega^\infty H\Q[X]
  \longrightarrow H\Q[X]
\]
is an equivalence.
Since $H\Q \otimes \sigma^\infty$ and $\omega^\infty$ preserve colimits and finite limits, and $\Mod_{H\Q}(k)^\eff$ is generated under these operations by the $H\Q[X]$, it follows that $\omega^\infty$ is fully faithful.
The displayed description of the unit shows that this localization is generated by $qfh$-equivalences.
As it also inverts $h$-equivalences, it agrees with $h$-localization.
\end{proof}

We also note the following class of local objects.
\begin{lemma} \label{lemm:NMon-Q-local-objects}
Let $k$ be unorderable.
Then $\SH^{S^1}(k)_\Q(2) \subset \SH^{S^1}(k)_\Q$ consists of local objects for the localization of Proposition \ref{prop:rational-vs-h}.
\end{lemma}
\begin{proof}
The category $\SH^{S^1}(k)_\Q(2)$ is generated under colimits by objects of the form $\Sigma^{n,2}\Sigma^\infty_{S^1}(X)_\Q$, where $X \in \Spc(k)_*$ and $n \in \Z$.
Such objects lie in the closure under limits of the image of $\omega^\infty: \SH(k)_\Q \to \SH^{S^1}(k)_\Q$ by \cite[Theorem~4.14(1)]{bachmann-moving2}, whose hypothesis $T_2(k)$ holds by the proof of \cite[Theorem~2.2.30]{freudenthal}.
Since also $\Mod_{H\Q}(k) \wequi \SH(k)_\Q$, $k$ being unorderable \cite[Corollary~16.2.14]{triangulated-mixed-motives}, we see that the generators are indeed local.

The general case follows from this, since the localization preserves colimits (the left adjoint preserving compact generators).
\end{proof}

Here is a version for normed monoids.
\begin{corollary} \label{cor:NMon-rat-loc}
The functor \[ \NMon(\Spc(k))_\Q \to \Mod_{H\Q}(k)^\veff \] left adjoint to $\Omega^\infty$ is a localization.
Additionally, it upgrades to a normed functor.
Finally, if $k$ is unorderable, very $2$-effective objects are local.
\end{corollary}
For the second statement, on the left hand side we mean the norm structure coming from Lemma \ref{lemm:norms-for-nhopf}, and on the right hand side the one obtained by restriction from $\Mod_{H\Q}(k)$.
\begin{proof}
We can factor the functor as \begin{gather*} \NMon(\Spc(k))_\Q \xrightarrow{(1)} \NMon(\Spc^\et(k))_\Q \stackrel{(2)}{\wequi} \CMon(\Spc^\et(k))_\Q \stackrel{(3)}{\wequi} \CMon(\Spc)_\Q \otimes \Spc^\et(k) \\ \xrightarrow{(4)} \CMon(\Spc)_\Q \otimes \Mod_{H\Q}(k)^\veff \wequi \Mod_{H\Q}(k)^\veff. \end{gather*}
Here $\Spc^\et(k)$ denotes the étale localized version of $\Spc(k)$.
Since $p_*$ is compatible with étale equivalences and motivic equivalences (use/argue as in \cite[Proposition 2.11]{bachmann-norms}), the localization $\Spc(k) \to \Spc^\et(k)_\Q$ is compatible with norms in the sense of \cite[Definition 6.14]{bachmann-norms} and consequently induces the localization (1) \cite[Proposition 6.16 and Lemmas D.3, D.6]{bachmann-norms}, which is moreover normed by Lemma \ref{lemm:norms-for-nhopf}.
The equivalence (2) is \cite[Corollary C.16]{bachmann-norms}.
The equivalence (3) is \cite[Theorem~4.6]{gepner2016universality}; here $\otimes$ means the Lurie tensor product applied sectionwise.
The map (4) is obtained by noting that $\Mod_{H\Q}(k)^\veff$ is étale local \cite[Theorem~14.3.4(1)]{triangulated-mixed-motives}, and hence $\Spc(k) \to \Mod_{H\Q}(k)^\veff$ factors through $\Spc^\et(k)$; then tensor with $\CMon(\Spc)_\Q$.
It remains to see that (4) is a localization, i.e., the right adjoint is fully faithful.
This follows from the facts that $\CMon(\Spc^\et(k)) \to \CMon(\Spc(k))$, $\SH^{S^1}(k)_{\ge 0} \to \CMon(\Spc(k))$ and $\Mod_{H\Q}(k)^\veff \to \SH^{S^1}(k)_{\ge 0}$ are fully faithful (the first one is clear, the second one follows from \cite[Proposition~3.1.9]{EHKSY}, and the third one is Proposition \ref{prop:rational-vs-h}).

For the ``additionally'' part, we have already proved that the first functor is normed.
Since all the other functors are between étale sheaves, they are normed if and only if they are symmetric monoidal (compatibly in the base) \cite[Corollary C.13]{bachmann-norms}.
This is clear from the constructions.

The last statement follows, under the unorderability assumption, from our description of the functors and Lemma \ref{lemm:NMon-Q-local-objects}.
\end{proof}

\section{Normed Hopf algebras: characteristic $p$ considerations} \label{sec:normed-hopf-charp}
We now specialize the theory from \S\ref{sec:norm-hopf-1} to positive characteristic motivic cohomology.
Thus we will be working over a (usually perfect) field $k$ of characteristic $p$.
Our aim is to show that, in this setting, suitable cellular subcategories of normed algebraic structures coincide with derived algebraic structures.

In the language of \S\ref{sec:norm-hopf-1}, we will mainly deal with $A = H\F_p \in \SH(k)$.
For technical reasons, we need to be slightly more general, so throughout we write $R$ for a chosen coefficient ring for homology, and work with $A=HR \in \SH(k)$.
The only case besides $R=\F_p$ we really need is $R=\Z/p^2$, although if $k$ is finite then sometimes $R=\Z_{(p)}$ is also possible.

\subsection{Normed and derived algebras} \label{subsec:normed-and-derived}
As a warm-up case, we show in this section that cellular normed algebras (with $\F_p$-coefficients in characteristic $p$) are the same as derived algebras.

Recall from \S\ref{subsub:charp} the equivalence $\Mod_{HR}(k)^\cell \wequi \Mod_R^\Z$.
Recall also the theory of derived algebras in a symmetric monoidal stable $\infty$-category with a $t$-structure as developed in \cite[\S4]{raksit2020hochschild}.
We denote the free derived algebra monad by $\LSym$.

\begin{remark} \label{remark:LSym-universal-property}
$\LSym$ is the underlying monad of a \emph{filtered monad} $\LSym_*$.
This is in fact the initial filtered monad under the usual filtered symmetric monad which preserves sifted colimits, is degreewise polynomial, and has the following property:
\begin{itemize}
\item Iterated application of the monad to a finite projective $R$-module (placed in homotopical degree zero) produces a coconnective $R$-module.
\end{itemize}
This is proved in \cite[Lemma 2.13]{bachmann-MHH} for $R=\F_2$ and with a slightly different formulation; the proof in fact yields the above more general assertion.
\end{remark}

\begin{proposition} \label{prop:NAlg-DAlg}
Let $k$ be perfect of characteristic $p$, and $R=\Z/p^\nu$ for some $\nu \ge 1$.
For $E \in \Mod_{HR}(k)^\cell$ , we have $\NSym(E) \in \Mod_{HR}(k)^\cell$.
Moreover \[ \NAlg(\Mod_{HR}(k)^\cell) \wequi \DAlg(\Mod_{R}^\Z). \]

Similar results apply with $R=\Z_{(p)}$, provided that $k$ is finite.
\end{proposition}
\begin{proof}
Under our assumptions on $k$, we have $\Mod_{HR}(k)^\cell \wequi \Mod_{R}^\Z$, as discussed above.

Since $\NSym \wequi \bigoplus_n D_n$ by \eqref{eq:formula-for-nsym}, the first assertion is Lemma \ref{lemm:motives-of-Dn}.
The functor $\NSym$ thus defines a monad on $\Mod_{HR}(k)^\cell$, with category of algebras $\NAlg(\Mod_{HR}(k)^\cell)$ (see \cite[Theorem~4.7.3.5]{lurie-ha}).

Write $\Mod_{HR}(k)^0 \subset \Mod_{HR}(k)$ for the subcategory on finite sums of $\Sigma^{2*,*}HR$, i.e., finite split Tate motives.
This is a symmetric monoidal $1$-category.
By inspection (using that $HR[B_\et\Sigma_n] \wequi HR$ \cite[Proposition 2.1(3)]{bachmann-linearity}), $D_n|_{\Mod_{HR}(k)^0} \wequi \Sym_n|_{\Mod_{HR}(k)^0} \wequi \LSym_n|_{\Mod_{HR}(k)^0}$.\NB{say more?}
We now extend this to an equivalence of monads, on all of $\Mod_{HR}(k)^\cell$.

By \cite[(proof of) Lemma 2.14]{bachmann-MHH}, the monad $\NSym$ upgrades to filtered monad $\NSym_*$.
We now claim that $\NSym_*$ satisfies the properties of Remark \ref{remark:LSym-universal-property}, i.e., that $D_n$ preserves sifted colimits, is polynomial, and iterated applications send finite sums of $\Sigma^{2n,n}HR$ to a Chow coconnective spectrum.
For the first two properties, see \cite[Example~5.12 and Corollary~5.33]{bachmann-colimits}.
But applying $\NSym$ (iteratively) to $\Mod_{HR}(k)^0$ is the same as applying $\LSym$, by what we said above.
So the required coconnectivity for $\NSym$ just follows from the one for $\LSym$.

We thus have a map of filtered monads $\LSym_* \to \NSym_*$.
To prove this is an equivalence, it suffices to do so after restricting to $\Mod_{HR}(k)^0$ \cite[Proposition 4.2.15]{raksit2020hochschild}.
This case was already dealt with.
\end{proof}

Let us immediately record a first application.
Recall that $J$ denotes the James construction, i.e., the free $\scr E_1$-monoid on a pointed motivic space, and $T$ denotes the free normed monoid on a pointed motivic space.
\begin{corollary} \label{cor:identify-J-T}
Let $k$ be of characteristic $p$.
Consider the canonical maps \[ S^{2n+1,n} \to T(S^{2n+1,n}) \quad\text{and}\quad JS^{2n,n} \to TS^{2n,n}. \]
\begin{enumerate}
\item The maps are inverted by $H\Z_{(p)}[\ph]$.
\item For $n \ge 2$, the maps are inverted by $(\ph)_{(p)}$.\NB{in fact $n \ge 1$ works}
\end{enumerate}
\end{corollary}
\begin{proof}
By essentially smooth base change\NB{details?}, we can assume that $k=\F_p$.

(1) Combining Proposition \ref{prop:NAlg-DAlg} and Lemma \ref{lemm:M-of-T} we see that $H\Z_{(p)}[TS^{2n+\epsilon,n}]$ has the homotopy of a free derived algebra on a generator in Chow degree $\epsilon$, i.e., for $\epsilon = 0$ we get a polynomial ring and for $\epsilon = 1$ we get an exterior algebra.
This is the same as the homology of, respectively, $JS^{2n,n}$ and $S^{2n+1,n}$.

(2) By Lemma \ref{lemm:limit-determined-M-test}(2), it suffices to show that the spaces are limit-determined.
This holds by Example \ref{ex:resolvable-spaces}, Example \ref{ex:J-effective} and Lemma \ref{lemm:Dn-effectivity-Spc*}(2).
\end{proof}

Let us also point out the following straightforward consequence.
\begin{lemma} \label{lemm:Dn-weights}
Same assumptions as for Proposition \ref{prop:NAlg-DAlg}.\NB{statement about $\scr C_i$ holds much more generally}
Write $\scr C_i, \scr D_i \subset \Mod_{HR}(k)^\cell$ for the subcategories generated under colimits by $\Sigma^{m,n}HR$, where $n \ge i$ respectively $n \le i$.
Then $D_n(\scr C_i) \subset \scr C_{ni}$ and $D_n(\scr D_i) \subset \scr D_{ni}$.
\end{lemma}
\begin{proof}
We discuss the analogous statement for $\Mod_R^\Z$, which is equivalent by the proposition.
We prove the statement about $\scr C_i$, the one about $\scr D_i$ is similar.
Write $D_n(E)_j \in \Mod_R$ for the weight $j$ part of $D_n(E)$.
We must show that for $E \in \scr C_i$, $D_n(E)_j = 0$ for $j < ni$.
Since $D_n(\ph)_j$ (for fixed $j$) is a polynomial functor preserving sifted colimits, by \cite[Proposition 4.2.15]{raksit2020hochschild}, it suffices to show the vanishing on projective generators of $\scr C_i$.
But if $X$ is such a generator, $D_n(X)$ is a summand of $X^{\otimes n}$, which indeed lies in $\scr C_{ni}$.
\end{proof}

\subsection{Accessing normed Hopf algebras}
We now want to analyze the theory of normed Hopf algebras over a perfect field of characteristic $p$.
We aim to produce similar results as we did for normed algebras in the previous subsection.
In this subsection we begin our study of normed Hopf algebras; comparison with derived Hopf algebras is relegated to the next subsection.
To simplify things a bit, we set $R=\Z/p^\nu$ and disregard the $\Z_{(p)}$ case.

Define \[ \Mod_{HR}(k)^{\cell+} = \Mod_{HR}(k)^{\cell} \cap \Mod_{HR}(k)^{\eff} \cap \Mod_{HR}(k)^{\fet+} \subset \Mod_{HR}(k). \]

\begin{remark} \label{ref:Mod-cell+-explicit}
Since $k$ is perfect and $R=\Z/p^\nu$, objects of $\Mod_{HR}(k)^{\cell+}$ can be written in the form $\bigoplus_i \Sigma^{2i,i} M_i$, where $M_i$ is a perfect $R$-module (placed here in weight zero) and $M_i = 0$ for $i < 0$.
\end{remark}

Set \[ \Mod_{HR}(k)^{\cell\cn+}_{/HR} = \{ C \to HR \mid C \in \Mod_{HR}(k)^{\cell+}, s_0(C) \xrightarrow{\wequi} HR \} \subset \Mod_{HR}(k)_{/HR}, \] and define $\Mod_{HR}(k)^{\cell\cn+}_{HR/}$ similarly.
Define \[ \NHopf(\Mod_{HR}(k))^{\cell\cn+} \subset \NHopf(\Mod_{HR}(k)) \text{ and }  \cNAlg(\Mod_{HR}(k))^{\cell\cn+} \subset \cNAlg(\Mod_{HR}(k)) \] by requiring the underlying object to lie in $\Mod_{HR}(k)^{\cell\cn+}_{/HR}$.
Define $(\ph)^{-}$ similarly: e.g. $\Mod_{HR}(k)^{\cell-}$ consists of objects of the form $\bigoplus_{i \le 0} \Sigma^{2i,i} M_i$, where each $M_i$ is a perfect $R$-module placed in weight zero.

\begin{example}
The unit of the symmetric monoidal category $\NHopf(\Mod_{HR}(k))$ is \[ HR[\Z] = \bigoplus_{n \in \Z} HR. \]
This does not lie in $\NHopf(\Mod_{HR}(k))^{\cell\cn+}$.
For this and related reasons, in what follows we will often consider \emph{nonunital} symmetric monoidal categories.
\end{example}

\begin{remark}
Write $\Mod_{HR}(k)^{\cell\cn+} \subset \Mod_{HR}(k)^{\cell+}$ for the subcategory on those modules $M$ concentrated in strictly positive weights.
Then the functor $(\ph)_+: \Mod_{HR}(k)^{\cell\cn+} \to \Mod_{HR}(k)^{\cell\cn+}_{HR/}$ is an equivalence of categories, with right adjoint (i.e., inverse) the reduced forgetful functor $M \mapsto \overline{M} = M/HR$.
Similarly $\Mod_{HR}(k)^{\cell\cn+} \wequi \Mod_{HR}(k)^{\cell\cn+}_{/HR}$.
\end{remark}

We note that $HR$ is weight-bounded in the sense of Definition \ref{def:wt-bounded} (see Example \ref{ex:charp-wt-bounded}).
Consequently we may specialize Proposition \ref{prop:duality} to deduce that passing to duals induces an antiequivalence between $\Mod_{HR}(k)^{\cell+}$ and $\Mod_{HR}(k)^{\cell-}$, and similarly for the other categories.
Thus the following definition is reasonable.
\begin{definition}
For $X \in \Mod_{HR}(k)^{\cell+}$ we set $D^n(X) = D(D_n(DX))$.
\end{definition}

We will apply the notion of quick convergence from Definition \ref{def:quick-conv} in $\NHopf(\Mod_{HR}(k)), \bNAlg(\Mod_{HR}(k))$, and so on.
\begin{example}
For $A, B \in \bNAlg(\Mod_{HR}(k))$ the coproduct is given by the tensor product of underlying modules (see Lemma \ref{lemm:colim-in-bNAlg}) and so the cofiber of the transition map relevant for the definition of quick convergence is $A \otimes \overline{B}$, where $\overline{B} = \cof(HR \to B)$.
In particular, a sequence $A_1, A_2, \dots \in \bNAlg(\Mod_{HR}(k)^\eff)$ is quickly converging if the effectivity of the $\overline{A_i}$ tends to infinity with $i$.
\end{example}

\begin{remark} \label{rmk:closure-under-quickly-conv}
Essentially by construction, the categories $\Mod_{HR}(k)^{\cell\pm} \subset \Mod_{HR}(k)^\cell$ are closed under finite limits, finite colimits, tensor products, and quickly converging sums/products.
Employing Remark \ref{ref:Mod-cell+-explicit}, we see that they are also closed under quickly converging $\N$-indexed limits and colimits.
\end{remark}

We now record some basic properties of $\NHopf(\Mod_{HR}(k))^{\cell\cn+}$ and related categories.
\begin{proposition} \label{prop:NHopf-basics}
\begin{enumerate}
\item Passing to duals induces an equivalence \[ \cNAlg(\Mod_{HR}(k))^{\cell\cn+} \wequi (\NAlg(\Mod_{HR}(k))^{\cell\cn-})^\op. \]
\item The forgetful functor $\NHopf(\Mod_{HR}(k))^{\cell\cn+} \to \bNAlg(\Mod_{HR}(k))^{\cell\cn+}$ is an equivalence.
\item The subcategory $\NHopf(\Mod_{HR}(k))^{\cell\cn+} \subset \NHopf(\Mod_{HR}(k))$ is closed under finite colimits and quickly converging $\N$-indexed colimits (e.g., quickly converging sums).
\end{enumerate}
\end{proposition}
\begin{proof}
(1) Immediate from Lemma \ref{lemm:duals-cNAlg}.

(2) We must show that any $M \in \bNAlg(\Mod_{HR}(k))^{\cell\cn+}$ is grouplike.
This can be checked after tensoring down from $R$ to $\F_p$\NB{details?}.
Under this assumption, $\pi_{**}M$ is a connected graded bialgebra and hence grouplike \cite[Proposition 8.2]{milnor-moore}; this implies what we want.

(3) By (2), we may replace $\NHopf$ by $\bNAlg$.
By Lemma \ref{lemm:colim-in-bNAlg}, colimits of normed bialgebras are computed in normed algebras, and so it will suffice to show that $\NAlg(\Mod_{HR}(k))^{\cell\cn+} \subset \NAlg(\Mod_{HR}(k))$ is closed under finite and quickly converging colimits, or equivalently, finite coproducts, pushouts, and quickly converging colimits.
Since the forgetful functor to $\Mod_{HR}(k)$ preserves sifted colimits, closure under quickly converging colimits follows from Remark \ref{rmk:closure-under-quickly-conv}.
Finite coproducts are given by tensor products and so are treated also in Remark \ref{rmk:closure-under-quickly-conv}.
To deal with a pushout $B_1 \leftarrow A \to B_2$, we express it as a geometric realization of tensor products $B_1 \otimes A^{\otimes \bullet} \otimes B_2$ \cite[Theorem~4.4.2.8]{lurie-ha}, commute the geometric realization with the forgetful functor (use Lemmas \ref{lemm:cNAlg-basics}(1) and \ref{lemm:NAlg-basics}) and then write the geometric realization as an $\N$-indexed colimit of finite colimits.
By \cite[Remark~1.2.4.3]{lurie-ha}, the $r$-th skeletal subquotient is $\Sigma^r(B_1 \otimes \overline{A}^{\otimes r} \otimes B_2)$: indeed, the degeneracies insert the unit of $A$, so the latching quotient replaces each factor of $A$ by $\overline{A}$.
This is quickly converging, since $\overline{A} \in \Sigma^{2,1}\Mod_{HR}(k)^\eff$; we conclude by Remark \ref{rmk:closure-under-quickly-conv}.
\end{proof}

\begin{proposition} \label{prop:NHopf-adj}
\begin{enumerate}
\item The reduced forgetful functor $\cNAlg(\Mod_{HR}(k))^{\cell\cn+} \to \Mod_{HR}(k)^{\cell\cn+}$ has a right adjoint $\Cofree$ satisfying \[ \Cofree(X) \wequi \bigoplus_{n \ge 0} D^n(X). \]
\item The forgetful functor $\NHopf(\Mod_{HR}(k))^{\cell\cn+} \to \cNAlg(\Mod_{HR}(k))^{\cell\cn+}_*$ has a left adjoint $\Free$ satisfying \[ \Free X \wequi \bigoplus_n D_n(\overline{X}). \]
  Here $\overline{X} = \cof_{\Mod_{HR}(k)}(HR \xrightarrow{*} X)$.
\end{enumerate}
\end{proposition}
\begin{proof}
(1) Passing to duals via Proposition \ref{prop:NHopf-basics}(1), it is equivalent to show: the forgetful functor $\NAlg(\Mod_{HR}(k))^{\cell\cn-} \to \Mod_{HR}(k)^{\cell\cn-}$ has a left adjoint $L$ satisfying $L(X) \wequi \bigoplus_n D_n(X)$ (use the equivalence of quickly converging sums and products from Lemma \ref{lemm:quickly-converging-sum}(4)).
This is clear without the size restrictions (from the formula \eqref{eq:formula-for-nsym} for free normed objects), and so we need only show that if $X \in \Mod_{HR}(k)^{\cell\cn-}$ indeed $\bigoplus_n D_n(X) \in \Mod_{HR}(k)^{\cell\cn-}$.
This follows from Lemma \ref{lemm:Dn-weights} (noting that, in the notation of that lemma, $X \in \scr D_{-1}$).

(2) We first prove the analogous statement for $\bNAlg(\Mod_{HR}(k)) \to \cNAlg(\Mod_{HR}(k))_*$.
For $X \in \cNAlg(\Mod_{HR}(k))_*$ we have the pushout square\NB{details?}
\begin{equation*}
\begin{CD}
HR_+ @>>> X_+ \\
@VVV     @VVV \\
HR @>>> X.
\end{CD}
\end{equation*}
Since the forgetful functor $\cNAlg(\Mod_{HR}(k)) \to \Mod_{HR}(k)$  is presentably normed (see Lemma \ref{lemm:cNAlg-basics} and Corollary \ref{cor:forgetful-normed}), the formula \eqref{eq:formula-for-nsym} for free normed objects shows that the forgetful functor commutes with the formation thereof (on unpointed objects).
Moreover, since $\bNAlg(\Mod_{HR}(k)) \to \NAlg(\Mod_{HR}(k))$ preserves colimits by Lemma \ref{lemm:colim-in-bNAlg} we obtain a pushout square in either category of the form
\begin{equation*}
\begin{CD}
\Free(HR_+) @>>> \Free(X_+) \\
@VVV     @VVV \\
\Free(HR) @>>> \Free(X).
\end{CD}
\end{equation*}
By inspection, this is $\NSym$ applied to the pushout square
\begin{equation*}
\begin{CD}
HR @>>> X \\
@VVV   @VVV \\
0  @>>> \overline{X}.
\end{CD}
\end{equation*}
We now conclude from the formula for free normed spectra.

It thus remains to show that if $X \in \cNAlg(\Mod_{HR}(k))^{\cell\cn+}_*$ then $\Free X \wequi \bigoplus_n D_n(\overline{X})$ lies in $\NHopf(\Mod_{HR}(k))^{\cell\cn+}$.
That $\Free X$ lies in $\Mod_{HR}(k)_{HR/}^{\cell\cn+}$ follows from Lemma \ref{lemm:Dn-weights}.
The grouplikeness now follows from Proposition \ref{prop:NHopf-basics}(2).
\end{proof}

\begin{remark} \label{rmk:Free-as-coalg}
The formulas in Proposition \ref{prop:NHopf-adj} are understood as objects of $\Mod_{HR}(k)$.
For (2), we can provide an expression that works as an object of $\cNAlg(\Mod_{HR}(k))$, namely \[ \colim_n D_n(X) \xrightarrow{\wequi} \Free X. \]
(Produce a map using universal properties, and check it is an equivalence by computing in $\Mod_{HR}(k)$.)
\end{remark}

We can deduce some stability properties.
\begin{proposition} \label{prop:NHopf-stability}
\begin{enumerate}
\item $\NHopf(\Mod_{HR}(k))^{\cell\cn+} \subset \NHopf(\Mod_{HR}(k))$ is closed under $\hotimes$, so is a nonunital symmetric monoidal subcategory.
\item $\NHopf(\Mod_{HR}(k))^{\cell\cn+}$ is stable.
\item Quickly converging coproducts in $\NHopf(\Mod_{HR}(k))^{\cell\cn+}$ are biproducts.
\end{enumerate}
\end{proposition}
\begin{proof}
(1) Call $A \in \NHopf(\Mod_{HR}(k))^{\cell\cn+}$ $n$-effective if $\overline{A} := \fib(A \to HR) \in \Mod_{HR}(k)$ is $n$-effective.
Call $A \in \Mod_{HR}(k)^{\cell}$ \emph{finite in weights $\le n$} if $f^{n+1} A$ is compact.
We will show the following.
In (b) and (c), we assume that all Hopf algebras involved have cellular underlying modules and $1$-effective augmentation ideals, without assuming finiteness in every weight.
\begin{enumerate}[(a)]
\item If $A, B \in \NHopf(\Mod_{HR}(k))^{\cell\cn+}$ with $A$ $n$-effective, also $A \hotimes B$ is $n$-effective.
\item If $A \in \NHopf(\Mod_{HR}(k))$ is finite in weights $\le n$, $B \in \NHopf(\Mod_{HR}(k))^{\cell\cn+}$, then for any map $A \to B$, $B \otimes_A HR$ is finite in weights $\le n$.
\item If $A \to B \to C$ is a cofiber sequence in $\NHopf(\Mod_{HR}(k))$ with $A, C$ finite in weights $\le n$, then also $B$ is finite in weights $\le n$.
\item If $A \in \NHopf(\Mod_{HR}(k))^{\cell\cn+}$ is $n$-effective then $A \otimes_{\Free A} HR =: C_A$ is $2n$-effective.
\end{enumerate}

Supposing this for now, we prove that if $A, B \in \NHopf(\Mod_{HR}(k))^{\cell\cn+}$, then $A \hotimes B$ is finite in weights $\le n$, for any $n$.
This will imply that $A \hotimes B \in \NHopf(\Mod_{HR}(k))^{\cell\cn+}$, as desired.
The cofiber sequence $\Free A \hotimes B \to A \hotimes B \to C_A \hotimes B$ together with (c) allows us to replace $A$ by $\Free A$ and $C_A$.
Repeating this process and applying (d), the term $C_A$ eventually becomes arbitrarily effective (note that every object in $\NHopf(\Mod_{HR}(k))^{\cell\cn+}$ is $1$-effective), whence the same is true for $C_A \hotimes B$ by (a), and so $C_A \hotimes B$ is finite in weights $\le n$.
We have thus reduced to the case $A = \Free A'$.
We can argue similarly on the other side, reducing to $B = \Free B'$.
But then $A \hotimes B = \Free(A' \wedge B')$, which is certainly finite in weights $\le n$.

(a) We can write $A \wequi |\Free^{\bullet+1} A|$ \cite[Proposition~4.7.3.14]{lurie-ha}, and similarly for $B$.
The formula for $\Free A$ from Proposition \ref{prop:NHopf-adj} shows that if $A$ is $n$-effective, so is $\Free A$.
Consequently $A \hotimes B$ is a colimit of objects of the form $\Free X \hotimes \Free Y \wequi \Free(X \wedge Y)$, where $X$ is $n$-effective in the sense of pointed objects (and $Y$ is $1$-effective).
Then $\Free(X \wedge Y)$ is $n$-effective, and we conclude since $n$-effective modules are stable under colimits.

(b) Compute $B\otimes_A HR$ by the two-sided bar construction.
Its skeletal filtration starts at $B$ and has successive cofibers $\Sigma^d(B\otimes\overline{A}^{\otimes d})$ for $d\ge1$.
Since $\overline{A}$ is $1$-effective and $B$ is effective, these cofibers have weights at least $d$, so only finitely many stages contribute in weights $\le n$.
In those weights each term is perfect weightwise, by the assumed finiteness of $A$ and $B$, and hence the same holds for the realization.

(c) We prove this by induction on $n\ge0$; the case $n=0$ follows because the augmentations identify the weight-zero parts with $HR$.
Assume the result proved for $n-1$, so that $B$ is finite in weights $\le n-1$.
The bar filtration computing $C=B\otimes_A HR$ starts at $B$ and has successive cofibers $\Sigma^d(B\otimes\overline{A}^{\otimes d})$ for $d\ge1$.
Through weight $n$, these terms involve only weights $\le n-1$ of $B$, because every factor of $\overline{A}$ has strictly positive weight.
They are therefore finite in weights $\le n$ by the induction hypothesis and the assumed finiteness of $A$.
Moreover, terms with $d>n$ contribute nothing in those weights, so $\cof(B\to C)$ is finite in weights $\le n$.
Since $C$ is finite in weights $\le n$, the same follows for $B$ by the fiber sequence $B\to C\to\cof(B\to C)$.

(d) As in (c), $C_A$ differs from $\cof(\Free A \to A)$ (computed in $HR$-modules) only by terms of the form $\overline{\Free A}^{\otimes d} \otimes A$ for $d > 1$, all of which are $2n$-effective.
That $\cof(\Free A \to A)$ is $2n$-effective follows from the formula of Proposition \ref{prop:NHopf-adj}(2).

(2) Pullbacks in $\NHopf(\Mod_{HR}(k))^{\cell\cn+}$ are computed in $\cNAlg(\Mod_{HR}(k))^{\cell\cn+}$, the forgetful functor admitting a left adjoint by Proposition \ref{prop:NHopf-adj}(2).
But the latter category is anti-equivalent to an appropriate subcategory of derived algebras, by Proposition \ref{prop:NHopf-basics}(1) and Proposition \ref{prop:NAlg-DAlg}.
Consequently pullbacks in $\NHopf(\Mod_{HR}(k))^{\cell\cn+}$ are computed via a double dualization and a bar construction.
Because of our connectivity assumption, this bar construction behaves degreewise like a finite colimit (see e.g. the discussion in \cite[\S2.2.2]{bachmann-wilson}), and so passes through the double dualization.
The upshot is that pullbacks are computed as a cobar construction.
Similarly pushouts are computed as a bar construction, by Lemma \ref{lemm:colim-in-bNAlg} and Proposition \ref{prop:NHopf-basics}(2,3).
That bar and cobar constructions are inverse in the connected setting is well-known; see e.g. \cite[Remark 2.17]{bachmann-wilson}.

(3) For a quickly converging coproduct $X = \bigoplus_i X_i$, we must show that the canonical maps $X \to X_i$ exhibit $X$ as the product.
Since the left adjoint $\Free$ has dense image, it suffices to check this after forgetting to $\cNAlg(\Mod_{HR}(k))^{\cell\cn+}$.
We thus have a quickly converging, split $\N$-indexed diagram of partial coproducts in $\cNAlg(\Mod_{HR}(k))^{\cell\cn+}$, and we seek to identify its colimit along the inclusions with its limit along the canonical retractions.
Dualizing via Proposition \ref{prop:NHopf-basics}(1) and Proposition \ref{prop:NAlg-DAlg}, we reduce to showing the same in derived algebras.
But for such, limits and filtered colimits are computed on underlying modules, so it remains to compare the colimit along inclusions with the limit along retractions in $\Mod_{HR}(k)^\cell$. These are respectively a quickly converging sum and product, which agree by Lemma \ref{lemm:quickly-converging-sum}(4).
\end{proof}

We have the following immediate corollary of our description of cofree normed Hopf algebras.
Let us call $C \in \cNAlg(\Mod_{HR}(k))^{\cell\cn+}$ \emph{copolynomial} if $C$ is split Tate, and $\pi_{**}(DC)$ is a polynomial ring over $R$ (necessarily on generators in Chow degree $0$, only finitely many in each weight).
\begin{lemma} \label{lemm:maps-into-copoly}
If $C \in \cNAlg(\Mod_{HR}(k))^{\cell\cn+}$ is copolynomial then there exists $C_0 \in \Mod_{HR}^{\cell+}$ with $C \wequi \Cofree C_0$.
\end{lemma}
\begin{proof}
Pick a summand $V$ of $\pi_{**}(DC)$ generating freely.
Let $C_0$ correspond to the dual module of $V$, which we can view as a summand of $C$.
To show that $C \to \Cofree C_0$ is an equivalence we may dualize (using Proposition \ref{prop:NHopf-basics}(1) and Proposition \ref{prop:NAlg-DAlg}).
But the dual assertion holds by construction.
\end{proof}

\begin{remark} \label{rmk:accessing-minus}
The discussion in this section focussed on $\Mod_{HR}(k)^{\cell\cn+}$.
Everything can be repeated with $\Mod_{HR}(k)^{\cell\cn-}$ in place thereof.
This is especially fruitful when combined with the duality equivalence \[ \NHopf(\Mod_{HR}(k))^{\cell\cn+} \wequi (\NHopf(\Mod_{HR}(k))^{\cell\cn-})^\op \] deduced from Lemma \ref{lemm:duals-cNAlg}.
\end{remark}

\begin{example} \label{ex:cofree-NHopf-alge}
We can deduce from Remark \ref{rmk:accessing-minus} and Proposition \ref{prop:NHopf-adj}(2) that the forgetful functor $\NHopf(\Mod_{HR}(k))^{\cell\cn+} \to \NAlg(\Mod_{HR}(k))^{\cell\cn+}_{/HR}$ admits a right adjoint.
In other words, cofree normed Hopf algebras exist and take the expected shape.
\end{example}

\subsection{Derived bipolynomial Hopf algebras}
We now define and study a well-behaved subcategory of $\NHopf(\Mod_{HR}(k))^{\cell\cn+}$ which we call the \emph{category of derived bipolynomial Hopf algebras}.

\subsubsection{Definition and first properties}
To begin with, write $\Mod_{HR}(k)^{\tate} \subset \Mod_{HR}(k)^{\cell}$ for the subcategory on sums of objects of the form $\Sigma^{2n,n}HR$ (for various $n$), and $\NHopf(\Mod_{HR}(k))^{\tate\cn+} \subset \NHopf(\Mod_{HR}(k))^{\cell\cn+}$ for the subcategory on those objects whose underlying module is in $\Mod_{HR}(k)^{\tate}$.
The functor $\pi_{2*,*}$ implements a symmetric monoidal equivalence of $\infty$-categories between $\Mod_{HR}(k)^{\tate}$ and the category of evenly graded free $R$-modules.
Consequently for $E \in \NHopf(\Mod_{HR}(k))^{\tate\cn+}$, $\pi_{2*,*}E$ is a Hopf algebra over $R$.
The functors $\Cofree$ and $\Free$ preserve the subcategories $(\ph)^\tate$, and $\pi_{2*,*}$ commutes with these functors (use Propositions \ref{prop:NAlg-DAlg} and \ref{prop:NHopf-basics}).

\begin{definition} \label{def:dhopfbp}
For $n>0$ we write $\A^1(n) = \Cofree \Sigma^{2n,n}HR \in \cNAlg(\Mod_{HR}(k))^{\tate\cn+}_*$.
Write \[ \Hopf_R^{bp} \subset \NHopf(\Mod_{HR}(k))^{\tate\cn+} \] for the category obtained as the closure under summands and quickly converging sums from objects of the form $\Free \A^1(n)$, for $n > 0$.
We denote by \[ \dHopf_R^{bp} \subset \NHopf(\Mod_{HR}(k))^{\cell\cn+} \] the smallest thick stable subcategory containing $\Hopf_R^{bp}$.
\end{definition}

\begin{example} \label{ex:Poo-A1}
We have $\A^1(1) \wequi HR[\P^\infty]$ (e.g. use Lemma \ref{lemm:maps-into-copoly}).
\end{example}

\begin{remark} \label{rmk:Hopf-scalar-ext}
The formation of $\dHopf_R^{bp}$ is functorial in $R$.
In other words, given $X \in \dHopf_R^{bp}$ and a ring homomorphism $R \to R'$, then $X \otimes_R R'$ (viewed as an object of $\NHopf(\Mod_{HR'}(k))$) lies in $\dHopf_{R'}^{bp}$.
We refer to the resulting functor $\dHopf_R^{bp} \to \dHopf_{R'}^{bp}$ as \emph{scalar extension}.
Indeed on the underlying modules, it is just computed as the usual scalar extension.
\end{remark}

\begin{lemma} \label{lemm:NHopf-comp-start}
\begin{enumerate}
\item Let $N \in \NHopf(\Mod_{HR}(k))^{\tate\cn+}$ be copolynomial and $M \in \cNAlg(\Mod_{HR}(k))^{\tate\cn+}_*$.
  Then $\Map_\NHopf(\Free M, N)$ is a discrete space, coinciding via $\pi_{2*,*}$ with the classical Hopf algebra maps.
\item The $\infty$-category $\Hopf_R^{bp}$ is a $1$-category, equivalent via $\pi_{2*,*}$ to the corresponding category of classical Hopf algebras.
\end{enumerate}
\end{lemma}
\begin{proof}
(1) Immediate from Lemma \ref{lemm:maps-into-copoly} and the commutation of $\pi_{2*,*}$ with $\Cofree$/$\Free$.

(2) In both categories, quickly converging sums and products coincide (see Proposition \ref{prop:NHopf-stability}(3) for $\Hopf_R^{bp}$; the case of classical Hopf algebras is straightforward), so we need only investigate maps between the generators.
By (1), it will be enough to show that $\Free\A^1(n)$ is copolynomial.
Using the description of $\Free\A^1(n)$ as a coalgebra from Remark \ref{rmk:Free-as-coalg}, we see that $\pi_{**}D\Free\A^1(n)$ is the inverse limit of rings $R[x_1, \dots, x_d]^{\Sigma_d}$ (with $|x_i| = (2n,n)$) along the maps sending the last variable to zero.
That this is a polynomial ring follows from the symmetric function theorem \cite[Chapter~I, \S2, (2.4)]{macdonald-symmetric-functions}.
\end{proof}

There is the following strengthening of part (1) of the above lemma, the proof of which we defer to the next subsection.
\begin{proposition} \label{prop:dHopfbp-disc}
For $X, Y \in \Hopf_R^{bp}$, the mapping spectrum $\map_{\dHopf_R^{bp}}(X, Y)$ is discrete.
\end{proposition}

\begin{remark} \label{rmk:bipoly-explanation}
The $1$-category $\Hopf_R^{bp}$ is additive, and so we may view it as a discrete spectral category.
Proposition \ref{prop:dHopfbp-disc} is equivalent to the assertion that \[ \dHopf_R^{bp} \wequi \Perf(\Hopf_R^{bp}). \]
When $R=\F_p$ (or $\Z_{(p)}$), by a theorem of Ravenel--Wilson \cite[Proposition~2.3]{ravenel-wilson-bipolynomial}, $\Hopf_R^{bp} \subset \Hopf_R$ consists of precisely the \emph{bipolynomial Hopf algebras}, i.e., Hopf algebras $H$ such that both $H$ and $DH$ are polynomial rings (on finitely many generators in each degree).
This explains our terminology.
\end{remark}

We also have the following, again deferring the proof to the next subsection.
\begin{lemma} \label{lemm:bipoly-monoidal}
Both $\dHopf_R^{bp}$ and $\Hopf_R^{bp}$ are closed under $\hotimes$ and $D_n^{\hotimes}$ for $n \ge 1$ (when viewed as a subcategory of $\NHopf(\Mod_{HR}(k))$).
The restriction of these functors to $\Hopf_R^{bp}$ coincides with the tensor product, respectively symmetric power, of Hopf algebras.
\end{lemma}

In particular $\Hopf_R^{bp}$ and $\dHopf_R^{bp}$ are nonunital symmetric monoidal categories, and the equivalence of Remark \ref{rmk:bipoly-explanation} is symmetric monoidal.

\begin{remark}
We have chosen to define $\dHopf_R^{bp}$ in terms of $\NHopf(\Mod_{HR}(k))^{\cell\cn+}$, i.e., using positive gradings.
We could have equivalently worked with $\NHopf(\Mod_{HR}(k))^{\cell\cn-}$, containing $\Free\A^1(n)$ for $n < 0$, and then obtaining an analogous category which we denote by $\dHopf_R^{bp-}$.
We defer the relationship between $\dHopf_R^{bp}$ and $\dHopf_R^{bp-}$ to the next section.
\end{remark}

\subsubsection{The Witt Hopf algebra and its filtrations} \label{subsub:witt}
Lemma \ref{lemm:NHopf-comp-start} provides us with access to part of the category $\dHopf_R^{bp}$.
In particular, it allows us to construct certain important objects that will feature in the sequel.
From now on we conflate classical and normed Hopf algebras, as long as they lie in $\Hopf_R^{bp}$.

Our first example is, of course, $\Free\A^1(n)$.
This is the coordinate Hopf algebra of the \emph{big Witt vectors}, with the generator placed in degree $n$ \cite[\S9]{hazewinkel-witt-vectors}.
There is a well-known summand of $\Free\A^1(n)$ denoted $W(n)$ and called the \emph{$p$-typical Witt vectors} on a generator in degree $n$ \cite[\S14.25]{hazewinkel-witt-vectors}.
In fact one finds that \cite[\S14.25]{hazewinkel-witt-vectors} \[ \Free\A^1(n) \wequi \bigotimes_{p \nmid m} W(nm). \]
As an algebra, one has $W(n) \wequi R[x_0, x_1, x_2, \dots]$ with $|x_i| = p^in$.
The primitives of $W(n)$ are freely spanned by the \emph{Witt polynomials} $w_i = x_0^{p^i} + px_1^{p^{i-1}} + \dots + p^i x_i$ \cite[Remark~(8.4)]{husemoller-hopf}.
There are two important Hopf algebra morphisms, the \emph{Verschiebung} $V: W(n) \to W(pn)$ and the \emph{Frobenius} $F: W(pn) \to W(n)$.
They are determined\footnote{In order for $F$ to be determined by this, we must work with $R=\Z_{(p)}$, or at least require the formulas for all $R=\Z/p^\nu$ at the same time.} by \cite[\S\S5.25, 5.27]{hazewinkel-witt-vectors} \[ V(x_i) = \begin{cases} x_{i-1} & i > 0 \\ 0 & i=0 \end{cases} \text{ and } F(w_i) = w_{i+1}. \]
Using these we define \[ W(n)_{\le d} = \fib(V^{d+1}: W(n) \to W(p^{d+1}n))  \in \dHopf_R^{bp} \] and \[ W(n)^{\le d} = \cof(F^{d+1}: W(p^{d+1}n) \to W(n)) \in \dHopf_R^{bp}. \]
We also write \[ \Ga(n) := W(n)^{\le 0}. \]

\begin{remark} \label{rmk:Witt-filtrations}
A formal consequence of the construction is that we have bifiber sequences \[ W(n)_{\le d} \to W(n)_{\le d+e+1} \to W(p^{d+1}n)_{\le e} \] and \[ W(p^{d+1}n)^{\le e} \to W(n)^{\le d+e+1} \to W(n)^{\le d} \] as well as diagrams \[ W(n)_{\le 0} \to W(n)_{\le 1} \to \dots \to W(n) \] and \[ W(n) \to \dots \to W(n)^{\le 1} \to W(n)^{\le 0}. \]
\end{remark}

\begin{remark} \label{rmk:W-scalar-ext}
The objects $\A^1(n)$, $W(n)$, and consequently also $W(n)_{\le d}$, $W(n)^{\le d}$, $\Ga(n)$ are stable under scalar extension in the sense of Remark \ref{rmk:Hopf-scalar-ext}.
For $\A^1(n)$ and $\Free \A^1(n)$ this follows from the formulas for $\Cofree$ and $\Free$ from Proposition \ref{prop:NHopf-adj}, the case of $W(n)$ follows by passing to a summand, and the remaining cases by taking fibers/cofibers.
\end{remark}

We have defined these as derived bipolynomial Hopf algebras, but they are in fact closely related to classical Hopf algebras of the same name.
\begin{lemma} \label{lemm:Witt-filtrations}
Both $W(n)_{\le d}$ and $W(n)^{\le d}$ lie in $\NHopf(\Mod_{HR}(k))^{\tate\cn+}$.
Their bigraded homotopy groups are given respectively by the Hopf algebra kernel of $V^{d+1}$ and cokernel of $F^{d+1}$.
Moreover we have \[ \colim_d W(n)_{\le d} \wequi W(n) \wequi \lim_d W(n)^{\le d} \in \dHopf_R^{bp}. \]
\end{lemma}
\begin{proof}
By Proposition \ref{prop:NHopf-adj}(2), the fiber defining $W(n)_{\le d}$ is computed in $\cNAlg(\Mod_{HR}(k))^{\cell\cn+}$, whence corresponds by duality (Proposition \ref{prop:NHopf-basics}(1) and Proposition \ref{prop:NAlg-DAlg}) to a pushout in $\DAlg(\Mod_R^{\Z_{\le 0}})$.
Using the self-duality of $W(n)$ \cite[Propositions~(7.2) and~(8.3)]{husemoller-hopf}, under which $F$ and $V$ are exchanged \cite[Theorem~13.48, equation~(13.56)]{hazewinkel-witt-vectors}, we see by explicit computation that $W(n)_{\le d}$ is concentrated in Chow degree zero and has homotopy given by the Hopf algebra kernel.\NB{details?}

The cofiber defining $W(n)^{\le d}$ is computed in $\NAlg(\Mod_{HR}(k))$ (by Lemma \ref{lemm:colim-in-bNAlg} and Proposition \ref{prop:NHopf-basics}) and thus (by the same computation as in the first paragraph) seen to be concentrated in Chow degree zero, with homotopy given by the Hopf algebra cokernel.

Similarly the colimit of the $W(n)_{\le d}$ is computed in $\NAlg(\Mod_{HR}(k))$ and thus immediately seen to be $W(n)$.
The limit of the $W(n)^{\le d}$ is up to dualization (as in the first paragraph) given by the same colimit, i.e., $W(n)$ again.
\end{proof}

\begin{example} \label{ex:sphere-bipolynomial}
$HR[T(S^{2n,n})] \in \NHopf(\Mod_{HR}(k))$ actually lies in $\dHopf_R^{bp}$ and coincides with $W(n)_{\le 0}$.
To show this, we must construct a map $HR[T(S^{2n,n})] \to W(n)$ and a nullhomotopy of the composite into $W(pn)$.
Since $HR[T(S^{2n,n})] \wequi \Free HR[S^{2n,n}]$, this reduces to a statement about pointed normed coalgebras, and so by duality to ordinary algebras.
It is straightforward to verify.
\end{example}

\begin{example} \label{ex:Poo-bipolynomial}
$HR[K(\Z(1),2)] \in \NHopf(\Mod_{HR}(k))$ actually lies in $\dHopf_R^{bp}$ and coincides with $\Ga(1)$.
Indeed since $\A^1(1) \wequi HR[K(\Z(1),2)]$ (Example \ref{ex:Poo-A1}) and $W(1)$ is a summand of $\Free\A^1$, we obtain a canonical map $W(1) \to HR[K(\Z(1),2)]$.
We must verify that this map factors through the cofiber of Frobenius and induces an equivalence $W(1)/F \wequi HR[K(\Z(1),2)]$.
This is straightforward.\NB{details?}
\end{example}

We want to briefly discuss duality for $\dHopf_R^{bp}$.
Given $X \in \NHopf(\Mod_{HR}(k))^{\cell\cn+}$, the dual $DX$ defines an object of $\NHopf(\Mod_{HR}(k))^{\cell\cn-}$ (see Remark \ref{rmk:accessing-minus}).
For $n<0$ and $d\ge0$, we use the reversed filtration notation $W(n)_{\ge -d} := \fib(V^{d+1}: W(n) \to W(p^{d+1}n))$ and $W(n)^{\ge -d} := \cof(F^{d+1}: W(p^{d+1}n) \to W(n))$.
\begin{lemma} \label{lemm:Witt-dual}
We have an equivalence $DW(n) \wequi W(-n)$, inducing $D(W(n)_{\le d}) \wequi W(-n)^{\ge -d}$ and $D(W(n)^{\le d}) \wequi W(-n)_{\ge -d}$.
\end{lemma}
\begin{proof}
Since $DDW(n) \wequi W(n)$, $DW(n)$ is copolynomial and so, in order to exhibit an equivalence $W(-n) \to DW(n)$ we need only exhibit a map of the underlying classical Hopf algebras (Lemma \ref{lemm:NHopf-comp-start}(1)) and show that it is an equivalence.
The classical self-duality is given by \cite[Propositions~(7.2) and~(8.3)]{husemoller-hopf}.
Under this duality $F$ and $V$ are exchanged \cite[Theorem~13.48, equation~(13.56)]{hazewinkel-witt-vectors}, so the other statements follow.
\end{proof}

\begin{corollary}
Passing to duals induces $\dHopf_R^{bp,\op} \wequi \dHopf_R^{bp-}$.
\end{corollary}
\begin{proof}
Immediate from Lemma \ref{lemm:Witt-dual} and Remark \ref{rmk:accessing-minus}.
\end{proof}

We can now supply the proofs deferred from the previous subsection.
\begin{remark} \label{rmk:hotimes-howto}
It may be helpful to recall how to compute tensor products of (some) Hopf algebras.
When working over $\Z_p$, there is a notion of Hopf algebras with a lift of Verschiebung.
Taking indecomposables then yields a graded $\Z_p$-module with a graded endomorphism $V$.
On appropriate subcategories, this functor is a symmetric monoidal equivalence.
The objects $W(n)_{\le i}$ and $W(n)$ all come from Hopf algebras with a lift of Verschiebung; consequently their tensor products may be readily computed.
More details are recalled, for example, in \cite[\S2.4.5]{bachmann-wilson}.
\end{remark}

\begin{proof}[Proof of Lemma \ref{lemm:bipoly-monoidal}.]
We begin by establishing the claims about $\hotimes$.
Let us momentarily denote the usual tensor product of Hopf algebras by $\hotimes^{cl}$.
Since both $\hotimes$ and $\hotimes^{cl}$ preserve colimits in each variable (see Lemma \ref{lemm:norms-for-nhopf} for the normed case), it will suffice to show that $W(n) \hotimes W(m)$ and $W(n) \hotimes^{cl} W(m)$ lie in $\Hopf_R^{bp}$ and agree.
In the classical case, there is a well-known formula to write $W(n) \hotimes^{cl} W(m)$ as a (quickly converging) sum of $W(d)$ for various $d$, which can be obtained for example from the method of Remark \ref{rmk:hotimes-howto}.
Let us write $\scr C \subset \dHopf_R^{bp}$ for the subcategory on objects which:
\begin{itemize}
\item are summands of objects of the form $\Free M$, where $M \in \cNAlg(\Mod_{HR}(k))^{\tate\cn+}_*$, and
\item are copolynomial.
\end{itemize}
By Lemma \ref{lemm:NHopf-comp-start}, $\scr C$ is a $1$-category which is equivalent (via $\pi_{2*,*}$) to a subcategory of $\Hopf_R$.
Observing that $\Free M \hotimes \Free N \wequi \Free(M \wedge N)$ (both classically and in the normed setting), we see that if $X, Y \in \scr C$ then $X \hotimes Y$ satisfies the first condition and $\pi_{2*,*}X \hotimes^{cl} \pi_{2*,*}Y \wequi \pi_{2*,*}(X \hotimes Y)$.
Apply this to $X=W(n)$, $Y=W(m)$.
Since $W(n) \hotimes^{cl} W(m)$ splits into $W(d)$-s, we see that $W(n) \hotimes W(m)$ is copolynomial and so lies in $\scr C$.
But since $\pi_{2*,*}$ is fully faithful on $\scr C$, we deduce that $W(n) \hotimes W(m)$ splits into the same sum as classically.
This concludes the proof for $\hotimes$.

Now we treat $D_n^\hotimes$.
The functor is obtained by restriction from a polynomial functor preserving sifted colimits, and hence commutes with finite totalizations \cite[Proposition 2.15]{k-theory-and-poly-functors}.
We also have $D_n^\hotimes(A \otimes B) \wequi \bigotimes_{n = i+j}D_i^\hotimes(A) \hotimes D_j^\hotimes(B)$ (by the formula \eqref{eq:formula-for-nsym} for free normed objects).
Thus if $\scr D \subset \dHopf_R^{bp}$ is the subcategory of objects $X$ with $D_n^\hotimes(X) \in \dHopf_R^{bp}$, then $\scr D$ is closed under finite sums, finite colimits, and finite limits.
It thus suffices to treat $\Hopf_R^{bp}$.
We will determine $D_n^\hotimes(W(a))$ and will see that it splits into a quickly converging sum of $W(b)$ for various $b \ge a$.
By induction, we can assume that, for $i<n$, $D_i^\hotimes$ preserves $a$-effective Hopf algebras for every $a$.
Using the previous computation, this shows that $D_n^\hotimes$ preserves quickly converging split colimits corresponding to sums of $W(a)$-s, and hence also preserves $a$-effectivity.

The only thing that remains to be done is to determine $D_n^\hotimes(W(a))$.
Using Remark \ref{rmk:smash-product-normed-str}, we see that $D_n^\hotimes(\Free\A^1(a)) \wequi \Free(D_n^\wedge(\A^1(a)))$.
For $X \in \cNAlg(\Mod_{HR}(k))_*$ let $\overline{X} = \cof(HR \to X)$ and consider the canonical map of $HR$-modules \[ \overline{D_n^\wedge(X)} \hookrightarrow D_n^\wedge(X) \to D_n(X) \to D_n(\overline{X}) \] (the middle map coming from the fact that $\cNAlg(\Mod_{HR}(k))_* \to \cNAlg(\Mod_{HR}(k))$ is lax normed).
This construction commutes with sifted colimits and is an equivalence on objects of the form $X_+$; hence it is an equivalence in general.
We conclude that $D_n^\wedge(\A^1(a))$ is concentrated in Chow degree zero and hence the same is true for $D_n^\hotimes(\Free\A^1(a))$, and its retract $D_n^\hotimes(W(a))$.
Moreover, we have the symmetric map $W(a)^{\hotimes n} \to D_n^\hotimes(W(a))$ which induces by applying $\pi_{2*,*}$ and passing to orbits the map $\Sym_n^\hotimes W(a) \to \pi_{2*,*} D_n^\hotimes W(a)$.
To show this is an equivalence, it suffices to prove the analogous statement for $\Free\A^1(a)$.
This follows from the classical identities $\Sym_n^\hotimes(\Free X) \wequi \Free(\Sym_n^\wedge X)$ and $\overline{\Sym_n^\wedge X} \wequi \Sym_n(\overline X)$, where $\Sym_n^\wedge$ denotes the symmetric power for the smash product of pointed coalgebras, together with the corresponding normed identities established above.
The classical Hopf symmetric power $\Sym_n^\hotimes(W(a))$ splits into a sum of $W(b)$-s \cite[Proposition~2.44, Example~2.46, and Remark~4.36]{bachmann-wilson}.
Consequently we have now shown that $D_n^\hotimes(W(a))$ is a retract of a free monoid and copolynomial, i.e., lies in the category $\scr C$ from above.
The equivalence of $\Sym^n(W(a))$ with a sum of $W(b)$ thus lifts to $\dHopf_R^{bp}$, concluding the proof.
\end{proof}

\begin{lemma} \label{lemm:Wled-rep}
The following mapping spectra are concentrated in degree zero and given by the usual Hopf algebra maps:
\begin{enumerate}
\item $\map_{\dHopf_R^{bp}}(W(n), W(m)^{\le d})$
\item $\map_{\dHopf_R^{bp}}(W(n)^{\le r}, \Ga(m))$, provided that $m$ is not of the form $p^en$ with $e > r$ (e.g., $m$ coprime to $p$)
\end{enumerate}
\end{lemma}
\begin{proof}
(1) Let us first treat the case $d=0$.

We may as well replace $W(n)$ by $\Free\A^1(n)$.
Noting that $\Ga$ is copolynomial (e.g. by Example \ref{ex:Poo-bipolynomial}), Lemma \ref{lemm:NHopf-comp-start}(1) shows that $\Map(\Free\A^1(n), \Ga(m))$ is discrete and given by the usual Hopf algebra maps.
It remains to show that $\Map(\Free\A^1(n), B^i \Ga(m))$ is $i$-connective.
Dualizing via Proposition \ref{prop:NHopf-basics}(1) and Proposition \ref{prop:NAlg-DAlg}, we must show that $\Map_{\DAlg_R^{aug}}(C^i D\Ga(m), D\A^1(n))$ is $i$-connective.
But $D\Ga(m)$ is not just a polynomial ring (i.e. $\Ga$ is copolynomial), the generator is even primitive (for degree reasons).
Thus $D\Ga(m)$ coincides with the free derived algebra as an $\scr E_\infty$-Hopf algebra.
It follows that $C^i D\Ga(m)$ is a free derived algebra on a generator in Chow degree $-i$ (see, e.g., \eqref{eq:LSym-CB}).
Since $D\A^1(n)$ is concentrated in homotopical degree $0$, the claim follows.

We now treat the general case.
The bifiber sequences of Remark \ref{rmk:Witt-filtrations} are exact sequences of Hopf algebras, and hence yield left exact sequences \[ 0 \to \Hom_{\Hopf_R}(W(n), \Ga(p^{d+1}m)) \to \Hom_{\Hopf_R}(W(n), W(m)^{\le d+1}) \xrightarrow{\alpha_d} \Hom_{\Hopf_R}(W(n), W(m)^{\le d}). \]
By induction on $d$, $\map_{\dHopf_R^{bp}}(W(n), \Ga(p^{d+1}m))$ and $\map_{\dHopf_R^{bp}}(W(n), W(m)^{\le d})$ are concentrated in degree zero and given by the usual Hopf algebra maps.
Using the bifiber sequence, this immediately implies that $\map_{\dHopf_R^{bp}}(W(n), W(m)^{\le d+1})$ is concentrated in degree zero and that \[ \pi_0\map_{\dHopf_R^{bp}}(W(n), W(m)^{\le d+1}) \to \pi_0\map_{\dHopf_R^{bp}}(W(n), W(m)^{\le d}) \] is surjective.
This implies that $\alpha_d$ is surjective and thus also $\pi_0 \map_{\dHopf_R^{bp}}(W(n), W(m)^{\le d+1})$ is given by the usual Hopf algebra maps.


(2)
Observe that $\Ga(m)$ is the cofree normed Hopf algebra, in the sense of Example \ref{ex:cofree-NHopf-alge}, on $R[y_m]/y_m^2$.
(Equivalently, $D\Ga(m) = W(-m)_{\ge 0}$ is a free normed Hopf algebra, for which see Example \ref{ex:sphere-bipolynomial}.)
Using this, we find $\Hom_{\Hopf}(W(p^{r+1}n), \Ga(m)) = 0$ under our assumption and so $(F^{r+1})^*: \Hom_{\Hopf}(W(n), \Ga(m)) \to \Hom_{\Hopf}(W(p^{r+1}n), \Ga(m))$ is surjective.
This implies that $\map_{\dHopf_R^{bp}}(W(n)^{\le r}, \Ga(m))$ is concentrated in degree zero and given by the kernel of $(F^{r+1})^*$, i.e., the classical Hopf algebra maps.
\end{proof}

\begin{proof}[Proof of Proposition \ref{prop:dHopfbp-disc}.]
Since $W(n) \wequi \lim_d W(n)^{\le d}$ both in $\Hopf_R$ and $\dHopf_R^{bp}$ (by Lemma \ref{lemm:Witt-filtrations}), this is immediate from Lemma \ref{lemm:Wled-rep} and its proof (the proof shows that certain maps are surjective, so that the $\lim^1$-term in the Milnor exact sequence vanishes).
\end{proof}

\begin{example} \label{ex:Ga-injective}
Lemma \ref{lemm:Wled-rep}(2) implies that for $m$ coprime to $p$, $[\Ga(n), B\Ga(m)] = 0$ and so any cofiber sequence \[ \Ga(m) \to X \to \Ga(n) \in \dHopf_R^{bp} \] is split.
\end{example}

\begin{remark}
Lemma \ref{lemm:Wled-rep} is a generalization of the fact that, over a perfect field, $W(n)^{\le r}$ is injective for $n$ coprime to $p$.
\end{remark}

\subsubsection{Derived Hopf algebras} \label{subsec:derived-hopf-alg}
Over a perfect field $k$ of characteristic $p$, the abelian category $\Hopf_k^{cn}$ of connected, evenly graded $k$-Hopf algebras has a derived category, which we denote by $\dHopf_k^{cn}$.
We elucidate the connection between $\dHopf_{\F_p}^{bp}$ and $\dHopf_{\F_p}^{cn}$.

\begin{lemma} \label{lemm:dHopf-bp-nobp}
Via identifying $W(n)$ in both categories, we obtain a fully faithful inclusion $\dHopf_{\F_p}^{bp} \hookrightarrow \dHopf_{\F_p}^{cn}$.
It is compatible with $\hotimes$ and converts $D_n^\hotimes$ to $\LSym_n^\hotimes$.
\end{lemma}
\begin{proof}
We must show that $\map_{\dHopf_{\F_p}^{cn}}(W(n), W(m))$ is discrete and given by the usual Hopf algebra maps.
The fully-faithfulness claim then follows since quickly converging sums are biproducts in both categories.

Suppose first that $n$ is coprime to $p$.
Writing $W(n) \wequi \colim_d W(n)_{\le d} \in \dHopf_{\F_p}^{cn}$, and using that $W(n)_{\le d}$ is projective \cite[Lemma~4.1]{goerss-dieudonne}, we get \[ \map_{\dHopf_{\F_p}^{cn}}(W(n), W(m)) \wequi \lim_d \map_{\dHopf_{\F_p}^{cn}}(W(n)_{\le d}, W(m)) \wequi \lim_d (\Hom_{\Hopf_{\F_p}}(W(n)_{\le d}, W(m))[0]). \]
By passing to duals in Lemma \ref{lemm:Wled-rep}(1) (using Lemma \ref{lemm:Witt-dual}), this is the same as \[ \lim_d \map_{\dHopf_{\F_p}^{bp}}(W(n)_{\le d}, W(m)) \wequi \map_{\dHopf_{\F_p}^{bp}}(W(n), W(m)), \] as needed.

Now let $n$ be general.
We can write $n = p^{i+1} n_0$ with $n_0$ coprime to $p$ and have a short exact sequence $W(n_0)_{\le i} \to W(n_0) \to W(n) \in \Hopf_{\F_p}^{cn}$.
This becomes a cofiber sequence in $\dHopf_{\F_p}^{cn}$, and also in $\dHopf_{\F_p}^{bp}$ by Lemma \ref{lemm:Witt-filtrations}.
It follows that \begin{gather*} \map_{\dHopf_{\F_p}^{cn}}(W(n), W(m)) \wequi \fib(\map_{\dHopf_{\F_p}^{cn}}(W(n_0), W(m)) \to \map_{\dHopf_{\F_p}^{cn}}(W(n_0)_{\le i}, W(m))) \\ \wequi \fib(\map_{\dHopf_{\F_p}^{bp}}(W(n_0), W(m)) \to \map_{\dHopf_{\F_p}^{bp}}(W(n_0)_{\le i}, W(m))) \wequi \map_{\dHopf_{\F_p}^{bp}}(W(n), W(m)) \end{gather*} has the desired form.

It remains to establish the last part.
All the relevant functors are obtained by restriction from sifted-cocontinuous polynomial functors, and hence are polynomial in the sense of \cite[Proposition~2.15]{k-theory-and-poly-functors}.
By \cite[Theorem~2.19]{k-theory-and-poly-functors} and the perfect-envelope description in Remark \ref{rmk:bipoly-explanation}, they are therefore determined by their restrictions to $\Hopf_R^{bp}$.
In light of Lemma \ref{lemm:bipoly-monoidal}, we must thus show that the derived tensor product of $W(n)$ and $W(m)$ is zero-truncated, and similarly for $\LSym_n(W(m))$.
Using that the derived indecomposables functor detects equivalences of connected graded derived algebras \cite[Lemma 2.11]{brantner2025formal}, for this it suffices to show that the relevant indecomposables are concentrated in degree $0$.\NB{details?}
But $Q$ is symmetric monoidal and commutes with $\LSym_n$ (both assertions are readily checked on the projective generators, e.g. using Remark \ref{rmk:hotimes-howto}), so it will suffice to show that $Q(W(n)) \otimes Q(W(m))$ and $\LSym_n(Q(W(m)))$ are concentrated in degree zero.
This is clear.
\end{proof}

\begin{remark}
The abelian category $\Hopf_{\F_p}^{cn}$ is projectively generated by the Hopf algebras $W(n)_{\le d}$ for $n$ coprime to $p$, which are also defined over any $\Z_{(p)}$-algebra $R$.
Consequently $\dHopf_{\F_p}^{cn}$ can be constructed by animating (and stabilizing) the subcategory of $\Hopf_{\F_p}^{cn}$ generated under finite sums by the $W(n)_{\le d}$, and one might propose to do the same over $R$.
This category will, however, not contain $\dHopf_R^{bp}$ (unless $R=\F_p$).
To see this, it suffices to exhibit $n, d, m, e$ (with $n$ coprime to $p$) such that $\map_{\dHopf_R^{bp}}(W(n)_{\le d}, W(m)_{\le e})$ fails to be discrete.
This is the case for $n=1$, $d=2$, $m=1$, $e=0$.
\end{remark}
\begin{proof}
Since $\map_{\dHopf_R^{bp}}(W(1)_{\le 2}, W(m))$ is discrete by passing to duals in Lemma \ref{lemm:Wled-rep}(1) and $W(1)_{\le 0} \wequi \fib(V: W(1) \to W(p))$, it suffices to show that $\Hom_{\Hopf_R}(W(1)_{\le 2}, W(1)) \xrightarrow{V} \Hom_{\Hopf_R}(W(1)_{\le 2}, W(p))$ fails to be surjective.
We claim that specifically the composite \[ u: W(1)_{\le 2} \to W(p^2)_{\le 0} \to W(p^2) \xrightarrow{F} W(p) \] cannot be factored through $V$.
To see this, we may put $R=\Z/p^2$ and write $W(1) = R[x_0, x_1, x_2, \dots]$, with $|x_i| = p^i$, and we use subsets of these variables in all the other Hopf algebras.
Then $V$ is the map $W(1) = R[x_0, x_1, \dots] \to R[x_1, x_2, \dots] = W(p)$ (modding out by $x_0$).
The map $F: W(p^2) \to W(p)$ is quite complicated, but satisfies $F(x_2) = x_1^p + px_2$ \cite[Proposition 1.8]{goerss-dieudonne}.
Thus $u$ is given by $u(x_0) = 0 = u(x_1)$, $u(x_2) = x_1^p + px_2$.
Suppose given $g: W(1)_{\le 2} \to W(1)$ with $Vg = u$.
Then $g(x_0) = ax_0$ (for degree reasons) and $g(x_1) = bx_0^p$ (for degree reasons, and since we need $g(x_1)$ to be divisible by $x_0$).
We have $\Delta(x_0) = x_0 \otimes 1 + 1 \otimes x_0$ and $\Delta(x_1) = x_1 \otimes 1 + 1 \otimes x_1 - \sum_{i=1}^{p-1}\frac{1}{p} {p \choose i} x_0^i \otimes x_0^{p-i}$ \cite[\S5.9, equation~(5.10)]{hazewinkel-witt-vectors}.
Note that the coefficient of $x_0 \otimes x_0^{p-1}$ in $\Delta(x_1)$ is $-1$.
Comparing these coefficients in $g(\Delta(x_1)) = \Delta(g(x_1))$ we see that $pb=-a^p$, which implies that both $a$ and $b$ are divisible by $p$.
It follows that the reduction mod $p$ of $g(x_2)$ is primitive, whence a scalar multiple of $w_2 \equiv x_0^{p^2} \pmod{p}$ \cite[Remark~(8.4)]{husemoller-hopf}.
But this is impossible since $g(x_2) \equiv x_1^p \pmod{x_0, p}$.
\end{proof}

\section{Main results} \label{sec:main}
In this section, we prove our main results.
First in \S\ref{subsec:cyclic-norm} we use the fact that in positive characteristic, $C_p$ acts freely on $\A^1$, to show that any $\P^1$-loop space affords a universal cyclic norm map (analogous to the map $D_n(X) \to X$ on a normed algebra $X$).
We use this in \S\ref{subsec:homology-Gm} to analyze the Toda fiber sequences of \S\ref{subsub:toda} and determine the homology of $\Omega_\Gm JS^{2n,n}$.
In \S\ref{subsec:more-dHopf} we perform some constructions with derived bipolynomial Hopf algebras.
Building on this, in \S\ref{subsec:P1-loops} we describe the homology of $\P^1$-loops on certain normed monoids in terms of derived primitives and indecomposables.
This allows us in \S\ref{subse:model} to prove our main theorem, i.e., to establish a new model of the motivic Eilenberg--Mac Lane spaces at the characteristic.

Throughout, $k$ is a perfect field of characteristic $p$.

\subsection{The cyclic norm} \label{subsec:cyclic-norm}
In this subsection we show that if $X \in \Spc(k)_*$, then $\Omega_{\P^1} X$ admits a normed cyclic multiplication, i.e., a map $D_{C_p}(\Omega_{\P^1}X) \to \Omega_{\P^1}X$.
Moreover, if $X$ was a normed monoid, we show that this cyclic multiplication coincides with the one coming from the norm structure.
This is in sharp contrast to classical homotopy theory (and motivic homotopy theory in characteristic $0$), where $\Omega^2 X$ in general only has an $\scr E_2$-multiplication.

We will build our cyclic norms out of conorms for $\P^1$.
The basis of this is the following, which holds for general schemes.
\begin{lemma} \label{lemm:conorm-basis}
Let $f: B \to A$ be a finite étale morphism of schemes and $X \in \Spc(A)_*$.
Then
\[
f_*\Omega_{\P^1_B}(f^*X)
\wequi \iMap_*(f_\vee\P^1_B,X).
\]
Here we write $f_\vee:\Spc(B)_*\to\Spc(A)_*$ for the left adjoint of $f^*$.
\end{lemma}
\begin{proof}
Via Yoneda and the adjunction $f_\vee\dashv f^*$, this follows from the pointed projection formula
\[
f_\vee(\P^1_B)\wedge T
\wequi f_\vee(\P^1_B\wedge f^*T)
\]
(which holds on pointed representables and extends by colimits).
\end{proof}

From now on we work in characteristic $p$ again.
Identify $C_p$ with the additive group of $\F_p$, and let it act on $\scr V=\scr O^{\oplus 2}$ by $a\cdot(u,v)=(u+av,v)$.
On $\P(\scr V)$, the chart $\scr L=\{v\ne0\}$ has coordinate $x=u/v$ and action $x\mapsto x+a$.
Thus $\scr L$ is an equivariant affine line (with a non-linear $C_p$-action).
Embed $\scr B:=C_p$ into $\scr L$ by $a\mapsto a$; its ideal is generated by the invariant polynomial $x^p-x$.
The point $\infty=[1:0] \in \P(\scr V)$ is fixed.
Near infinity the coordinate $z=v/u$ transforms as $z\mapsto z/(1+az)$, so its class modulo $z^2$ gives an invariant trivialization of the conormal line at infinity, and hence of $N_{\infty/\P(\scr V)}$.
Likewise, the class of $x^p-x$ in the conormal line of $\scr B$ gives an equivariant trivialization $N_{\scr B/\P(\scr V)}\simeq\scr B\times\A^1$, with trivial action on the second factor.

Now suppose given a scheme $A/\F_p$ and a cyclic étale cover $B \to A$, i.e., a map of stacks $A \to B_\et C_p$.
Then $B = A \times_{B_\et C_p} \scr B$.
Pulling back all the other data from $B_\et C_p$, we obtain:
\begin{itemize}
\item a rank-two vector bundle $V$ with a section $\infty_A$ of $\P(V)$ and an embedding of $B$ into $L=\P(V)\setminus\infty_A$,\footnote{$L$ is an affine-line torsor, but we do not really need this.}
\item a trivialization of $N_{\infty_A/\P(V)}$, and
\item a trivialization of $N_{B/\P(V)}$.
\end{itemize}
Moreover, all this data is functorial in the cyclic cover.

\begin{construction}
Write $f:B\to A$ and point $\P(V)$ at $\infty_A$.
In the quotient $\P(V)/L$, let $s$ denote the image of $\infty_A$, rather than the usual quotient basepoint.
Purity identifies $s$ with the zero section of $\Th(N_{\infty_A/\P(V)})$.
Using the chosen normal trivialization and inversion on $\P^1$, we obtain an equivalence $(\P^1_A,\infty)\wequi(\Th(N_{\infty_A/\P(V)}),0)$.
The \emph{co-norm corresponding to $B/A$} is the composite \begin{gather*} c_{B/A}: (\P^1_A,\infty) \wequi (\Th(N_{\infty_A/\P(V)}),0) \wequi (\P(V)/L,s) \wequi (\P(V),\infty_A) \\ \to \P(V)/(\P(V)\setminus B) \wequi \Th(N_{B/\P(V)}) \wequi B_+\wedge\P^1_A \wequi f_\vee\P^1_B. \end{gather*}
Here the third equivalence is the inverse of the quotient map, which is a motivic equivalence because $L$ is contractible.\NB{It's the complement of $\infty$ in $\P(V)$... Also it's an affine line torsor, so Zariski-locally isomorphic to $\A^1$.}
The collapse map is pointed since $\infty_A$ lies in $\P(V)\setminus B$, and the remaining equivalences use purity \cite[\S3.2]{A1-homotopy-theory} and the chosen trivialization of $N_{B/\P(V)}$.
\end{construction}
Note that this construction is functorial in $B/A$, since the data entering it was.
Combining with Lemma \ref{lemm:conorm-basis}, for any $X \in \Spc(A)_*$ we obtain the cyclic norm \[ N_{B/A}: f_* \Omega_\P^1(f^*X) \wequi \iMap_*(f_\vee \P^1, X) \xrightarrow{c_{B/A}^*} \iMap_*(\P^1, X) \wequi \Omega_{\P^1} X. \]
Now recall that, since we are using the cartesian norm structure, the cyclic extended power of $X \in \Spc(S)$ is defined as \[ D_{C_p}(X) = \colim_{B/A \in B_\et C_p} (A \to S)_\sharp(B \to A)_*(X|_B). \]

\begin{definition} \label{def:cyclic-norm}
For $X \in \Spc(S)_*$, $S$ a scheme of characteristic $p$, the \emph{universal cyclic norm} on $\Omega_{\P^1}(X)$ is the map \[ N: D_{C_p}(\Omega_{\P^1} X) \wequi \colim_{B/A \in B_\et C_p} (A \to S)_\sharp \iMap_*((B \to A)_\vee \P^1, X|_A) \to \Omega_{\P^1} X \in \Spc(S) \] induced via the functoriality of the individual cyclic norms.
\end{definition}

Let $Y \in \NMon(\Spc(S))$.
Then $Y$ has norms for any finite étale map, and so in particular, given a cyclic cover $f: B \to A$, we obtain the norm map $f_*(Y|_B) \to Y|_A$.
Taking the colimit over all such, we obtain the universal cyclic norm on $Y$, $N_Y: D_{C_p}(Y) \to Y \in \Spc(S)$.

\begin{proposition} \label{prop:compare-norms}
Let $S$ be a scheme of characteristic $p$ and $Y \in \NMon(\Spc(S))$, so that also $\Omega_{\P^1} Y \in \NMon(\Spc(S))$.
Then the two universal cyclic norms on $\Omega_{\P^1} Y$ above coincide.
\end{proposition}
\begin{proof}
It suffices to show the claim for the individual cyclic norms, functorially in $B/A$.
Thus let $f: B \to A$ be a cyclic cover, and write $N_{B/A}, N'_{B/A}: f_*f^* \Omega_{\P^1} Y \to \Omega_{\P^1} Y$ respectively for the cyclic norm constructed above and the one coming from the normed monoid structure.
Then $N_{B/A}$ is induced by the map $T(c_{B/A}): T(\P^1) \to T(f_\vee \P^1) \wequi f_\sharp T(\P^1)$, whereas $N'_{B/A}$ is induced by the map $T(\P^1) \to f_\sharp T(\P^1)$ dual to $f$.
It suffices to show that these two agree.
Recall that $\NMon(\Spc(S))$ is $\fet$-semiadditive, in the sense that $f_\sharp \wequi f_*$ on this category \cite[Proposition~6.8]{bachmann-MGM}, and under this equivalence, the map $T(\P^1) \to f_\sharp T(\P^1) \wequi f_* T(\P^1)$ dual to $f$ is just the unit of adjunction.
It thus suffices to observe that the composite $\P^1 \xrightarrow{c_{B/A}} f_\vee \P^1 \to f_* \P^1$ (the last map being obtained from the base point of $\P^1$ via the ambidexterity construction \cite[\S5.3, (5.20)]{hoyois-equivariant}) is also the unit of adjunction.\NB{details...?}
\end{proof}

\subsection{Homology of some $\Gm$-loop spaces} \label{subsec:homology-Gm}
Our goal now is to describe $H\F_p[\Omega^{1,1}JS^{2n,n}]$.
We know quite a lot about this object, from a variety of perspectives.
Here we summarize the main properties we will need for our applications.

\begin{theorem}\label{thm:loops2S3-new}
Let $k$ be a perfect field of characteristic $p$ and $n \ge 2$.
\begin{enumerate}
\item $H\Z_{(p)}[\Omega^{1,1}JS^{2n,n}]$ is cellular.
\item There are classes \[ a_0, a_1, \dots, b_1, b_2, \dots \in \pi_{**}H\F_p[\Omega^{1,1}JS^{2n,n}] \] with \[ |a_i| = (2np^i-1,np^i-1), |b_i|=(2np^i-2,np^i-1) \] such that \[ \pi_{**}H\F_p[\Omega^{1,1}JS^{2n,n}] \cong \F_p[a_0, a_1, \dots, b_1, b_2, \dots]/(a_0^2,a_1^2,\dots) \] as primitively generated Hopf algebras.
\item We have $\beta a_i = b_i$ for $i \ge 1$.
\item The co-unit map $\Sigma^{1,1}\Omega^{1,1}JS^{2n,n} \to JS^{2n,n}$ sends $\Sigma^{1,1}a_i$ to $x^{p^i}$, where $x \in H_{2n,n}(JS^{2n,n}, \Z_{(p)})$ is the bottom cell.
\end{enumerate}
\end{theorem}

\begin{remark}
We strongly suspect that this theorem remains true for other fields (with the coefficient ring $\pi_{**} H\F_p$ inserted appropriately where necessary).
We concentrate on the stated case since this is what we will need later.
\end{remark}

\begin{remark}
Consider the sequence of maps \[ \Omega^{1,1} JS^{2n,n} \to \Omega^{1,1} JS^{2np,np} \to \Omega^{1,1} JS^{2np^2,np^2} \to \dots \] given by the Hopf maps from Construction \ref{cons:H-2}.
Proposition \ref{prop:Toda-fibrations} identifies the ($p$-local) fibers as $\Omega^{1,1}J_{p-1}S^{2np^i,np^i}$.
Theorem \ref{thm:loops2S3-new} thus asserts that these fiber sequences split in homology.
\end{remark}

We now embark on the proof of Theorem \ref{thm:loops2S3-new}.
We first dispense with the cellularity claim.
\begin{lemma} \label{lemm:loops2S3-cellular}
$H\Z_{(p)}[\Omega^{1,1}J_{p-1}S^{2n,n}]$ and $H\Z_{(p)}[\Omega^{1,1}JS^{2n,n}]$ are cellular.
\end{lemma}
\begin{proof}\NB{perhaps more details?}
Noting that $\Omega^{1,1}J_{p-1}S^{2n,n}$ is connected, we have \[ \Omega^{1,1}J_{p-1}S^{2n,n} \wequi B \Omega^{2,1} J_{p-1} S^{2n,n}, \] and so in order to show that $H\Z_{(p)}[\Omega^{1,1}J_{p-1}S^{2n,n}]$ is cellular it suffices to show the same for $H\Z_{(p)}[\Omega^{2,1}J_{p-1}S^{2n,n}]$.
Rotating the fiber sequence of Proposition \ref{prop:Toda-fibrations} (noting that $\Omega^{2,1} S^{2np-1,np}$ is connected), for this it suffices to show that $H\Z_{(p)}[\Omega^{1,1}S^{2n-1,n}]$ and $H\Z_{(p)}[\Omega^{3,1}S^{2np-1,np}]$ are cellular.
For the first term, see Corollary \ref{cor:Omega11-S2n-1}.
For the second, set $X = \Omega^{3,1}S^{2np-1,np}$.
Then $B^2 X \wequi \Omega^{1,1}S^{2np-1,np}$, whence also $B^2 H\Z_{(p)}[X]$ is cellular.
This implies that $H\Z_{(p)}[X]$ is cellular, by e.g. \cite[Lemma 5.9]{bachmann-wilson}.

We will argue somewhat differently for $\Omega^{1,1}JS^{2n,n}$.
Write $X \mapsto X^\cell$ for the right adjoint of the inclusion of cellular objects in $\Mod_{H\Z_{(p)}}(k)$, i.e., the cellularization functor.
This is lax symmetric monoidal.
We call the map $X^\cell \to X$ the \emph{cellularization map for $X$}.
For any $N$ we shall prove: for all $i \ge 0$ the cellularization map for $H\Z_{(p)}[\Omega^{i+1,1} J S^{2np^i,np^i}]$ is $N$-connective (i.e., its cofiber $R$ satisfies $\ul\pi_j(R)_* = 0$ for $j < N$).
Noting that $JS^{2np^i,np^i}$ is $np^i$-connective, we see that the assertion is vacuous for $i$ sufficiently large.
We can now prove the claim by descending induction on $i$.
Employing the second Toda fibration of Proposition \ref{prop:Toda-fibrations}, we can write \[ H\Z_{(p)}[\Omega^{i+1,1}JS^{2np^i,np^i}] \wequi H\Z_{(p)}[\Omega^{i+1,1}J_{p-1}S^{2np^i,np^i}] \sslash H\Z_{(p)}[\Omega^{i+2,1} JS^{2np^{i+1},np^{i+1}}]. \]
Note that $H\Z_{(p)}[\Omega^{i+1,1}J_{p-1}S^{2np^i,np^i}]$ is cellular: as above (i.e., by \cite[Lemma 5.9]{bachmann-wilson}) it suffices to prove this after applying arbitrarily many bar constructions, but this reduces to the cellularity of \[ B^i H\Z_{(p)}[\Omega^{i+1,1}J_{p-1}S^{2np^i,np^i}] \wequi H\Z_{(p)}[B^i \Omega^{i+1,1}J_{p-1}S^{2np^i,np^i}] \wequi H\Z_{(p)}[\Omega^{1,1}J_{p-1}S^{2np^i,np^i}], \] which was already established.
We conclude by the observation that if $X, Y \in \Mod_{H\Z_{(p)}}(k)$ are connective and the cellularization maps of $X$ and $Y$ are $N$-connective, then the same is true for $X \otimes Y$.
(Indeed $X^\cell \to X$ is $N$-connective if and only if there exists a cellular spectrum $X'$ and a cofiber sequence $X' \to X \to R$ with $R$ $N$-connective---in the one direction take $X'=X^\cell$, and in the other note that $X/X^\cell \wequi R/R^\cell$ is $N$-connective.
But now $X^\cell \otimes Y^\cell \to X \otimes Y$ is $N$-connective if $X^\cell \to X$ and $Y^\cell \to Y$ are.)
\end{proof}

From now on, we thus need only make assertions about motives with $\F_p$-coefficients, and only for objects which are cellular and effective.
As discussed in \S\ref{subsub:charp}, we have \[ \Mod_{H\F_p}(k)^{\cell\eff} \wequi \Mod_{\F_p}^\N, \] i.e., we are really just working with graded $\F_p$-module spectra.
We call the external degrees \emph{weights}.
Given an object of $\Mod_{\F_p}^\N$ which is augmented or pointed (such as an algebra or a coalgebra), we call it \emph{connected} if the only term in weight $0$ is the (co)unit.
Algebraic operations on such objects have very favorable properties which we will utilize.
Here are some examples.
\begin{itemize}
\item Let $E \in \cAlg_{\scr E_1}(\Mod_{\F_p}^\N)^\cn_{\F_p/}$.
  Then the cobar construction $CE = \CB(\F_p, E, \F_p)$ is, in every individual weight, given by a finite limit.
  This implies, for example, that $E$ is recovered as the bar construction $BCE$.
\item More generally, given $0 \le m,n \le \infty$, bar--cobar duality restricts to an equivalence on connected objects \cite[Theorem~A.6 and its proof]{burklund2022multiplicative}
\[
B: \Alg_{\scr E_{m+1}}(\cAlg_{\scr E_n})^\cn
\adj
\Alg_{\scr E_m}(\cAlg_{\scr E_{n+1}})^\cn : C.
\]
\item Also relative Bar-Cobar duality works well.
  That is, given $E \in \cAlg_{\scr E_1}(\Mod_{\F_p}^\N)^\cn_{\F_p/}$ and $M \in \cMod_E$, $\CB(M, E, \F_p)$ is a $CE$-module (recall from the previous point that $CE$ is an algebra).
  Conversely, given a $CE$-module $N$, we can form the two-sided bar construction $B(N, CE, \F_p) = N \otimes_{CE} \F_p$ and obtain a comodule.
  These constructions upgrade to an equivalence \cite[Theorem~4.18]{matsuoka2014koszul} \[ B(\ph, CE, \F_p): \Mod_{CE} \adj \cMod_E: C(\ph, E, \F_p). \]
\end{itemize}

\begin{remark} \label{rmk:automatic-augmentation}
For $E \in \Alg(\Mod_{\F_p}^\N)^\cn$, we have $\Map(E, \F_p) \wequi *$ (use that the inclusion $\Mod_{\F_p} \to \Mod_{\F_p}^\N$ of the degree zero part has a symmetric monoidal right adjoint), i.e., every connected algebra has a unique augmentation.
Similarly, every connected coalgebra has a unique pointing.
We may thus sometimes drop the augmentation/pointing from the notation.
\end{remark}

\begin{example} \label{ex:rotate-fiber-seq}
Let $E \to B \in \cCAlg(\Mod_{\F_p}^\N)^\cn$.
Then $F := \CB(E, B, \F_p)$ upgrades to a coalgebra, namely the fiber of $E \to B$ in $\cCAlg(\Mod_{\F_p}^\N)^\cn$.
By the third point above, $F$ acquires an action by $CB$, and $E \wequi F \sslash CB = B(F, CB, \F_p) = F \otimes_{CB} \F_p$.
Note that also $CB \wequi \fib(F \to E)$.
(The dual statement for algebras is more apparent, and the automatic finiteness of limits involved implies that the dual argument works.)
\end{example}

In order to facilitate our computations, we will express answers in terms of \emph{derived algebras}.
Recall from \cite[\S4.2]{raksit2020hochschild} the free derived algebra functor \[ \LSym: \Mod_{\F_p}^\N \to \Mod_{\F_p}^\N. \]
It has the following properties.
\begin{itemize}
\item It upgrades to a symmetric monoidal functor from $(\Mod_{\F_p}^\N, \oplus)$ to $(\Mod_{\F_p}^\N, \otimes)$.
  Since the first symmetric monoidal structure is additive, every object $X$ is uniquely an algebra and a coalgebra, i.e., an $\scr E_\infty$-Hopf algebra (for $\oplus$).
  Consequently $\LSym(M)$ is an $\scr E_\infty$-Hopf algebra for $\otimes$.
\item If $M \in \Mod_{\F_p}^\N$ is connected in the sense that the weight zero part vanishes, then $\LSym(M)$ is connected in the sense of pointed objects (i.e., the only part in weight zero is the unit).
\item Since $\LSym$ converts $\oplus$ into $\otimes$, we have \begin{equation} \label{eq:LSym-CB} \LSym(\Sigma M) \wequi B\LSym(M) \text{ and } \LSym(\Sigma^{-1} M) \wequi C\LSym(M) \end{equation} (provided $M$ is connected).
\item If $M$ is concentrated in Chow degree zero, then so is $\LSym(M)$, and it is given by the usual symmetric algebra on $M$, with the usual (primitively generated) coalgebra structure.
\end{itemize}

\begin{remark}
Given a generator $x$ in bidegree $(m,n)$, we write $\LSym(x)$ for $\LSym(\Sigma^{m,n}\F_p)$.
Similarly for more than one generator, we use notation like $\LSym(x,y,z)$.
\end{remark}

\begin{example} \label{ex:OmegaGmS2n-1}
Since $\Omega_\Gm S^{2n-1,n}$ is a space, $H\F_p[\Omega_\Gm S^{2n-1,n}]$ is an $\scr E_\infty$-coalgebra.
In fact it follows from Corollary \ref{cor:Omega11-S2n-1} that \[ H\F_p[\Omega_\Gm S^{2n-1,n}] \wequi \LSym(b_{2n-2,n-1}) \in \cCAlg(\Mod_{\F_p}^\N). \]
(Indeed the corollary shows that the motive is cellular, concentrated in Chow degree zero, and has homotopy groups given by the right hand side.
But taking homotopy groups is an equivalence of categories on Chow degree zero.)
\end{example}

\begin{lemma} \label{lemm:cAlg-fiber-splitting}
Let $B \in \cAlg_{\scr E_1}(\Mod_{\F_p}^\N)^\cn$ and $E \in \cMod_B$.
Set $F = \CB(E, B, \F_p)$.
Given any splitting $E \to F \in \Mod_{\F_p}^\N$, the induced map $E \to F \otimes B$ (to a cofree comodule) is an equivalence.

Dually, given $A \in \Alg_{\scr E_1}(\Mod_{\F_p}^\N)^\cn$ and $E \in \Mod_A$, any section $F := E \otimes_A \F_p \to E \in \Mod_{\F_p}^\N$ induces an equivalence $F \otimes A \to E \in \Mod_A$.
\end{lemma}
\begin{proof}
Indeed we have $\CB(F \otimes B, B, \F_p) \wequi F$ and so applying relative cobar duality to $E \to F \otimes B$ yields an equivalence, as needed.
The dual statement is proved in the same way.
\end{proof}

We now analyze the Toda fiber sequences in this language.
\begin{lemma} \label{lemm:identify-Jp-1}
We have an equivalence \[ H\F_p[\Omega^{2,1} J_{p-1} S^{2np^i,np^i}] \wequi C\LSym(a_i, b_{i+1}) \in \cCAlg(\Mod_{\F_p}^\N). \]
Here $a_i$ denotes the formal suspended generator, so that $Ca_i$ is the actual bottom class in the displayed motive.
(The $a_i, b_i$ will turn out to be the classes from Theorem \ref{thm:loops2S3-new}.)
\end{lemma}
\begin{proof}
Applying $H\F_p[\Omega^{1,1}(\ph)]$ to the first fiber sequence from Proposition \ref{prop:Toda-fibrations} yields by Corollary \ref{cor:EM} (as well as the discussion in \S\ref{subsec:p-loc}) a fiber sequence of coalgebras.
Using Example \ref{ex:OmegaGmS2n-1}, we see that this fiber sequence takes the form \[ C\LSym(a_i) \to H\F_p[\Omega^{2,1} J_{p-1} S^{2np^i,np^i}] \to C\LSym(b_{i+1}). \]
Applying Lemma \ref{lemm:cAlg-fiber-splitting}, it suffices to produce a splitting $H\F_p[\Omega^{2,1} J_{p-1} S^{2np^i,np^i}] \to C\LSym(a_i) \in \cCAlg(\Mod_{\F_p}^\N)$.
Note that by rotating the fiber sequence (see Example \ref{ex:rotate-fiber-seq}) and using that bar and cobar constructions in our setting are given by finite limits and colimits in each weight, we see that $H\F_p[\Omega^{2,1} J_{p-1} S^{2np^i,np^i}]$ is bidegreewise finitely generated.
The same is true for the other terms, as well.
Consequently it will suffice to pass to duals and construct a section on cochains.

Recall that for any $X \in \Spc(k)$, $H\F_p[X]$ is a normed coalgebra, and its dual is a normed algebra.
But cellular normed $H\F_p$-algebras are derived algebras: see Proposition \ref{prop:NAlg-DAlg}.
We are thus in the following situation.
We are given a cofiber sequence \[ B\Gamma_{\F_p}(b_{i+1}') \to DH\F_p[\Omega^{2,1} J_{p-1} S^{2np^i,np^i}] \to \Gamma_{\F_p}(a_i') \in \DAlg(\Mod_{\F_p}^\Z), \] and we are trying to construct a section of the right hand map.
(Indeed the dual of a primitively generated polynomial ring is a divided power algebra, so $D\LSym(x) \wequi \Gamma_{\F_p}(x')$ for $x$ in Chow degree zero.)
Now $\Gamma_{\F_p}(x) \wequi \F_p[x, x^{(p)}, x^{(p^2)}, \dots]/((x^{(p^i)})^p)$ and so $B\Gamma_{\F_p}(b_{i+1}') \wequi \LSym(c_0, d_0, c_1, d_1, \dots)$ where the $c_i$ have Chow degree $1$ and the $d_i$ have Chow degree $2$.
It follows that $\pi_{2*,*} DH\F_p[\Omega^{2,1} J_{p-1} S^{2np^i,np^i}] \wequi \Gamma_{\F_p}(a_i')$ (i.e., the cofiber sequence only modifies positive Chow degrees) and so we have the desired section on homotopy groups.
But this already induces a section as derived algebras, as follows from the presentation of $\Gamma_{\F_p}(x)$ in terms of free derived algebras given above (i.e., a derived algebra map $\Gamma_{\F_p}(x) \to A$ is the same as a family of elements $y_i$ in $\pi_{2*,*} A$ together with nullhomotopies of the $y_i^p$, and then $x^{(p^i)} \mapsto y_i$).
\end{proof}

Before we analyze the second Toda fiber sequence, we need some further preparation.
We must now determine multiplicative structures, about which Lemma \ref{lemm:identify-Jp-1} says nothing.
These will arise in the following ways:
\begin{itemize}
\item If $X$ is an $\scr E_\infty$-coalgebra, then $C^n X$ is an $\scr E_n$-algebra for $n \ge 1$.
  In particular, $\pi_{**} C^n X$ is a ring, carrying $\scr E_n$-operations.
\item If $X \in \Spc(k)_*$, then $\Omega^n X$ is an $\scr E_n$-monoid, and so $H\F_p[\Omega^n X]$ is also an $\scr E_n$-algebra.
  Note that if $X$ is $(n+1)$-connected, then $H\F_p[\Omega^n X] \wequi C^n H\F_p[X]$ and the two constructions coincide.
\item If $X$ is an $\scr E_\infty$-algebra (e.g., a derived algebra, e.g., $H\F_p[M]$ where $M$ is a normed monoid), then $X$ carries all $\scr E_\infty$-operations.
\item If $X \in \Spc(k)_*$, then $\Omega^{2,1} X$ carries the universal cyclic norm of Definition~\ref{def:cyclic-norm}, inducing an operation $D_{C_p}(H\F_p[\Omega^{2,1}X]) \to H\F_p[\Omega^{2,1}X]$.
  This also provides $\pi_{**} H\F_p[\Omega^{2,1}X]$ with all $\scr E_\infty$-operations (and a ring structure, $\Omega^{2,1}X$ being a loop space).
\item If $X \in \NMon(\Spc(k))$, then the universal cyclic norm on $\Omega^{2,1}X$ agrees with the one induced by its normed-monoid structure by Proposition~\ref{prop:compare-norms}, and hence the resulting power operations on $\pi_{**} H\F_p[\Omega^{2,1}X]$ agree.
\end{itemize}

In terms of these operations, we summarize the following about the structure of free derived algebras.
\begin{itemize}
\item For $d \le 0$, the free derived algebra $\LSym(x_d)$ has homotopy groups a graded polynomial ring (i.e., in characteristic $2$, just a polynomial ring, and in characteristic $p\ne 2$, a tensor product of a polynomial ring on even generators and an exterior algebra on odd generators), freely generated by appropriate operations (in fact, by $\scr E_{1-d}$-operations) \cite[Theorem~4.1.1]{priddy1973derived}.
\item In particular, $\pi_{**}\LSym(x_d)$ is a primitively generated, graded polynomial ring: its generators are obtained from the primitive class $x_d$ by power operations, which preserve primitive classes by naturality and additivity.
  The primitives are given by $(Px_d)^{p^r}$ for $r \ge 0$, where $Px_d$ denote the power operation generators.
\end{itemize}

\begin{remark} \label{rmk:loops2S3-derived-algebra}
The main part of Theorem \ref{thm:loops2S3-new} thus asserts that \[ H\F_p[\Omega^{1,1}JS^{2n,n}] \wequi \LSym(a_0, a_1, \dots, b_1, b_2, \dots) \in \DAlg(\Mod_{\F_p}^\N), \] and that at least on $\pi_{**}$, this equivalence is compatible with the coalgebra structures.
(Indeed the left hand side is a derived algebra by Corollary~\ref{cor:identify-J-T} and Proposition~\ref{prop:NAlg-DAlg}, so receives a map from the right hand side by universal properties, and this map is then an equivalence by computing the homotopy rings.)
This suggests that also \[ C^i H\F_p[\Omega^{1,1}JS^{2n,n}] \wequi C^i \LSym(a_0, a_1, \dots, b_1, b_2, \dots) \] for $i \ge 1$.
In fact, such an equivalence for any $i$ implies it for $i' < i$ by applying iterated bar constructions.
\end{remark}

\begin{remark} \label{rmk:cobar-cyclic-norm}
Let $D_{C_p}^\naive: \cCAlg(\Mod_{\F_p}^\N) \to \cCAlg(\Mod_{\F_p}^\N)$ be the naive extended power construction $D_{C_p}^\naive(X) = (X^{\otimes p})_{hC_p}$.
By a coalgebra with naive cyclic multiplication, we mean a pair $(X, m)$ with $X \in \cCAlg(\Mod_{\F_p}^\N)_*$ and $m: D_{C_p}^\naive(X) \to X \in \cCAlg(\Mod_{\F_p}^\N)_*$.
Such objects are stable under cobar: the canonical exchange transformation for a limit yields \[ D_{C_p}^\naive(CX) \to C(D_{C_p}^\naive(X)) \to CX. \]
This construction is functorial in the pair $(X, m)$.
\end{remark}

\begin{example} \label{ex:cobar-cyc-mult}
If $X$ comes from a normed coalgebra with a normed cyclic multiplication (for example, $X = H\F_p[\Omega^{2,1}Y]$) then the canonical map $D_{C_p}^\naive(X) \to D_{C_p}(X)$ equips $X$ with a naive cyclic multiplication, and hence the same is true for cobars on $X$.
\end{example}

\begin{remark} \label{rmk:cobar-cyclic-commute}
Note that the comparison
\[
H\F_p[\Omega\Omega_{\P^1}Y]\longrightarrow CH\F_p[\Omega_{\P^1}Y]
\]
respects cyclic multiplication, where on the right hand side we use the construction of Remark \ref{rmk:cobar-cyclic-norm}.
Indeed, the geometric multiplication is precomposition with the conorm of $\P^1$, which commutes with taking loops in $Y$.
Moreover, $H\F_p[\ph]$ commutes with naive extended powers, and the comparison from the motive of a limit to the limit of the motives respects the limit-exchange construction of Example~\ref{ex:cobar-cyc-mult}.
\end{remark}

Now consider again $H\F_p[\Omega^{2,1}J_{p-1}S^{2n,n}]$, which in Lemma \ref{lemm:identify-Jp-1} we identified with $C\LSym(a_0, b_1)$ \emph{as coalgebras}.
By the above discussion, this is in fact an $\scr E_1$-ring with cyclic multiplication, and so in particular its homotopy is a ring with power operations.
In much of what follows, we will be working to determine this structure.
We begin by determining a small part of that structure.
\begin{lemma}
The map $C\LSym(a_i) \to H\F_p[\Omega^{2,1}J_{p-1}S^{2np^i,np^i}]$ lifts to $\Alg_{\scr E_1}(\cCAlg(\Mod_{\F_p}^\N))$.
\end{lemma}
\begin{proof}
The functor $\tau_{c \ge 0}$ is symmetric monoidal on coconnective $\F_p$-module spectra, and so \[ \tau_{c \ge 0} H\F_p[\Omega^{2,1}J_{p-1}S^{2np^i,np^i}] \to H\F_p[\Omega^{2,1}J_{p-1}S^{2np^i,np^i}] \in \Alg_{\scr E_1}(\cCAlg(\Mod_{\F_p}^\N)). \]
The object on the left hand side is discrete, and so, it will suffice to show that the iterated loop multiplication $z^{*k}$ of $z=Ca_i$ with itself yields $z^k$.

To do this, we may work over $k=\F_p$, in which case we can consider the conormed coalgebra $H\Z_{(p)}[\Omega^{2,1}J_{p-1}S^{2np^i,np^i}]$.
We can repeat the discussion in Lemma \ref{lemm:identify-Jp-1} integrally (the key points being that cellular normed $\Z_{(p)}$-algebras over $\F_p$ are derived algebras by Proposition~\ref{prop:NAlg-DAlg}, and that $\Gamma_{\Z_{(p)}}(a_i)$ is the quotient of a polynomial algebra by a regular sequence) to learn that \[ H\Z_{(p)}[\Omega^{2,1}J_{p-1}S^{2np^i,np^i}] \wequi C\LSym_{\Z_{(p)}}(a_i, b_{i+1}) \] as coalgebras.
Using that the loop multiplication is a coalgebra map, from this one deduces that $\Delta(z^{*k}) = \Delta(z)^{*k}$.
For degree reasons $z^{*k} = c_k z^k$ for some $c_k \in \Z_{(p)}$.
Comparing coefficients on both sides and using that $\Z_{(p)}$ is torsion-free, we inductively see that $c_k = 1$ for all $k$.
\end{proof}

\begin{remark} \label{rmk:a-family}
Since $C\LSym(a_i)$ is in fact free as an $\scr E_1$-algebra, it follows that $C^r \LSym(a_i)$ is free as an $\scr E_r$-algebra (for $r \ge 1$).
Suppose now given a map $f: \Omega^{1,1}J_{p-1}S^{2np^i,np^i} \to M$ where $M$ is a normed monoid.
Then $C^r H\F_p[f]$ is an $\scr E_r$-algebra map, and so is compatible with all the $\scr E_r$-operations.
On $C^r\LSym(a_i)$, these operations generate the homotopy of the free model, and on $C^r H\F_p[M]$, these are just the derived algebra operations.
The conclusion is that, when mapping $\Omega^{1,1}J_{p-1}S^{2np^i,np^i}$ into a normed monoid and then cobaring down, we understand all the $\scr E_r$-operations on the classes generated by $a_i$.
\end{remark}

We want to do something similar for $b_{i+1}$.
This is a bit more complicated.
\begin{lemma} \label{lemm:b-family}
Suppose given a map $\alpha: \Omega J_{p-1}S^{2np^i,np^i} \to M \in \Spc(k)_*$, where $M \in \NMon(\Spc(k))$ is simply-connected.
Then for $r \ge 1$, $\pi_{**} C^r H\F_p[\Omega^{1,1} \alpha]$ is compatible with all power operations on $C^{r+1} b_{i+1}$ coming from the $\P^1$-loop structure on $\Omega^{2,1} J_{p-1}S^{2np^i,np^i}$ (not the $\scr E_r$-operations obtained by cobaring down).
\end{lemma}
\begin{proof}
Replacing $n$ by $np^i$, we may assume $i=0$.
For $n>2$, both $\Omega J_{p-1}S^{2n,n}$ and $M$ are simply-connected, as are their $\Gm$-loop spaces by contraction.
It follows (by Corollary~\ref{cor:EM}) that $CH\F_p[\Omega^{1,1} \alpha] \wequi H\F_p[\Omega^{2,1} \alpha]$ is a $\P^1$-loop map, so compatible with the operations.
This property is retained by further cobars by Remark \ref{rmk:cobar-cyclic-commute}.

A problem arises when $n=2$, because then $\Omega J_{p-1}S^{2n,n}$ is no longer simply connected, and we can no longer appeal to Eilenberg--Moore.\NB{The comparison map is, in fact, not an equivalence.}
Put $J=J_{p-1}S^{4,2}$, $V=\Omega^{2,1}J$, and $B_b=C\LSym(b_1)$.
It suffices to find a $2$-connected space $X \to J$ and a coalgebra map \[ B_b \to H\F_p[\Omega^{2,1}X] \] whose composite to $H\F_p[V]$ preserves $b$.
Indeed, in \[ C^rH\F_p[\Omega^{2,1}X]\xrightarrow{f_1} C^rH\F_p[V]\xrightarrow{f_2} C^rH\F_p[\Omega^{1,1}M] \] both $f_1$ and $f_2f_1$ respect operations (the former being induced by a $\P^1$-loop map, the latter by arguing as above), and by our construction of the map from $B_b$, $C^{r+1}b_1$ lifts to the source.
It thus remains to construct $X$ and the auxiliary data.

\emph{If $p=2$,} take the Hopf construction of multiplication on $\SL_2\simeq S^{3,2}$:
\[
\nu:X=S^{7,4}\longrightarrow S^{4,2}=J
\]
(see \cite[\S4.6]{dugger-isaksen-hopf}).
By Example~\ref{ex:OmegaGmS2n-1} and Corollary~\ref{cor:EM}, $H\F_2[\Omega^{2,1}X] \wequi B_b$, giving the desired map.
We just need to verify that our maps preserve $b$.
This is exactly saying that $\nu$ has Hopf invariant $1$, and is verified by repeating the classical argument (see, e.g., \cite[\S4B]{hatcher606algebraic}).

\emph{If $p$ is odd,} work $p$-locally and let
\[
X=\fib\bigl(J\longrightarrow K(\ul K_2^{MW},2)\bigr)
\]
be the $2$-connected cover.
Here the higher James cells do not change $\ul\pi_2$, and $\ul\pi_2(S^{4,2})=\ul K_2^{MW}$ by \cite[Theorem~6.40]{A1-alg-top}.
Since $(\ul K_2^{MW})_{-1}=\ul K_1^{MW}$ and the Witt-theoretic kernel of $\ul K_1^{MW}\to\Gm$ is $2$-primary over an unorderable field, we obtain
\[
\widetilde V:=\Omega^{2,1}X\longrightarrow V
\xrightarrow{\ell}B\Gm_{(p)}.
\]
Write $\scr C(W)=DH\F_p[W]$, $\scr C(B\Gm) = \F_p[c]$ and $\scr L$ for the line bundle corresponding to $\widetilde V$.
We then have \[ \scr C(\widetilde V)\wequi \scr C(V) \otimes_{\F_p[c]}\F_p \] (in other words, $\scr C(\widetilde V)$ is the cofiber of multiplication by the first Chern class of $\scr L$ on $\scr C(V)$, which follows via the Thom isomorphism from the cofiber sequence $\widetilde V \to \scr L \to \Th(\scr L)$).

The bottom map $\P^1\to V$ is adjoint to the bottom inclusion $S^{4,2}\to J$, which induces the identity on $\ul K_2^{MW}$.
After contraction and passage to $\Gm_{(p)}$, its composite with $\ell$ therefore has first Chern class a generator.
Thus the cochain splitting in Lemma~\ref{lemm:identify-Jp-1} may be chosen under $\F_p[c]$, with
\[
\scr C(V)\simeq\Gamma_{\F_p}(x)\otimes DB_b,\qquad c\longmapsto x.
\]
It follows that
\[
\scr C(\widetilde V)\simeq
\bigl(\Gamma_{\F_p}(x)\otimes_{\F_p[c]}\F_p\bigr)\otimes DB_b.
\]
Augmenting the first factor and dualizing gives the required map $B_b\to H\F_p[\widetilde V]$; its composite with $H\F_p[V]$ is the second-factor inclusion of the chosen splitting, as needed.
\end{proof}

\begin{lemma} \label{lemm:loops2S3-main}
For $r \ge 0$, we have $\pi_{**} C^r H\F_p[\Omega^{1,1}JS^{2n,n}] \wequi \pi_{**} C^r \LSym(a_0, a_1, \dots, b_1, b_2, \dots)$ as Hopf algebras with power operations.
\end{lemma}
\begin{proof}
Let us put $A_i = \LSym(a_i, a_{i+1}, \dots, b_{i+1}, b_{i+2}, \dots)$.
Fix $N \ge 0$.
We shall construct, for each $r \ge 2$ and each $i$ a map of derived algebras $C^r A_i \to C^r H\F_p[\Omega^{1,1}JS^{2np^i,np^i}]$ which is an equivalence in weights $\le N$.
From this the result follows easily.\NB{obtain a sequence of generators, map in a free derived algebra, bar up as much as needed}

From now on we work in weights $\le N$, without explicitly stating so.
Since $\Omega^{2,1}JS^{2np^i,np^i}$ is $(np^i-1)$-effective (by Example~\ref{ex:J-effective} and \cite[Proposition~4.2.1(3)]{freudenthal}), we see by Corollary \ref{cor:module-cellular-criterion}(4) that $\overline{H\F_p[\Omega^{2,1}JS^{2np^i,np^i}]}$ is concentrated in weights $\ge np^i-1$, and so, the claim is vacuous for $i$ sufficiently large.
Let us assume that it is proved for $i+1$ and all $r$; we deduce the same for $i$.
We take the second Toda fiber sequence, apply $C^{r-1} H\F_p[\Omega^{2,1}(\ph)]$ and rotate, to obtain \[ C^{r+1}H\F_p[\Omega^{1,1}JS^{2np^{i+1},np^{i+1}}] \xrightarrow{\delta} C^rH\F_p[\Omega^{1,1}J_{p-1}S^{2np^i,np^i}] \to C^r H\F_p[\Omega^{1,1}JS^{2np^i,np^i}] \xrightarrow{\gamma} C^r H\F_p[\Omega^{1,1}JS^{2np^{i+1},np^{i+1}}]. \]
Since the result is proved for $i+1$ and all $r$, we know the source of $\delta$ as a Hopf algebra with power operations.
In particular, it is generated by geometric power operations on the primitives $a_*, b_*$.
We claim that $\pi_{**}\delta$ factors through the augmentation.
By the above generation statement, we need only verify that $\delta$ annihilates the primitives $a_*, b_*$ (observe that $\delta$ is a $\P^1$-loop map, so preserves geometric power operations).
But the only possible primitives in $C^r \LSym(a_i, b_{i+1})$ of weight not divisible by $p$ are $C^ra_i, C^rb_{i+1}$.
Comparing bidegrees, we thus see that $\delta(C^{r+1}a_j) = 0 = \delta(C^{r+1}b_j)$ for $j \ge i+2$, whereas $\delta(C^{r+1}a_{i+1}) = \lambda C^rb_{i+1}$ for some constant $\lambda$.
Tracing through the constructions, we see that $\lambda$ is precisely the degree of the bottom cell attaching map and hence vanishes mod $p$ by Lemma \ref{lem:boundary-p-map}.
We thus have proved the claim that $\pi_{**} \delta$ factors through the augmentation.
Via Lemma \ref{lemm:Toda-bar-lifts} below, this implies that the classes $C^ra_j, C^rb_j$ for $j>i$ lift to $\pi_{**} C^r H\F_p[\Omega^{1,1}JS^{2np^i,np^i}]$.
Applying Lemma \ref{lemm:cAlg-fiber-splitting} now shows that \[ \pi_{**} C^r H\F_p[\Omega^{1,1}JS^{2np^i,np^i}] \wequi \pi_{**} C^rH\F_p[\Omega^{1,1}J_{p-1}S^{2np^i,np^i}] \otimes \pi_{**} C^r A_{i+1}. \]
To conclude the proof, it will suffice to show that the factor $\pi_{**} C^r H\F_p[\Omega^{1,1}J_{p-1}S^{2np^i,np^i}]$ is a graded polynomial ring on appropriate power operations on $a_i$ and $b_{i+1}$.
The $a_i$-factor is a polynomial ring on $\scr E_r$-operations on $a_i$, by Remark \ref{rmk:a-family}, as needed.
On the other hand, since we know $CH\F_p[\Omega^{1,1}J_{p-1}S^{2np^i,np^i}]$ as a coalgebra by Corollary~\ref{cor:EM} and Lemma~\ref{lemm:identify-Jp-1}, we know that $\pi_{**} C^r H\F_p[\Omega^{1,1}J_{p-1}S^{2np^i,np^i}]$ is a graded polynomial ring generated by its subspace of primitives.
Consider the map \[ \Omega J_{p-1}S^{2np^i,np^i} \xrightarrow{(1)} \Omega S^{2np^{i+1}-1,np^{i+1}} \xrightarrow{(2)} \Omega^2 J_{p-1} S^{2np^{i+1},np^{i+1}} \xrightarrow{(3)} \Omega^2 JS^{2np^{i+1},np^{i+1}}. \]
(Here (1) is the projection from the first Toda fiber sequence, (2) is the inclusion in a \emph{different} first Toda fiber sequence, and (3) is the inclusion in a second Toda fiber sequence.
Crucially, by induction we have already analyzed (3).)
By Lemma \ref{lemm:b-family}, the induced map on motives, after the appropriate algebraic cobars, is compatible with (geometric) operations on $C^rb_{i+1}$.
On the other hand, by our induction hypothesis, the target contains a free derived algebra on the image of $C^rb_{i+1}$ (corresponding now to the $a$-family!).
Also, by construction, the map annihilates all classes corresponding to $a_i$ (and its power operations).
This shows that, starting from $b_{i+1}$ in $\pi_{**}C^r H\F_p[\Omega^{1,1}J_{p-1}S^{2np^i,np^i}]$ and applying geometric power operations, the space of primitives we build is big enough to exhaust all of the required generators.\NB{details?}
\end{proof}

We have used the following well-known calculation.
\begin{lemma} \label{lemm:Toda-bar-lifts}
Let $D\to E$ be a map of connected augmented $\scr E_1$-algebras, with $D$ a derived algebra, and put $T=B(E,D,\F_p)$.
Suppose the $\pi_{**}D$-action on $\pi_{**}E$ factors through the augmentation and $D\simeq\LSym(V)$ in the weights under consideration.
Then the classes of $\Sigma V\subset \pi_{**}BD$ lift to $\pi_{**}T$.
\end{lemma}
\begin{proof}
The relative bar spectral sequence has \[ E^2_{s,t,w}=\bigl(\pi_{**} E\otimes \operatorname{Tor}^{\pi_{**}D}_{s,*,*}(\F_p,\F_p)\bigr)_{t,w}, \]
Here $s$ is bar length, $t$ is internal degree, and $w$ is weight.
Augmentation gives a filtered map to the absolute bar, inducing augmentation on the first factor.
A named generator of $V$ defines a nonzero length-one class on $E^2$.
No $d_r$, $r\ge2$, can leave it, since $d_r$ lowers bar length by $r$.
Nor can a differential hit it: its image in the absolute bar is permanent, since symmetric degree one splits off as $\Sigma V$ in $B\LSym(V)\simeq\LSym(\Sigma V)$.
Naturality would send any differential hitting the relative class to one hitting this nonzero absolute class.
Positive weights make the filtration finite in each weight, so the permanent class gives the required lift.
\end{proof}

We can now put everything together.
\begin{proof}[Proof of Theorem \ref{thm:loops2S3-new}.]
(1) is Lemma \ref{lemm:loops2S3-cellular}.
(2) is proved in Lemma \ref{lemm:loops2S3-main}.
For (3) and (4), the case $a_1$ respectively $a_0$ is immediate (see Lemma \ref{lem:boundary-p-map} for (3)).
Now consider the commutative diagram
\begin{equation*}
\begin{CD}
\Sigma^{1,1}\Omega^{1,1} JS^{2n,n} @>>> \Sigma^{1,1}\Omega^{1,1} JS^{2np,np} \\
@VVV         @VVV \\
JS^{2n,n} @>>> JS^{2np,np} \\
\end{CD}
\end{equation*}
with the horizontal maps given by the Hopf maps from Construction \ref{cons:H-2}.
The top map annihilates $a_0$, but preserves the other $a_i$ (essentially by construction).
The bottom map sends $x_{2n,n}^{p^{i+1}}$ to $x_{2np,np}^{p^i}$.\footnote{The argument is the same as classically: in the notation of the proof of Proposition \ref{prop:Toda-fibrations}, passing to cohomology we see that $y_{2np,np}^{(p^i)}$ maps to $(x_{2n,n}^{(p)})^{(p^i)} = \frac{p^{i+1}!}{(p!)^{p^i}p^i!}x_{2n,n}^{(p^{i+1})}$ and we claim that this coefficient is $1$ mod $p$. Indeed the expression equals $\prod_{j=1}^{p^i}\binom{jp-1}{p-1}$, whose factors are congruent to $1$ mod $p$ by Lucas' theorem.}
Thus to prove the claim (4) for $a_1$ in $\Omega^{1,1}JS^{2n,n}$ it suffices to prove it for the corresponding bottom class $a_0$ in $\Omega^{1,1}JS^{2np,np}$, where it again holds by construction.
Continuing in this way treats all the $a_i$.
Assertion (3) is proved similarly.
\end{proof}

\subsection{More about $\dHopf_R^{bp}$} \label{subsec:more-dHopf}
Let $R=\Z/p^\nu$.
Recall the nonunital, stable, symmetric monoidal $\infty$-category $\dHopf_R^{bp}$ of \emph{derived bipolynomial Hopf algebras} from Definition \ref{def:dhopfbp}.
We will without further comment use the objects $W(n)_{\le d}$, $W(n)^{\le d}$ and $\Ga(n)$ ($=W(n)^{\le 0}$) from \S\ref{subsub:witt}.

\subsubsection{Algebraic loops}
Since $\dHopf_R^{bp}$ is (nonunital) symmetric monoidal, we have the functor $W(1)_{\le 0} \hotimes (\ph)$.
It has a partially defined right adjoint which we call $\hLoops^{2,1}$.

\begin{lemma} \label{lemm:Ga-loops}
For $m > 1$, $\hLoops^{2,1}\Ga(m)$ exists and is given by $\Ga(m-1)$.
Also $\hLoops^{2,1} \Ga(1) = R$.
\end{lemma}
\begin{proof}
We begin by constructing the desired map.
Recalling that $W(1) \hotimes W(m-1) \wequi W(m) \oplus \dots$ (see the proof of Lemma \ref{lemm:bipoly-monoidal} and Remark \ref{rmk:hotimes-howto}), we obtain a map \[ W(1)_{\le 0} \hotimes W(m-1) \to W(1) \hotimes W(m-1) \to W(m) \to \Ga(m). \]
We claim that this factors through the cofiber of Frobenius on $W(m-1)$, and in fact this is already true for $W(1) \hotimes W(m-1) \to \Ga(m)$.
By Lemma \ref{lemm:Wled-rep} this reduces to a statement about classical Hopf algebras, which holds just because filtration by Frobenius is lax symmetric monoidal (see, e.g., \cite[\S2.5.1]{bachmann-wilson}).
Consequently we obtain by adjunction the desired map $W(1)_{\le 0} \hotimes \Ga(m-1) \to \Ga(m)$.

To show that this map exhibits $\Ga(m-1)$ as $\hLoops^{2,1}\Ga(m)$, we need only test the induced adjunction map on generators, i.e., we must show that (for $j \ge 0$, $n>0$) \[ \Map_{\NHopf}(C^j \Free \A^1(n), \Ga(m-1)) \wequi \Map_{\NHopf}(C^j \Free \A^1(n) \hotimes W(1)_{\le 0}, \Ga(m)). \]
Using that $W(1)_{\le 0} \wequi \Free HR[S^{2,1}]$ (Example \ref{ex:sphere-bipolynomial}), this translates into the following: the canonical map \[ \Map_{\cNAlg_*}(\A^1(n), B^j \Ga(m-1)) \to \Map_{\cNAlg_*}(\A^1(n) \wedge HR[S^{2,1}], B^j \Ga(m)) \] is an equivalence.
Dualizing via Proposition \ref{prop:NHopf-basics}(1) and Proposition \ref{prop:NAlg-DAlg}, we equivalently need to show that \[ \Map_{\DAlg^{aug}}(C^j D\Ga(m-1), D\A^1(n)) \wequi \Map_{\DAlg^{aug}}(C^j D\Ga(m), D(\A^1(n) \wedge HR[S^{2,1}])). \]
As in the proof of Lemma \ref{lemm:Wled-rep}, $C^j D\Ga(m-1)$ is the free derived algebra on a generator in Chow degree $-j$ and weight $1-m$; thus we reduce to showing that $D(\A^1(n) \wedge HR[S^{2,1}])$ is concentrated in Chow degree zero and that the degree $1-m$ part of $D\A^1(n)$ identifies with the degree $-m$ part of $D(\A^1(n) \wedge HR[S^{2,1}])$.
Equivalently, $\A^1(n) \wedge HR[S^{2,1}]$ is concentrated in Chow degree zero and the degree $m-1$ part of $\A^1(n)$ identifies with the degree $m$ part of $\A^1(n) \wedge HR[S^{2,1}]$.
This is straightforward.\NB{details?}

For the final statement, we may again check on generators, i.e., verify that $\Map_{\NHopf}(\Free \A^1(n) \hotimes W(1)_{\le 0}, B^j \Ga(1)) = 0$ for all $n > 0$, $j \ge 0$.
Passing to coalgebras and dualizing as above, this translates into $\Map_{\DAlg^{aug}}(C^j D\Ga(1), D(\A^1(n) \wedge HR[S^{2,1}])) = *$.
In other words we must show that $\A^1(n) \wedge HR[S^{2,1}]$ has no component in weight $1$, which is again straightforward.
\end{proof}

\begin{proposition} \label{prop:alg-loops-defined}
The functor $\hLoops^{2,1}$ is defined on all of $\dHopf_R^{bp}$.
It is compatible with scalar extension, i.e., for $X \in \dHopf_R^{bp}$ we have \[ \hLoops^{2,1}(X \otimes_R R') \wequi \hLoops^{2,1}(X) \otimes_R R' \in \dHopf_{R'}^{bp}. \]
\end{proposition}
\begin{proof}
The subcategory of objects on which either claim (existence/compatibility with scalar extension) holds is thick and closed under quickly converging biproducts (as in Proposition \ref{prop:NHopf-stability}).
It consequently suffices to prove the claims for $W(n)$.
The case $W(n)^{\le 0} = \Ga(n)$ has been dealt with in Lemma \ref{lemm:Ga-loops} and Remark \ref{rmk:W-scalar-ext}.
Using the fiber sequences from Remark \ref{rmk:Witt-filtrations}, inductively we deal with $W(n)^{\le r}$.
In fact, repeated application of Example \ref{ex:Ga-injective} shows that the fiber sequences \[ \Ga(p^{r+1}n - 1) \wequi \hLoops^{2,1} \Ga(p^{r+1}n) \to \hLoops^{2,1} W(n)^{\le r+1} \to \hLoops^{2,1} W(n)^{\le r} \] split ($p^{r+1}n-1$ being coprime to $p$).
The system $\{\hLoops^{2,1}W(n)^{\le d}\}_d$ is thus quickly converging and split.
Its limit exists in the larger category $\NHopf(\Mod_{HR}(k))^{\cell\cn+}$ and is there given by a certain infinite sum of $\Ga$.
This can be written as a cofiber of an analogous infinite sum of $W$, and so lies in $\dHopf_R^{bp}$.
We deduce that $\hLoops^{2,1}W(n)$ exists and its formation commutes with scalar extension (by our explicit computation).
\end{proof}

\begin{remark}
It may be worth spelling out the end of the above proof: we showed that \[ \hLoops^{2,1}W(n) \wequi \bigotimes_{d \ge 0} \Ga(p^dn-1) \] (when $n=1$, the term $\Ga(0)$ should be dropped).
In other words, $\hLoops^{2,1}$ turns the classical indecomposable bipolynomial Hopf algebra $W(n)$ into a large divided-power Hopf algebra.
\end{remark}

\subsubsection{Indecomposables and primitives}
Given any ordinary augmented $R$-algebra $A$, we have the \emph{indecomposables} $QA = I/I^2$.
Here $I$ is the augmentation ideal, and $QA$ is an $R$-module.
We obtain by restriction an additive, symmetric monoidal functor \[ Q: \Hopf_R^{bp} \to \FFree_R^{\N}, \] where the right hand side denotes degreewise finitely generated free $R$-modules.
Passing to $\Perf$ and including into all modules, we obtain the stably nonunital symmetric monoidal functor \[ Q: \dHopf_R^{bp} \to \Mod_R^\N \wequi \Mod_{HR}(k)^{\cell\eff}. \]
We call this the \emph{derived indecomposables}.
We define the \emph{derived primitives} $P$ similarly (but beware that this functor is not symmetric monoidal).

\begin{remark} \label{rmk:P-Q-explanation}
Let $E \in \dHopf_R^{bp}$.
Then, shifting up sufficiently much if necessary, $E$ can be built using finite colimits out of (quickly converging sums of) $W(n)$-s.
This yields a resolution of the underlying derived algebra of $E$ by polynomial algebras.
We thus see that $QE$ coincides with the usual derived indecomposables of the derived algebra underlying $E$.
In particular, if $E$ is concentrated in Chow degree zero, we see the usual left derived functor of indecomposables (which is a right exact functor) of the underlying algebra.

Dually we can also build $E$ by first shifting down and then using finite limits of $W(n)$-s.
At least when $R = \F_p$, this corresponds to resolving a Hopf algebra \emph{injectively}, and so computes the right derived functor of primitives (which is indeed left exact, at least in this situation).
\end{remark}

\begin{remark} \label{rmk:P-Q-compatible-scalar-ext}
By construction, the functors $P$ and $Q$ commute with scalar extension (see Remark \ref{rmk:Hopf-scalar-ext} for scalar extension of derived Hopf algebras).
\end{remark}

Now let $E \in \dHopf_R^{bp}$.
We have the co-unit of adjunction $c: HR[TS^{2,1}] \hotimes\hLoops^{2,1}E \to E$.
Our next goal is to show that, analogously with the classical situation, this factors through primitives/indecomposables.

We begin with the case of discrete bipolynomial Hopf algebras.
\begin{lemma}\label{lemm:hopfbp-nat-trans}
Write $I: \Hopf_R^{bp} \to \Mod_R^\N$ for the functor sending $E$ to its underlying module, with the degree zero part removed (i.e., the augmentation ideal).
There is a unique natural transformation \[ \Sigma^{2,1}Q\hLoops^{2,1} \to P \in \Fun(\Hopf_R^{bp}, \Mod_R^\N) \] such that the following diagram commutes
\begin{equation*}
\begin{CD}
\Sigma^{2,1}Q\hLoops^{2,1} @>>> P \\
@AAA                       @VVV \\
I(HR[TS^{2,1}] \hotimes\hLoops^{2,1}) @>Ic>> I.
\end{CD}
\end{equation*}
\end{lemma}
\begin{proof}
By construction, the left hand map is surjective (on $\pi_{2*,*}$) and the right hand map is injective, so uniqueness is clear.
For existence, as usual it suffices to treat $W(n)$; extension to quickly converging sums is straightforward.
Applied to $W(n)$, by definition $P$ is concentrated in homotopical degree zero (and just given by the usual primitives).
On the other hand, applying $Q$ to $\hLoops^{2,1}W(n)$ may yield higher derived functors.
Nonetheless, since $\hLoops^{2,1}W(n)$ is concentrated in Chow degrees $\ge 0$ (see the proof of Proposition \ref{prop:alg-loops-defined}), so is $Q\hLoops^{2,1}W(n)$, and its $\pi_0$ is given by the ordinary indecomposables (see Remark \ref{rmk:P-Q-explanation}).

To verify that we land in primitives, consider the commutative diagram
\begin{equation*}
\begin{CD}
\Sigma^{2,1}Q\hLoops^{2,1}W(n) @>>> IW(n) \\
@AAA                                    @AAA \\
\Sigma^{2,1}I\hLoops^{2,1}W(n) @>>> I(HR[TS^{2,1}] \hotimes\hLoops^{2,1}W(n)).
\end{CD}
\end{equation*}
The left hand map is surjective on $\pi_0$, and in the bottom left hand coalgebra, every element is primitive (by Eckmann--Hilton, for any pointed normed coalgebra $C$, the reduced comultiplication in $\Sigma C \wequi HR[S^1] \wedge C$ comes from the reduced diagonal $S^1 \to S^1 \wedge S^1$, which is null\NB{details?}).

That the map factors through indecomposables is verified similarly.
In the end we want to show that $\Sigma^{2,1}I\hLoops^{2,1}W(n) \to W(n)$ annihilates products.
For this we observe that the composite $\Sigma^{2,1}(\hLoops^{2,1}W(n) \times \hLoops^{2,1}W(n)) \to \Sigma^{2,1} \hLoops^{2,1}W(n) \to W(n)$ of the loop multiplication and the counit is given by the sum of the two projections, both of which annihilate products of positive-degree classes from the two different factors.
\end{proof}

Consider now the functor \[ \Hopf_R^{bp} \to \Fun(\Delta^1 \times \Delta^1, \Mod_R^\N) \] constructed in Lemma \ref{lemm:hopfbp-nat-trans}.
The two top corners are additive functors.
The two bottom corner functors are degreewise polynomial (this is easy to check for the bottom right hand corner, i.e., the forgetful functor, and it follows for the bottom left hand functor since it is obtained by composition with additive functors).
\begin{construction}\label{cons:dhopfbp-nat-trans}
Using \cite[Theorem 2.19]{k-theory-and-poly-functors}, we may thus extend the above functor to $\Perf(\Hopf_R^{bp}) \wequi \dHopf_R^{bp}$.
By uniqueness of extension applied to three of the sides, we see that we have constructed for $E \in \Perf(\Hopf_R^{bp})$ a natural commutative square
\begin{equation*}
\begin{CD}
\Sigma^{2,1}Q\hLoops^{2,1}E @>>> PE \\
@AAA                       @VVV \\
I(HR[TS^{2,1}] \hotimes\hLoops^{2,1})E @>Ic>> IE,
\end{CD}
\end{equation*}
where the vertical maps are the canonical ones.
\end{construction}

\begin{remark} \label{rmk:c-compat-scalar-extn}
All functors involved in the above transformations are compatible with scalar extension, by Proposition \ref{prop:alg-loops-defined} and Remark \ref{rmk:P-Q-compatible-scalar-ext}.
The construction shows that the same is true for the natural transformations.
\end{remark}

\begin{remark}
Consider the case $R=\F_p$.
We can extend the functors $P$ and $Q$ to $\dHopf^{cn}_{\F_p}$ by first restricting to the thick subcategory of $\dHopf_{\F_p}^{bp}$ generated by the $W(n)_{\le r}$ and then extending by filtered colimits.
(See the discussion in \S\ref{subsec:derived-hopf-alg} for the relationship between $\dHopf_{\F_p}^{bp}$ and $\dHopf^{cn}_{\F_p}$.)
Since the functor $Q$ on derived algebras preserves colimits, our ``new'' functor on $\dHopf^{cn}_{\F_p}$ is still just given by the derived indecomposables.
The new functor $P$ still coincides with the derived primitives on $\dHopf_{\F_p}^{bp} \subset \dHopf^{cn}_{\F_p}$, as we can check by explicit computation.\NB{details?}
In general, however, it is given by some kind of ind-extension of the derived primitives.

The natural transformations of Construction \ref{cons:dhopfbp-nat-trans} also extend, with similar caveats.
\end{remark}

\subsection{Homology of $\P^1$-loops} \label{subsec:P1-loops}
\subsubsection{}
Recall from Definition \ref{def:cellular-categories} the categories $\NMon(\SH^{S^1}(k))^{\cell i \eff}$ and $\NMon(\Spc(k))^{\cell i \eff}_{\ge 1}$.
We introduce some convenient variants.

\begin{definition}
Let $i \ge 1$.
We define \[ \NMon(\Spc(k))^{\cell i \eff \omega}_{\ge 1} \subset \NMon(\Spc(k))^{\cell i \eff}_{\ge 1} \text{ and } \NMon(\SH^{S^1}(k))^{\cell i \eff \omega} \subset \NMon(\SH^{S^1}(k))^{\cell i \eff} \] as generated by $T(S^{n+1,n})$ for $n \ge i$, respectively
under finite colimits, extensions and retracts or as a thick subcategory.
We write $\cn$ in place of $1\eff$.
\end{definition}
It follows from Lemma \ref{lemm:NMon-prestable} that \[ \NMon(\SH^{S^1}(k))^{\cell i \eff} \wequi \Sp(\NMon(\Spc(k))^{\cell i \eff}_{\ge 1}). \]
In particular, the $T(S^{n+1,n})$ are compact in the latter category and so \cite[\S5.5.7]{lurie-htt} \[ \NMon(\SH^{S^1}(k))^{\cell i \eff} \wequi \Ind(\NMon(\SH^{S^1}(k))^{\cell i \eff \omega}). \]
Similarly \[ \NMon(\SH^{S^1}(k))^{\cell i \eff \omega} \wequi \SW(\NMon(\Spc(k))^{\cell i \eff \omega}_{\ge 1}). \]

\subsubsection{}
Now consider the functor $H\F_p[\ph]: \NMon(\Spc(k))^\gp \to \NHopf(\Mod_{H\F_p}(k))$ from Example \ref{ex:homology-hopf-algebra}.
It is cocontinuous and normed.
By Example \ref{ex:sphere-bipolynomial}, we have $H\F_p[T(S^{2n,n})] \in \dHopf_{\F_p}^{bp} \subset \NHopf(\Mod_{H\F_p}(k))$.
Consequently we obtain a functor \[ H\F_p[\ph]: \NMon(\Spc(k))^{\cell \cn \omega}_{\ge 1} \to \dHopf_{\F_p}^{bp}. \]
The latter category being stable by Proposition \ref{prop:NHopf-stability}(2), we can extend over $\SW$ to obtain \[ \NMon(\SH^{S^1}(k))^{\cell \cn \omega} \to \dHopf_{\F_p}^{bp}. \]
Extending over $\Ind$-objects yields \[ \Ind(H\F_p)[\ph]: \NMon(\SH^{S^1}(k))^{\cell\cn} \to \Ind(\dHopf_{\F_p}^{bp}). \]

\subsubsection{}
Let us write $\hLoops^{2,1}$ for the functor on $\NMon(\SH^{S^1}(k))^{\cell\cn}$ right adjoint to tensoring with $T(S^{2,1})$.

\begin{construction} \label{cons:main-trafo}
Using the transformation from Construction \ref{cons:dhopfbp-nat-trans}, as well as the canonical exchange transformation for $\hLoops^{2,1}$ and $\Ind(H\F_p)[\ph]$, we obtain a natural transformation \[ \mu: \Sigma^{2,1} \Ind(Q) \Ind(H\F_p)[\hLoops^{2,1}(\ph)] \to \Ind(P) \Ind(H\F_p)[\ph] \in \Fun(\NMon(\SH^{S^1}(k))^{\cell\cn}, \Mod_{\F_p}^\N). \]
\end{construction}

\begin{remark} \label{rmk:trafo-scalar-ext}
We can repeat the entire discussion with $R=\Z/p^\nu$ in place of $\F_p$.
In the end this yields a natural transformation \[ \mu_R: \Sigma^{2,1} \Ind(Q) \Ind(HR)[\hLoops^{2,1}(\ph)] \to \Ind(P) \Ind(HR)[\ph] \in \Fun(\NMon(\SH^{S^1}(k))^{\cell\cn}, \Mod_{R}^\N). \]
It follows from Proposition \ref{prop:alg-loops-defined} and Remark \ref{rmk:c-compat-scalar-extn} that the resulting natural transformation is even compatible with scalar extension, i.e., $\mu_{R'} \wequi \mu_R \otimes_R R'$.
\end{remark}

\begin{lemma} \label{lemm:cell-cn-loops}
If $E \in \NMon(\SH^{S^1}(k))^{\cell 2\eff}_{(p)} \subset \NMon(\SH^{S^1}(k))^{\cell\cn}$, then $\hLoops^{2,1} E \wequi \Omega^{2,1} E$.
\end{lemma}
\begin{proof}
It suffices to show that if $X \in \NMon(\SH^{S^1})_{(p)}^{\cell 2\eff} \subset \NMon(\SH^{S^1})$, then $\Omega^{2,1} X \in \NMon(\SH^{S^1})^{\cell\cn}$.
This reduces to the case $X = T(S^{2n,n})_{(p)}$.
Instead of showing that $\Omega^{2,1}X$ is cellular, we may as well prove this for $\Omega^{1,1}X$, cellular objects being stable under shifting by definition.
Suppose first that $n \ge 3$.
Then $T(S^{2n,n})$ is very $3$-effective by Lemma \ref{lemm:Dn-effectivity-Spc*}, and so $\Omega^{1,1} T(S^{2n,n})$ is $2$-effective by \cite[Proposition~4.2.1(2)]{freudenthal}.
Consequently by Lemma \ref{lem:cellularity} (and Remark \ref{rmk:cell-p-loc}), it suffices to show that $H\Z_{(p)}[\Omega^{1,1}T(S^{2n,n})]$ is cellular.
Since $T(S^{2n,n})_{(p)} \wequi J(S^{2n,n})_{(p)}$ by Corollary \ref{cor:identify-J-T}, this is dealt with in Theorem \ref{thm:loops2S3-new}(1).

It remains to handle the case $n=2$.
The argument works entirely the same way, except that, unfortunately, $\Omega^{1,1} T(S^{4,2})$ is not $2$-effective, and so we cannot directly apply Lemma \ref{lem:cellularity}.
Instead we shall show that $\cof(TS^{2,1} \to \Omega^{2,1} T(S^{4,2}))_{(p)}$ is $2$-effective.
The same is then true for $\cof(TS^{3,1} \to \Omega^{1,1} T(S^{4,2}))_{(p)}$, and we can verify its cellularity via Lemma \ref{lem:cellularity}.
By appeal to Lemma \ref{lemm:2-eff-trick} (as in Example \ref{ex:K1-cellular}), it will thus suffice to show that the canonical map $S^{2,1} \to \Omega^{2,1}J(S^{4,2})_{(p)}$ has cofiber with $p$-localization which is simply-connected and $2$-effective.
As a preliminary observation, if we produce a factorization $S^{2,1} \to X \to \Omega^{2,1}J(S^{4,2})_{(p)}$ such that $X \to \Omega^{2,1}J(S^{4,2})_{(p)}$ induces a surjection $\ul\pi_1$ and an equivalence on $S^{3,2}$-nullifications, and $X$ is connected, then it suffices to prove that $\cof(S^{2,1} \to X)_{(p)}$ is simply connected and $2$-effective.
We have $JS^{4p,2p} \in O(S^{8,4})$ by Example~\ref{ex:J-effective}, and hence so is its $p$-localization by \cite[Proposition 4.3.4]{freudenthal}, and hence $\Omega\Omega^{2,1} JS^{4p,2p}_{(p)} \in O(S^{5,3})$ \cite[Proposition 4.2.1]{freudenthal}.
The second fiber sequence of Proposition \ref{prop:Toda-fibrations} together with \cite[Proposition 3.1.23]{freudenthal} thus allows us to replace (via the preliminary observation) $\Omega^{2,1}J(S^{4,2})_{(p)}$ by $\Omega^{2,1} J_{p-1}S^{4,2}_{(p)}$.
Similarly the first fiber sequence allows us to replace $\Omega^{2,1} J_{p-1}S^{4,2}_{(p)}$ by $\Omega^{1,1}S^{3,2}_{(p)}$.
By using \cite[Proposition 4.3.4]{freudenthal} once more, we may look at the actual cofiber.
Since $S^{3,2} \wequi \A^2 \setminus 0 \wequi \SL_2$ \cite[p.~3, Notation/Preliminaries]{asok2014algebraic}, by the affine Grassmannian model from Proposition \ref{prop:aff-grass}, this cofiber may be built (even before $p$-localization)\NB{details?} by starting with $S^{4,2}$ ($\wequi \Sigma \A^2 \setminus 0$) and attaching cells of the form $\A^n \setminus 0$ for $n > 2$, which thus lie in $O(S^{4,2})$ as needed.
\end{proof}

\begin{remark}
Recall from \S\ref{subsec:derived-hopf-alg} that $\Ind(\dHopf_{\F_p}^{bp})$ contains $\dHopf^{cn}_{\F_p}$ as a full subcategory.
Since the objects $H\F_p[T(S^{2n,n})] \wequi W(n)_{\le 0} \in \dHopf^{cn}_{\F_p}$ are in fact compact objects in $\dHopf^{cn}_{\F_p}$, we see that $\Ind(H\F_p)[\ph]$ actually lands in $\dHopf^{cn}_{\F_p}$.
We thus obtain a (cocontinuous, nonunital symmetric monoidal) functor \[ \dHFp[\ph]: \NMon(\SH^{S^1}(k))^{\cell\cn} \to \dHopf^{cn}_{\F_p}. \]
Since $H\F_p[\ph]$ factors through $p$-localization (for connected normed monoids this agrees with localization of the underlying nilpotent motivic space, as explained in \S\ref{subsub:mon-ploc}, and such localization preserves mod-$p$ homology), for $E \in \NMon(\SH^{S^1}(k))^{\cell 2\eff}$, the natural transformation of Construction \ref{cons:main-trafo} takes the form \[ \mu: \Sigma^{2,1}Q\dHFp[\Omega^{2,1}E] \to P\dHFp[E]. \]
\end{remark}

\subsubsection{}
We now describe the homology of $\P^1$-loops.
We begin with the following Example, illustrating that a more naive construction is insufficient.

\begin{example}
Instead of the transformation $\mu$, there is also the more straightforward map \[ \dHFp[\Omega^{2,1}E] \to \hLoops^{2,1}\dHFp[E]. \]
It is typically not an equivalence, even for very well-behaved $E$.
For example, consider $E = B\SL$, so that $\Omega^{2,1} E \wequi B\GL$.
One finds that $\hLoops^{2,1}\dHFp[B\SL]$ is concentrated in Chow degree zero, so from now on we need only contemplate classical Hopf algebras.
The functor $\hLoops^{2,1}$ is already defined for coalgebras, where a priori it only yields a coalgebra.
However, since $S^{2,1}$ is a cocommutative cogroup, $\hLoops^{2,1}C$ is a Hopf algebra for any coalgebra C.
By an Eckmann--Hilton argument, if $C$ was already a Hopf algebra, then the two resulting Hopf algebra structures agree.
Noting that $\dHFp[B\SL]$ is cofree as a coalgebra (equivalently, the cohomology of $B\SL$ is polynomial), we see that as a coalgebra, $\hLoops^{2,1}\dHFp[B\SL]$ is also cofree, on shifted generators.
This already shows that $\dHFp[B\GL] \to \hLoops^{2,1}\dHFp[B\SL]$ is not an equivalence, since in the source, the shifted generators are not primitive, whereas they become so in the target.
In fact $\hLoops^{2,1}$ on a cofree Hopf algebra produces a \emph{divided power Hopf algebra} (use Lemma \ref{lemm:Ga-loops}, i.e., $\hLoops^{2,1} \A^1(n+1) \wequi \Ga(n)$), and the map $\dHFp[B\GL] \to L^{2,1}\dHFp[B\SL]$ identifies with \[ \F_p[b_1, b_2, \dots] \to \Gamma[b_1,b_2,\dots] .\]
\end{example}

Here is our main result.
\begin{theorem} \label{thm:homology-of-loops}
Let $k$ be a perfect field of characteristic $p$ and $R=\Z/p^\nu$.
For $X \in \NMon(\SH^{S^1}(k))^{\cell 2\eff}$, the map \[ \mu_R: \Sigma^{2,1}\Ind(Q)\Ind(HR)[\Omega^{2,1}(X)] \to \Ind(P) \Ind(HR)[X] \] from Construction \ref{cons:main-trafo} is an equivalence.
\end{theorem}
\begin{proof}
Since $\otimes_R \F_p$ is conservative and $\mu_R$ is compatible with scalar extension (Remark \ref{rmk:trafo-scalar-ext}), we reduce to the case $R=\F_p$.
Since both sides preserve colimits and desuspensions, we reduce to the case $X = T(S^{2n,n})$ (with $n \ge 2$).
We can compute the right hand side by using the resolution $\dHFp[TS^{2n,n}] \to W(n) \xrightarrow{V} W(pn)$ (see Example \ref{ex:sphere-bipolynomial}).
The primitives in $W(n)$ are spanned by $\{w_{np^i}\}_{i \ge 0}$; here $w_{np^i}$ is given by applying the Witt polynomial to the polynomial generators $\{x_{np^i}\}_{i \ge 0}$ \cite[p.4]{goerss-dieudonne}.
Since $V(x_{np^i}) = x_{np^i}$ for $i>0$ and $V(x_n) = 0$ \cite[Proposition 2.1]{goerss-dieudonne} (beware that our indexing is slightly different from Goerss') we see that integrally we have $Vw_{np^i} = pw_{np^i}$ for $i > 0$, and $Vw_{np^i} = 0$ for $i=0$.
In particular, with $\F_p$-coefficients the map induced by $V$ on primitives is just $0$ and we get \[ P \dHFp[TS^{2n,n}] \wequi \F_p\{w_{np^i}, \Sigma^{-1}w_{np^j} \mid i \ge 0, j \ge 1\} \text{ with } \beta w_{np^j} = \Sigma^{-1}w_{np^j}. \]
(In fact, with $\F_p$-coefficients, $w_{np^i} = x_n^{p^i}$.
But from this alone we would not have been able to deduce the Bockstein action.)

On the other hand, Theorem \ref{thm:loops2S3-new}(2), Remark \ref{rmk:loops2S3-derived-algebra}, and Corollary \ref{cor:identify-J-T} show that $H\F_p[\Omega^{1,1}TS^{2n,n}] \wequi \NSym(a_0, a_1, \dots, b_1,b_2, \dots)$ and so $\Sigma^{2,1}Q\dHFp[\Omega^{2,1}TS^{2n,n}] \wequi \Sigma^{1,1} Q\dHFp[\Omega^{1,1}TS^{2n,n}]$ is abstractly isomorphic to $P \dHFp[TS^{2n,n}]$.
It remains to check that the map is the obvious one.
Recalling that our map is obtained by scalar extension from a map defined for $\Z/p^2$, we see that Bockstein acts on both sides in a compatible way.
It thus suffices (by Theorem \ref{thm:loops2S3-new}(3)) to show that $a_i$ maps to $w_{np^i}$.
By construction (see Construction \ref{cons:dhopfbp-nat-trans}), we have the following commutative diagram
\begin{equation*}
\begin{CD}
I \pi_{*-2,*'-1} H\F_p[\Omega^{2,1}TS^{2n,n}] @>>> \pi_{*-2,*'-1} Q\dHFp[\Omega^{2,1}TS^{2n,n}] \\
@VcVV @V{\mu}VV \\
\pi_{*,*'} H\F_p[TS^{2n,n}] @<u<< \pi_{*,*'} P\dHFp[TS^{2n,n}].
\end{CD}
\end{equation*}
Here $I$ means the augmentation ideal, $\mu$ is the map we are studying, the top horizontal map is the obvious quotient map, and $u$ is induced by the inclusion of primitives into a coalgebra (via taking derived functors).
Thus $u(w_{np^i}) = x_n^{p^i}$.
Also from Theorem \ref{thm:loops2S3-new}(4) we deduce that $c(a_i) = x_{2n,n}^{p^i}$.
The result follows.
\end{proof}

\subsection{A model for $K(\Z_{(p)}(n), 2n)$} \label{subse:model}
We now utilize Theorem \ref{thm:homology-of-loops} to implement the strategy outlined in \S\ref{subsec:into-summary}.

\begin{theorem} \label{thm:mainmain}
Let $k$ be a field of characteristic $p$ and $n \ge 1$.
The canonical map (see \S\ref{subsec:intro-main}) \[ D_n^{\otimes} K(\Z_{(p)}(1), 2) \to K(\Z_{(p)}(n),2n) \in \NMon(\Spc(k)) \] is an equivalence.
\end{theorem}
\begin{proof}
By essentially smooth base change, we may assume that $k$ is perfect (and in fact $k=\F_p$).

Let us put \[ E^\bullet = \NSym^{\Z\otimes}(K(\Z_{(p)}(1),2),1)^\bullet \wequi D_\bullet^\otimes(K(\Z_{(p)}(1),2)) \in \NAlg(\NMon(\Spc(k))^\Z). \]
The map $T(\P^1) \to K(\Z_{(p)}(1), 2) = E^1$ induces \[ \P^1 \wedge E^n \to E^1 \otimes E^n \to E^{n+1}, \] turning $E^\bullet$ into a $\P^1$-prespectrum.
By Proposition \ref{prop:Dn-tens-eff} and Example \ref{ex:K1-very-1-eff}, $E^n$ is very $n$-effective (for $n \ge 1$), and hence the spectrification \[ sE = \colim_n \Sigma^{-2n,-n} \Sigma^\infty E^n \] of $E$ is very effective.
We shall show that \begin{equation} \label{eq:key} \Omega^{2,1}E^{n+1} \wequi E^n \text{ for } n \ge 1. \end{equation}
Assuming this, we see that \[ (sE)^0 \wequi \Omega^{2,1} E^1 \wequi \Omega^{2,1} K(\Z_{(p)}(1),2) \wequi \Z_{(p)}. \]
Hence the map $\Omega^\infty(sE \to H\Z_{(p)})$ is an equivalence, and thus so is $sE \to H\Z_{(p)}$, $H\Z_{(p)}$ also being very effective \cite[Lemma~13.7]{bachmann-norms}, and $\Omega^\infty$ being conservative on very effective spectra \cite[Corollary~3.5.15]{EHKSY}.

It remains to prove \eqref{eq:key}.
It is clear that $E^1$ is limit-determined (in the sense of Definition \ref{def:limit-determined}), being an infinite $\P^1$-loop space.
Also $E^n$ is limit-determined by Example \ref{ex:resolvable-spaces} for $n \ge 2$, being $2$-effective (as argued above).
Consequently also $\Omega^{2,1}E^n$ is limit-determined for all $n$, by Example \ref{ex:loops-limit-determined}.
Hence by Lemma \ref{lemm:limit-determined-M-test} (which is applicable by Lemma \ref{lemm:Dnotimes-plocal}), in order to show that $E^n \to \Omega^{2,1}E^{n+1}$ is an equivalence, it suffices to show this after applying $H\Z_{(p)}[\ph]$.
From this it follows that we need only show it for $H\F_p[\ph]$ and $H\Q[\ph]$ (see e.g. \cite[Lemma~2.16]{bachmann-eta}), and for the latter it suffices to show this after unstable rationalization \cite[Theorem~4.3.9(2)]{AFHLocalization}.

We first treat the rational case.
Note that (unstable) rationalization of simply-connected connected normed monoids commutes with $\Omega^{2,1}$ (see the discussion in \S\ref{subsec:p-loc}) with $D_n^\otimes$ by Lemma \ref{lemm:Dnotimes-plocal}.
By Lemma \ref{lemm:key1} below we have $D_n^\otimes(K(\Q(1),2)) \xrightarrow{\wequi} K(\Q(n),2n)$.
Consequently \[ (\Omega^{2,1} E^{n+1})_\Q \wequi \Omega^{2,1} D_{n+1}^\otimes K(\Q(1),2) \wequi \Omega^{2,1} K(\Q(n+1),2n+2) \wequi K(\Q(n),2n) \wequi E^n_\Q, \] as needed.

Now we treat the case of $\F_p$-homology.
We know that $E^n$ and $\Omega^{2,1}E^{n+1}$ both lie in $\NMon(\SH^{S^1}(k))^{\cell\cn}$, by Proposition \ref{prop:Dn-cellular} and Example \ref{ex:K1-cellular} for the former and by Lemma \ref{lemm:cell-cn-loops} for the latter.
It is thus enough to show that $\dHFp[E^n] \to \dHFp[\Omega^{2,1}E^{n+1}]$ is an equivalence.
By \cite[Lemma 2.11]{brantner2025formal}\footnote{This reference deals with $\Z$-algebras instead of $\F_p$-algebras, but the same proof works in the latter case.}, it suffices to check this on derived indecomposables $Q$.
By Theorem \ref{thm:homology-of-loops}, we have \[ \Sigma^{2,1}Q\dHFp[\Omega^{2,1}E^{n+1}] \wequi P\dHFp[E^{n+1}]. \]
Since $\dHFp[K(\Z_{(p)}(1),2)] \wequi \Ga(1)$ (Example \ref{ex:Poo-bipolynomial}) and $\dHFp[\ph]$ turns $D_n^\otimes$ into $\LSym_n^\hotimes$ (see Lemma \ref{lemm:dHopf-bp-nobp}), we conclude by Lemma \ref{lemm:key2} below.
\end{proof}

\begin{lemma} \label{lemm:key1}
The canonical map \[ D_n^\otimes(K(\Q(1),2)) \to K(\Q(n),2n) \] is an equivalence.
\end{lemma}
\begin{proof}
Consider the adjunction \[ L: \NMon(\Spc(k))^{\gp\A^1}_\Q \adj \Mod_{H\Q}(k)^\veff: \Omega^\infty. \]
Then $L$ is a normed functor, and so sends $D_n^\otimes$ to $D_n$.
By Corollary \ref{cor:NMon-rat-loc}, $L$ is a localization and each $D_n^\otimes K(\Q(1),2)$ is a local object.
The same is also true for each $K(\Q(n),2n)$, being in the image of the right adjoint.
It thus suffices to check that the map becomes an equivalence after applying $L$.
When varying $n$, both sides turn into graded rings.
In degree one, they agree by construction, and are given by the free $H\Q$-module on the tautological class $t$.
We shall argue that in degree $n$, both sides are given by the free $H\Q$-module on the tautological class $t^n$.
Since $L\Omega^\infty \wequi \id$, the claims about the right hand side are clear.
For the left hand side, we have \[ L D_n^\otimes(K(\Q(1),2)) \wequi D_n L K(\Q(1),2) \wequi D_n \Sigma^{2,1} H\Q. \]
By the Thom isomorphism \cite[Proposition~5.37]{bachmann-colimits} we have $D_n \Sigma^{2,1}H\Q \wequi \Sigma^{2n,n} D_n H\Q$.
We know that $D_n H\Q \wequi H\Q[B_\et\Sigma_n]$ \cite[Example~5.8 and Proposition~5.14]{bachmann-colimits}.
By \cite[Proposition~5.39]{bachmann-colimits}, the transfer for the trivial subgroup exhibits $H\Q[B_\et\Sigma_n]$ as a retract of $H\Q$, since $n!$ is invertible in $\Q$.
Thus the unit $H\Q \to H\Q[B_\et\Sigma_n]$ is an equivalence.
Under the Thom isomorphism, this unit corresponds to $t^n$.
\end{proof}

\begin{lemma} \label{lemm:key2}
For $n \ge 1$ the canonical map \[ \Sigma^{2,1} Q\LSym_n^\hotimes \Ga(1) \to P\LSym_{n+1}^\hotimes \Ga(1) \] is an equivalence.
\end{lemma}
\begin{proof}
Let us write \[ RW_{(p)}^\bullet = \LSym^\hotimes_{L_{(p)}}(W(1)) \hotimes_{\LSym^\hotimes_{L_{(p)}}(W(p))} L_{(p)} \in \DAlg(\dHopf_{\F_p}) \] for the ($p$-local) \emph{derived Ravenel--Wilson Hopf ring} of Definition \ref{def:rw-bar}, viewed as a derived Hopf ring via Remark \ref{rmk:rw-bar-derived}. 
Here we move classical rings into $\dHopf_{\F_p}$ by implicitly applying $\F_p[\ph]$.
The map $\LSym^\hotimes_{L_{(p)}}(W(p)) \to \LSym^\hotimes_{L_{(p)}}(W(1))$, denoted $f_p$ in the appendix, has the property that after tensoring down from $L_{(p)}$ to $\Z_{(p)}$, it just becomes $\LSym^\hotimes(F: W(p) \to W(1))$.
In particular \[ RW_{(p)} \otimes_L \Z \wequi \LSym^\hotimes(W(1)/W(p)), \] and so (recalling that $\Ga(1) = W(1)/W(p)$) \[ \LSym_n^\hotimes \Ga(1) \wequi (RW_{(p)} \otimes_L \Z)^n. \]
Since tensoring down from $L$ to $\Z$ is just a big colimit\NB{explain?}, which is preserved by the functors $Q$ and $P$, it will thus be enough to show that the canonical map \[ \Sigma^{2,1}QRW_{(p)}^n \to PRW_{(p)}^{n+1} \] is an equivalence.
We shall identify it in more familiar terms.
First note that, by the fundamental theorem about the derived Ravenel--Wilson Hopf ring, i.e., Theorem \ref{thm:rw-bar-exact}, $RW_{(p)}$ is actually concentrated in Chow degree zero, i.e., corresponds to a classical Hopf ring, namely the classical Ravenel--Wilson Hopf ring.\NB{the following references are really about $BP$ instead of $MU_{(p)}$...}
Second, the underlying Hopf algebras of this (in each degree) are bipolynomial, and hence, sums of Witt Hopf algebras \cite[Proposition 3.6]{WSWThesis}.
Consequently, the derived indecomposables/primitives in our statements are just given by the classical indecomposables and primitives.
Third, ``the canonical map'' in the statement of the lemma arises by construction from multiplication by a certain class in $\pi_{2,1} RW_{(p)}^1$, namely the restriction of the orientation class (corresponding to $\P^\infty \to \Sigma^{2,1}\MGL$) to $\P^1$; Ravenel--Wilson call this $b_{(0)}$.\NB{justify last assertion?}
Our assertion thus translates into the claim that \[ b_{(0)} \circ: \Sigma^{2,1} QRW_{(p)}^n \to PRW_{(p)}^{n+1} \] is an isomorphism (for $n \ge 1$).
This is \cite[Theorem 5.3]{ravenel-wilson}.
\end{proof}

We can now compute the homology of motivic Eilenberg--Mac Lane spaces at the characteristic.
\begin{corollary} \label{cor:EM-homology}
Let $k$ be a field of characteristic $p$ and $n \ge 1$.
There is a canonical equivalence \[ H\F_p[K(\Z_{(p)}(n),2n)] \wequi \LSym_n^\hotimes H\F_p[K(\Z_{(p)}(1),2)] \in \dHopf_{\F_p}^{bp} \subset \NHopf(\Mod_{H\F_p}(k)). \]
\end{corollary}
\begin{proof}
Immediate from Theorem \ref{thm:mainmain} and the fact that $H\F_p[\ph]$ transforms $D_n^\otimes$ into $\LSym_n^\hotimes$, i.e., Lemma \ref{lemm:dHopf-bp-nobp}.
\end{proof}

\section{Applications} \label{sec:applications} \todo{most proofs somewhat sketchy}

\begin{theorem} \label{thm:dual-steenrod}
Let $S$ be any scheme.
The motivic dual Steenrod algebra over $S$ has the expected form.
In other words, the canonical map \[ \bigoplus_{I \in B}\Sigma^{p(I),q(I)} H\F_p \to H\F_p \otimes H\F_p \in \SH(S) \] constructed in \cite[Theorem 10.26]{spitzweck2012commutative} is an equivalence.
\end{theorem}
\begin{proof}
By construction, the maps are compatible for all $S$, so it suffices to treat $S=\Spec(\Z)$, and thus by \cite[Proposition~B.3]{bachmann-norms} it suffices to treat $S=\Q$, $S=\F_\ell$ and $S=\F_p$.
For the first case see \cite{voevodsky2003reduced, voevodsky2010motivic}, and for the second see \cite[Theorem~1.1]{motivic-steenrod}.
In the final case, it was proved in \cite[Theorem~1.2]{frankland-spitzweck} that the desired map is split injective.
As we showed in the introduction (Corollary \ref{cor:intro-steenrod}), the target is cellular (as is the source, of course), and under $\pi_{**}$ they are degreewise finite dimensional of the same dimension.
Thus the map is an equivalence.
\end{proof}

\begin{theorem} \label{thm:hopkins-morel}
Let $S$ be any scheme.
Then $\MGL/(t_1, t_2, \dots) \wequi H\Z$.
\end{theorem}
\begin{proof}
As before, the statement need only be proved over prime fields.
The proof of \cite[\S7]{hoyois-algebraic-cobordism}, which so far worked only over fields of characteristic zero, now goes through over all fields, thanks to Theorem \ref{thm:dual-steenrod} (which was the only missing ingredient in positive characteristic).
\end{proof}

\begin{theorem}
Let $S$ be any scheme.
We have $s_* \MGL \wequi L_* \otimes H\Z$.
\end{theorem}
\begin{proof}
As explained in \cite[Theorem 4.6]{spitzweck2010relations}, this follows from Theorem \ref{thm:hopkins-morel} and knowing that $s_0H\Z \wequi H\Z$, for which see \cite[Theorem~1.6]{bachmann-elmanto-morrow}.
\end{proof}

\begin{theorem} \label{thm:EM} \NB{can this be generalized to all schemes?}
Let $S$ be a Dedekind scheme.
Then $H\F_p[K(\F_p(a), b)] \in \Mod_{H\F_p}(S)$ takes the expected form (see \cite[Theorem 3.16]{bachmann-wilson}), for $b \ge 2a$.
\end{theorem}
\begin{proof}
The argument is the same as in Theorem \ref{thm:dual-steenrod} (we have a map in from the expected answer over any base, and we can verify it is an equivalence over fields, the only new case following from Corollary \ref{cor:EM-homology}), once we prove that $K(\F_p(a), b)$ is stable under base change among Dedekind schemes.
This follows from \cite[Proposition 7.4]{bachmann-elmanto-morrow}.
\end{proof}

\begin{theorem}
Let $S$ be a Dedekind scheme.
Then $\MGL[\MGL^\bullet] \in \Mod_\MGL(S)$ takes the expected form, i.e., is canonically up to homotopy equivalent to the Ravenel--Wilson Hopf ring.
Here $\MGL^n := \Omega^\infty \Sigma^{2n,n}\MGL$.
\end{theorem}
\begin{proof}
The construction of \cite[\S5.1]{bachmann-wilson} gives a canonical map, up to homotopy, from the Ravenel--Wilson Hopf ring to $\MGL[\MGL^\bullet]$.
We claim that $\MGL^\bullet$ is stable under base change among Dedekind schemes.
We prove this at the end.
Assuming this for now, both sides of the canonical Ravenel--Wilson comparison map commute with base change, and so by \cite[Proposition~B.3]{bachmann-norms}, it suffices to check that the comparison map is an equivalence over fields.
There the argument of \cite[\S5]{bachmann-wilson} applies integrally: Theorem \ref{thm:hopkins-morel} supplies the integral Hopkins--Morel equivalence and slice filtration, and Corollary \ref{cor:EM-homology} supplies the previously missing Eilenberg--Mac Lane homology calculation in characteristic $p$.

It remains to establish the claim.
Write $F^n = f_n\MGL$ and $Q_n = \MGL/F^n$, for $n \ge 0$.
Theorem \ref{thm:hopkins-morel} and the construction in the proof of \cite[Theorem~4.6]{spitzweck2010relations} give
\[
F^0 = \MGL,
\qquad
F^n/F^{n+1} \wequi A_n := \Sigma^{2n,n}H\Z \otimes L_n,
\qquad
F^n \in \Sigma^{2n,n}\SH(S)^{\mathrm{veff}}.
\]
Moreover, this filtration commutes with base change: the cited proof constructs $F^n$ as a colimit of the diagram of Lazard monomials of degree at least $n$, with entries the corresponding Tate suspensions of $\MGL$, which is independent of the base.
Let $f:T \to S$ be a morphism of Dedekind schemes.
Instead of pulling back motivic spaces, we prove the stronger statement about pulling back Nisnevich sheaves.
Note that $f^*$ preserves colimits and finite products, and hence commutes with delooping commutative monoids.
It follows that, for $E \in \Shv_\Nis(\Sm_S,\mathrm{Sp})_{\ge 0}$, $f^*\Omega^\infty E \wequi \Omega^\infty f^*E$.
(Use that pullback of spectra is modeled by levelwise pullback of prespectra followed by spectrification, $E$ is represented by $(X,BX,B^2X,\dots)$, where $X=\Omega^\infty E$, and $(f^*X, f^*BX, f^* B^2 X, \dots)$ is already a spectrum by what we said above.)

Write $\omega^\infty$ for the functor taking a motivic spectrum to its underlying Nisnevich sheaf of spectra.
By the above, it suffices to show that $\omega^\infty\Sigma^{2r,r}\MGL$ is connective and that its formation commutes with Nisnevich pullback.
We have cofiber sequences
\[
\omega^\infty\Sigma^{2r,r}A_n
\longrightarrow \omega^\infty\Sigma^{2r,r}Q_{n+1}
\longrightarrow \omega^\infty\Sigma^{2r,r}Q_n.
\]
The first terms are connective by \cite[Corollary~6.24 and Remark~7.2]{bachmann-elmanto-morrow} and stable under Nisnevich pullback by \cite[Proposition~7.4]{bachmann-elmanto-morrow}.
Since $Q_0=0$ and pullback of spectral sheaves is exact, induction shows that $\omega^\infty\Sigma^{2r,r}Q_n$ is connective and stable under Nisnevich pullback.
Moreover, the fiber of
\[
\omega^\infty\Sigma^{2r,r}\MGL
\longrightarrow \omega^\infty\Sigma^{2r,r}Q_n
\]
is $\omega^\infty\Sigma^{2r,r}F^n$, which is $(n+r-1)$-connective by \cite[Proposition~3.44(1)]{bachmann-elmanto-morrow}.
Since pullback of spectral sheaves preserves connectivity, these approximations show that $\omega^\infty\Sigma^{2r,r}\MGL$ is also connective and stable under Nisnevich pullback (i.e., $f^*$ commutes with the limit of this rapidly converging tower), as needed.
\end{proof}

\begin{theorem}
Let $S$ be the spectrum of a field.
Then $\pi_{2*,*}\MGL \wequi L_*$.
\end{theorem}
\begin{proof}
The proof of \cite[Proposition 8.2]{hoyois-algebraic-cobordism} goes through unchanged.
\end{proof}

\begin{theorem} \label{thm:MGL-htpy}
Let $k$ be a field of characteristic $p$.
Then $\pi_{**} \MGL_p^\comp \wequi K_*^M(k)_p^\comp \otimes L_*$.
\end{theorem}
\begin{proof}
The argument in \cite[Lemma 2.40]{bachmann-stems} shows that the slice spectral sequence for $\MGL_k$ converges strongly.
It follows from \cite[Theorem~1.1]{geisser2000k} that $\pi_{**} H\Z_p^\comp \wequi K_*^M(k)_p^\comp$.
Consequently the spectral sequence collapses as indicated.
\end{proof}

\begin{theorem} \label{thm:stable-stems}
Let $k$ be a field of characteristic $p$.
Then $\pi_{**} \1_p^\comp \wequi \Ext(BP_{2*}BP) \otimes_{\Z[\eta]} K_*^{MW}(k)_p^\comp$.
\end{theorem}
\begin{proof}
This follows from the convergence of the motivic Adams--Novikov spectral sequence (see \cite[\S5.1 and Theorem~7.3.4]{mantovani2018localizations} and \cite[Lemma~6.43]{bachmann-resolution}) and Theorem \ref{thm:MGL-htpy}.\NB{details}
\end{proof}

%

\appendix
\section{The derived Ravenel--Wilson Hopf ring}\label{app:rw-derived}
\providecommand{\gpt}[1]{\todo[color=cyan!25,linecolor=cyan!60!black,bordercolor=cyan!60!black]{\textcolor{teal!80!black}{#1}}}
\providecommand{\Free}{\operatorname{Free}}
\providecommand{\Ga}{{\mathbb G}_a}
\providecommand{\RW}{\mathrm{RW}}
\providecommand{\dRW}{\scr{RW}}
\providecommand{\Sym}{\mathrm{Sym}}
\providecommand{\LSym}{\operatorname{LSym}}
\providecommand{\Vect}{\mathrm{Vect}}
\providecommand{\hotimes}{\odot}
\providecommand{\htimes}{\otimes}
\providecommand{\hH}{\mathfrak H}
\providecommand{\frob}{\mathrm f}
\providecommand{\ver}{\mathrm v}
\providecommand{\wittprojection}{\mathrm \pi}
\providecommand{\HCAlg}{\mathrm{HCAlg}}
\providecommand{\Hopf}{\mathrm{Hopf}}
\providecommand{\MU}{\mathrm{MU}}
\providecommand{\Set}{\mathrm{Set}}
\providecommand{\Ab}{\mathrm{Ab}}
\providecommand{\cn}{\mathrm{cn}}
\providecommand{\cl}{\pi}
\providecommand{\A}{\mathbb A}
\providecommand{\cof}{\operatorname{cof}}
\providecommand{\comp}{\wedge}
\providecommand{\wequi}{\simeq}
\providecommand{\ph}{\mathord-}
\providecommand{\fr}{\mathrm{fr}}
\providecommand{\Hom}{\operatorname{Hom}}
\providecommand{\CAlg}{\mathrm{CAlg}}
\providecommand{\cCAlg}{\mathrm{cCAlg}}
\providecommand{\bCAlg}{\mathrm{bCAlg}}
\providecommand{\1}{\mathbbm{1}}

In this section we introduce the \emph{derived Ravenel--Wilson Hopf ring} and show that it coincides with the ordinary Ravenel--Wilson Hopf ring of \cite{ravenel-wilson}.
The material has been extracted from our forthcoming article \cite{bachmann-wilson} in order to make the current manuscript more self-contained.

We fix a prime $p$ and work over $\F_p$.
We write $L=L_{(p)}$ for the $p$-local Lazard ring.  We write $\Hopf_{\F_p}$ for the category of evenly graded bicommutative Hopf algebras over $\F_p$. We write $H\hotimes K$ for the tensor product of Hopf algebras, so maps $H\hotimes K\to T$ correspond to bilinear coalgebra maps $H\htimes K\to T$. A \emph{Hopf ring} is a commutative ring object for $\hotimes$, with a separate external grading.
The Hopf ring associated with a graded commutative ring $A$ is denoted $\F_p[A]$.

If $H$ is a Hopf algebra, $\Sym(H)$ denotes the free commutative Hopf ring on $H$ placed in external degree one; its component of external degree $n$ is $H^{\hotimes n}/\Sigma_n$. Similarly $\Lambda(H)$ denotes the exterior algebra for $\hotimes$, obtained by killing the image of the diagonal. For a graded coefficient ring $A$, we write $\Sym_A(H)=\F_p[A]\hotimes\Sym(H)$.

We regard a simplicial Hopf algebra as a simplicial object in the additive category of graded Hopf algebras. Its homology, written $\hH_i$, is the homology of the normalized chain complex under Dold--Kan. The same convention applies degreewise to simplicial Hopf rings.

\subsection{Frobenius, Verschiebung, and Dieudonn\'e theory}\label{app:rw-witt}

For each $n>0$ we have the functor $\phi_n^*: \Vect_{\F_p}^\N \to \Vect_{\F_p}^\N$ defined by $(\phi_n^*V)_{ni} = V_i$ and $(\phi_n^*V)_j = 0$ for $j$ not a multiple of $n$.
The functor $\phi_n^*$ admits an evident symmetric monoidal structure and hence extends canonically to $\cCAlg(\Vect_{\F_p}^\N)$, and so on.
Moreover, $\phi_n^*$ is symmetric monoidal for $\hotimes$.

Given $A \in \CAlg(\Vect_{\F_p})$ we have the Frobenius endomorphism $F: A \to A, a \mapsto a^p$.
This is in fact a symmetric monoidal endotransformation of the identity functor.
When working with graded algebras, the Frobenius changes grading, and so for $A \in \CAlg(\Vect_{\F_p}^\N)$ takes the form $F: \phi_p^* A \to A$.
Being symmetric monoidal, this upgrades to \[ F:\phi_p^* \to \id \in \Fun(\Hopf_{\F_p}, \Hopf_{\F_p}). \]

We shall also need the \emph{Verschiebung} map, which arises in a dual fashion.
For a (graded) vector space $C$, denote by $C^{tp}$ the quotient of $(C^{\otimes p})^{C_p}$ by the subspace generated by vectors of the form $\bigoplus_{\sigma \in C_p} \sigma v$.
Then there is a canonical equivalence of (oplax symmetric monoidal) endofunctors $(\ph)^{tp} \wequi \phi_p^*$.\NB{ref?}
Consequently for a graded coalgebra we obtain a map \[ V: C \xrightarrow{\Delta^p} (C^{\otimes p})^{h\Sigma_p} \to (C^{\otimes p})^{tp} \wequi \phi_p^* C \] called the \emph{Verschiebung}.
Again this is a symmetric monoidal endotransformation and hence induces \[ V: \id \to \phi_p^*  \Fun(\Hopf_{\F_p}, \Hopf_{\F_p}). \]

We recall the following important facts.
\begin{lemma} \label{lemm:FV-basics}
\begin{enumerate}
\item We have $FV = [p] = VF$.
\item Given $H, K \in \Hopf_{\F_p}$ we have $V_{H \hotimes K} = V_H \hotimes V_K$ and the following diagram commutes
\begin{equation*}
\begin{CD}
\phi_p^*H \hotimes \phi_p^*K @>{\wequi}>> \phi_p^*(H \hotimes K) \\
@A{\id_H \hotimes V_K}AA                   @V{F_{H \hotimes K}}VV \\
\phi_p^* H \hotimes K @>{F_H \hotimes \id_K}>> H \hotimes K.
\end{CD}
\end{equation*}
\end{enumerate}
\end{lemma}
\begin{proof}
(1) That $FV = [p]$ is essentially the point of the definition, and follows from the fact that $F$ vanishes on elements of $(H^{\otimes p})^{C_p}$ which are \emph{induced}, i.e., in the image of $(H^{\otimes p})_{C_p}$.
Naturality of $V$ applied to the map $F: \phi_p^* H \to H$ shows that $V_H F_H = F_{\phi_p^* H} V_{\phi_p^* H} = [p]_{\phi_p^* H}$, as needed.

(2) It suffices to check the first equality on free Hopf algebras, where it reduces to the analogous claim for coalgebras and the usual tensor product, which is clear (by construction).
The second statement follows from the distributivity of $\hotimes$ over $\htimes$.
(In slightly more detail, given $x \in H_{i/p}$  and $y \in K_j$ we can express $x^p \hotimes y$ in terms of the iterated diagonal $\Delta^p(y)$, writing out this formula precisely yields the desired commutative square.)
\end{proof}

Given a Hopf algebra $H$, applying indecomposables to Frobenius and Verschiebung yields maps \[ QF: \phi_p^* QH \to QH \quad\text{and}\quad QV: QH \to \phi_p^* QH. \]
Clearly $QF = 0$ (i.e., $a^p$ is always decomposable).
This is not so for $QV$ and so the graded vector space $QH$ acquires an additional structure, namely the Verschiebung endomorphism.
Let us write $\Vect_{\F_p}^{\N,V}$ for the category of such vector spaces, i.e., graded vector spaces $X$ together with a map $V: X \to \phi_p^* X$.

\begin{proposition} \label{prop:indecomp-charp}
The above construction upgrades to a symmetric monoidal functor \[ Q: \Hopf_{\F_p} \to \Vect_{\F_p}^{\N,V} \] which restricts to an equivalence \[ \Hopf_{\F_p}^{\cn,F=0} \wequi (\Vect_{\F_p}^{\N,V})^\cn. \]
Here by $\Hopf_{\F_p}^{\cn,F=0} \subset \Hopf_{\F_p}^{\cn}$ we denote the subcategory of those Hopf algebras whose Frobenius endomorphism vanishes, i.e., $H$ such that $F: \phi_p^* H \to H$ factors through $\F_p$.
\end{proposition}
\begin{proof}
A symmetric monoidal structure is inherited from the universal property of $\hotimes$ and Lemma \ref{lemm:FV-basics}(2).
The fact that the restricted functor is an equivalence is part of Dieudonné theory: there is an equivalence $D$ from $\Hopf_{\F_p}^{\cn}$ to a category of graded modules over $\Z_p$ with an action by $F$ and $V$ such that $FV = VF = p$.
Moreover there is a canonical short exact sequence $D\phi_p^* H \xrightarrow{DF} DH \to QH \to 0$.
See, e.g., \cite[Theorem 4.7, Lemma 4.4]{goerss-dieudonne}.
The result follows.
\end{proof}

\subsubsection{Additive group} \label{subsec:Ga}
An important Hopf algebra for us will be $\Ga$, the Hopf algebra corresponding to the additive group over $\F_p$.
As a ring it is given by \[ \Ga = \F_p[a_1, a_2, a_3, \dots]/(a_i a_j = {i+j \choose j} a_{i+j}), \text{ where } |a_i| = i. \]
The coalgebra structure is determined by \[ \Delta(a_i) = \sum_{i = j + k} a_j \otimes a_k. \]
The underlying coalgebra of $\Ga$ is the base change to $\F_p$ of the $L$-coalgebra $\A^1$ defined in Definition~\ref{def:A1}. We use the same symbol for its base change whenever the coefficient ring is clear.

\begin{example} \label{ex:Ga-hotimes-Ga}
We now find $a_i^p = 0$\NB{details?} and \[ \Ga \wequi \F_p[a_{(0)}, a_{(1)}, \dots]/a_{(i)}^p, \text{ where } a_{(i)} = a_{p^i}. \]
It follows that $F=0$.
One may check that $V: \Ga \to \phi_p^* \Ga$ sends $a_{(0)}$ to $0$, whereas for $i>0$ it sends $a_{(i)}$ to an element which might be denoted $\phi_p^*a_{(i-1)}$.
In particular $V$ is surjective and we tend to denote the generators of $\phi_p^* \Ga$ by their names coming from $\Ga$, i.e., we tend to think of the element $\phi_p^*a_{(i-1)}$ as $Va_{(i)}$ and even denote it by just $a_{(i)}$.
We thus have \[ Q\phi_{p^n}^*\Ga \wequi \F_p\{a_{(n)}, a_{(n+1)}, \dots\} \text{ with } Va_{(i)} = a_{(i-1)} \text{ and } a_{(i)} = 0 \text{ for } i<n. \]
Using the equivalence between connected Hopf algebras and graded vector spaces with Verschiebung (Proposition \ref{prop:indecomp-charp}) we see that, for example, \begin{gather*} \Ga \hotimes \Ga \wequi \F_p[a_{(i)} \circ a_{(j)}, i,j \ge 0]/(a_{(i)} \circ a_{(j)})^p = 0 \text{ with } V(a_{(i)} \circ a_{(j)}) = a_{(i-1)} \circ a_{(j-1)} \\ \text{ where } a_{(i)} \circ a_{(j)} = 0 \text{ whenever $i<0$ or $j<0$}. \end{gather*}
Equivalently \[ \Ga \hotimes \Ga \wequi \bigotimes_{i > 0} \phi_{p^i+1}^*\Ga \otimes \phi_2^*\Ga \otimes \bigotimes_{j > 0} \phi_{p^j+1}^*\Ga, \] with the copies of $\Ga$ corresponding respectively to $a_{(i)} \circ a_{(0)}$, $a_{(0)} \circ a_{(0)}$ and $a_{(0)} \circ a_{(j)}$.
Many similar computations are possible and straightforward, e.g., \[ \Sym^2 \Ga \wequi \F_p[a_{(i)} \circ a_{(j)}, i,j \ge 0]/(a_{(i)} \circ a_{(j)})^p = 0,a_{(i)} \circ a_{(j)} = a_{(j)} \circ a_{(i)}. \]
\end{example}

\subsubsection{Hopf algebras with a lift of Verschiebung} \label{subsec:lift-of-v}
Following Goerss, by a \emph{Hopf algebra with a lift of Verschiebung} we mean a pair of
\begin{itemize}
\item an object $H \in \HCAlg(\Mod_{\Z_p}^{\heart\N})$, and
\item a Hopf algebra map $V: H \to \phi_p^* H$ such that $V \otimes \F_p = V_{H \otimes \F_p}$.
\end{itemize}
We denote the category of these objects by $\HCAlg_{\Z_p}^V$.
It carries a symmetric monoidal structure, i.e., a version of $\hotimes$.
Indeed given two Hopf algebras over $\Z_p$ $H_1, H_2$, with lifts of Verschiebung $V_1, V_2$, the map $V_1 \hotimes V_2$ lifts the Verschiebung on $H_1 \hotimes H_2$ (Lemma \ref{lemm:FV-basics}(2)).

Let us write $\Mod_{\Z_p}^{\N,V}$ for the category of graded $\Z_p$-modules together with an action by $V$.
Then as before we have the functor of indecomposables $Q: \HCAlg_{\Z_p}^V \to \Mod_{\Z_p}^{\N,V}$.
\begin{proposition}[Goerss] \label{prop:Hopf-with-V}
The functor $Q$ is cocontinuous, symmetric monoidal and restricts to an equivalence of categories between (on either side) the connected, torsion-free objects.
\end{proposition}
\begin{proof}
Symmetric monoidality is proved as before.
For the equivalence, see \cite[Theorem 2.7]{goerss-dieudonne}.
\end{proof}

\begin{example} \label{ex:M-versch}
For $n \in \N$ we have the object $M_n \in \Mod_{\Z_p}^{\N,V}$ which is given by $\Z_p$ in degrees of the form $n/p^*$ and by $0$ else, with $V$ the identity whenever this makes sense.
Applying the inverse of $Q$, we obtain important Hopf algebras with lift of Verschiebung, still denoted $M_n$.
\end{example}

\begin{example} \label{ex:Mp-hotimes-Mp}
Computations of tensor products of Hopf algebras with lift of Verschiebung are straightforward using the above equivalence.
For example, let us write $M_{p^\infty} := \colim_i M_{p^i}$.
Then \[ M_{p^\infty} \hotimes M_{p^\infty} \wequi \bigoplus_{i > 0} \phi_{p^i + 1}^* M_{p^\infty} \oplus \phi_2^* M_{p^\infty} \oplus \bigoplus_{j > 0} \phi_{p^j + 1}^* M_{p^\infty}. \]
Similarly $M_n \hotimes M_m$ is easily seen to be given by a finite sum of copies of $M_i$, for various $i$.
Symmetric and exterior powers are treated similarly.
\end{example}

\begin{remark}
The similarities between Examples \ref{ex:Mp-hotimes-Mp} and \ref{ex:Ga-hotimes-Ga} are not just superficial.
It is in fact the case that $(M_{p^\infty} \otimes \F_p)/F \wequi \Ga$; see Lemma \ref{lemm:W-mod-F}.
\end{remark}

\subsection{Homomorphic curves}\label{subsec:RW-moduli}
Let $\scr C$ be an additive, symmetric monoidal $1$-category (in which the tensor product is additive).
By a Hopf ring in $\scr C$, we mean a commutative graded ring object in the cartesian symmetric monoidal category $\cCAlg(\scr C)$.
Equivalently, giving $R \in \cCAlg(\scr C)$ the structure of a Hopf ring is to supply a factorization of \[ \Hom_{\cCAlg(\scr C)}(\ph, R): \cCAlg(\scr C)^\op \to \Set \] through commutative graded rings.
Similarly, given an ordinary commutative graded ring $A$, we can define a Hopf ring under $A$ in $\scr C$ by requiring a factorization of the represented functor through commutative graded $A$-algebras.

\begin{example}
By definition, if $R$ is a Hopf ring under $A$ then the set of (graded) coalgebra maps $\1 \to R$ is an $A$-algebra, and thus contains elements $[a]$ for any $a \in A$.
In the situation of the above remark, elements $a \in A$ (i.e., maps $a: * \to A$) induce the elements $[a] := \1_\scr{C}[a] \in \1_\scr{C}[A]$ (i.e., maps $[a]: \1_\scr{C} \to \1_\scr{C}[A]$) via the functor $\1_\scr{C}[\ph]$.
\end{example}

\begin{example} \label{ex:hopfring-frob}
Consider the case $A=\F_p$, i.e., Hopf rings with $[p]=0$.
Any functor $\cCAlg(\scr C)^\op \to \CAlg_{\F_p/}^\heart$ has a canonical endomorphism, given sectionwise by the Frobenius.
This induces a canonical endotransformation of the identity functor on Hopf rings under $\F_p$.
In more down to earth terms, if $R$ is a Hopf ring with $[p]=0$, then the coalgebra map $R \xrightarrow{\Delta_p} R^{\otimes p} \xrightarrow{m_p} R$ is a Hopf algebra map, where $m$ means the Hopf ring multiplication.
\end{example}

\begin{definition} \label{def:ModLfr}
We write $L$ for the $p$-local Lazard ring, viewed as a commutative graded ring, and $\Mod_L^\fr$ for the symmetric monoidal $1$-category of free graded $L$-modules on generators in non-negative degrees.
For a multiplicative subset $S \subset \Z$, define $\Mod_{S^{-1}L}^\fr$ similarly.
\end{definition}
\begin{remark}\label{rmk:lazard-ring}
We shall use frequently and without further comment that $L$ (as defined above) is also the universal ring carrying a (commutative, 1-dimensional) formal group law \cite{quillen-MU}.
In other words there exist elements $a_{i,j} \in L_{1+i+j}$ such that the formal power series $\sum_{i,j} a_{i,j} x^i y^j \in L\fpsr{x,y}$ is a formal group law with a certain universal property.
\end{remark}

\begin{definition} \label{def:A1}
By standard computations in homotopy theory, the $p$-local complex cobordism homology
\[\A^1=\MU_{2*}\C\P^\infty\otimes_{\MU_{2*}}L\]
defines an object of $\cCAlg(\Mod_L^\fr)$.
For $r\ge 0$, write
\[A_r=\A^1_{\le r}=\MU_{2*}\C\P^r\otimes_{\MU_{2*}}L\subset\A^1.\]
\end{definition}
\begin{remark}
The multiplication map $\mu: \C\P^\infty \times \C\P^\infty \to \C\P^\infty$ and the pointing $* \to \C\P^\infty$ in fact propagate $\A^1$ to a Hopf algebra in $\Mod_L^\fr$.
We will however usually only think of $\A^1$ as a coalgebra, and keep track of the multiplication and pointing explicitly.
\end{remark}
\begin{remark}
We can define the above objects without reference to homotopy theory.
The coalgebra $\A^1$ is free over $L$, with basis $\{b_i\}_{i \ge 0}$ where $|b_i| = i$, $\Delta(b_i) = \sum_{i = k+j} b_k \otimes b_j$, and $b_i b_j = {i+j \choose j} b_{i+j}$.
Then the algebra dual to $\A^1$ is the formal power series ring $L\fpsr{c}$ with $b_i$ dual to $c^i$.
The map $\mu$ is most easily described on duals, where $\mu^*: L\fpsr{c} \to L\fpsr{x,y}$ is induced by sending $c$ to the universal formal group law $\sum_{ij} a_{ij} x^i y^j$ (see Remark \ref{rmk:lazard-ring}).
\end{remark}

In what follows, we will routinely find ourselves in situations where certain seemingly infinite sums can (and must) be formed.
The basis of this is the following.
\begin{definition} \label{def:conil}
We call a pointed coalgebra $C \in \cCAlg(\scr C)_{\1/}$ \emph{conilpotent} if there is $N \gg 0$ such that the iterated reduced diagonal \[ C \xrightarrow{\Delta^N} C^{\otimes N} \to (C/\1)^{\otimes N} \] is null.
We call a pointed coalgebra $C \in \cCAlg(\scr C)_{\1/}$ \emph{ind-conilpotent} if we can present $C = \colim_i C_i$ as a filtered colimit in pointed coalgebras, where each $C_i$ is conilpotent.
\end{definition}
\begin{example}
We have $\A^1 \wequi \colim_r A_r$, exhibiting $\A^1$ as ind-conilpotent.
(Indeed $\MU^{2*}\C\P^r\otimes_{\MU_{2*}}L\wequi L[c]/c^{r+1}$, from which conilpotence of $A_r$ easily follows.)
\end{example}

\begin{remark} \label{rmk:pointed-vanishing}
Let $R \in \CAlg(\HCAlg(\scr C))$ and $C = \colim_i C_i$ an ind-conilpotent coalgebra, with each $C_i$ conilpotent.
For any pointed map $x: C_i \to R$ ($R$ pointed by $[0]$), $x^n$ factors through the reduced $n$-fold diagonal (because the multiplication in $R$ induces $R/\1 \otimes R/1 \to R/1$, the pointing being by $[0]$); thus $x^n = [0]$ for $n \gg 0$.
Since also $\Hom(C, R) \wequi \lim_i \Hom(C_i, R)$, for $x \in \Hom(C, R)$ we can make sense of infinite sums like $\sum_i x^i$ (since they are finite when restricted to any $C_i$).
\end{remark}

We will now work towards defining the Ravenel--Wilson Hopf rings, via the moduli problem they corepresent.
This is done in the context of the following data:
\begin{itemize}
\item A commutative graded formal group law $L \to A$.
\item A strongly additive symmetric monoidal functor $F: \Mod_L^\fr \to \scr C$ (where $\scr C$ is some strongly additive symmetric monoidal $1$-category).
\end{itemize}
\begin{remark} \label{rmk:FA1-indconil}
Note that $F(\A^1) \wequi \colim_i F(A_i)$, since the colimit on underlying objects (ignoring the coalgebra structures) is an infinite sum (and thus preserved by $F$, by assumption).
In other words, $F(\A^1)$ is ind-conilpotent.
\end{remark}
\begin{definition} \label{def:homomorphic-curve}
In the above situation, for any Hopf ring $R$ in $\scr C$ under $A$ we define the \emph{set of homomorphic curves in $R$} as the set of all coalgebra maps $\gamma: F(\A^1) \to R^1$ such that:
\begin{itemize}
\item $\gamma$ is pointed.
\item Writing $\mu: \A^1 \times \A^1 \to \A^1$ for the multiplication map and $p_i: \A^1 \times \A^1 \to \A^1$ for the projections, we have \[ \mu^*(\gamma) = \sum_{i,j} [a_{ij}] p_1^*(\gamma)^i p_2^*(\gamma)^j. \]
  Here on the right hand side, the multiplication uses the multiplication in $R$, the sum uses the addition in $R$, and the sum is defined by Remarks \ref{rmk:pointed-vanishing} and \ref{rmk:FA1-indconil}.
\end{itemize}

If this functor is corepresentable by a Hopf ring in $\scr C$ under $A$ (e.g., $\scr C$ presentable), we denote by $\RW(\scr C, A)$ the corepresenting object.
\end{definition}
\begin{remark} \label{rmk:construct-RW}
In other words, there exists at most one Hopf ring $\RW(\scr C, A)$ under $A$ in $\scr C$ (up to unique isomorphism) such that $A$-algebra homomorphisms from $\RW(\scr C, A)$ to $R$ are the same as homomorphic curves in $R$.
In fact, if $\scr C$ has enough colimits, we can build a pushout of the form
\begin{equation*}
\begin{CD}
F(\A^2) @>{\alpha}>> F(\A^1) \\
@VVV           @VVV  \\
\1_\scr{C}[L] @>>> \RW(\scr C, A).
\end{CD}
\end{equation*}
Here $F(\ph)$ denotes the free Hopf ring on a pointed coalgebra placed in degree $1$, and the map $\alpha$ corepresents the map \[ \gamma \mapsto \mu^*(\gamma) - \sum_{i,j} [a_{ij}] p_1^*(\gamma)^i p_2^*(\gamma)^j. \]
Then the bottom right corner corepresents homomorphic curves, essentially by construction.
\end{remark}

\begin{remark}[Additive curves] \label{rmk:additive-curves}
Let $F: \Mod_L^\fr \to \scr C$ be an additive symmetric monoidal functor.
Then the Hopf algebra $(\A^1, \mu)$ is sent by $F$ to another Hopf algebra $(F\A^1, F\mu)$.
It is possible that $(F\A^1, F\mu)$ is isomorphic to $F\Ga$.
(E.g., take $L \to \F_p$ and the induced functor $F: \Mod_L^\fr \to \Vect_{\F_p}$.)
In this case homomorphic curves in $\scr C$ are usually called \emph{additive curves}.
\end{remark}

\subsection{Curves in formal groups}\label{subsec:curves}
In preparation for the next section, here we discuss some further structures in $\Hopf_{\F_p}$.
We recalled above the functors $\phi_n^*$ and the endotransformations $F,V$.

\begin{definition}
By a \emph{formal group} we shall mean a functor \[ X: \cCAlg(\Vect_{\F_p}^\N)^{\mathrm{fd},\cn,\op}_{\F_p/} \to \Ab \] from (the opposite of) pointed, connected, finite dimensional (over $\F_p$) coalgebras to abelian groups, such that there exists a pointed coalgebra $\tilde X$ and an isomorphism \[ X|_{\Fun(\cCAlg(\Vect_{\F_p}^\N)^{\mathrm{fd},\cn, \op}_{\F_p/}, \Set)} \wequi \Hom(\ph, \tilde X). \]
\end{definition}

One class of examples arises as follows.
\begin{example} \label{ex:fg-1}
Let $H \in \Hopf_{\F_p}$.
Define $H^\comp(C) = \Hom_{\cCAlg(\Vect_{\F_p}^\N)_{\F_p/}}(C, H)$.
This inherits an abelian group structure from $H$, and so is a formal group.
\end{example}

Now let $X$ be a formal group.
We call \[ X_n = X(\phi_n^* \A^1) \] the set of \emph{pointed curves of degree $n$ in $X$}.

\subsubsection{Formal groups and Hopf algebras}
Our definition of formal groups has mostly psychological utility, as the following result shows.
\begin{lemma} \label{lemm:HCAlg-vs-FG}
The functor \[ (\ph)^\comp: \Hopf_{\F_p}^\cn \to (\text{(formal groups)}) \] is an equivalence of categories.
\end{lemma}
\begin{proof}
Let us call a functor
\[X: \cCAlg(\Vect_{\F_p}^\N)^{\mathrm{fd},\cn,\op}_{\F_p/} \to \Set\]
representable by a pointed coalgebra $\tilde X$ a \emph{pointed formal scheme}.
We have
\[(\text{formal groups}) = \Ab(\text{pointed formal schemes}).\]
By the fundamental theorem of coalgebras \cite[Corollary 2.2.2]{sweedler-hopf} we have
\[\cCAlg(\Vect_{\F_p}^\N)^\cn_{\F_p/}
 \wequi \Ind(\cCAlg(\Vect_{\F_p}^\N)^{\mathrm{fd},\cn}_{\F_p/}).\]
Consequently the restricted Yoneda functor
\[\cCAlg(\Vect_{\F_p}^\N)^\cn \to (\text{formal schemes})\]
is fully faithful.
It will suffice to prove that it is essentially surjective.
I.e., given a pointed coalgebra $\tilde X$, we need to find a pointed connected coalgebra $X$ such that $X^\comp \wequi \tilde X^\comp$.
We can take $X$ to be the image of \[ \bigvee_{C \to \tilde X} C \to \tilde X, \] where the wedge is over all morphisms of pointed coalgebras $C \to \tilde X$ with $C$ finite dimensional.\NB{details?}
\end{proof}

Denote by $\Free$ the left adjoint of the forgetful functor $\Hopf_{\F_p} \to \cCAlg(\Vect_{\F_p}^\N)_{\F_p/}$.
The Hopf algebra $\Free\A^1$ is called the \emph{big Witt Hopf algebra}.
As a ring we have \[ \Free\A^1 = \F_p[a_1, a_2, \dots]; \] in particular $\Free\A^1$ is connected.
In fact, $(\Free\A^1)^\comp$ is the free formal group on $\A^1$, in the following sense.

\begin{corollary} \label{cor:Free-FG}
Let $X$ be a formal group and $\gamma \in X_n$.
Then $\gamma$ extends uniquely to a morphism of formal groups $(\Free\phi_n^*\A^1)^\comp \to X$.
\end{corollary}
\begin{proof}
Immediate from Lemma \ref{lemm:HCAlg-vs-FG} and connectivity of $\Free\phi_n^*\A^1$.
\end{proof}

\begin{corollary} \label{cor:fg-p-local}
All formal groups are $p$-local, i.e., for any $m$ with $(m,p)=1$, the multiplication by $m$ map is an isomorphism on any formal group.
\end{corollary}
\begin{proof}
By Lemma~\ref{lemm:HCAlg-vs-FG}, it suffices to consider a connected graded Hopf algebra. In any fixed positive degree, a sufficiently high iterate of $V$ vanishes for degree reasons. Thus $[p]=FV$ acts locally nilpotently by Lemma~\ref{lemm:FV-basics}, and hence multiplication by any integer prime to $p$ is invertible.
\end{proof}

\subsubsection{Frobenius and Verschiebung for curves}
We will now construct the Frobenius and Verschiebung operators on the sets $X_n$.

\begin{definition}
Let $X$ be a formal group and $\gamma: \phi_n^* \A^1 \to X$ a pointed curve of degree $n$.
Form \[ \phi_n^* \A^m \xrightarrow{p_1^*\gamma + \dots + p_m^*\gamma} X, \] which is symmetric in the $m$ coordinates, and hence factors uniquely through the quotient coalgebra $\phi_n^*\A^m/\Sigma_m$.\footnote{The existence of $\tilde X$ implies both that the factorization exists and is unique.}
Define the \emph{Frobenius operator} $\frob_m(\gamma) \in X_{mn}$ as the composite \[ \phi_{mn}^*\A^1 \to \phi_n^*\A^1 \htimes \dots \htimes \phi_{nm}^* \A^1 \wequi \phi_n^*\A^m/\Sigma_m \xrightarrow{p_1^*\gamma + \dots + p_m^*\gamma} X, \] where the equivalence is induced by the symmetric function theorem.\footnote{In other words, passing to duals, an isomorphism of coalgebras $\phi_n^*\A^1 \htimes \dots \htimes \phi_{nm}^* \A^1 \wequi \phi_n^*\A^m/\Sigma_m$ is the same as an isomorphism of rings $\F_p[X_n, \dots, X_{mn}] \wequi \F_p[Y_n^{(1)}, \dots, Y_n^{(m)}]^{\Sigma_m}$, which is supplied by the classical symmetric function theorem.}

We also define the \emph{Verschiebung operators} given by \[ \ver_m: X_{mn} \to X_{n}, (\phi_{mn}^*\A^1 \xrightarrow{\gamma} X) \mapsto (\phi_n^* \A^1 \to \phi_{mn}^* \A^1 \xrightarrow{\gamma} X). \]
Here the map $\phi_n^* \A^1 \to \phi_{mn}^* \A^1$ is dual to $\F_p[t_{nm}] \to \F_p[t_n], t_{mn} \mapsto t_n^m$.
\end{definition}

\begin{remark}
Essentially by definition, $\phi_n^* \Free\A^1$ corepresents the functor \[ \Hopf_{\F_p} \to \Ab, \quad H \mapsto H^\comp_n. \]
The operations $\frob_m$, $\ver_m$ on pointed curves thus correspond to operations \[ F_m: \phi_{mn}^* \Free\A^1 \to \phi_n^*\Free\A^1 \text{ and } V_m: \phi_n^*\Free\A^1 \to \phi_{mn}^* \Free\A^1. \]
\end{remark}

\begin{lemma} \label{lemm:FV-on-W}
On $\Free\A^1$ have $F_p = F$ and $V_p = V$.
\end{lemma}
\begin{proof}
The map $F_p$ is obtained as the composite \[ \Free \phi_p^* \A^1 \to \Free(\A^p/\Sigma_p) \to (\Free \A^1)^{\htimes p}/\Sigma_p \to \Free \A^1, \] where the first map comes from the symmetric function theorem, the second map is a categorical comparison map, and the third is the multiplication.
Unwinding the definitions, this is precisely $F$.\NB{details?}
Since $V$ is already defined for coalgebras and functorial, in order to prove that $V=V_p$ for $\Free\A^1$, it suffices to prove the same for $\A^1$ itself.
In other words, we must show that $V$ on $\A^1$ is dual to the Frobenius.
This is true essentially by definition.
\end{proof}

\begin{example} \label{ex:lift-FreeA1}
It follows that $\Free\A^1$ lifts to a torsion-free Hopf algebra over $\Z_p$, with a lift of Verschiebung given by $V_p$ (constructed the same way over $\Z_p$ as over $\F_p$).
We can analyze this lifted Hopf algebra as in \S\ref{subsec:lift-of-v}, i.e., in terms of the action of $V$ on its indecomposables.
Using \cite[(13.4)]{hazewinkel-witt-vectors}, we see that the indecomposables are given, in the notation of Example \ref{ex:M-versch}, by $M_\infty := \colim_n M_n$.
\end{example}

\begin{remark} \label{rmk:symm-hazewinkel}
The Hopf algebra $\Free\A^1$ is denoted by $\mathrm{Symm}$ in \cite{hazewinkel-witt-vectors}.
(See in particular (9.4) and (9.5) in \emph{loc. cit.})
Moreover, unwinding the definitions, the maps $F_m$ and $V_m$ we have defined coincide with Hazewinkel's \cite[\S13]{hazewinkel-witt-vectors}.
This is useful because \emph{loc. cit.} establishes numerous properties of these operations.
(Hazewinkel does not keep track of gradings the way we do, but most desired properties can be deduced nonetheless.)
\end{remark}

These operators satisfy the following relations.
\begin{lemma} \label{lemm:frob-ver-identities}
The operations $\frob_m$ and $\ver_m$ are natural in $X$ and satisfy
\begin{align*}
\frob_1 &= \ver_1 = 1 \\
\frob_m\ver_m &= m \\
\frob_m\frob_m' &= \frob_{m m'} \\
\ver_m\ver_{m'} &= \ver_{m m'} \\
\frob_m\ver_{m'} &= \ver_{m'}\frob_m \quad (m',m)=1 \\
\ver_p\frob_p &= p,
\end{align*}
whenever they are defined.
\end{lemma}
Note that $\ver_m$ can only be applied to curves of degree divisible by $m$.
So the second equation above claims that the endomorphism $\frob_m \ver_m$ of $X_{nm}$ is given by multiplication by $m$ for every $X, n$.
If $m$ does not divide $n$, then $\ver_m$ on $X_n$ is not defined, not even as zero, and we are not claiming that multiplication by $m$ on $X_n$ is zero.
\begin{proof}
By Corollary~\ref{cor:Free-FG}, this reduces to the universal case $X=(\Free\A^1)^\comp$, i.e., the big Witt vectors.
The relations in this context are proved in \cite[Theorem 13.48]{hazewinkel-witt-vectors}, except for the last one, which is a consequence of Lemma~\ref{lemm:FV-on-W}.
\end{proof}

\begin{definition}
We call a pointed curve $\gamma \in X_n$ \emph{$p$-typical} if $\frob_\ell \gamma = 0$ for all primes $\ell \ne p$.
\end{definition}

With this at hand, we can now construct the \emph{$p$-typicalization operator}.
\begin{lemma} \label{lemm:p-typicalization}
Let $\mu$ be the Möbius function.\footnote{Thus $\mu(n) = 0$ if $n$ is divisible by a square $>1$ and $\mu(n) = (-1)^k$ if $n$ is a product of $k$ distinct primes.}
Suppose that $X$ is a formal group.
The operator
\[
\wittprojection_p = \sum_{(m,p)=1} \frac{\mu(m)}{m} \ver_m \frob_m
\]
is well-defined and has the following properties:
\begin{enumerate}
\item For all $\gamma$, the curve $\wittprojection_p \gamma$ is $p$-typical;
\item If $\gamma$ is $p$-typical then $\wittprojection_p \gamma=\gamma$;
\item $\wittprojection_p \wittprojection_p =\wittprojection_p$.
\end{enumerate}
\end{lemma}
\begin{proof}
We first check that the operator is well-defined.
The division by $m$ expressions make sense by Corollary \ref{cor:fg-p-local}.
Let $A_r \subset \A^1$ be the subcoalgebra spanned by $b_0,\dots,b_r$, so that $\A^1 = \colim_r A_r$.
We must show that when restricting along $\phi_m^* A_r \to \phi_m^* \A^1$, the sum defining the operator becomes finite.
For this, it will be enough to show that given $r$ there exists $N$ such that for all $n>N$ and any pointed curve $\gamma: \phi_{mn}^* \A^1 \to X$, the restriction of $\ver_n \gamma$ to $\phi_m^* A_r$ is zero.
Equivalently, using the definition of $\ver_n$, we want that the composite $\phi_m^* A_r \to \phi_m^* \A^1 \to \phi_{mn}^* \A^1$ (where the last map is the $n$-th power) is null.
This is dual to the composite $\F_p[X_{mn}] \xrightarrow{X_{mn} \mapsto X_m^n} \F_p[X_m]/X_m^{r+1}$, so $N=r$ works.

(1) Fix a prime $\ell \ne p$; we shall show that $\frob_\ell \wittprojection_p = 0$.
Let $S = \{m \mid (m,p) = 1 = (m,\ell)\}$.
Then \[ \wittprojection_p = \sum_{m \in S} \frac{\mu(m)}{m} \ver_m \frob_m + \sum_{m \in S} \frac{\mu(\ell m)}{\ell m} \ver_{\ell m} \frob_{\ell m}. \]
Using that $\mu(\ell m) = -\mu(m)$ for $m \in S$ and the identities in Lemma \ref{lemm:frob-ver-identities}, the claim follows quickly.

(2) Immediate.

(3) Clear from (1) and (2).
\end{proof}

\subsubsection{Witt Hopf algebra} \label{app:rw-witt-hopf}
Let $H\in\Hopf_{\F_p}$. Applying the $p$-typicalization operator of Lemma~\ref{lemm:p-typicalization} to $H^\comp_1$ defines a retract of the functor corepresented by $\Free\A^1$.
\begin{definition}\label{def:W}
The corresponding summand of $\Free\A^1$ is denoted by $W$ and called the \emph{Witt Hopf algebra}. As an algebra it is polynomial, with generators $u_i$ in degrees $p^i$ (after normalizing the first degree to one).
\end{definition}

\begin{remark} \label{rmk:W-lift}
Comparing \cite[Theorem 14.21]{hazewinkel-witt-vectors} with our construction of $p$-typicalization, we see that our definition of $W$ agrees with the usual (Hopf algebra representing) $p$-adic Witt vectors.
As in Example \ref{ex:lift-FreeA1}, $W$ lifts to $\HCAlg_{\Z_p}^V$ and, being torsion-free, the lift is determined by its indecomposables with their $V$-action.
In the notation of Example \ref{ex:Mp-hotimes-Mp}, it is given by $M_{p^\infty}$ (use \cite[\S14.8, (13.4)]{hazewinkel-witt-vectors}).
\end{remark}

\begin{lemma} \label{lemm:W-mod-F}
$F: \phi_p^* W \to W$ is injective and $W/F \wequi \Ga$, via the canonical map.
\end{lemma}
\begin{proof}
That $F$ is injective on $W$ follows from the same statement for $\Free\A^1$, which is clear because $\Free\A^1$ is a polynomial ring.
Since $F=0$ on $W/F$ as well as $\Ga$, by Proposition \ref{prop:indecomp-charp} it suffices to check that $W/F \to \Ga$ induces an isomorphism on indecomposables.
This is true by Remark \ref{rmk:W-lift} and Example \ref{ex:Ga-hotimes-Ga}.
\end{proof}

\begin{remark}\label{app:rw-projective}
The Hopf algebras $M_n\otimes \F_p$ above are projective generators of the additive category of connected graded Hopf algebras \cite{goerss-dieudonne,MR281709}. Hence $W$ is a filtered colimit of projectives, and its symmetric powers are direct sums of twists of $W$ (compare Example~\ref{ex:Mp-hotimes-Mp}). In particular, the symmetric powers of $W$ have flat underlying objects for $\hotimes$.
\end{remark}
Let us also observe the following.
\begin{lemma} \label{lemm:additive-curves-frob}
A pointed curve $\gamma \in X_n$ corresponds to a morphism of formal groups $\phi_n^* \Ga^\comp \to X$ if and only if $\frob_q \gamma = 0$ for all primes $q$.
\end{lemma}
\begin{proof}
In light of Lemma \ref{lemm:HCAlg-vs-FG} (identifying formal groups and connected Hopf algebras), an equivalent formulation is that the sequence \[ \bigotimes_q \Free \phi_q^* \A^1 \xrightarrow{F_{(\ph)}} \Free\A^1 \to \Ga \] is right exact.
For the $p$-typical part, this is Lemma \ref{lemm:W-mod-F} above.
\end{proof}

\subsection{An alternative presentation of the Ravenel--Wilson Hopf ring}\label{subsec:RW-p}
In this section we consider the case $\scr C = \Vect_{\F_p}^\N$, $A = L$.
We put $\RW_{(p)} := \RW(\Vect_{\F_p}^\N, L)$, and we shall construct an alternative presentation for this ring.

Here is a new example of a formal group.
\begin{example}\label{ex:fg-2}
Let $R \in \CAlg(\Hopf_{\F_p}^\N)_{\F_p[L]/}$.
Then the set $(R^1)^\comp(C)$ defined in Example \ref{ex:fg-1} acquires a new abelian group structure denoted $+^F$, where \[ \gamma_1 +^F \gamma_2 = \sum_{ij} [a_{ij}] \gamma_1^i \gamma_2^j. \]
This sum is well-defined by conilpotence (see Remark \ref{rmk:pointed-vanishing}; note that a connected finite-dimensional coalgebra is always conilpotent).
\end{example}
Via Examples \ref{ex:fg-1} and \ref{ex:fg-2}, we have thus found \emph{two} formal group structures on the set-valued functor $(R^1)^\comp$.
To distinguish them, we denote the operations by $+$ (Example \ref{ex:fg-1}) and $+^F$ (Example \ref{ex:fg-2}).
Correspondingly we obtain two notions of Frobenius operators, which we write $\frob_q$ and $\frob_q^F$, and we speak of $p$-typical and $p^F$-typical curves.

\begin{theorem} \label{thm:alternative-presentation}
\begin{enumerate}
\item There is a canonical isomorphism between the two functors \[ \CAlg(\Hopf_{\F_p}^\Z)_{L/} \to \Set \] sending $R$ to the $p$-typical curves, respectively $p^F$-typical curves (in $(R^1)^\comp$).
\item There is a pushout square in $\CAlg(\Hopf_{\F_p}^\Z)_{L/}$ of the form
\begin{equation*}
\begin{CD}
\Sym_L \phi_p^* W @>>> L_{(p)} \\
  @V\frob_p^FVV        @VVV \\
\Sym_L W @>>> \RW_{(p)}.
\end{CD}
\end{equation*}
  Here $W$ is placed in degree $1$ and $\Sym_L W$ corepresents, via (1), the functor of $p^F$-typical curves.
  Under this isomorphism, the map $\frob_p^F$ corepresents the operator of the same name.
\end{enumerate}
\end{theorem}
Let us elaborate slightly on part (2).
The functor $\Hopf_{\F_p} \to \Ab$, $H \mapsto \wittprojection_p H^\comp_1$ is corepresented by $W$, essentially by definition.
Consequently the Hopf ring $\Sym_L W$ corepresents the functor sending $R$ to the $p$-typical curves in $R^1$.
By part (1), this functor is isomorphic (but not equal) to the functor sending $R$ to the $p^F$-typical curves in $R^1$.
The operation $\frob_p^F$ on $p^F$-typical curves thus corresponds to an endomorphism (up to twisting) of the functor corepresented by $\Sym_L W$, whence an endomorphism (up to twisting) of this object.
This endomorphism is what appears in the left hand vertical map of the diagram.
\begin{proof}
Assuming (1), (2) is immediate from Lemma \ref{lemm:additive-curves-frob}.
We thus prove (1).
On $X=(R^1)^\comp$, we have two formal group structures denoted $+$ and $+^F$.
We consequently obtain two notions of $p$-typicalization, denoted $\wittprojection_p$ and $\wittprojection_p^F$.
We claim that more precisely, $\wittprojection_p$ induces a bijection from $p^F$-typical curves to $p$-typical curves.
Write $I \subset L_{(p)}$ for the augmentation ideal.
The two representing Hopf rings are positive, as follows from their presentations by free Hopf rings on the curve coalgebra. Since the generators of $I$ have strictly negative external degrees, in each fixed bidegree sufficiently high powers of $I$ vanish on connected components. Thus the two representing Hopf rings are $I$-adically complete.\gpt{Replaces references to abstract positivity and completion results.}
It thus suffices to prove the claim for $I$-adically complete Hopf rings (under $L_{(p)}$).
Let $R$ be one such.

\emph{Injectivity}:
Let $\gamma_1, \gamma_2: \A^1 \to R^1$ be $p^F$-typical curves with $\wittprojection_p \gamma_1 = \wittprojection_p \gamma_2$.
We need to show that $\gamma_1 = \gamma_2$, and by completeness, it will suffice to show that $\gamma_1$ and $\gamma_2$ have the same image in $R/I^n R$ for every $n$.
We do this by induction, the case $n=0$ being clear.
We have an exact sequence \[ \F_p \to I^nR/I^{n+1}R \to R/I^{n+1}R \to R/I^nR \to \F_p \in \Mod_{\F_p[L]}(\Hopf_{\F_p}^\Z). \]
Applying the (trivially!) left exact functor $\Hom(\A^1, (\ph)^1)$ we get \[ \begin{aligned} 0 &\to \Hom(\A^1,(I^nR/I^{n+1}R)^1) \\ &\to \Hom(\A^1,(R/I^{n+1}R)^1) \to \Hom(\A^1,(R/I^n R)^1). \end{aligned} \]
By induction, the images of $\gamma_1, \gamma_2$ in the right hand group agree, and consequently we find $\delta \in \Hom(\A^1, (I^nR/I^{n+1}R)^1)$ with \[ \overline{\gamma_1} = \overline{\gamma_2} + \delta \in \Hom(\A^1, (R/I^{n+1}R)^1). \]
Note that $R/I^{n+1}R$ is still a ring, and so curves in $R/I^{n+1}$ can still be added using the formal group law (producing an alternative group structure on $\Hom(\A^1, (R/I^{n+1}R)^1)$ denoted $+^F$).
We then have $\overline{\gamma_2} +^F \delta = \overline{\gamma_2} + \delta$ since $[a_{ij}]\delta = 0$ for any $a_{ij} \in I$.
Consequently $\delta = \overline{\gamma_1} -^F \overline{\gamma_2}$ is $p^F$-typical (by construction, $\wittprojection_p^F$ is a homomorphism for $+^F$).
On the other hand $\wittprojection_p(\delta) = \wittprojection_p(\overline{\gamma_1} - \overline{\gamma_2}) = 0$, since $\wittprojection_p$ is a homomorphism for the usual $+$ and $\gamma_1, \gamma_2$ are supposed to have the same $p$-typicalization.
We claim that $\wittprojection_p(\delta) = \wittprojection_p^F(\delta)$.
Assuming this, since $\delta$ was shown to be $p^F$-typical, we deduce that $\delta=0$ as needed.
To prove the claim, recall that the set-valued functor $((R/I^{n+1})^1)^\comp$ carries two formal group structures, one coming from the usual Hopf algebra structure and one involving the $[a_{ij}]$.
This functor has a subset, namely $((I^nR/I^{n+1}R)^1)^\comp$.
This is a subgroup for either formal group structure (namely the kernel of $R/I^{n+1}R \to R/I^nR$), and the two formal group structures on the subfunctor in fact coincide (because, as above, given $x: C \to I^nR/I^{n+1}R$ we have $[a_{ij}]x = 0$ for any $a_{ij} \in I$).
This proves that $\wittprojection_p = \wittprojection_p^F$ on the subfunctor, as needed.

\emph{Surjectivity}:
Let $\gamma: \A^1 \to R^1$ be $p$-typical.
To simplify notation, given $a, b: \A^1 \to R^1$, we shall say that $a=b \pmod{I^n}$ if $a$ and $b$ have the same image in $\Hom(\A^1, (R/I^nR)^1)$.
We shall construct a sequence of $p^F$-typical curves $\tilde \gamma_n$ such that for each $n$, (1) $\tilde\gamma_n = \tilde \gamma_{n-1} \pmod{I^{n-1}}$  and (2) $\wittprojection_p \tilde\gamma_n = \gamma \pmod{I^n}$.
Assuming this, the $\tilde \gamma_n$ define an element \[ \tilde \gamma \in \Hom(\A^1, R^1) \wequi \lim_n \Hom(\A^1, (R/I^nR)^1). \]
By construction, $\wittprojection_p^F(\tilde \gamma) = \tilde \gamma \pmod{I^n}$ for every $n$, i.e. $\wittprojection_p^F(\tilde \gamma)$ and $\tilde \gamma$ have the same image in $\Hom(\A^1, (R/I^nR)^1)$ for every $n$, so agree; in other words, $\tilde\gamma$ is $p^F$-typical.
Also by construction, $\wittprojection_p(\tilde \gamma) = \gamma \pmod{I^n}$ for every $n$, so $\wittprojection_p(\tilde \gamma) = \gamma$.
It thus suffices to construct inductively the $\tilde\gamma_n$; we set $\tilde\gamma_0 = [0]$.
We put \[ \begin{aligned} \tilde\gamma_{n+1} &= \wittprojection_p^F\bigl(\tilde \gamma_n -^F \\ &\qquad \wittprojection_p^F \wittprojection_p (\tilde \gamma_n - \gamma)\bigr). \end{aligned} \] and verify (1) and (2).
By induction, $\wittprojection_p(\tilde \gamma_n - \gamma) = \wittprojection_p(\tilde \gamma_n) - \wittprojection_p(\gamma) = 0 \pmod{I^n}$.
This implies that $\tilde\gamma_{n+1} = \wittprojection_p^F(\tilde \gamma_n -^F [0]) = \tilde \gamma_n \pmod{I^n}$, which is (1).
As in the proof of injectivity, this also implies that \[ (*) \quad \wittprojection_p^F \wittprojection_p(\tilde \gamma_n - \gamma) = \wittprojection_p \wittprojection_p(\tilde \gamma_n - \gamma) \pmod{I^{n+1}} \] and \[ (**) \quad a -^F \wittprojection_p(\tilde \gamma_{n}-\gamma) = a - \wittprojection_p(\tilde \gamma_{n}-\gamma) \pmod{I^{n+1}} \] for any $a: \A^1 \to R^1$.
Now we compute
\begin{align*}
\wittprojection_p\tilde\gamma_{n+1} &= \wittprojection_p\wittprojection_p^F\big( \tilde\gamma_{n} -^F \wittprojection_p^F\wittprojection_p(\tilde \gamma_{n}-\gamma)\big) \pmod{I^{n+1}} \\
&=\wittprojection_p(\wittprojection_p^F \tilde\gamma_{n} -^F \wittprojection_p^F\wittprojection_p^F\wittprojection_p(\tilde \gamma_{n}-\gamma)\big) \pmod{I^{n+1}} \\
&=\wittprojection_p(\tilde\gamma_{n} -^F \wittprojection_p^F\wittprojection_p(\tilde \gamma_{n}-\gamma)\big) \pmod{I^{n+1}} \\
&\stackrel{(*)}{=}\wittprojection_p(\tilde\gamma_{n} -^F \wittprojection_p\wittprojection_p(\tilde \gamma_{n}-\gamma)\big) \pmod{I^{n+1}} \\
&=\wittprojection_p(\tilde\gamma_{n} -^F \wittprojection_p(\tilde \gamma_{n}-\gamma)\big) \pmod{I^{n+1}} \\
&\stackrel{(**)}{=}\wittprojection_p(\tilde\gamma_{n} - \wittprojection_p(\tilde \gamma_{n}-\gamma)\big) \pmod{I^{n+1}} \\
&=\wittprojection_p\tilde\gamma_{n} - (\wittprojection_p\tilde \gamma_{n}-\gamma) \pmod{I^{n+1}} \\
&=\gamma \pmod{I^{n+1}}.
\end{align*}
\end{proof}

\subsection{The derived Ravenel--Wilson Hopf ring}\label{app:rw-bar}
We recall the definition of the map to be used in the bar construction. For any Hopf ring $T$ under $L$, the pointed curves in $T^1$ have two formal group laws: the one induced by the Hopf-algebra addition, and the one obtained by applying the universal formal group law $\mathcal F(X,Y)=\sum a_{ij}X^iY^j$ using the Hopf-ring product. We denote the associated Frobenius operations by $\frob_p$ and $\frob_p^{\mathcal F}$, respectively. The operation $\frob_p^{\mathcal F}$ on the $p^{\mathcal F}$-typical universal curve determines a morphism
\[
 f=\frob_p^{\mathcal F}:\Sym_L(\phi_p^*W)\longrightarrow \Sym_L(W).
\]
The augmentation $\epsilon:\Sym_L(\phi_p^*W)\to L$ kills the positive external degrees. Put
\[
 A=\Sym_L(\phi_p^*W),\qquad B=\Sym_L(W),\qquad
 R=B\hotimes_A L.
\]
By Theorem~\ref{thm:alternative-presentation}, we have $R\simeq\RW_{(p)}$.

\begin{definition}\label{def:rw-bar}
Define the derived Ravenel--Wilson Hopf ring $\dRW_{(p)}$ to be the simplicial Hopf ring given by the two-sided bar construction
\[
 (\dRW_{(p)})_n=B\hotimes_L A^{\hotimes_L n}\hotimes_L L.
\]
The first and last face maps are induced by $f$ and $\epsilon$; the remaining face maps multiply adjacent $A$-factors. Degeneracies insert the unit. There is a canonical augmentation $\dRW_{(p)}\to R$.
\end{definition}

\begin{theorem}\label{thm:rw-bar-exact}
The augmentation $\dRW_{(p)}\to R$ induces an isomorphism on zeroth homology and $\hH_i(\dRW_{(p)})=0$ for $i>0$.
\end{theorem}

\begin{remark}\label{rmk:rw-bar-derived}
The degree-zero statement follows formally from the bar faces: $\hH_0(\dRW_{(p)})=B\hotimes_A L$.
In the language of derived Hopf rings, the above bar construction computes the derived pushout
\[
\LSym_L^\hotimes(W)\hotimes_{\LSym_L^\hotimes(\phi_p^*W)}L.
\]
\end{remark}

The remainder of this appendix is dedicated to proving Theorem \ref{thm:rw-bar-exact}.

\subsection{The special fiber}\label{app:rw-fiber}

\begin{lemma} \label{lemm:filter-derived-RW}
\begin{enumerate}
\item Write $\tilde y: \A^1 \to W \to (\Sym_{L_{(p)}} W)^1 \to \dRW_{(p)}^1$ for the canonical ($[p]$-typical) coordinate.
  Then the composite $\A^1 \xrightarrow{\tilde y} \dRW_{(p)}^1 \to (\dRW_{(p)} \hotimes_L \Z_{(p)})^1$ is an equivalence of coalgebras (so in particular the right hand side is concentrated in homological degree zero) compatibly with addition, i.e., of Hopf algebras.
\item Consequently $(\dRW_{(p)} \hotimes_L \F_p)^1 \wequi \cof(\Ga \xrightarrow{p} \Ga) \wequi \Ga \htimes \Ga[1]$.
  Denote the copy of $\Ga$ in degree zero by $y$ and the one in homological degree $1$ by $x$.
\item We have \[ \hH_* (\dRW_{(p)} \hotimes_L \F_p) \wequi \Sym(\Ga) \hotimes \Lambda(\Ga), \] where the generator $\Ga \wequi \Sym^1(\Ga)$ is $y$ and the generator $\Ga \wequi \Lambda^1(\Ga)$ is $x$.
\end{enumerate}
\end{lemma}
\begin{proof}
By Definition~\ref{def:rw-bar}, the specialization $\dRW_{(p)} \hotimes_L \Z_{(p)}$ is obtained as the bar construction corresponding to the pushout
\begin{equation*}
\begin{CD}
\Sym_{\Z_{(p)}} \phi_p^* W @>>> \Z_{(p)} \\
  @V\frob_p^FVV        @VVV \\
\Sym_{\Z_{(p)}} W @>>> \dRW_{(p)} \hotimes_L \Z_{(p)}.
\end{CD}
\end{equation*}
Because the formal group law on $\Sym_{\Z_{(p)}} W$ is additive, the map $\frob_p^F$ is just $\Sym(\frob_p)$\NB{details?}, i.e., $\Sym(F)$ (Lemma \ref{lemm:FV-on-W}).
Note that $\Sym^0_{\Z_{(p)}}(X) = \F_p[\Z_{(p)}]$ and $\Sym^i_{\Z_{(p)}}(X) = \Sym^i(X)$ whenever $X$ is connected.

(1) Restricting the bar construction defining $\dRW_{(p)} \hotimes_L \Z_{(p)}$ to degree $1$ produces a bar construction modelling $\cof(F: \phi_p^* W \to W)$.
This is $\Ga$ by Lemma \ref{lemm:W-mod-F}.

(2) Immediate.

(3)
We use the filtration by iterates of Frobenius on the normalized simplicial Hopf-algebra complex. Its associated graded has a spectral sequence, converging on connected components by the positive internal grading.\gpt{Termwise Frobenius filtration replaces abstract derived filtration.}
Moreover, this construction is multiplicative.
The degree zero part of $\dRW_{(p)} \hotimes_L \F_p$ is clearly just $\F_p[\F_p]$, so we need only determine the connected part.
We will thus determine $\hH_*(\dRW_{(p)} \hotimes_L \F_p)$ by studying the spectral sequence of the filtration by Frobenius on $\dRW_{(p)} \hotimes_L \F_p$.

\emph{Step 1. Determining the associated graded.}
Since $F$ commutes with $\phi_p^*$, the associated graded in each (filtration) degree is just a twist of the associated graded in degree zero.
It thus suffices to determine associated graded in degree zero, i.e. to determine $\hH_*(\dRW_{(p)} \hotimes_L \F_p/F)$.
We will first do this for $\hH_*(\dRW_{(p)} \hotimes_L \Z_{(p)}/F)$ and then mod out by $[p]$ in a second step.
Note that $W$ is a filtered colimit of projective generators (Remark~\ref{app:rw-projective}), so its derived tensor products are computed by ordinary tensor products.
The symmetric powers of $W$ split into sums of twists of $W$ (Remark~\ref{app:rw-projective}).
Consequently derived tensor products of these objects are also zero-truncated.
We find that the bar construction modelling $\dRW_{(p)} \hotimes_L \Z_{(p)}$ is a simplicial object with entries ordinary (underived) Hopf algebras, all of which are in fact sums of twists of $W$.
To determine the cofiber of $F$ on $\dRW_{(p)} \hotimes_L \Z_{(p)}$, we may take the cofiber of $F$ termwise in the bar construction.
By Lemma \ref{lemm:W-mod-F}, this turns every summand $W$ into $\Ga$.
It follows that the simplicial Hopf ring $\dRW_{(p)} \hotimes_L \Z_{(p)}/F$ is modeled by the simplicial object which is the \emph{underived} bar construction corresponding to 
\begin{equation*}
\begin{CD}
\Sym_{\F_p} \phi_p^* \Ga @>>> \F_p \\
  @V\Sym(0)VV        @VVV \\
\Sym_{\F_p} \Ga @>>> (\dRW_{(p)} \hotimes_L \Z_{(p)})/F.
\end{CD}
\end{equation*}
Using the equivalence between Hopf algebras with $F=0$ and vector spaces with Verschiebung (Proposition \ref{prop:indecomp-charp}), Dold--Kan identifies this simplicial object with $\Sym(\Ga) \hotimes \Lambda(\phi_p^* \Ga)$.
Here and below, $\Sym$ and $\Lambda$ without subscripts means over $\F_p[\F_p]$ (which only affects the degree zero part).
As usual, $\Sym(\Ga)$ is placed in chain degree zero, and $\Lambda^1(\Ga)$ in degree one.
The differential is zero (since the maps inducing the above square are), and so we have found that \[ \hH_*(\dRW_{(p)} \hotimes_L \Z_{(p)}/F) \wequi \Sym(\Ga) \hotimes \Lambda(\phi_p^* \Ga). \]
In particular, $[p]=0$ on these Hopf algebras, from which we deduce that \[ \hH_*(\dRW_{(p)} \hotimes_L \F_p/F) \wequi \Sym(\Ga) \hotimes \Lambda(\phi_p^* \Ga) \hotimes \Lambda(\epsilon), \] where $\epsilon$ is placed in chain degree one and internal degree zero.

\emph{Step 2. Determining the $E_1$-page.}
In what follows, almost all our manipulations involve Hopf algebras with $F=0$, and so we will mainly be thinking about indecomposables.
We write $\alpha_{(0)}, \alpha_{(1)}$ etc. for the generators of $\Sym^1(\Ga)$.
We will systematically need to deal with also $\phi_p^*\Sym^1(\Ga)$, $\phi_p^{2*}\Sym^1(\Ga)$, and so on.
As in Example \ref{ex:Ga-hotimes-Ga}, we will use the (surjective but not injective) Verschiebung maps to identify generators of the twists.
Thus e.g. $\phi_p^* \Ga$ has generators $\alpha_{(1)}, \alpha_{(2)}$, and so on.
We symbolically denote this correspondence as $\phi_p^* \Ga \wequi \Ga/\alpha_{(0)}$, and so on.
We denote the generators of $\Lambda^1(\phi_p^* \Ga)$ by $\beta_{(1)}, \beta_{(2)}$, and so on, and similarly for $\Lambda^1 \phi_p^{2*}\Ga$, etc.
Let us denote by $G^*$ the associated graded of the (homology of the) filtration by Frobenius on $\dRW_{(p)} \hotimes_L \F_p$.
Since $VF = [p] = FV$, Verschiebung induces a map $G^i \to G^{i+1}$ (not to be confused with the internal Verschiebung on each $G^i$!).
We will view $G^*$ as an $\F_p[t]$-module via this operator.
As such, we identify\NB{details?} \[ G^* \wequi \Sym(\Ga) \hotimes \Lambda(\Ga) \hotimes \Lambda(\epsilon) \hotimes \F_p[t] / (\beta_{(0)}, t\alpha_{(0)}, t\beta_{(1)}, t^2\alpha_{(1)}, \dots, t^i\alpha_{(i-1)},t^i \beta_{(i)}, \dots). \]
This is the $E_1$-page of our spectral sequence, converging strongly to $\hH_*(\dRW_{(p)} \hotimes_L \F_p)$ (at least on connected parts).
Its $d_1$-differential is induced by the filtration from the original differentials in the bar construction.
Using again that $Fa \hotimes b = F(a \hotimes Vb)$ one checks\NB{details} that the differential is determined by \[ d(\beta_{(i)}) = t\alpha_{(i)}, \quad d(\alpha_{(i)}) = 0, \quad d(\epsilon) = t. \]

\emph{Step 3. Studying the spectral sequence.}
Write \[ \tilde G^* = \Sym(\Ga) \hotimes \Lambda(\Ga) \hotimes \Lambda(\epsilon) \hotimes \F_p[t] \] with differential as above, so that there is an evident surjection $\tilde G^* \to G^*$.
We shall show at the end that $\tilde G^* \to G^*$ is a quasi-isomorphism.
Assuming this, write \[ \tilde G_0^* = \Sym(\Ga) \hotimes \Lambda(\Ga)\hotimes \F_p[t], \] also with the above differential.
Then we see that \[ \tilde G^* \cong \tilde G_0^* \otimes^L_{\F_p[t]} \F_p \wequi \tilde G_0^* \otimes_{\F_p[t]} \F_p, \] the latter equivalence since $\tilde G_0^*$ consists of free $\F_p[t]$-modules.
But \[ G_0^* \otimes_{\F_p[t]} \F_p \wequi \Sym(\Ga) \hotimes \Lambda(\Ga) \] placed entirely in $*=0$, so in particular with zero differential.
This proves what we want.

\emph{Step 4. Showing that $\tilde G^* \to G^*$ is a quasi-isomorphism.}
We work entirely on indecomposables.
Let $A_0 \subset \tilde G^*$ be the submodule spanned by $\beta_{(0)}$, so that \[ A_0 = \beta_{(0)} \wedge \Sym(\alpha_{(i)}, i \ge 0) \otimes \Lambda(\beta_{(j)}, j \ge 1) \otimes \Lambda(\epsilon) \otimes \F_p[t]. \]
For $\beta_{(0)} \wedge F \in A_0$ we have $d(\beta_{(0)} \wedge F) = t\alpha_{(0)} F \pm \beta_{(0)} \wedge dF$ which vanishes only if $F=0$.
Hence $d: A_0 \to dA_0$ is an isomorphism and $\tilde G^* \to \tilde G^*/(\beta_{(0)}, t\alpha_{(0)})$ is a quasi-isomorphism.
Similarly we can let $A_1 \subset \tilde G^*/(A_0, dA_0)$ be the submodule spanned by $t\beta_{(1)}$ and argue as before\NB{details?} to see that $\tilde G^*/(\beta_{(0)}, t\alpha_{(0)}) \to \tilde G^*/(\beta_{(0)}, t\alpha_{(1)}, t\beta_{(1)}, t^2\alpha_{(1)})$ is a quasi-isomorphism.
Continuing this way (using that quasi-isomorphism are stable under filtered colimits) we obtain the desired result.
\end{proof}

\subsection{A regular-sequence calculation}\label{app:rw-koszul}

We now discuss an exactness property.
We discuss it purely as a piece of algebra.
Consider the Hopf ring \[ \scr K = \Lambda(\Ga) \hotimes \Sym(\Ga) \hotimes \F_p[\F_p[a_0, a_1, a_2, \dots]] \in \Alg(\Hopf_{\F_p}) \] (which is not commutative because of the exterior algebra part).
There are canonical ``coordinates'' \[ x: \Ga \wequi \Lambda^1 \Ga \to \scr K \quad\text{and}\quad y: \Ga \wequi \Sym^1 \Ga \to \scr K. \]
We shall promote $\scr K$ to a differential graded Hopf ring, so in particular, a chain complex in $\Hopf_{\F_p}$.
To this end, we declare that $x$ has degree $1$ and everything else has degree $0$, so that \[ \scr K_n = \Lambda^n(\Ga) \hotimes \Sym(\Ga) \hotimes \F_p[\F_p[a_0, a_1, a_2, \dots]] \in \Hopf_{\F_p}. \]
The differential $d: \scr K_* \to \scr K_{*-1}$ is induced by \[ d(x) = \sum_i a_i y^{p^i}. \]
In other words, we are looking at the Koszul complex of the ``element'' $d(x)$.

\begin{remark} \label{rmk:build-complex}
Some explanations are in order.
First note that $d(x)$ is well-defined, since $y: \Ga \to \Sym(\Ga) \hotimes \F_p[\F_p[a_0, a_1, a_2, \dots]]$ is pointed and $\Ga$ is ind-conilpotent (using ind-conilpotence of $\A^1$); it is a map of Hopf algebras by the Frobenius identity for Hopf-ring multiplication.
Now we define a map
\begin{gather*}
\tilde d_n: \Ga^{\htimes n} \longrightarrow
\Lambda(\Ga) \hotimes \Sym(\Ga) \hotimes \F_p[\F_p[a_0,a_1,\dots]],\\
\tilde d_n = \sum_i [(-1)^{i-1}]
 p_1\cdot\dots\cdot p_{i-1}\cdot (d_1\circ p_i)
 \cdot p_{i+1}\cdot\dots\cdot p_n.
\end{gather*}
Here $p_i: \Ga^{\htimes n} \to \Ga$ is projection to the $i$-th factor, $d_1$ is $d(x)$, and the multiplication uses the ring structure of $\scr K$.
Then $\tilde d_n$ is multilinear, so passes to $\Ga^{\hotimes n}$, and moreover satisfies ``$d(a \hotimes a) = 0$'', so passes to the exterior power (by definition of the exterior power).
Extending linearly over the remaining tensor factors, we obtain $d_n: \scr K_n \to \scr K_{n-1}$.
This is easily verified to be a differential (i.e. square to zero).
\end{remark}

\begin{lemma} \label{lemm:regular-sequence}
The complex $\scr K$ is acyclic.
In other words, for $n > 0$, the sequence of Hopf algebras \[ \scr K_{n+1} \xrightarrow{d} \scr K_n \xrightarrow{d} \scr K_{n-1} \in \Hopf_{\F_p} \] is exact.
\end{lemma}
\begin{proof}
The internal degree zero part of $\scr K$ is concentrated in chain degree zero, and so acyclic.
It thus suffices to consider the connected part of $\scr K$.
Since $F=0$ on $\Ga$, the same is true on $\Ga^{\hotimes n}$ and thus on $\Lambda(\Ga), \Sym(\Ga)$ (see Example \ref{ex:Ga-hotimes-Ga}).
In other words, $F=0$ on all of $\scr K^\cn$.
To check exactness, we may thus apply the symmetric monoidal functor of indecomposables (Proposition \ref{prop:indecomp-charp}).
Now \[ Q(\Sym(\Ga) \hotimes \F_p[\F_p[a_0, a_1, a_2, \dots]]) \wequi \F_p[\alpha_{(0)}, \alpha_{(1)}, \dots] \otimes \F_p[a_0, a_1, a_2, \dots] \] and $Q(\scr K)$ is the Koszul complex of the sequence $r_1, r_2, \dots$, where \[ r_j = \sum_{i} a_i F^i(\alpha_{(j)}). \]
Here $F^i(\alpha_{(j)})$ means the $i$-fold Hopf ring Frobenius applied to $\alpha_{(j)}$.
To determine $F(\alpha_{(j)})$, we must consider $\Delta_p(\alpha_{(j)}) \in \Ga^{\hotimes p}$ and its image in $Q(\Sym^p \Ga)$.
Using the formula for the coproduct on $\Ga$ (see Example \ref{ex:Ga-hotimes-Ga}) we see that $F(\alpha_{(j)})=\alpha_{(j-1)}^p$ and thus $F^i(\alpha_{(j)}) = \alpha_{(j-i)}^{p^i}$.
(In particular the sum is finite, as we knew it had to be.)
Regularity of this sequence is proved by Ravenel--Wilson \cite[Lemma 4.15(b)]{ravenel-wilson}.\NB{They do not have $a_0$, but this is irrelevant since our $r_0$ is just $a_0$...}
\end{proof}

\subsection{The comparison}\label{app:rw-comparison}

Put $I=\ker(L\to\Z_{(p)})$ and $J=\ker(L\to\F_p)=(p,I)$.

\begin{lemma}\label{lemm:coefficients}
Write $[p]_{\mathcal F}(T)=\sum_{m\ge1}c_mT^m\in L\fpsr{T}$. All $c_m$ lie in $J$, so $[p]_{\mathcal F}(T)\in J\fpsr{T}$, and modulo $J^2$ we have, in $(J/J^2)\fpsr{T}$,
\[ [p]_{\mathcal F}(T)=\sum_{i\ge0}a_iT^{p^i},\qquad a_0=\bar p. \]
The elements $a_i\in J/J^2$ are linearly independent and can be completed to a homogeneous basis. Consequently
\[\operatorname{gr}_J L\simeq \F_p[a_0,a_1,\ldots,b_1,b_2,\ldots].\]
\end{lemma}
\begin{proof}
This follows from standard formal group law theory; see, e.g., \cite[Appendix A2.1--A2.2]{ravenel1986complex}.
\end{proof}

In order to prove our main theorem, we filter $\dRW_{(p)}$ by powers of $J$ and study the resulting spectral sequence.
The coefficient exact sequences \[0\to J^n/J^{n+1}\to L/J^{n+1}\to L/J^n\to0\] give a Bockstein exact couple, and we now determine the first differential.
\begin{lemma}\label{lemm:bockstein}
Put $S=\dRW_{(p)}\hotimes_L \F_p$, and let $y$ and $x$ be the generators of $\hH_0(S)^1$ and $\hH_1(S)^1$ from Lemma~\ref{lemm:filter-derived-RW}. The first differential in the above exact couple is the derivation determined by
\[d_1(x)=\sum_{i\ge0}a_i y^{p^i},\qquad d_1(y)=0,\qquad d_1(a)=0\quad(a\in\operatorname{gr}_J L).\]
\end{lemma}
\begin{proof}
Let $\tilde y:\A^1\to \dRW_{(p)}^1$ be the canonical $[p]$-typical coordinate. The defining relation gives $\frob_p^{\mathcal F}(\tilde y)=0$. Applying $\ver_p^{\mathcal F}$ yields
\[ [p]_{\mathcal F}(\tilde y)=0\quad\text{in }\hH_0(\dRW_{(p)}). \]
Since every coefficient of the $p$-series belongs to $J$, the expression
\[r=\sum_m c_m\tilde y^m\in\hH_0(J\hotimes_L \dRW_{(p)})\]
is defined by conilpotence and maps to zero in $\hH_0(\dRW_{(p)})$. The coefficient fiber sequence $J\hotimes_L \dRW_{(p)}\to \dRW_{(p)}\to S$ thus gives a class $z\in\hH_1(S)^1$ with $\partial_J(z)=r$.

We identify $z$ by passing to the additive specialization $L\to\Z_{(p)}$. Write $Z=\dRW_{(p)}\hotimes_L\Z_{(p)}$. By Lemma~\ref{lemm:filter-derived-RW}(1), $Z^1\simeq\Ga$, with the image of $\tilde y$ the identity curve. The connecting homomorphism for multiplication by $p$ is therefore an isomorphism
\[ \partial_p:\hH_1(S)^1\xrightarrow{\simeq}\hH_0((p)\hotimes_{\Z_{(p)}}Z)^1,\]
and is normalized by $\partial_p(x)=p\otimes y$. Naturality of the connecting homomorphism identifies the image of $r$ with $p\otimes y$, since $[p]_{\mathcal F}(T)$ becomes $pT$. Thus $z=x$ and
\[\partial_J(x)=\sum_m c_m\tilde y^m.\]
Passing to $J/J^2$, Lemma~\ref{lemm:coefficients} gives the asserted formula for $d_1(x)$. The class $y$ lifts to $\dRW_{(p)}$, and the coefficient classes lift from $L$, so their differentials vanish. Multiplicativity of the exact couple determines the differential on all of its first page.
\end{proof}

\begin{lemma}\label{lemm:rw-positivity}
The normalized bar complex $\dRW_{(p)}$ is positive: $ (\dRW_{(p)}^r/\F_p)(s)=0$ for $s<r$. For a positive complex, its derived $I$-adic quotient tower ($I=\ker(L\to\Z_{(p)})$) is degreewise eventually constant. The connected part of a connective graded Hopf-algebra complex over $\F_p$ is $[p]$-complete.
\end{lemma}
\begin{proof}
For the first two assertions, positivity holds on $L$ and on $W$ placed in external degree one, and is preserved by tensor products, symmetric powers, and the bar realization. In any fixed bidegree, a sufficiently high product of elements of $I$ comes from external degree exceeding the internal degree; the same bound applies to a bar resolution of a derived $I$-quotient. For the last assertion one reduces to graded Hopf algebras and uses $[p]=FV=VF$; on the connected part, iterations of Verschiebung lower the positive internal degree until they vanish. Hence the localization at $[p]$ is zero and the connected part is $[p]$-complete. Compare \cite[Proposition 9.1]{goerss-dieudonne}.
\end{proof}

\begin{lemma}\label{lemm:joint-completeness}

The connected part $\dRW_{(p)}^{\cn}$ is derived $J$-complete: the canonical map
\[\dRW_{(p)}^{\cn}\longrightarrow\lim_n(\dRW_{(p)}\hotimes_L L/J^n)^{\cn}\]
is an equivalence.
\end{lemma}
\begin{proof}
The defining bar construction shows, by Lemma~\ref{lemm:rw-positivity}, that $\dRW_{(p)}$ is positive:
\[\big(\dRW_{(p)}^r/\F_p\big)(s)=0\qquad(s<r).\]
Put $K_n=(p^n,I^n)$. The inclusions
\[J^{2n}\subseteq K_n\subseteq J^n\]
show that the $K_n$-adic and $J$-adic quotient towers are pro-equivalent.

Choose homogeneous polynomial generators $t_j$ for $I$, all in strictly negative external degrees. The degree bound of Lemma~\ref{lemm:rw-positivity} applies uniformly to their products: in any fixed external and internal bidegree, the contribution of $I^n$ to a positive object vanishes for $n\gg0$. The same degree bound applies to derived $I$-quotients, obtained from the bar construction, and holds uniformly after reduction modulo $p^n$. Thus the tower
\[(\dRW_{(p)}/p^n)^{\cn}\longrightarrow ((\dRW_{(p)}/p^n)\hotimes_L L/I^n)^{\cn}\]
is degreewise eventually an equivalence as $n\to\infty$. Since $L/I^n$ is free over $\Z_{(p)}$, derived base change identifies the right hand side with $(\dRW_{(p)}\hotimes_L L/K_n)^{\cn}$. Consequently
\[\lim_n(\dRW_{(p)}\hotimes_L L/J^n)^{\cn}\simeq\lim_n(\dRW_{(p)}/p^n)^{\cn}.\]
The connected part of $\dRW_{(p)}$ is $[p]$-complete by Lemma~\ref{lemm:rw-positivity}, applied to the connected graded Hopf algebras underlying its external components (and passage to normalized simplicial complexes). This identifies the last limit with $\dRW_{(p)}^{\cn}$.
\end{proof}

We can now establish our main result.

\begin{proof}[Proof of Theorem \ref{thm:rw-bar-exact}.]
With the notation of Lemma~\ref{lemm:bockstein}, let
\[P_n=\dRW_{(p)}\hotimes_L L/J^{n+1},\qquad V_n=(J^n/J^{n+1})\hotimes_{\F_p} S.\]
There are fiber sequences $V_n\to P_n\to P_{n-1}$, and hence a multiplicative exact couple whose first page is
\[E_1^{s,t}= (J^s/J^{s+1})\otimes_{\F_p}\hH_t(S).\]
By Lemma~\ref{lemm:filter-derived-RW}, this is
\[\operatorname{gr}_J L\otimes_{\F_p}\Sym(\Ga)\hotimes\Lambda(\Ga).\]
The differential of Lemma~\ref{lemm:bockstein}, together with the choice of generators in Lemma~\ref{lemm:coefficients}, identifies this complex with $\scr K\otimes_{\F_p} \F_p[b_1,b_2,\ldots]$. Lemma~\ref{lemm:regular-sequence} therefore gives
\[E_2^{s,t}=0\qquad(t>0).\]

We explain why this implies vanishing of positive homology, without assuming convergence for the uncompleted ring. By Lemma~\ref{lemm:joint-completeness}, the connected part of the tower is complete and separated. Represent it by a complete filtered complex. If $z$ is a cycle in positive homological degree and filtration $s$, its leading term is a $d_1$-cycle. The vanishing above allows us to subtract a boundary to move $z$ into filtration $s+1$. Iterating gives a sequence of corrections in strictly increasing filtration whose sum converges, and whose differential is $z$. Hence $\hH_t(\dRW_{(p)}^{\cn})=0$ for $t>0$.

In internal degree zero the defining bar construction is just the coefficient ring $L$, and is already discrete. We conclude that $\hH_t(\dRW_{(p)})=0$ for $t>0$. Since $\hH_0(\dRW_{(p)})=R$, the augmentation is a homology equivalence.
\end{proof}

\bibliographystyle{alpha}
\bibliography{bibliography}

\end{document}